# Composite local observables of sine–Gordon via stochastic analysis

by Massimiliano Gubinelli, Peter Paulovics
University of Oxford
*Email:* {massimiliano.gubinelli, peter.paulovics}@maths.ox.ac.uk

**Abstract**

We develop a rigorous framework for the description of general local composite observables in Euclidean quantum field theories. The approach is nonperturbative and combines renormalisation group ideas with techniques of stochastic analysis. Via renormalisation dictated by sufficiently precise approximate solution theory of corresponding flow equations, we construct various local observables of the sine–Gordon model on the full space up to $\beta^2 < 6\pi$, including highly singular ones such as vertex fields $\exp(i\alpha\phi(y))$ in the optimal Wick regime $|\Re e(\alpha)| < 4\pi / \beta$ and derivatives of arbitrarily high order, i.e. $\partial^\mu\phi(y)$ for any multiindex $\mu \in \mathbb{N}_0^2$. Applications of the approach include a direct proof of smoothness of observable-inserted correlation functions, and verification that the constructed objects satisfy a version of the prefactorization algebra axioms adapted to our setting.

## Table of contents

# 1 Introduction

The purpose of this paper is to develop an understanding of observables of Euclidean (quantum) field theories, i.e. to provide a nonperturbative, probabilistic framework for the construction and analysis of nonlinear functionals of certain generalized random fields. Euclidean field theories are notoriously difficult to construct and describe, in part owing to the need for renormalisation: heuristically, their laws are given by the Radon–Nikodym derivative

$$“\mathrm{d}\nu(\phi) \propto e\left(-\lambda\int_{\mathbb{R}^d}\tilde{V}(\phi(x))\mathrm{d}x\right)\mathrm{d}\mu(\phi)”$$

with $\mu$ a measure on $\mathcal{S}'(\mathbb{R}^d)$ that is not supported on any space of bona fide functions, thus making the expression $\tilde{V}(\phi(x))$ for general nonlinear $\tilde{V}$ a priori ill-defined. To give rigorous meaning to such measures $\nu$, one therefore needs to understand how the nonlinearities give rise to singularities, and remove them in a way that does not lead to trivial objects.

After a random field $\varphi_\infty$ with $\mathrm{Law}(\varphi_\infty) = \nu$ has been constructed, one would like to study expressions of the form $\int \prod_j q^{(j)}(\phi(y_j))\mathrm{d}\nu(\phi)$ with $y_j \in \mathbb{R}^d$ and $q^{(j)}$ nonlinear functions. One then finds themselves facing the above renormalisation problem *again*. What is even worse, the $q^{(j)}$ are allowed to be much more singular than $\tilde{V}$, in fact so much that for *any* field $\varphi_\infty$, it is not guaranteed that after renormalisation, the "correct" object $[\![q^{(j)}(\varphi_\infty(y_j))]\!]$ exists as an $\mathcal{S}'(\mathbb{R}^d)$-valued random variable (take, for example, monomials of sufficiently high order derivatives of $\varphi_\infty$).

Our approach is based on the scale-by-scale coupling of the free field with the interacting EQFT [5]. In this framework, the target random field $\varphi_\infty$ is realised as the terminal value of a field-valued Markov process $(\varphi_t)_{t\in\mathbb{R}_+}$. The interpretation of the parameter $t$ is energy scale or resolution; the law of $\varphi_t$ is the Wilsonian effective description of the interacting field theory. The infinitesimal generator $\mathscr{L}^V$ of the process is determined by the effective potential $(V_t)_{t\in\mathbb{R}_+}$ of the theory, which in turn is given by the solution of the renormalisation group flow equation.

To renormalise a nonlinear functional $q$ of the field $\varphi_\infty$, one then needs to study the kernel of the unbounded operator $\mathscr{L}^V$. Since precise control of the Polchinski equation is in general difficult without a priori knowledge about the theory, we construct scale-dependent functionals $(q_t)_{t\in\mathbb{R}_+}$ such that $\mathscr{L}^V q \approx 0$ *approximately* in some appropriate topology. To obtain the full nonperturbative correlation functions, we then employ a probabilistic correction. Indeed, by Ito's formula, the process $\mathfrak{Q} := (\mathfrak{Q}_t)_{t\in\mathbb{R}_+}$, where

$$\mathfrak{Q}_t := q_t(\varphi_t) + \mathbb{E}\left[\int_t^\infty \mathscr{L}^V q_s(\varphi_s)\mathrm{d}s \middle| \varphi_t\right], \tag{1.1}$$

is a martingale in the filtration of $\varphi$. Whenever the martingale $\mathfrak{Q}$ is closed with $\sigma(\varphi_\infty)$-measurable terminal value, the random variable $\mathfrak{Q}_\infty$ is the functional $q_\infty(\varphi_\infty)$ of the UV-complete field one is interested in. However, the above martingale, and hence the expectation $\langle\mathfrak{Q}\rangle := \mathbb{E}[\mathfrak{Q}_0]$, is well-defined even if $q_\infty(\varphi_\infty)$ fails to exist as a random element of $\mathcal{S}'(\mathbb{R}^d)$. This latter failure of existence is indeed the case for sufficiently high Euclidean dimension $d$ or sufficiently singular functional $q_\infty$ (e.g. Wick square of the free field in $d = 4$, cf. [2], see also [1]). We define the martingales $\mathfrak{Q}$ as the observables of the theory $\varphi_\infty$, and denote by $[\![q(\varphi)]\!] := \mathfrak{Q}$ where $q$ is the "classical core" of the functional they arise as renormalisations of; see Section 2 for precise definitions.

We carry out the above program for various observables of the infinite volume massive sine–Gordon model, which is defined as follows. Given a $T < \infty$, $\rho \in \mathcal{C}_c^\infty(\mathbb{R}^2)$ with $\rho(x) \in [0, 1]$, $\lambda \in \mathbb{R}$, $\beta > 0$, denote by

$$\nu_T^\rho(\phi) := Z_T^{-1}\exp\left(-\lambda\int_{\mathbb{R}^2}\rho(x){:}\cos(\beta\phi(x)){:}_T\,\mathrm{d}x - c_T\right)\mathrm{d}\mu_T(\phi) \tag{1.2}$$

where $\mu_T$ is the centered Gaussian measure with covariance $\int_0^T s^{-2}\exp(-s^{-1}(-\Delta_{\mathbb{R}^2} + m^2))\mathrm{d}s$, $m > 0$, ${:}\cdot{:}_T$ denotes Wick ordering with respect to $\mu_T$, $c_T$ is a diverging constant corresponding to vacuum energy renormalisation, and $Z_T$ is a normalisation constant to keep $\nu_T$ a probability measure. A limit of $\{\nu_T^\rho(\phi)\}_{T,\rho}$ as $T \to \infty$ and $\rho \to \mathbb{1}_{\mathbb{R}^2}$ is called an infinite-volume sine–Gordon$(\beta^2, m, \lambda)$ measure. The construction of the interacting scale interpolation process $(\varphi_t)_{t\in\mathbb{R}_+}$, which we take as starting point, was carried out for this model in [29] in the regime $\beta^2 < 6\pi$. Our first main result is the following

**Theorem 1.1.** *Let $\beta^2 < 6\pi$, $m > 0$, $\lambda \in \mathbb{R}$, fix $y \in \mathbb{R}^2$. There exists a process $(\varphi_t)_{t\in\mathbb{R}_+}$ such that $\nu_T^\rho(\phi) \to \mathrm{Law}(\varphi_\infty)$ weakly in $H^{-\varepsilon}(\langle x\rangle^{-n})$ for every $\varepsilon > 0$, $n > 2$ as $T \to \infty$ and $\rho \to \mathbb{1}_{\mathbb{R}^2}$, possibly along a subsequence. The following local observables associated to $\varphi$ exist:*

- $[\![\varphi^k(y)]\!]$ *for any* $k \in \mathbb{N}_0$
- $[\![\partial^\mu\varphi(y)]\!]$ *for any multiindex* $\mu \in \mathbb{N}_0^2$
- $[\![\partial^\mu e^{i\alpha\varphi(y)}]\!]$ *for any* $\alpha \in \mathbb{C}$ *with real part* $|\Re e(\alpha)| < 4\pi/\beta$ *and multiindex* $\mu \in \mathbb{N}_0^2$.

*In each of these cases, there exist scale-dependent functionals $q_t$ with $\lim_{t\to\infty}\|q_t(\varphi_t) - \mathfrak{Q}_t\|_{L^1} = 0$.*

We emphasize that our method does not require smallness of the coupling constant $\lambda$ for any purposes.

The condition $|\Re e(\alpha)| < 4\pi/\beta$ on the Wick renormalised vertex field $[\![e^{i\alpha\varphi(y)}]\!]$ is conjectured to be optimal [35]. The way we achieve this condition is through a novel screening mechanism developed for external vertex insertion, i.e. sharp $L^\infty$ bounds on the propagator of the corresponding linearised Polchinski equation; see Section 5.2 and in particular Lemma 5.4 below. This can be considered the main model-specific contribution of the paper.

In our framework, the product of observables $\mathfrak{Q}^{(1)}$, $\mathfrak{Q}^{(2)}$ is again a martingale, denoted by $\mathfrak{Q}^{(1)}\star\mathfrak{Q}^{(2)}$, that satisfies the approximation property

$$\lim_{t\to\infty} \left\| (\mathfrak{Q}^{(1)}\star\mathfrak{Q}^{(2)})_t - \mathfrak{Q}_t^{(1)}\mathfrak{Q}_t^{(2)} \right\|_{L^1} = 0.$$

As a consequence, whenever the martingales are closed and therefore their terminal value exists as a random variable, we have

$$\langle \mathfrak{O}^{(1)}\star\mathfrak{O}^{(2)} \rangle = \mathbb{E}[\mathfrak{O}_\infty^{(1)}\mathfrak{O}_\infty^{(2)}],$$

a property that any sensible notion of product should satisfy. The details of the definition of the product $\star$ are given in Section 2.

Denoting the linear span of above observables considered in Theorem 1.1 by $\mathcal{O}_W^{\mathrm{loc}}$, the following result shows that we can construct any finite products of elements of $\mathcal{O}_W^{\mathrm{loc}}$, provided they are local at distinct points. Moreover, we prove that $\mathcal{O}_W^{\mathrm{loc}}$ admits a structure similar in spirit to the prefactorization algebras of Costello–Gwilliam [18],[19].

**Theorem 1.2.** *Suppose $(\beta^2, m, \lambda)$ satisfy the conditions of Theorem 1.1. For any $\mathfrak{Q}^{(1)}, \ldots, \mathfrak{Q}^{(n)} \in \mathcal{O}_W^{\mathrm{loc}}$ localized at pairwise distinct points $y_1, \ldots, y_n \in \mathbb{R}^2$, their product $\mathfrak{Q}^{(1)}\star\ldots\star\mathfrak{Q}^{(n)}$ exists. Moreover, the assignment*

$$(\mathfrak{Q}^{(1)}, \ldots, \mathfrak{Q}^{(n)}) \mapsto \mathfrak{Q}^{(1)}\star\ldots\star\mathfrak{Q}^{(n)}$$

*is associative, commutative, $n$-linear, and admits identity.*

In particular, given disjoint finite index sets $J, L, M \subseteq \mathbb{N}$, $\{k_j\}_{j\in J} \subseteq \mathbb{N}_0$, multiindices $\{\mu_l\}_{l\in L} \subseteq \mathbb{N}_0^2$, $\{\alpha_m\}_{m\in M} \subseteq \mathbb{C}$ with $|\Re e(\alpha_m)| < 4\pi/\beta$, we prove that the correlation functions

$$\left\langle \prod_{j\in J} [\![\varphi^{k_j}(y_j)]\!] \prod_{l\in L} [\![\partial^{\mu_l}\varphi(y_l)]\!] \prod_{m\in M} [\![e^{i\alpha_m\varphi(y_m)}]\!] \right\rangle \tag{1.3}$$

exist for any $\{y_i\}_{i\in J\cup L\cup M} \subseteq \mathbb{R}^2$ such that $y_{i_1} \neq y_{i_2}$ whenever $i_1 \neq i_2$. It is natural to ask about the smoothness of such correlation functions in the insertion points $y_i$. As the integral kernel of $(-\Delta + m^2)^{-1}$ is smooth outside the diagonal, we expect the same property to be inherited for the observable-inserted correlation functions. Given an $n \in \mathbb{N}$, the big diagonal $\mathrm{Diag}_n$ is defined as

$$\mathrm{Diag}_n := \{(x_1, \ldots, x_n) \in (\mathbb{R}^2)^n : \exists i \neq j \in \{1, \ldots, n\} \; s.t. \; \|x_i - x_j\| = 0\}.$$

We then have the following

**Theorem 1.3.** *Suppose $(\beta^2, m, \lambda)$ satisfy the conditions of Theorem 1.1. For any $\mathfrak{Q}^{(1)}, \ldots, \mathfrak{Q}^{(n)} \in \mathcal{O}_W^{\mathrm{loc}}$, the correlation functions*

$$(y_1, \ldots, y_n) \mapsto \langle \mathfrak{Q}^{(1)}(y_1)\star\ldots\star\mathfrak{Q}^{(n)}(y_n) \rangle$$

*are $\mathcal{C}^\infty$ on $(\mathbb{R}^2)^n \setminus \mathrm{Diag}_n$.*

We note that our proof strategy is general and relies only on the translation invariance of the Gaussian free field.

## 1.1 Related literature

Composite observables are important not only from the viewpoint of physical applications, but also from that of pure mathematics: studying the structure of observables in quantum field theory has inspired progress in functional analysis [30],[3] and various subfields of algebra and geometry [13],[33],[43]. We highlight, in particular, the work of Costello–Gwilliam [18],[19], which in turn served as an inspiration for us to set up the framework of Section 2.4.

As we have mentioned, our approach builds on the scale decomposition process $\varphi$ of the interacting field theory. This point of view can be traced back to [7]; we refer to [46], [24] for this perspective in the context of fermionic theories. As we have mentioned in the introduction, in the present work we rely on the construction of the sine–Gordon scale decomposition process obtained in [29]. See also the preprint [5] for a conjectural framework for the description of Euclidean field theories based around the process $\varphi$.

The sine–Gordon model is an important example of a scalar bosonic interacting quantum field theory. The parameter $\beta$ in (1.2) governs the singularity of the perturbation, and hence the difficulty of the renormalisation problem; in our convention, naive Wick renormalisation of the potential is sufficient when $\beta^2 < 4\pi$, and the theory is superrenormalisable/subcritical in the regime $\beta^2 < 8\pi$. It has a long history of mathematically rigorous study with various methods in various parameter ranges, see e.g. [27], [23], [42], [4], [6], [14], [36] and references therein.

We would like to highlight the continuous-scale renormalisation group approach to the model that relies on convergence of perturbation theory of the effective potential, initiated by the breakthrough work of Brydges–Kennedy [15] in the parameter range $\beta^2 < 16\pi/3$ and extended to $\beta^2 < 6\pi$ recently by Bauerschmidt–Bodineau [7]. Such perturbation theory is in fact nontrivial owing to the non-polynomial nature of the cosine,, see the series of papers [12], [38], [37], the recent work [36], and for similar problems in the context of regularity structures, [32] and [16]. Many of the technical difficulties we face in this paper stem from reconciling how this "Mayer-type" perturbation theory behaves in the presence of the external field insertions of Theorem 1.1.

Rigorous construction of composite observables of the sine–Gordon model has been of great interest, in part because of the Coleman correspondence [17], a conjectured duality between massless sine–Gordon and massive Thirring model observables. For partial verification, see [11], [22], [10] and references therein. In the Wick regime $\beta^2 < 4\pi$ and in the presence of an adiabatic/finite volume cutoff, [26] and [34] studied the insertions of products of first-order derivatives, $[\![e^{i\alpha\phi(y)}]\!]$ for $|\alpha| = \beta < \sqrt{4\pi}$, and the stress-energy tensor. For the massless model in infinite volume at the free fermion point $\beta^2 = 4\pi$, the vertex fields $[\![e^{i\alpha\phi(y)}]\!]$ for $|\alpha| < \sqrt{4\pi}$ have been constructed recently in [9] with imaginary multiplicative chaos techniques, and the corresponding one-point function has been shown to satisfy the Lukyanov–Zamolodchikov formula [35]. See Remark 6.3 for a discussion of parameter regimes of $(\beta, \alpha)$ in relation to our construction.

For general scalar bosonic EQFTs, the Osterwalder–Schrader axioms ([40], [41]) imply real analyticity of the Schwinger functions at non-coinciding points (see [28], Chapter 19). As for smoothness of observable-inserted correlation functions, the only proof we are aware of is for the Seiberg bound-constrained vertex operators of Liouville conformal field theory [39], relying on the Gaussian Multiplicative Chaos construction of the model [21]. Our proof is based on a direct flow equation argument, applies for all composite observable insertions we consider, and relies only on translation invariance of the free field. Moreover, the argument itself does not rely on smallness assumption of the coupling $\lambda$, and hence, in principle can be employed in cases where full verification of the OS axioms is not yet known.

## 1.2 Notation

We write $\mathbb{N} = \{1, 2, \dots\}$, $\mathbb{N}_0 = \mathbb{N} \cup \{0\}$, and $\mathbb{R}_+ = [0, \infty)$. Given an $\ell \in \mathbb{N}_0$, denote by

$$1{:}\ell := \{1, \dots, \ell\}, \qquad 1{:}0 := \emptyset.$$

For index sets $I_1, I_2 \subseteq 1{:}\ell$, we denote their union by $I_1 \uplus I_2$ instead of $I_1 \cup I_2$ when $I_1 \cap I_2 = \emptyset$.

Throughout the paper, unless stated otherwise, $\rho$ denotes an element of $\mathcal{C}_b^\infty(\mathbb{R}^d)$ with $\rho(x) \in [0, 1]$. We either take $\rho \equiv 1$ or $\rho(x) := \rho_R(x) := \varrho(x/L)$ for some $L > 0$ and $\varrho$ fixed compactly supported test function.

Given an $x\in\mathbb{R}^d$ and $t\geq 0$, set

$$\langle x\rangle := \sqrt{1+|x|^2}, \qquad \mathfrak{l}_t := 1+\log(1\vee t).$$

We do not keep track of logarithmic factors throughout the paper, hence the exponent $c$ of $\mathfrak{l}_t^c$ may change between occurrences.

Write $A\lesssim B$ if $A\leq CB$ for some constant $C>0$ independent of the variables and cutoffs under consideration, and $A\asymp B$ if $A\lesssim B$ and $B\lesssim A$. Dependence of $C$ on parameters is sometimes indicated by subscripts.

For a multiindex $\mu\in\mathbb{N}_0^d$, denote

$$|\mu| := \sum_{k=1}^d \mu_k, \qquad \partial^\mu := \prod_{k=1}^d \partial_{z_k}^{\mu_k}$$

For a Banach space $E$, we write $L^p(\Omega;E)$ for the corresponding Bochner-Lebesgue space.

Given a real Hilbert space $E$ and a bounded linear functional $f\in L(E;\mathbb{R})$, we denote by $f^\sharp\in E$ the unique vector induced by the Riesz–Frechet isomorphism: $\langle f^\sharp, e\rangle_E = f(e)$ for all $e\in E$. For a Frechet differentiable functional $f\colon E\to\mathbb{R}$, we denote its derivative at $\phi\in E$ by $\mathrm{D}f|_\phi\in L(E;\mathbb{R})$, and abbreviate $(\mathrm{D}f|_\phi)^\sharp$ to $\mathrm{D}f^\sharp|_\phi$ or $\mathrm{D}f|_\phi^\sharp$.

Once an arbitrary filtration $(\mathscr{F}_t)_{t\geq 0}$ is fixed, we write $\mathbb{E}_t := \mathbb{E}_t^{\mathscr{F}} := \mathbb{E}[\cdot|\mathscr{F}_t]$. Conditional expectation with respect to the Brownian filtration $(\mathscr{F}_t^B)_{t\geq 0}$ is denoted by $\mathbb{E}_t^B$.

**Acknowledgments.** M.G. is supported by the UKRI Frontier Research Grant "StochFields" (EP/Z534328/1). P.P. is supported by the EPSRC Centre for Doctoral Training in Mathematics of Random Systems: Analysis, Modelling and Simulation (EPSRC Grant EP/S023925/1).

P.P. would like to thank Christopher Beem, Tom Klose, Jaka Pelaič, Andrea Pitrone and Martin Peev for useful discussions and comments on earlier versions of the manuscript and related topics.

Special thanks go to Sarah–Jean Meyer for countless insightful conversations.

This paper has been written in TeXmacs (www.texmacs.org). All colons and semicolons of this article were generated by P.P.

# 2 General theory of stochastic observables

## 2.1 Motivation

The basis of our approach to local composite operators is that from a renormalisation group point of view, the natural notion of an observable of a field theory is a martingale with respect to the filtration generated by the multiscale description of the theory, instead of a functional of the interacting field. In this subsection we sketch the motivation behind this idea; the discussion here is semi-informal and rigorous mathematics starts in Section 2.2.

Let $(\dot{G}_u)_{u\in\mathbb{R}_+}$ be such that $\int_0^\infty \dot{G}_u \mathrm{d}u = (-\Delta+m^2)^{-1}$, and denote by $\mu_{s,t}$ the centered Gaussian measure with covariance $G_{s,t} := \int_s^t \dot{G}_u \mathrm{d}u$, $s\leq t\in[0,\infty]$. Suppose we have access to the effective potentials $(V_t)_{t\in[0,\infty]}$, a family of functionals satisfying the RG equation

$$V_s(\phi) = -\log\int_{\mathcal{S}'} e^{-V_t(\phi+\xi)}\mathrm{d}\mu_{s,t}(\xi) \tag{2.1}$$

for all $s\leq t\in[0,\infty]$ (see e.g. [8] for a review for mathematicians). We then may define the renormalised/effective measure $\nu_t$ as

$$\mathrm{d}\nu_t(\phi) := \mathcal{Z}_t^{-1} e^{-V_t(\phi)}\mathrm{d}\mu_{0,t}(\phi), \tag{2.2}$$

interpreted to be the part of the model that is visible to a detector of resolution at most $t$. Here the normalisation constant $\mathcal{Z}_t$ is chosen so that $\nu_t$ is a probability measure.

One way to think about the RG equation (2.1) is that it gives rise to the scale description process $\varphi_t$ we mentioned in the introduction, and encodes the Markov property. To see this, define for bounded Borel functions $h: \mathcal{S}'(\mathbb{R}^d) \to \mathbb{R}$ the linear operator

$$\mathbf{P}_{s,t} h(\phi) := e^{V_s(\phi)} \int_{\mathcal{S}'} e^{-V_t(\phi+\xi)} h(\phi+\xi) \mathrm{d}\mu_{s,t}(\xi).$$

Note that we can express taking expectation with respect to the effective measures in terms of the action of $\mathbf{P}_{s,t}$. Indeed, using (2.1) and (2.2),

$$\begin{aligned} \int_{\mathcal{S}'} h(\xi) \mathrm{d}\nu_t(\xi) &= \frac{1}{\int_{\mathcal{S}'} e^{-V_t(\xi)} \mathrm{d}\mu_{0,t}(\xi)} \int_{\mathcal{S}'} e^{-V_t(\xi)} h(\xi) \mathrm{d}\mu_{0,t}(\xi) \\ &= \frac{1}{e^{-V_0(0)}} \int_{\mathcal{S}'} e^{-V_t(0+\xi)} h(0+\xi) \mathrm{d}\mu_{0,t}(0+\xi) \\ &= \mathbf{P}_{0,t} h(0). \end{aligned}$$

We also see that $\mathbf{P}_{s,t} h \geq 0$ whenever $h \geq 0$, $\mathbf{P}_{s,t} \mathbb{1} = \mathbb{1}$ and, as $h$ is bounded, $V_\infty$ is assumed to be bounded from below and $\mu_{s,t}$ is a probability measure, by dominated convergence $\mathbf{P}_{s,t}$ preserves pointwise convergence. Moreover, as $\lim_{s \to t} \mu_{s,t} = \delta_{\xi=0}$ weakly we obtain $\mathbf{P}_{t,t} = \mathrm{id}$. Finally, the Chapman–Kolmogorov identity is a consequence of the convolution semigroup property of Gaussian measures: for all $s \leq u \leq t$ we have $\mu_{s,t} = \mu_{s,u} * \mu_{u,t}$, thus

$$\begin{aligned} \mathbf{P}_{s,u} \mathbf{P}_{u,t} h(\phi) &= e^{V_s(\phi)} \int_{\mathcal{S}'} e^{-V_u(\phi+\zeta)} \mathbf{P}_{u,t} h(\phi+\zeta) \mathrm{d}\mu_{s,u}(\zeta) \\ &= e^{V_s(\phi)} \int_{\mathcal{S}'} e^{-V_u(\phi+\zeta)} e^{V_u(\phi+\zeta)} \int_{\mathcal{S}'} e^{-V_t(\phi+\zeta+\xi)} h(\phi+\zeta+\xi) \mathrm{d}\mu_{u,t}(\xi) \mathrm{d}\mu_{s,u}(\zeta) \\ &= e^{V_s(\phi)} \int_{\mathcal{S}'} e^{-V_t(\phi+\eta)} h(\phi+\eta) \mathrm{d}\mu_{s,t}(\eta) \\ &= \mathbf{P}_{s,t} h(\phi). \end{aligned}$$

Thus $\mathbf{P}_{s,t}(\phi, A) := \mathbf{P}_{s,t} \mathbb{1}_A(\phi)$, $\phi \in \mathcal{S}'$, $A \in \mathscr{B}(\mathcal{S}')$ defines a Markovian transition kernel describing the evolution of the measures $\nu_t$:

$$\begin{aligned} \mathbf{P}_{s,t}^* \nu_s(A) &:= \int_{\mathcal{S}'} \mathbf{P}_{s,t} \mathbb{1}_A(\phi) \mathrm{d}\nu_s(\phi) \\ &= \int_{\mathcal{S}'} e^{V_s(\phi)} \int_{\mathcal{S}'} e^{-V_t(\phi+\xi)} \mathbb{1}_A(\phi+\xi) \mathrm{d}\mu_{s,t}(\xi) \mathrm{d}\nu_s(\phi) \\ &= \frac{1}{e^{-V_0(0)}} \int_{\mathcal{S}'} e^{-V_s(\phi)} e^{V_s(\phi)} \int_{\mathcal{S}'} e^{-V_t(\phi+\xi)} \mathbb{1}_A(\phi+\xi) \mathrm{d}\mu_{s,t}(\xi) \mathrm{d}\mu_{0,s}(\phi) \\ &= \frac{1}{e^{-V_0(0)}} \int_{\mathcal{S}'} e^{-V_t(\psi)} \mathbb{1}_A(\psi) \mathrm{d}\mu_{0,t}(\psi) \\ &= \nu_t(A). \end{aligned} \tag{2.3}$$

Hence, by Kolmogorov's consistency theorem[2.1] there exists a probability space $(\Omega, \mathscr{A}, \mathbb{P})$ with associated Markov process $(\varphi_t)_{t \in \mathbb{R}_+}$ in its canonical filtration $\mathscr{F}_t := \sigma(\varphi_s: s \leq t)$ such that $\mathrm{Law}(\varphi_t) = \nu_t$.

In the usual probabilistic approach to EQFTs ([28], [44]), observables are defined as real-valued functionals $f_\infty$ of $\varphi_\infty$, or equivalently by a Doob–Dynkin argument as $\sigma(\varphi_\infty)$-measurable random variables $F_\infty: \Omega \to \mathbb{R}$. Given such an observable $F_\infty$, we may define $F_t$ as conditional expectation with respect to $\mathscr{F}_t$. Explicitly, setting $F_t := f_t(\varphi_t)$ with

$$f_t(\phi) := \mathbf{P}_{t,\infty} f_\infty(\phi) = e^{V_t(\phi)} \int_{\mathcal{S}'} e^{-V_\infty(\phi+\xi)} f_\infty(\phi+\xi) \mathrm{d}\mu_{t,\infty}(\xi),$$

2.1. All the EQFTs we are interested in are expected to be supported on some weighted Sobolev space of sufficiently low regularity, hence we may replace $\mathcal{S}'$ with these spaces so that the assumption of Polishness on the state space is satisfied.

by (2.3) we see that for all $s \le t$,

$$\mathbb{E}[f_t(\varphi_t)|\mathscr{F}_s] = \mathbb{E}[\mathbf{P}_{t,\infty} f_\infty(\varphi_t)|\varphi_s] = \mathbf{P}_{s,t}\mathbf{P}_{t,\infty} f_\infty(\varphi_s) = \mathbf{P}_{s,\infty} f_\infty(\varphi_s) = f_s(\varphi_s).$$

In particular, expectation is preserved along the flow:

$$\int_{\mathcal{S}'} f_t(\phi)\mathrm{d}\nu_t(\phi) = \int_{\mathcal{S}'} \mathbf{P}_{t,\infty} f_\infty(\phi)\mathrm{d}\nu_t(\phi) = \int_{\mathcal{S}'} f_\infty(\phi)\mathrm{d}\mathbf{P}^*_{t,\infty}\nu_t(\phi) = \int_{\mathcal{S}'} f_\infty(\phi)\mathrm{d}\nu_\infty(\phi)$$

That is, we constructed a martingale $F = (F_t)_{t\in\mathbb{R}_+}$ with $L^1(\Omega)$-limit $F_\infty$ such that $F_t$ is the effective description of the original observable $F_\infty$: it contains details only up to resolution $t$ and averaging it with respect to the effective measure $\nu_t$ yields the same result as averaging the original observable $F_\infty$ with respect to the UV-complete EQFT $\nu_\infty$. Thus observables, in the usual sense, may be identified through their multiscale effective description with a subset of martingales with respect to the scale decomposition filtration $\mathscr{F}$.

It is a well-established fact that the random variable language is insufficient for a proper description of composite observables: it has been shown in [2] that the Wick square of the Gaussian free field, defined through analytic continuation of the Wick square of the relativistic free field, cannot exist as an $\mathcal{S}'(\mathbb{R}^d)$-valued random variable for $d \ge 4$. That is, $\langle :\varphi^2:, g\rangle_{\mathcal{S}',\mathcal{S}}$, $g \in \mathcal{S}(\mathbb{R}^d)$ is not a real-valued random variable, i.e. the problem is not only "spatial" but also "probabilistic" regularity. One expects a similar result for high enough Wick powers also in $d = 3$. The approach taken in [2] is to regard such objects as $\mathcal{S}'(\mathbb{R}^d)$-valued bona fide Hida distributions; see [1] for further discussion on this.

To circumvent the above issue, we thus define observables to be all $\mathscr{F}$-martingales. This is a proper generalization to the random variable language: those martingales that are not closed (equivalently, by Doob's martingale convergence theorem, not uniformly integrable) do not correspond to an observable $F_\infty$ in the traditional picture. Intuitively, this generalization has, in some sense, similar flavour to that of classical versus generalized functions: the notion of observable is "completed" by getting rid of uncomfortable convergence issues in the $t \to \infty$ limit; yet the rigidity in the structure of a martingale preserves enough information so that important computations can still be carried out. Moreover, as we will see, the multiscale point of view inherent in the martingale description is actually helpful in proving properties of observables.

## 2.2 Definitions and setting

Denote by $\mathcal{H}$ a Hilbert space densely and continuously embedded into $\mathcal{S}'(\mathbb{R}^d)$, where the latter is endowed with the strong topology. For practical purposes, it is convenient to take $\mathcal{H} := H^{k,\ell}(\mathbb{R}^d) = B^k_{2,2}(\langle x\rangle^{-\ell})$ with $k, \ell > 0$ a polynomially weighted Sobolev space of positive regularity. The parameter $k$ will be chosen based on the particular observables we are interested in, e.g. we will usually assume that $\mathcal{H} \hookrightarrow \mathcal{C}$ so that point evaluations $\delta_x$ are continuous functionals on $\mathcal{H}$; this can be achieved with $k > d/2$ by Morrey's inequality.

Let us fix a probability space $(\Omega, \mathscr{A}, \mathbb{P})$. Consider a filtration $\mathscr{F} = (\mathscr{F}_t)_{t\in\mathbb{R}_+}$ on $(\Omega, \mathscr{A}, \mathbb{P})$ and a stochastic process $\varphi = (\varphi_t)_{t\in\mathbb{R}_+}$ adapted to $\mathscr{F}$ such that $\varphi_t\colon \Omega \to \mathcal{H}$ for every finite $t$. We assume that $\varphi_t$ converges at least $\mathbb{P}$-a.s. in the weak* topology of $\mathcal{S}'$ to a random variable $\varphi_\infty\colon \Omega \to \mathcal{S}'$. Denote by $\mathbb{E}_t[\cdot] := \mathbb{E}[\cdot|\mathscr{F}_t]$.

The interpretation of this setting is as follows. $\mathrm{Law}(\varphi_\infty)$ is the probability measure we are interested in (usually an EQFT, but note that the framework is more general), the process $\varphi$ is the scale interpolation process of $\varphi_\infty$ associated to the coarse-graining procedure encoded in $\mathscr{F}$; that is, $\varphi_t$ for a $t < \infty$ is the effective field theory accessible at resolution smaller than $t$. Finally, $\mathcal{H}$ is simply the state space of coarse-grained fields.

In the traditional probabilistic setting, observables are real-valued random variables measurable with respect to $\sigma(\varphi_\infty)$ and integrable under $\mathrm{Law}(\varphi_\infty)$; equivalently, by the Doob–Dynkin representation, we can think of them as functionals of the random field $\varphi_\infty$. We call such objects *probabilistic* observables. In contrast, we have the following

**Definition 2.1.** A real-valued[2.2] $\mathscr{F}$-martingale is called a *stochastic observable*.

The set of stochastic observables is denoted by $\mathscr{O} = \mathscr{O}(\mathbb{R}^d, \mathscr{F}, \mathbb{P})$; it naturally admits a vector space structure under $t$-wise addition and multiplication by scalars. The space of probabilistic observables is denoted by $\mathscr{O}_p$; associating to each probabilistic observable $X$ its Doob martingale $\mathfrak{X}$, i.e.

$$\mathfrak{X}_t := \mathbb{E}[X | \mathscr{F}_t],$$

it is clear that $\mathscr{O}_p \subseteq \mathscr{O}$. When no confusion may arise, we shall call stochastic observables simply observables.

Just as probabilistic observables can be thought of functionals of $\varphi_\infty$, we may ask if an observable $\mathfrak{Q}$ is such that $\mathfrak{Q}_t$ is also measurable with respect to $\varphi_t$, so that by Doob–Dynkin we can think of it as a functional of the coarse-grained field $\varphi_t$ (consider, e.g., the Doob martingales associated to probabilistic observables when $\varphi$ is a time-homogeneous diffusion and $\mathscr{F}$ is its natural filtration). This leads us to try to parametrize the space of observables with families of scale-dependent functionals of $\mathcal{H}$.

**Definition 2.2.** A function $q: \mathbb{R}_+ \times \mathcal{H} \to \mathbb{R}$ is said to be a *germ* for an observable $\mathfrak{Q}$ if

$$\lim_{t \to \infty} \| \mathfrak{Q}_t - q_t(\varphi_t) \|_{L^1(\Omega)} = 0. \tag{2.4}$$

Based on this definition, there is a natural notion of equivalence relation on the space of scale-dependent functionals: given $q, \tilde{q}: \mathbb{R}_+ \times \mathcal{H} \to \mathbb{R}$, we write $q \sim_\varphi \tilde{q}$ if

$$\lim_{t \to \infty} \| \tilde{q}_t(\varphi_t) - q_t(\varphi_t) \|_{L^1(\Omega)} = 0.$$

Denoting by $[q]_\varphi$ the corresponding equivalence classes, i.e. $[q]_\varphi = q + [0]_\varphi$, where

$$[0]_\varphi := \left\{ r: \mathbb{R}_+ \times \mathcal{H} \to \mathbb{R}: \lim_{t \to \infty} \| r_t(\varphi_t) \|_{L^1(\Omega)} = 0 \right\},$$

the elements of the quotient space are represented by observables: if $q$ is a germ of an observable $\mathfrak{Q}$, then all $\tilde{q} \in [q]_\varphi$ are also germs for $\mathfrak{Q}$, since

$$\| \mathfrak{Q}_t - \tilde{q}_t(\varphi_t) \|_{L^1(\Omega)} \leq \| \mathfrak{Q}_t - q(\varphi_t) \|_{L^1(\Omega)} + \| \tilde{q}(\varphi_t) - q(\varphi_t) \|_{L^1(\Omega)} \to 0; \tag{2.5}$$

conversely, if $q$ is a germ for both $\mathfrak{Q}, \mathfrak{P} \in \mathscr{O}$, then for all $t \in \mathbb{R}_+$

$$\begin{aligned} \| \mathfrak{Q}_t - \mathfrak{P}_t \|_{L^1(\Omega)} &= \| \mathbb{E}_t[\mathfrak{Q}_T - \mathfrak{P}_T] \|_{L^1(\Omega)} \\ &\leq \| \mathfrak{Q}_T - \mathfrak{P}_T \|_{L^1(\Omega)} \\ &\leq \| \mathfrak{Q}_T - q_T(\varphi_T) \|_{L^1(\Omega)} + \| \mathfrak{P}_T - q_T(\varphi_T) \|_{L^1(\Omega)} \end{aligned} \tag{2.6}$$

vanishing as $T \to \infty$, i.e. $\mathfrak{Q} = \mathfrak{P}$.

A collection of observables $\mathcal{Q} := \{ \mathfrak{Q}_t^{(x)} \}_{x \in \mathbb{R}^d}$ is called an observable field. Note that while in our definitions we always took $\mathbb{R}^d$, it is not hard to generalize the discussion to the case where spacetime is a Riemannian manifold $M$.

As the notations $\mathscr{O}(\mathbb{R}^d, \mathscr{F}, \mathbb{P})$ and $\sim_\varphi$, $[q]_\varphi$ suggest, the space of observables and the equivalence classes of germs depend on the filtration $\mathscr{F}$ and the scale decomposition process $\varphi$. This might seem strange at first, but it is in fact natural: the physically relevant information of an observable is encoded in the $t = 0$ scale.

2.2. Strictly speaking then, vertex fields correspond to a pair of observables.

**Definition 2.3.** A *state* is the linear functional $\langle\cdot\rangle:\mathscr{O}\to\mathbb{R}$ given by $\langle\mathfrak{O}\rangle:=\mathbb{E}[\mathfrak{O}_0]$.

One might think of the situation as follows: we are given a fixed set of abstract symbols $\mathcal{O}\mathrm{bs}$, containing, for example $\{\phi,\phi^n,(\partial^\mu\phi)^n,\exp(\phi),\ldots\}$; the choice of $(\Omega,\mathscr{A},\mathbb{P},\mathscr{F},\varphi)$ then induces a partial function $\pi^{\mathscr{F},\varphi}:\mathcal{O}\mathrm{bs}\to\mathscr{O}$ with a corresponding germ structure. What we are really interested in is the set $\pi^{\mathscr{F},\varphi}(\mathcal{O}\mathrm{bs})$ of physically relevant observables. Given a different choice $(\tilde{\mathscr{F}},\tilde{\varphi})$, as long as both $\pi^{\mathscr{F},\varphi}$ and $\pi^{\tilde{\mathscr{F}},\tilde{\varphi}}$ are injective and their domain is the same, there exists a bijection between the physically relevant observables under the different representation, i.e. $\pi^{\mathscr{F},\varphi}(\mathcal{O}\mathrm{bs})\simeq\pi^{\tilde{\mathscr{F}},\tilde{\varphi}}(\mathcal{O}\mathrm{bs})$.

From now on, we will work in the following representation. Consider the heat kernel decomposition of the covariance operator $(-\Delta+m^2)^{-1}$ of the massive Gaussian free field given by

$$Q_t:=\sqrt{\frac{1}{t^2}e^{-\left(\frac{-\Delta+m^2}{t}\right)}},$$

defined via functional calculus or Fourier multipliers. Moreover, denote by $\dot{G}_t:=Q_t^2$. By a standard computation, the $L^2(\mathbb{R}^d)$ integral kernel $k_t^Q$ of $Q_t$ and $\dot{\mathcal{G}}_t(x-y)=k_t^{Q^2}$ of $\dot{G}_t=Q_t^2$ is given by the formulae

$$k_t^Q(x,y)=\frac{t^{d/2-1}}{(2\pi)^{d/2}}e^{-\frac{t}{2}|x-y|^2}e^{-\frac{m^2}{2t}},\qquad\text{and}\qquad\dot{\mathcal{G}}_t(x-y)=k_t^{Q^2}(x,y)=\frac{t^{d/2-2}}{(4\pi)^{d/2}}e^{-\frac{t}{4}|x-y|^2}e^{-\frac{m^2}{t}}\tag{2.7}$$

We can regard $Q_t\in L(L^2(\mathbb{R}^d);\mathcal{H})$ (see e.g. Lemma A.7), denote by $Q_t^*:\mathcal{H}\to L^2(\mathbb{R}^d)$ its adjoint, and define the scale interpolation of the free field, an $\mathcal{H}$-valued Gaussian process

$$B_t:=\int_0^t Q_s\mathrm{d}W_s$$

where $W$ is a cylindrical Brownian motion on $L^2(\mathbb{R}^d)$. The fact that the stochastic integral is well-defined and converges to the Gaussian free field $B_\infty$ is well-known, but to be self-contained we prove it in Appendix A along with other basic heat kernel estimates.

For reasons that will be clear soon, we require the following property from the filtration $\mathscr{F}$.

**Definition 2.4.** A filtration $\mathscr{F}$ is called *admissible* if $\varphi$ is $\mathscr{F}$-adapted and $B$ is an $\mathscr{F}$-martingale.

In the latter part of this paper, we take the scale decomposition process $\varphi$ to be that of the sine–Gordon model; $\varphi$ and the admissible filtration $\mathscr{F}$ is constructed in Theorem 4.4 below.

In general, the basic structure one requires from $\varphi$ is that it can be represented as a scale-dependent perturbation of the free field process $B$, i.e. it is an Ito diffusion of the form $\varphi_t=Z_t+B_t$ where the drift process $Z$ is

$$Z_t=-\int_0^t\dot{G}_s\mathrm{D}V_s(\varphi_s)\mathrm{d}s\tag{2.8}$$

where $V:\mathbb{R}_+\times\mathcal{H}\to\mathbb{R}$ is a sufficiently nice scale-dependent functional on the field configuration space $\mathcal{H}$. The parabolic generator of $\varphi$, denoted by $\mathscr{L}^V$, is the operator

$$\mathscr{L}^Vf_t(\phi):=\partial_tf_t(\phi)+\langle\mathrm{D}f_t|_\phi,-\dot{G}_t\mathrm{D}V_t|_\phi\rangle+\frac{1}{2}\mathrm{Tr}[Q_tQ_t^*\mathrm{D}^2f_t|_\phi]\tag{2.9}$$

defined for sufficiently nice $f:\mathbb{R}_+\times\mathcal{H}\to\mathbb{R}$. We will frequently abuse notation to write

$$\mathrm{Tr}[Q_tQ_t^*\mathrm{D}^2f_t|_\phi]=\mathrm{Tr}[\dot{G}_t\mathrm{D}^2f_t|_\phi]=\Delta_{\dot{G}_t}f_t(\phi),$$

see Lemma B.1 and Notation B.2 for justification.

## 2.3 Coherence

Germs serve to parametrize the space of observables. In some cases, the correspondence between observables and germs is bijective: the germ itself completely determines the martingale. The canonical example is when $q$ solves the backward Kolmogorov PDE: indeed, in this case $\mathfrak{Q}_t := q_t(\varphi_t)$ is a local martingale. It is convenient to ask for germs that solve the backward equation only *approximately* in some suitable sense so that the corresponding martingale can still be obtained from the knowledge of $q$ alone. Such germs are called coherent and this section is concerned with explaining how to "reconstruct" the unique observable from a coherent germ and establish precise conditions that guarantee this can be done. The reconstruction is done via a simple application of Ito's formula, we thus start with recalling some notations from infinite-dimensional stochastic calculus; see e.g. [20] for more details.

Recall that given processes $Q, f$ satisfying sufficiently strong measurability conditions with $Q_t \in \mathcal{L}_2(L^2(\mathbb{R}^d); \mathcal{H})$ and $f_t \in \mathcal{H}$ such that $\mathbb{P}\big(\int_0^T \|Q_t^* f_t\|_{L^2(\mathbb{R}^d)}^2 \mathrm{d}t < \infty\big) = 1$, and setting, for an $L^2(\mathbb{R}^d)$-cylindrical Brownian motion $W$,

$$\int_0^T \langle f_t, Q_t \mathrm{d}W_t \rangle := \int_0^T \tilde{Q}_t^f \mathrm{d}W_t \qquad \text{with } \tilde{Q}_t^f : L^2(\mathbb{R}^d) \to \mathbb{R} \quad \tilde{Q}_t^f(u) := \langle f_t, Q_t u \rangle_{\mathcal{H}},$$

then $\|\tilde{Q}_t^f\|_{\mathcal{L}_2(L^2(\mathbb{R}^d);\mathbb{R})}^2 = \|Q_t^* f_t\|_{L^2(\mathbb{R}^d)}^2$ and hence the stochastic integral $\int_0^T \langle f_t, Q_t \mathrm{d}W_t \rangle$ is a well-defined continuous $\mathbb{R}$-valued stochastic process.

Also recall that given an $\mathcal{L}_2(L^2(\mathbb{R}^d); \mathcal{H})$-valued predictable integrand process $\Phi$, in order for a stochastic integral $\int_0^T \Phi_t \mathrm{d}W_t$ to be a $\mathcal{H}$-martingale, it suffices that $\Phi_t \in \mathbb{H}_W^2([0,T]; \mathcal{L}_2(L^2(\mathbb{R}^d); \mathcal{H}))$, i.e. that

$$\mathbb{E}\int_0^T \|\Phi_t\|_{\mathcal{L}_2(L^2(\mathbb{R}^d);\mathcal{H})}^2 \mathrm{d}t < \infty.$$

Thus for the stochastic integral $\int_0^T \langle f_t, Q_t \mathrm{d}W_t \rangle$ to be a real-valued martingale, it suffices that

$$\mathbb{E}\int_0^T \|\tilde{Q}_t^f\|_{\mathcal{L}_2(L^2(\mathbb{R}^d);\mathbb{R})}^2 \mathrm{d}t = \mathbb{E}\int_0^T \|Q_t^* f_t\|_{L^2(\mathbb{R}^d)}^2 \mathrm{d}t < \infty.$$

Given a $T \in \mathbb{R}_+$, we denote by $q^T$ the restriction of $q: \mathbb{R}_+ \times \mathcal{H} \to \mathbb{R}$ to $[0,T] \times \mathcal{H} \to \mathbb{R}$.

**Definition 2.5.** A functional $q: \mathbb{R}_+ \times \mathcal{H} \to \mathbb{R}$ is called *coherent* if

1. For each $T \in \mathbb{R}_+$, the functional $q^T$ is continuous on $[0,T] \times \mathcal{H}$, and the partial derivatives $\partial_t$, D, $\mathrm{D}^2$ are bounded and continuous on bounded subsets of $[0,T] \times \mathcal{H}$.
2. For each $T \in \mathbb{R}_+$, we have $\mathbb{E}\int_0^T \|Q_t^* \mathrm{D} q_t^T|_{\varphi_t}^\sharp\|_{L^2(\mathbb{R}^d)}^2 \mathrm{d}t < \infty$, where $Q$ is the operator in the heat kernel interpolation.
3. It holds that

$$\int_0^\infty \|\mathscr{L}^V q_s(\varphi_s)\|_{L^1(\Omega)} \mathrm{d}s < \infty. \tag{2.10}$$

**Lemma 2.6. (reconstruction)** *Suppose $q$ is a coherent functional and $\varphi$ is an Ito diffusion with generator $\mathscr{L}^V$. Then*

1. *For each $T \in \mathbb{R}_+$ there exists a unique observable $\mathfrak{Q}^{\{T\}}$ defined as $\mathfrak{Q}_t^{\{T\}} := q_t^T(\varphi_t) + r_{t,T}^q$ for $t \in [0,T]$ and extended to $\mathbb{R}_+$ by $\mathfrak{Q}_T^{\{T\}}$, where*

$$r_{t,T}^q := \mathbb{E}_t \int_t^T \mathscr{L}^V q_s(\varphi_s) \mathrm{d}s. \tag{2.11}$$

2. *There exists a unique observable $\mathfrak{Q}$ with germ $q$ defined as $\mathfrak{Q}_t := q_t(\varphi_t) + r_{t,\infty}^q$ such that*

$$\lim_{T \to \infty} \|\mathfrak{Q}_t - \mathfrak{Q}_t^{\{T\}}\|_{L^1(\Omega)} = 0$$

*for each* $t \in \mathbb{R}_+$.

**Proof.** As we have already indicated, this is a simple application of Ito's formula.

1. Indeed, fixing $[s,t] \subseteq [0,T]$, by condition 1. we can apply Ito's formula to obtain

$$q_t^T(\varphi_t) = q_s^T(\varphi_s) + \int_s^t \mathscr{L}^V q_u^T(\varphi_u)\mathrm{d}u + \int_s^t \langle \mathrm{D}q_u^T|_{\varphi_u}^\sharp, Q_u \mathrm{d}W_u\rangle.$$

By condition 2. and admissibility, the last term is a martingale, thus we have that

$$\mathbb{E}_s \mathfrak{Q}_t^{\{T\}} = \mathbb{E}_s\left[q_s^T(\varphi_s) + \int_s^t \mathscr{L}^V q_u^T(\varphi_u)\mathrm{d}u + \int_s^t \langle \mathrm{D}q_u^T|_{\varphi_u}^\sharp, Q_u \mathrm{d}W_u\rangle\right] + \mathbb{E}_s\mathbb{E}_t \int_t^T \mathscr{L}^V q_u^T(\varphi_u)\mathrm{d}u =$$

$$q_s^T(\varphi_s) + \mathbb{E}_s \int_s^t \mathscr{L}^V q_u^T(\varphi_u)\mathrm{d}u + \mathbb{E}_s\mathbb{E}_t \int_t^T \mathscr{L}^V q_u^T(\varphi_u)\mathrm{d}u = \mathfrak{Q}_s^{\{T\}}$$

so $\mathfrak{Q}^{\{T\}}$ is a martingale. The condition (2.4) is satisfied trivially.

2. For each $t \in \mathbb{R}_+$, for all $T \geq t$

$$\mathfrak{Q}_t - \mathfrak{Q}_t^{\{T\}} = \mathbb{E}_t \int_T^\infty \mathscr{L}^V q_s(\varphi_s)\mathrm{d}s$$

so that (2.10) yields by conditional Jensen inequality, Fubini and dominated convergence that

$$\|\mathfrak{Q}_t - \mathfrak{Q}_t^{\{T\}}\|_{L^1(\Omega)} \leq \left\|\int_T^\infty \mathscr{L}^V q_s(\varphi_s)\mathrm{d}s\right\|_{L^1(\Omega)} \leq \int_{\mathbb{R}_+} \mathbb{1}_{[T,\infty)} \|\mathscr{L}^V q_s(\varphi_s)\|_{L^1(\Omega)}\mathrm{d}s \to 0$$

as $T \to \infty$. Since the $\mathfrak{Q}^{\{T\}}$ are martingales, this implies that $\mathfrak{Q}$ is also a martingale. Similarly, (2.4) follows from (2.10). □

Let us fix some terminology for objects we see in the construction for future reference.

**Definition 2.7.** The observable $\mathfrak{Q}$ of the previous lemma is called the *reconstruction* of $q$.

It is clear from the linearity of $\mathscr{L}^V$ and the triangle inequality that coherent functionals form a vector space. Moreover, the corresponding linear structure induced on a subset of observables via the Ito reconstruction "lift" is compatible with the linear structure of $\mathscr{O}$ (i.e. $t$-wise addition and scalar multiplication). More precisely, we make the following

**Definition 2.8.** Denote by $\hat{\mathscr{O}}$ the vector space of coherent functionals. Define the *reconstruction map* $\mathfrak{R}$ as

$$\mathfrak{R}: \hat{\mathscr{O}} \to \mathscr{O}, \qquad (\mathfrak{R}(q))_t := q_t(\varphi_t) + r_{t,\infty}^q.$$

The reconstruction map is linear: given $q^{(a)}, q^{(b)} \in \hat{\mathscr{O}}$, linearity of the generator $\mathscr{L}^V$ yields

$$\begin{aligned}(\mathfrak{R}(q^{(a)} + q^{(b)}))_t &= (q^{(a)} + q^{(b)})_t(\varphi_t) + r_{t,\infty}^{q^{(a)}+q^{(b)}} \\ &= q_t^{(a)}(\varphi_t) + r_{t,\infty}^{q^{(a)}} + q_t^{(b)}(\varphi_t) + r_{t,\infty}^{q^{(b)}} \\ &= (\mathfrak{R}(q^{(a)}))_t + (\mathfrak{R}(q^{(b)}))_t \\ &= (\mathfrak{R}(q^{(a)}) + \mathfrak{R}(q^{(b)}))_t\end{aligned}$$

so that $\mathfrak{R}(q^{(a)} + q^{(b)}) = \mathfrak{R}(q^{(a)}) + \mathfrak{R}(q^{(b)})$; homogeneity is also trivial.

**Definition 2.9.** Let $\mathscr{O}_c := \mathrm{Ran}(\mathfrak{R}) \subseteq \mathscr{O}$ be the linear subspace of $\mathscr{O}$ consisting of those observables that are Ito reconstructions of coherent germs.

The above reconstruction reduces the existence of composite observables to the construction of coherent functionals $q$, which in turn depends largely on estimates of the scale decomposition process $\varphi$ and its generator $\mathscr{L}^V$. Once one has a sufficient control over these objects, the construction of coherent germs can be achieved in a systematic manner based on ideas of renormalisation group flow equations, see Section 2.6 below.

This shift in perspective of putting coherent functionals in the center of attention validates the following slight abuse of terminology: from now on, we will often call a scale-dependent functional $q: \mathbb{R}_+ \times \mathcal{H} \to \mathbb{R}$ a germ (without specifying the corresponding martingale $\mathfrak{Q}$ for which (2.4) holds); the tacit assumption is that such a functional will turn out to be coherent, and thus indeed a germ (in the strict sense of Definition 2.2) of its reconstructed observable.

**Remark 2.10. (Connection to RG flow equations)** Heuristically, it is easy to see how our framework connects to Wilsonian ideas. Suppose we are searching for the family of pairs $(V_t, \tilde{q}_t)_{t\in\mathbb{R}_+}$ such that the perturbed interaction

$$V_t^\theta := V_t + \theta \tilde{q}_t \tag{2.12}$$

satisfies the RG equation (2.1). One obtains from the resulting flow of perturbed renormalised densities $(e^{-V_t^\theta})_{t\in\mathbb{R}_+}$ the evolution of the effective description of the observable: indeed,

$$-\partial_\theta e^{-V_t^\theta}|_{\theta=0} = \tilde{q}_t e^{-V_t},$$

which yields the expectation of $q_t$ under the effective measure, cf. (2.2). By a standard argument, we get from (2.1) the observable-perturbed Polchinski equation

$$\partial_t V_t^\theta + \frac{1}{2}\Delta_{\dot{G}_t} V_t^\theta = \frac{1}{2}\langle \mathrm{D} V_t^\theta, \dot{G}_t \mathrm{D} V_t^\theta \rangle.$$

Taking $\theta = 0$ we obtain the usual Polchinski equation for the potential $V_t$, while taking $\partial_\theta$ at $\theta = 0$ gives the observable flow equation

$$\partial_t \tilde{q}_t + \frac{1}{2}\Delta_{\dot{G}_t} \tilde{q}_t = \langle \mathrm{D}\tilde{q}_t, \dot{G}_t \mathrm{D} V_t \rangle,$$

which is precisely the backward Kolmogorov PDE $\mathscr{L}^V \tilde{q} = 0$.

In light of (2.12), the functional $\tilde{q}_t$ has the geometric interpretation of a tangent vector at $V_t$ to the (formal) manifold of effective interactions $\mathcal{M}$. The Polchinski equation gives rise to a vector field $\beta$ on $\mathcal{M}$:

$$\beta(V_t) := -\frac{1}{2}\Delta_{\dot{G}_t} V_t + \frac{1}{2}\langle \mathrm{D} V_t, \dot{G}_t \mathrm{D} V_t \rangle.$$

The family $(V_t)_{t\in\mathbb{R}_+}$ is then an integral curve of this vector field, and the evolution of $\tilde{q}_t \in T_{V_t}\mathcal{M}$ describes how a tangent vector is transported along such a curve by $\beta$. Our martingales $\mathfrak{Q}$ contain the information of such a flow in an abstract way, without having to give meaning to $\mathcal{M}$: indeed, for each $t$, we can employ a Doob–Dynkin representation to obtain a functional $\mathfrak{Q}_t = \tilde{q}_t(\varphi_t)$. Allowing for coherent functionals $(q_t)_{t\in\mathbb{R}_+}$ in this picture has the interpretation of a family of approximate tangent vectors that are sufficiently close to the actual $\tilde{q}$ such that the defect can be recovered via the correction procedure of Lemma 2.6.

## 2.4 Observable product

Observables of a quantum field theory do not form an algebra: we do not expect to be able to multiply any two observables. From a mathematical point of view, in the Euclidean setting this is ultimately a consequence of the singularity of the covariance of the random field $\varphi_\infty$ on the diagonal in $d \geq 2$. From a physics perspective, it is related to the fact that the effect of a measurement (i.e. coupling of a measurement device with the physical system) on a quantum system is considerable, hence taking two measurements of a quantum particle simultaneously yields incoherent results. Thus in most axiomatic frameworks of quantum field theory (see e.g. Haag–Kastler–Araki axioms [30],[3] factorization algebras [18],[19]), there is a notion of locality: "combining" measurements should be possible only for observables that are supported on disjoint spacetime regions.

It is well-known that the space of martingales is not stable under taking pointwise (with respect to $t$) products. To have a well-defined notion of observable product in the stochastic observable language, we leverage the fact that germs parametrize the space of observables: given two observables, we can multiply their respective germs. The idea is then to use the reconstruction map $\mathfrak{R}$ to push forward the product to the space of observables; we note the analogy with theory of regularity structures [31].

When the product germ happens to be coherent, by the Ito reconstruction lemma it gives rise to a unique observable, which we take to be our definition of the product observable provided it approximates the pointwise products of observables in $L^1(\Omega)$ in the UV. In our setting, the failure of the algebra structure is then encoded in the nontrivial condition that the product germ is coherent.

**Definition 2.11.** Let $\{\mathfrak{O}^{(j)}\}_{j=1}^n$ be observables with coherent germs $\{q^{(j)}\}_{j=1}^n$, respectively, and write $q^I := \prod_{i\in I} q^{(i)}$. If for every nonempty $I \subseteq 1{:}\,n$, we have that

1. $q^I$ is coherent,
2. for any partition $P$ of $I$,

$$\lim_{t\to\infty} \left\| \mathfrak{R}(q^I)_t - \prod_{J\in P} \mathfrak{R}(q^J)_t \right\|_{L^1(\Omega)} = 0, \tag{2.13}$$

then we say that the tuple $(q^{(1)}, \ldots, q^{(n)})$ is *compatible* and call the martingale $\mathfrak{O}^{(1)} \star \ldots \star \mathfrak{O}^{(n)} := \mathfrak{R}(q^{1:n})$ the *observable product*.

**Remark 2.12.** The UV $L^1(\Omega)$ approximation property (2.13) is what intuition expects the product should behave like: we would like to have, for all $t\in\mathbb{R}_+$, that

$$(\mathfrak{O}^{(1)} \star \mathfrak{O}^{(2)})_t = \lim_{T\to\infty} \mathbb{E}_t\big[\mathfrak{O}_T^{(1)} \mathfrak{O}_T^{(2)}\big] \tag{2.14}$$

in $L^1(\Omega)$, and when the martingales $\mathfrak{O}^{(1)}, \mathfrak{O}^{(2)}, \mathfrak{O}^{(1)} \star \mathfrak{O}^{(2)}$ are closed, to have

$$\langle \mathfrak{O}^{(1)} \star \mathfrak{O}^{(2)} \rangle = \mathbb{E}\big[\mathfrak{O}_\infty^{(1)} \mathfrak{O}_\infty^{(2)}\big]. \tag{2.15}$$

The statements (2.14) and (2.15) follow from (2.13) trivially: indeed,

$$\begin{aligned} \big\|\mathbb{E}_t\big[\mathfrak{O}_T^{(1)}\mathfrak{O}_T^{(2)}\big] - (\mathfrak{O}^{(1)}\star\mathfrak{O}^{(2)})_t\big\|_{L^1(\Omega)} &= \big\|\mathbb{E}_t\big[\mathfrak{O}_T^{(1)}\mathfrak{O}_T^{(2)} - (\mathfrak{O}^{(1)}\star\mathfrak{O}^{(2)})_T\big]\big\|_{L^1(\Omega)} \\ &\le \big\|\mathfrak{O}_T^{(1)}\mathfrak{O}_T^{(2)} - (\mathfrak{O}^{(1)}\star\mathfrak{O}^{(2)})_T\big\|_{L^1(\Omega)} \end{aligned}$$

vanishing as $T\to\infty$; closedness of $\mathfrak{O}^{(1)}\star\mathfrak{O}^{(2)}$ implies

$$\begin{aligned} \big\|(\mathfrak{O}^{(1)}\star\mathfrak{O}^{(2)})_\infty - \mathfrak{O}_t^{(1)}\mathfrak{O}_t^{(2)}\big\|_{L^1(\Omega)} &\le \|(\mathfrak{O}^{(1)}\star\mathfrak{O}^{(2)})_\infty - (\mathfrak{O}^{(1)}\star\mathfrak{O}^{(2)})_t\|_{L^1(\Omega)} \\ &+ \big\|(\mathfrak{O}^{(1)}\star\mathfrak{O}^{(2)})_t - \mathfrak{O}_t^{(1)}\mathfrak{O}_t^{(2)}\big\|_{L^1(\Omega)} \end{aligned}$$

vanishing as $t\to\infty$, which, coupled with $\mathfrak{O}_t^{(1)}\mathfrak{O}_t^{(2)} \to \mathfrak{O}_\infty^{(1)}\mathfrak{O}_\infty^{(2)}$ in probability (implied by convergence of $\mathfrak{O}^{(j)}$), yields that $(\mathfrak{O}^{(1)}\star\mathfrak{O}^{(2)})_\infty = \mathfrak{O}_\infty^{(1)}\mathfrak{O}_\infty^{(2)}$ almost surely.

The formulae of the above remark are not the only reason for requiring the condition (2.13). Let $q^{(a)} \sim_\varphi p^{(a)}$ coherent germs for $\mathfrak{O}^{(a)}$ and $q^{(b)} \sim_\varphi p^{(b)}$ are coherent germs for $\mathfrak{O}^{(b)}$ such that $q := q^{(a)}q^{(b)}$ and $p := p^{(a)}p^{(b)}$ are coherent. To have a well-defined notion of the observable product, we then want that $p$ and $q$ reconstruct the same martingale, i.e. that $\mathfrak{R}(q^{(a)}q^{(b)}) = \mathfrak{R}(p^{(a)}p^{(b)})$. Leaving condition (2.13) out of the definition of compatibility leaves us unable to prove this consistency result in general: there is no guarantee that $q^{(a)}q^{(b)} \sim_\varphi p^{(a)}p^{(b)}$.

**Proposition 2.13.** *If $q^{(a)} \sim_{\varphi} p^{(a)}$ are coherent germs for $\mathfrak{Q}^{(a)}$ and $q^{(b)} \sim_{\varphi} p^{(b)}$ are coherent germs for $\mathfrak{Q}^{(b)}$ such that $(q^{(a)}, q^{(b)})$ and $(p^{(a)}, p^{(b)})$ are compatible, then $q^{(a)} q^{(b)} \sim_{\varphi} p^{(a)} p^{(b)}$. Consequently, $\mathfrak{R}(q^{(a)} q^{(b)}) = \mathfrak{R}(p^{(a)} p^{(b)})$.*

**Proof.** Indeed,

$$\begin{aligned}
& \|q_t^{(a)} q_t^{(b)}(\varphi_t) - p_t^{(a)} p_t^{(b)}(\varphi_t)\|_{L^1(\Omega)} \\
\leq & \ \|q_t^{(a)} q_t^{(b)}(\varphi_t) - \mathfrak{R}(q^{(a)} q^{(b)})_t\|_{L^1(\Omega)} + \|\mathfrak{R}(q^{(a)} q^{(b)})_t - \mathfrak{Q}_t^{(a)} \mathfrak{Q}_t^{(b)}\|_{L^1(\Omega)} \\
+ & \ \|p_t^{(a)} p_t^{(b)}(\varphi_t) - \mathfrak{R}(p^{(a)} p^{(b)})_t\|_{L^1(\Omega)} + \|\mathfrak{R}(p^{(a)} p^{(b)})_t - \mathfrak{Q}_t^{(a)} \mathfrak{Q}_t^{(b)}\|_{L^1(\Omega)}
\end{aligned}$$

vanishes as $t \to \infty$ by coherence and compatibility. Combining the fact that $q^{(a)} q^{(b)} \sim_{\varphi} p^{(a)} p^{(b)}$ with (2.5) and (2.6) yields $\mathfrak{R}(q^{(a)} q^{(b)}) = \mathfrak{R}(p^{(a)} p^{(b)})$. □

In light of the discussion above, we do not expect the product $\star$ to be a well-defined binary operation on the space $\mathscr{O}_c$ of observables possessing coherent germs. Regardless, it is well-defined as a partial function and is easily seen to satisfy the following properties.

**Proposition 2.14.** *The partial binary operation $\star \colon \mathscr{O}_c \times \mathscr{O}_c \rightharpoonup \mathscr{O}_c$ is*

1. *commutative in the following sense: let $\mathfrak{O}^{(a)}, \mathfrak{O}^{(b)} \in \mathscr{O}_c$ with coherent germs $q^{(a)}, q^{(b)}$, respectively. Then $(q^{(a)}, q^{(b)})$ is compatible if and only if $(q^{(b)}, q^{(a)})$ is compatible, and when they are it holds that*
$$\mathfrak{O}^{(a)} \star \mathfrak{O}^{(b)} = \mathfrak{O}^{(b)} \star \mathfrak{O}^{(a)}.$$
2. *distributive in the following sense: if $\mathfrak{O}^{(a)}, \mathfrak{O}^{(b)}, \mathfrak{O}^{(c)} \in \mathscr{O}_c$ are such that $(q^{(a)}, q^{(b)})$ and $(q^{(a)}, q^{(c)})$ are compatible, then $\mathfrak{O}^{(a)} \star \mathfrak{O}^{(b)}$, $\mathfrak{O}^{(a)} \star \mathfrak{O}^{(c)}$ and $\mathfrak{O}^{(a)} \star (\mathfrak{O}^{(b)} + \mathfrak{O}^{(c)})$ exist and*
$$\mathfrak{O}^{(a)} \star (\mathfrak{O}^{(b)} + \mathfrak{O}^{(c)}) = \mathfrak{O}^{(a)} \star \mathfrak{O}^{(b)} + \mathfrak{O}^{(a)} \star \mathfrak{O}^{(c)}.$$
3. *associative in the following sense: if $\mathfrak{O}^{(a)}, \mathfrak{O}^{(b)}, \mathfrak{O}^{(c)} \in \mathscr{O}_c$, are such that $(q^{(a)}, q^{(b)}, q^{(c)})$ is compatible, then*
$$(\mathfrak{O}^{(a)} \star \mathfrak{O}^{(b)}) \star \mathfrak{O}^{(c)} = \mathfrak{O}^{(a)} \star \mathfrak{O}^{(b)} \star \mathfrak{O}^{(c)} = \mathfrak{O}^{(a)} \star (\mathfrak{O}^{(b)} \star \mathfrak{O}^{(c)}).$$
4. *unital in the following sense: there exists a unique observable $\mathtt{1} = (\mathtt{1}_t)_{t \in \mathbb{R}_+} \in \mathscr{O}_c$ given by $\mathtt{1}_t := 1$ such that $\mathfrak{O} \star \mathtt{1}$ exists for all $\mathfrak{O} \in \mathscr{O}_c$ and $\mathfrak{O} \star \mathtt{1} = \mathfrak{O}$.*

**Proof.** Commutativity, associativity and unitality are trivial.

For distributivity, note that by linearity of $\mathscr{L}^V$ and the triangle inequality, we have that $q^{(a)}(q^{(b)} + q^{(c)})$ is coherent. Hence

$$\begin{aligned}
& \|\mathfrak{R}(q^{(a)}(q^{(b)} + q^{(c)}))_t - \mathfrak{O}_t^{(a)}(\mathfrak{O}_t^{(b)} + \mathfrak{O}_t^{(c)})\|_{L^1(\Omega)} \\
\leq & \ \|\mathfrak{R}(q^{(a)} q^{(b)})_t - \mathfrak{O}_t^{(a)} \mathfrak{O}_t^{(b)}\|_{L^1(\Omega)} + \|\mathfrak{R}(q^{(a)} q^{(c)})_t - \mathfrak{O}_t^{(a)} \mathfrak{O}_t^{(c)}\|_{L^1(\Omega)}
\end{aligned}$$

vanishing as $t \to \infty$, yielding existence of $\mathfrak{O}^{(a)} \star (\mathfrak{O}^{(b)} + \mathfrak{O}^{(c)})$. □

**Remark 2.15.** The reverse implication for distributivity is not true, as seen by the following simple example. Let $\varphi_t := B_t$ be the scale interpolation of the free field and let $x \in \mathbb{R}^d$, $d \geq 2$. Consider the germs $q^{(a)} = q^{(b)} := \delta_x$ and $q^{(c)} := -\delta_x$. These are coherent and the corresponding observables are $\mathfrak{Q}^{(a)} = \mathfrak{Q}^{(b)} = B_t(x)$ and $\mathfrak{Q}^{(c)} = -B_t(x)$. The germ $q^{(a)}(q^{(b)} + q^{(c)}) = 0$ is coherent, the observable $\mathfrak{O}^{(a)} \star (\mathfrak{O}^{(b)} + \mathfrak{O}^{(c)})$ exists and is the constant zero martingale. However, $\mathfrak{O}^{(a)} \star \mathfrak{O}^{(b)}$ and $\mathfrak{O}^{(a)} \star \mathfrak{O}^{(c)}$ do not exist, in the sense that e.g. $q^{(a)} q^{(b)} = \delta_x^2$ is not coherent; renormalisation is needed to define such a product.

## 2.5 Germ subspaces

It is important to know when the product of observables exists. The previous definitions suggest to introduce subspaces of germs that help determining whether they are compatible. The rest of this section is concerned with setting up such machinery.

Given $M\in\mathbb{N}$, $\boldsymbol{\alpha}=(\alpha_1,\ldots,\alpha_M)\in\mathbb{C}^M$, $\boldsymbol{y}=(y_1,\ldots,y_M)\in(\mathbb{R}^d)^M$, denote by

$$\langle\boldsymbol{\alpha},\phi(\boldsymbol{y})\rangle:=\sum_{j=1}^M\alpha_j\phi(y_j).$$

Moreover, given $\boldsymbol{\alpha}_a\in\mathbb{C}^{M_a}$, $\boldsymbol{\alpha}_b\in\mathbb{C}^{M_b}$ we define their join by

$$\boldsymbol{\alpha}_a\veebar\boldsymbol{\alpha}_b:=(\alpha_{a,1},\ldots,\alpha_{a,M_a},\alpha_{b,1},\ldots,\alpha_{b,M_b})\in\mathbb{C}^{M_a+M_b}$$

and similarly $\boldsymbol{y}_a\veebar\boldsymbol{y}_b$. Obviously

$$\langle\boldsymbol{\alpha}_a\veebar\boldsymbol{\alpha}_b,\phi(\boldsymbol{y}_a\veebar\boldsymbol{y}_b)\rangle=\langle\boldsymbol{\alpha}_a,\phi(\boldsymbol{y}_a)\rangle+\langle\boldsymbol{\alpha}_b,\phi(\boldsymbol{y}_b)\rangle. \tag{2.16}$$

**Definition 2.16.** Let $\kappa,\sigma\in\mathbb{R}$, $c\geq 0$, $p\in[1,\infty]$, $M\in\mathbb{N}$, $\boldsymbol{\alpha}\in\mathbb{C}^M$, $\boldsymbol{y}\in(\mathbb{R}^d)^M$. Consider the seminorm

$$\|q\|_{\hat{\mathscr{O}}^{\kappa,\sigma,c}_{p,\boldsymbol{\alpha},\boldsymbol{y}}}:=\sup_{t\in\mathbb{R}_+}\langle t\rangle^{-\kappa}[_t^{-c}\|e^{-\langle\boldsymbol{\alpha},\varphi_t(\boldsymbol{y})\rangle}q_t(\varphi_t)\|_{L^p(\Omega)}+\sup_{t\in\mathbb{R}_+}\langle t\rangle^{\sigma+1}[_t^{-c}\|e^{-\langle\boldsymbol{\alpha},\varphi_t(\boldsymbol{y})\rangle}\mathscr{L}^Vq_t(\varphi_t)\|_{L^p(\Omega)}$$

for a germ $q:\mathbb{R}_+\times\mathcal{H}\to\mathbb{R}$ such that $\mathscr{L}^Vq$ is well-defined. Define the associated subspaces of germs

$$\hat{\mathscr{O}}^{\kappa,\sigma,c}_{p,\boldsymbol{\alpha},\boldsymbol{y}}:=\left\{q:\mathbb{R}_+\times\mathcal{H}\to\mathbb{R}:\|q\|_{\hat{\mathscr{O}}^{\kappa,\sigma,c}_{p,\boldsymbol{\alpha},\boldsymbol{y}}}<\infty\right\}$$

We call the parameter $\kappa\in\mathbb{R}$ the *superficial degree of divergence*, and $\sigma>0$ the *accuracy of perturbative approximation.*

The $\boldsymbol{\alpha},\boldsymbol{y}$ in the notation is frequently suppressed when this does not lead to confusion. Clearly, for any $\varepsilon>0$, we have

$$\hat{\mathscr{O}}^{\kappa,\sigma,c}_{p+\varepsilon}\subseteq\hat{\mathscr{O}}^{\kappa,\sigma,c}_{p},\qquad\hat{\mathscr{O}}^{\kappa,\sigma,c}_{p}\subseteq\hat{\mathscr{O}}^{\kappa+\varepsilon,\sigma,c}_{p},\qquad\hat{\mathscr{O}}^{\kappa,\sigma,c}_{p}\subseteq\hat{\mathscr{O}}^{\kappa,\sigma-\varepsilon,c}_{p},\qquad\hat{\mathscr{O}}^{\kappa,\sigma,c}_{p}\subseteq\hat{\mathscr{O}}^{\kappa,\sigma,c+\varepsilon}_{p}. \tag{2.17}$$

Denote

$$\hat{\mathscr{O}}^{\kappa,\sigma}_{p}:=\bigcup_{c\geq 0}\hat{\mathscr{O}}^{\kappa,\sigma,c}_{p},\qquad\hat{\mathscr{O}}^{\kappa,\sigma,c}_{\infty-}:=\bigcap_{p\in[1,\infty)}\hat{\mathscr{O}}^{\kappa,\sigma,c}_{p},\qquad\hat{\mathscr{O}}^{\kappa,\sigma}_{\infty-}:=\bigcap_{p\in[1,\infty)}\bigcup_{c\geq 0}\hat{\mathscr{O}}^{\kappa,\sigma,c}_{p}.$$

Note that this definition of the $\hat{\mathscr{O}}^{\kappa,\sigma}_{\infty-}$ space allows for $p$-dependent logarithmic corrections, as opposed to the smaller space $\bigcup_{c\geq 0}\bigcap_{p\in[1,\infty)}\hat{\mathscr{O}}^{\kappa,\sigma,c}_{p}$. The applications we consider in this paper do not require this stronger property, hence for our purposes the weaker space suffices, even though the germs we will construct in fact lie in the stronger one.

The purpose of the objects introduced is the quantitative description of compatibility.

**Lemma 2.17. (factorization lemma)** *Let $\kappa_j,\sigma_j,p_j,M_j,\boldsymbol{\alpha}_j,\boldsymbol{y}_j$, $j\in 1{:}n$, be as in Definition 2.16, and denote by $\boldsymbol{\alpha}=\boldsymbol{\alpha}_1\veebar\ldots\veebar\boldsymbol{\alpha}_n$, $\boldsymbol{y}=\boldsymbol{y}_1\veebar\ldots\veebar\boldsymbol{y}_n$. Suppose $q^{\langle j\rangle}\in\hat{\mathscr{O}}^{\kappa_j,\sigma_j}_{p_j,\boldsymbol{\alpha}_j,\boldsymbol{y}_j}$, and for all $j_1\neq j_2\in 1{:}n$ we have*

$$\sup_{t\in\mathbb{R}_+}\langle t\rangle^{\sigma+1}\left\|e^{-\langle\boldsymbol{\alpha},\varphi_t(\boldsymbol{y})\rangle}\mathrm{Tr}\big[\dot{G}_t\mathrm{D}q_t^{\langle j_1\rangle}(\varphi_t)\mathrm{D}q_t^{\langle j_2\rangle}(\varphi_t)\big]\prod_{j\in 1:n\setminus\{j_1,j_2\}}q_t^{\langle j\rangle}(\varphi_t)\right\|_{L^p(\Omega)}<\infty$$

*for a $p\in[1,\infty]$, $\sigma\in\mathbb{R}$ satisfying*

$$\sum_{j\in 1:n}p_j^{-1}\leq p^{-1},\qquad\sigma\leq\bigwedge_{j\in 1:n}\left(\sigma_j-\sum_{m\in(1:n)\setminus\{j\}}\kappa_m\right).$$

*Then* $\prod_{j=1}^n q^{\langle j\rangle} \in \hat{\mathscr{O}}_{p,\boldsymbol{\alpha},\boldsymbol{y}}^{\kappa,\sigma}$ *for any* $\kappa$ *with* $\sum_{j\in 1:n}\kappa_j \leq \kappa$.

**Proof.** Clearly

$$\langle t\rangle^{-\sum_{j=1}^n \kappa_j} \left\| \prod_{j\in 1:n} q_t^{\langle j\rangle}(\varphi_t)\right\|_{L^p(\Omega)} \leq \prod_{j\in 1:n} \langle t\rangle^{-\kappa_j} \| q_t^{\langle j\rangle}(\varphi_t)\|_{L^{p_j}(\Omega)}.$$

The generator of the product can be bounded as

$$\begin{aligned}
&\langle t\rangle^{\sigma+1} \left\| \mathscr{L}^V \left( \prod_{j\in 1:n} q_t^{\langle j\rangle}\right)(\varphi_t)\right\|_{L^p(\Omega)} \\
\leq\ & \sum_{j=1}^n \langle t\rangle^{\sigma_j+1} \| \mathscr{L}^V q_t^{\langle j\rangle}(\varphi_t)\|_{L^{p_j}(\Omega)} \prod_{m\neq j} \langle t\rangle^{-\kappa_m} \| q_t^{\langle m\rangle}(\varphi_t)\|_{L^{p_m}(\Omega)} \\
+\ & \sum_{j<k} \langle t\rangle^{\sigma+1} \left\| \mathrm{Tr}\big[\dot{G}_t \mathrm{D} q_t^{\langle j\rangle} \mathrm{D} q_t^{\langle k\rangle}\big] \prod_{m\notin\{j,k\}} q_t^{\langle m\rangle}(\varphi_t)\right\|_{L^p(\Omega)}.
\end{aligned}$$

Coupled with (2.16) and (2.17), the statement follows. □

**Remark 2.18.** The reason for the exponential weights $e^{-\langle \boldsymbol{\alpha}, \varphi_t(\boldsymbol{y})\rangle}$ in Definition 2.16 is to allow for 2 dimensional theories in the general framework, where observables of the form $\exp(\alpha\phi(y))$ play an important role. In turn, such exponentials do not behave well with the Hölder inequalities of the above lemma (consider e.g. the $d=2$ Gaussian free field) and require slightly different treatment. Since we are going to apply the general theory on such a 2 dimensional model, we will see how the general case of allowing exponentials work and closes the arguments; for theories and observables that are not of this generality, one may always take $\boldsymbol{\alpha}=0$ as a special case.

## 2.6 The strategy of the construction

In order to find $q$ such that (2.10) holds, we need to have a good understanding of the solution theory of the observable flow equation $\mathscr{L}^V q=0$, i.e.

$$\partial_s q_s + \frac{1}{2}\Delta_{\dot{G}_s} q_s = \langle \mathrm{D}q_s, \dot{G}_s \mathrm{D}V_s\rangle \tag{2.18}$$

To do this, we construct $\ell^*$-th order approximate germs with a Picard iteration scheme: we let

$$q^{[\leq \ell^*]} := \sum_{\ell\leq\ell^*} q^{[\ell]}$$

where $q^{[0]}\equiv 0$, and $q^{[\ell]}$ for $\ell\geq 1$ is defined to be the solution of

$$\partial_s q_s^{[\ell]} + \frac{1}{2}\mathrm{Tr}\big[\dot{G}_s \mathrm{D}^2 q_s^{[\ell]}\big] = \sum_{\ell'+\ell''=\ell} \big\langle \mathrm{D}q_s^{[\ell']}, \dot{G}_s \mathrm{D}V_s^{[\ell'']}\big\rangle \tag{2.19}$$

with appropriate boundary conditions (see the next section), where the $V^{[\ell]}$ are iterates of a Picard scheme similar to (2.19) of the Polchinski equation for the potential $V$; namely, $V^{[0]}\equiv 0$ and for $\ell\geq 1$, $V^{[\ell]}$ solves

$$\partial_s V_s^{[\ell]} + \frac{1}{2}\mathrm{Tr}\big[\dot{G}_s \mathrm{D}^2 V_s^{[\ell]}\big] = \frac{1}{2}\sum_{\ell'+\ell''=\ell} \big\langle \mathrm{D}V_s^{[\ell']}, \dot{G}_s \mathrm{D}V_s^{[\ell'']}\big\rangle. \tag{2.20}$$

Set $V^{[\leq\ell^*]} := \sum_{\ell\leq\ell^*} V^{[\ell]}$ and write $F^{[\ell]} := \mathrm{D}V^{[\ell]}$.

We are now going to show semi-heuristically that solving (2.19) up to high enough order $\ell^*$ gives a good approximating germ, i.e. for which $\mathscr{L}^V q^{[\leq\ell^*]} \approx_{\ell^*} 0$ in a suitable sense.

Suppose there exists an $\tilde{\ell}\in\mathbb{N}$ such that for any $\ell'\geq\tilde{\ell}$, the canonical Ito remainder associated to the germ $F^{[\leq\ell']}(y)$, denoted by $r^{F^{[\leq\ell']}}$, exists, and

$$\lim_{t\to\infty}\big\|F_t^{[\leq\ell']}(\varphi_t)(y)-F^{[\leq\tilde{\ell}]}(\varphi_t)(y)\big\|_{L^1(\Omega)}=0. \tag{2.21}$$

Denote by $\tilde{\ell}$ the smallest of such integers. Combining (2.21) with (2.5) yields that all truncations with $\ell'\geq\tilde{\ell}$ reconstruct the same observable (field), which we denote by $\{F_s(y)\}_{s\geq 0}$ and call the effective force.

For any $\ell\in\mathbb{N}_0$, write

$$\dot{G}R^{[>\ell]}:=\begin{cases}\dot{G}r^{F^{[\leq\ell]}} & \text{if } \ell\geq\tilde{\ell}\\ \displaystyle\sum_{\ell''=\ell+1}^{\tilde{\ell}}\dot{G}F^{[\ell'']}+\dot{G}r^{F^{[\leq\tilde{\ell}]}} & \text{if } \ell<\tilde{\ell}.\end{cases} \tag{2.22}$$

**Lemma 2.19.** *Letting $q^{[\leq\ell^*]}:=\sum_{\ell\leq\ell^*}q^{[\ell]}$ where $q^{[0]}\equiv 0$, and $q^{[\ell]}$ for $\ell\geq 1$ is defined via (2.19), we have*

$$\mathscr{L}^V q^{[\leq\ell^*]}=-\sum_{\ell\leq\ell^*}\langle \mathrm{D}q^{[\ell]},\dot{G}R^{[>(\ell^*-\ell)]}\rangle. \tag{2.23}$$

**Proof.** To simplify notation, write $F_s^{[\ell]}:=F_s^{[\ell]}(\varphi_s)$ and similarly for $q^{[\ell]}$. Then Ito reconstruction and (2.22) gives that for arbitrary $\ell^*\in\mathbb{N}_0$ we have

$$\dot{G}F=\dot{G}F^{[\leq\ell^*]}+\dot{G}R^{[>\ell^*]}. \tag{2.24}$$

Thus compute with (2.9) and (2.19) that

$$\begin{aligned}\mathscr{L}^V q^{[\leq\ell^*]} &= \left(\partial_s+\frac{1}{2}\Delta_{\dot{G}}\right)q^{[\leq\ell^*]}-\langle\mathrm{D}q^{[\leq\ell^*]},\dot{G}F^{[\leq\ell^*]}\rangle-\langle\mathrm{D}q^{[\leq\ell^*]},\dot{G}R^{[>\ell^*]}\rangle\\ &= \sum_{\ell\leq\ell^*}\sum_{\substack{\ell'+\ell''=\ell\\ \ell,\ell',\ell''\in\mathbb{N}_0}}\langle\mathrm{D}q^{[\ell']},\dot{G}F^{[\ell'']}\rangle-\left\langle\sum_{\ell\leq\ell^*}\mathrm{D}q^{[\ell]},\sum_{\ell\leq\ell^*}\dot{G}F^{[\ell]}\right\rangle-\langle\mathrm{D}q^{[\leq\ell^*]},\dot{G}R^{[>\ell^*]}\rangle\\ &= -\sum_{\substack{\ell',\ell''\leq\ell^*\\ \ell'+\ell''\geq\ell^*+1}}\langle\mathrm{D}q^{[\ell']},\dot{G}F^{[\ell'']}\rangle-\left\langle\sum_{\ell\leq\ell^*}\mathrm{D}q^{[\ell]},\dot{G}R^{[>\ell^*]}\right\rangle.\end{aligned}$$

Pick an $\ell'\leq\ell^*$. Then the first sum runs over all the $\dot{G}F^{[\ell'']}$ with

$$\left\{\begin{array}{lcl}\ell''+\ell' & \geq & \ell^*+1\\ \ell'' & \leq & \ell^*\end{array}\right\}\Leftrightarrow \ell^*-\ell'+1\leq\ell''\leq\ell^*,$$

thus

$$\begin{aligned}\sum_{\ell''=\ell^*-\ell'+1}^{\ell^*}\dot{G}F^{[\ell'']}+\dot{G}R^{[>\ell^*]} &= \sum_{\ell''=\ell^*-\ell'+1}^{\ell^*}\dot{G}F^{[\ell'']}+(\dot{G}F-\dot{G}F^{[\leq\ell^*]})\\ &= \dot{G}F-\dot{G}F^{[\leq(\ell^*-\ell')]}\\ &= \dot{G}R^{[>(\ell^*-\ell')]}\end{aligned}$$

where we used (2.24) twice. □

Within the rest of this subsection, we tune the discussion to the $d=2$ sine–Gordon model in regime $\beta^2<6\pi$; nevertheless, the ideas introduced are general at their core and should extend to other models mutatis mutandis.

Suppose that there exists a fixed $\kappa \geq 0$, the superficial degree of divergence of the observable (cf. Definition 2.16), for which we can prove the formal estimate

$$\left\| \mathrm{D} q_s^{[\ell]} \right\|_{L^p(\Omega)} \lesssim \langle s \rangle^{\kappa - (\ell - 1)\delta} \tag{2.25}$$

where $\delta := 1 - \frac{\beta^2}{8\pi} \in (0, 1]$ is the distance from criticality, and that we also have the bound

$$\left\| \dot{G}_s R_s^{[>\ell]} \right\|_{L^\infty(\Omega; L^\infty)} \lesssim \langle s \rangle^{-\delta - \ell\delta - 1} \tag{2.26}$$

modulo logarithmic corrections.

Plugging these into (2.23), using Holder and Young's inequalities and $\|\dot{\mathcal{G}}_s\|_{L^1} \lesssim \langle s \rangle^{-2}$, we formally obtain

$$\begin{aligned} \left\| \mathscr{L}^V q_s^{[\leq \ell^*]} \right\|_{L^1(\Omega)} &\leq \sum_{\ell \leq \ell^*} \left\| \left\langle \mathrm{D} q_s^{[\ell]}, \dot{G}_s R_s^{[>(\ell^* - \ell)]} \right\rangle \right\|_{L^1(\Omega)} \\ &\lesssim \sum_{\ell \leq \ell^*} \langle s \rangle^{\kappa - (\ell - 1)\delta} \langle s \rangle^{-\delta - (\ell^* - \ell)\delta - 1} \\ &\lesssim \langle s \rangle^{\kappa - \ell^* \delta - 1}. \end{aligned} \tag{2.27}$$

Since $\delta > 0$, choosing $\ell^*$ large enough then yields an arbitrarily precise approximate solution $q_s^{[\leq \ell^*]}$ of the equation $\mathscr{L}^V q = 0$, i.e. we can construct germs $q_s^{[\leq \ell^*]} \in \hat{\mathscr{O}}_{\infty-}^{\kappa, \sigma}$ with arbitrarily high accuracy of perturbative approximation $\sigma$ (cf. Definition 2.16). The construction of observables and observable products will then follow from the general theory developed so far, in particular the factorization lemma (Lemma 2.17) and Ito reconstruction (Lemma 2.6).

Note that in the above argument we relied on the postulated bounds (2.21) and (2.26) which of course need to be proven first before we can move on to the construction of various germs $q$. It is not too difficult to see that in fact these can be obtained in a similar vein to what we have explained just now for general observables, by proving (2.25) for the germ of the force $F$ with superficial degree of divergence $\kappa := 1 - \delta$.

To construct the effective force and hence the sine-Gordon measure, in the regime $\beta^2 < 6\pi$, it is sufficient to develop the equation (2.20) up to order $\ell^* = 3$, and in general, to some $\ell^*(\beta^2)$ that blows up as we approach criticality $\beta^2 \to 8\pi$. However, since our observables can have arbitrarily high superficial degree of divergence $\kappa$, we are forced to exploit the high order cancellations in Lemma 2.19, and hence need to develop the iterative scheme to arbitrarily high $\ell^*$ for *any* $\beta^2 > 0$; indeed, cf. (2.27).

The following sections carry out the above roadmap rigorously. Section 3 is concerned with studying (2.20); as a consequence of the bounds derived there, Section 4 constructs the effective force, and verifies (2.21) and (2.26). Section 5 constructs $q^{[\ell]}$ for various observables by solving the equation (2.19). Finally, Section 6 finishes the construction and proves Theorems 1.1 and 1.2.

## 2.7 Renormalisation via flow equations

Constructing observables of a field theory requires additional renormalisation on top of the renormalisation of the model itself. Indeed, taking general nonlinear functionals $Q(\varphi_\infty)$ of an already-renormalised EQFT $\varphi_\infty$ might introduce new divergences that need to be accounted for. An interesting aspect of these additional singularities is that they can be especially severe: in the language of previous sections, the superficial degree of divergence $\kappa \geq 0$ can be arbitrarily large. In contrast, power counting arguments suggest that the renormalisation of the model itself needs to take care of only relatively mild singularities.

Curing the aforementioned arbitrarily severe singularities, in general, requires introducing arbitrarily many (but finite) renormalisation counterterms; moreover, the form of such counterterms for complicated nonlinearities may be quite non-obvious. Hence in order to get a complete theory of composite observable renormalisation, it is desirable to approach the introduction of counterterms in a systematic and fully automated manner. The requirement of propagating the inductive bounds (2.25) in $\ell$ through the flow equation (2.19) not only offers such a way, but in fact forces it upon us. The counterterms are generated through an appropriate choice of boundary conditions on the flow, which we now explain.

We will construct the germs by making the ansatz

$$q_t^{[\ell],T}(\phi) = \sum_{\boldsymbol{a} \in \mathfrak{M}} \left\langle \gamma_{t,T}^{(\boldsymbol{a})}, f^{(\boldsymbol{a})}(\phi) \right\rangle \tag{2.28}$$

where $\gamma_{t,T}^{(\boldsymbol{a})}$ are $\phi$-independent kernels, $f^{(\boldsymbol{a})}(\phi)$ are functionals of the field, and $\boldsymbol{a}$ are multiindices carrying information about the expansion, containing e.g. reference to $\ell$. We require that only finitely many kernels of the sum are nonzero.

Plugging the above ansatz into (2.19), matching the kernels corresponding to the same basis element $f^{(\boldsymbol{a})}(\phi)$ we obtain a hierarchy of linear ODEs of the form

$$\partial_t \gamma_{t,T}^{(\boldsymbol{a})} + \dot{\mathbb{U}}_t \gamma_{t,T}^{(\boldsymbol{a})} = -\mathscr{S}_{t,T}^{(\boldsymbol{a})}$$

for some $t$-dependent coefficient function $\dot{\mathbb{U}}_t$ and source term $\mathscr{S}_{t,T}^{(\boldsymbol{a})}$. The latter is a sum over $n$-linear functionals of kernels $\gamma_{t,T}^{(\boldsymbol{d})}$. Using the Duhamel principle, we rewrite the equation into mild form

$$\gamma_{t,T}^{(\boldsymbol{a})} = e^{\mathbb{U}_{t,T}} \gamma_{T,T}^{(\boldsymbol{a})} + \int_t^T e^{\mathbb{U}_{t,s}} \mathscr{S}_{s,T}^{(\boldsymbol{a})} \mathrm{d}s \tag{2.29}$$

In order to show an estimate of the form (2.25), we set up the ansatz (2.28) in a way so that the hierarchy (2.29) is triangular in $\boldsymbol{a}$; we then solve the system inductively and try to propagate a bound of the form

$$|||\gamma_{t,T}^{(\boldsymbol{a})}||| \lesssim \langle t \rangle^{[\![\boldsymbol{a}]\!]} \mathfrak{l}_t^c \tag{2.30}$$

with $\mathfrak{l}_t := 1 + \log(1 \vee t)$, where $c$ is a constant that we do not track and may change from line to line, and $|||\cdot|||$ is a suitable norm compatible with the structure of the flow equation. To this end, we will aim to show that the integrand of (2.29) satisfies

$$|||e^{\mathbb{U}_{t,s}} \mathscr{S}_{s,T}^{(\boldsymbol{a})}||| \lesssim \langle s \rangle^{[\![\boldsymbol{a}]\!]} \mathfrak{l}_s^c \langle s \rangle^{-1} \tag{2.31}$$

Kernels with $[\![\boldsymbol{a}]\!] < 0$ are called irrelevant: they are expected to remain bounded under the flow and hence get a trivial UV boundary condition. In contrast, kernels with $[\![\boldsymbol{a}]\!] \geq 0$ need explicit renormalisation and are called relevant. Thus in general for $\ell \geq 2$ we aim to impose the boundary conditions

$$\begin{aligned} \gamma_{T,T}^{(\boldsymbol{a})} &= 0 \quad \text{if} \quad [\![\boldsymbol{a}]\!] < 0 \\ \gamma_{0,T}^{(\boldsymbol{a})} &= 0 \quad \text{if} \quad [\![\boldsymbol{a}]\!] \geq 0. \end{aligned} \tag{2.32}$$

It is clear that the case $\ell = 1$ is slightly different, since it has no contributions from lower iteration orders, as $q^{[0]} \equiv 0$.

When $[\![\boldsymbol{a}]\!] < 0$, we flow *backward*: the boundary condition at $t = T$ yields the equation

$$\gamma_{t,T}^{(\boldsymbol{a})} = \int_t^T e^{\mathbb{U}_{t,s}} \mathscr{S}_{s,T}^{(\boldsymbol{a})} \, \mathrm{d}s \tag{2.33}$$

which propagates the inductive bound (2.30) using (2.31):

$$|||\gamma_{t,T}^{(\boldsymbol{a})}||| \lesssim \int_t^\infty \langle s \rangle^{[\![\boldsymbol{a}]\!]} \mathfrak{l}_s^c \langle s \rangle^{-1} \mathrm{d}s \lesssim \langle t \rangle^{[\![\boldsymbol{a}]\!]} \mathfrak{l}_t^c.$$

In the $[\![\boldsymbol{a}]\!] \geq 0$ case we flow *forward*: plugging into (2.29) the relevant boundary condition (2.32) at $t = 0$ we obtain

$$e^{\mathbb{U}_{0,T}} \gamma_{T,T}^{(\boldsymbol{a})} = -\int_0^T e^{\mathbb{U}_{0,s}} \mathscr{S}_{s,T}^{(\boldsymbol{a})} \, \mathrm{d}s,$$

which, since $-\mathbb{U}_{0,T}+\mathbb{U}_{0,s}=-\mathbb{U}_{s,T}$, means that $\gamma_{T,T}^{(\boldsymbol{a})}=-\int_0^T e^{-\mathbb{U}_{s,T}}\mathscr{S}_{s,T}^{(\boldsymbol{a})}\,\mathrm{d}s$. Plugging this back into (2.29) yields, as $-\mathbb{U}_{s,T}+\mathbb{U}_{t,T}=\mathbb{U}_{t,s}$,

$$\gamma_{t,T}^{(\boldsymbol{a})}=-\int_0^t e^{\mathbb{U}_{t,s}}\mathscr{S}_{s,T}^{(\boldsymbol{a})}\,\mathrm{d}s. \tag{2.34}$$

which again propagates the inductive bound (2.30), via (2.31):

$$|||\gamma_{t,T}^{(\boldsymbol{a})}|||\lesssim\int_0^t \langle s\rangle^{[\![\boldsymbol{a}]\!]}\mathfrak{l}_s^c\langle s\rangle^{-1}\mathrm{d}s\lesssim\langle t\rangle^{[\![\boldsymbol{a}]\!]}\mathfrak{l}_t^c.$$

In the relevant case $[\![\boldsymbol{a}]\!]\geq 0$, we expect from (2.34) the kernel $\gamma_{t,T}^{(\boldsymbol{a})}$ to blow up as $t\to\infty$; coupled with the corresponding basis element, $\langle\gamma_{t,T}^{(\boldsymbol{a})},f^{(\boldsymbol{a})}(\phi)\rangle$ is a renormalisation counterterm. Moreover, the uniform-in-$T$ bounds allow for inductively taking the limit $T\to\infty$ in the source terms $\mathscr{S}_{s,T}^{(\boldsymbol{a})}$, rendering the counterterm $\langle\gamma_{t,\infty}^{(\boldsymbol{a})},f^{(\boldsymbol{a})}(\phi)\rangle$ also independent of the UV cutoff $T$.

# 3 The Polchinski equation of sine–Gordon

In this section we study the Polchinski equation (2.20). The aim is twofold: first, since the observable flow equation (2.19) contains $V^{[\ell]}$, we expect that proving the source term bound (2.31) for kernels of the germ $q$ requires a good understanding of kernels of the potential germ $V$. Second, as we have already noted, closing the argument in Section 2.6 requires access to the estimates (2.21) and (2.26) which need to be proven first. While the analysis of the Polchinski equation of sine–Gordon in the regime $\beta^2<6\pi$ is not novel, see e.g. [7], [29], we include a fully self-contained treatment here partly to fix notation, and partly because we use it to develop ideas, techniques and motivation that we will use in analysis of the observable flow equations later on.

As a warmup, we first solve the equation (2.20) in the Wick regime $\beta^2<4\pi$ by choosing the ansatz (2.28) to be the usual Brydges–Kennedy Mayer expansion. We then explain how the regime $\beta^2\geq 4\pi$, because of the presence of relevant kernels, requires localisation and flowing forward, and hence induces a generalised kernel representation. Doing so, we prove a sharp screening estimate on the covariance matrix of the Mayer representation that will turn out to be crucial in the construction of vertex field observables.

## 3.1 Warmup: the Wick regime $\beta^2<4\pi$

Given a field configuration $\phi\in\mathcal{H}$, we denote by

$$\psi_\xi:=\psi_{(\sigma,x)}:=\exp(i\beta\sigma\phi(x))\qquad\text{where}\quad \xi:=(\sigma,x)\in\{-1,1\}\times\mathbb{R}^2$$

a vertex field with charge $\sigma\in\{-1,1\}$ at position $x\in\mathbb{R}^2$. Given a finite index set $I\subseteq\mathbb{N}$, we denote by $\sigma_I:=(\sigma_i)_{i\in I}$, $x_I:=(x_i)_{i\in I}$, $\xi_I:=(\sigma_I,x_I)$ and define

$$\chi(\xi_I):=\sum_{j\in I}\sigma_j$$

the total charge of the ensemble $\xi_I$. We write $1{:}\ell:=\{1,\ldots,\ell\}$. The representation (2.28) for the potential $V^{[\ell],T}$ is, at first, the Mayer expansion ansatz

$$\begin{aligned}V_t^{[\ell],T}(\phi) &:= \int v_{t,T}^{[\ell]}(\xi_{1:\ell})\psi_{1:\ell}\rho(x_{1:\ell})\mathrm{d}\xi_{1:\ell}\\ &:= \sum_{\sigma_{1:\ell}\in\{-1,1\}^\ell}\int_{(\mathbb{R}^2)^\ell} v_{t,T}^{[\ell]}(x_{1:\ell},\sigma_{1:\ell})\prod_{j=1}^{\ell}e^{i\beta\sigma_j\phi(x_j)}\prod_{j=1}^{\ell}\rho(x_j)\mathrm{d}x_{1:\ell}.\end{aligned} \tag{3.1}$$

where $\rho\in\mathcal{C}^\infty(\mathbb{R}^2)$, $\rho(x)\in[0,1]$ serves as a volume cutoff. The UV cutoff $T$ in the notation of the kernels $v_{t,T}^{[\ell]}$ will be frequently suppressed when it does not lead to confusion. Unless explicitly stated otherwise, all estimates will be uniform in $\rho$ and $T$.

Plugging the above ansatz into the Polchinski equation (2.20), by the Duhamel principle we obtain a mild equation for the kernels of the representation. That is, denoting by

$$\dot{\mathbb{G}}_t(\xi_{I_1}, \xi_{I_2}) := -\frac{\beta^2}{2} \sum_{j \in I_1} \sum_{k \in I_2} \sigma_j \sigma_k \dot{\mathcal{G}}_t(x_j - x_k), \qquad \text{and} \qquad \dot{\mathbb{G}}_t(\xi_I) := \dot{\mathbb{G}}_t(\xi_I, \xi_I), \tag{3.2}$$

for index sets $I_1, I_2, I \subseteq \mathbb{N}$, we have the following simple

**Lemma 3.1.** *Suppose that the kernels $v_{t,T}^{[\ell]}(\xi_{1:\ell})$ satisfy the equation*

$$\begin{aligned} v_{t,T}^{[\ell]}(\xi_{1:\ell}) &= e^{\mathbb{G}_{t,T}(\xi_{1:\ell})} v_{T,T}^{[\ell]}(\xi_{1:\ell}) \\ &- \int_t^T e^{\mathbb{G}_{t,s}(\xi_{1:\ell})} \sum_{I_1 \uplus I_2 = 1:\ell} C_{|I_1|,|I_2|} v_{s,T}^{[|I_1|]}(\xi_{I_1}) \dot{\mathbb{G}}_s(\xi_{I_1}, \xi_{I_2}) v_{s,T}^{[|I_2|]}(\xi_{I_2}) \mathrm{d}s \end{aligned} \tag{3.3}$$

*for some positive combinatorial coefficients $C_{|I_1|,|I_2|}$, where the sum runs over all $I_1, I_2 \subseteq 1{:}\ell$ such that $I_1 \cap I_2 = \emptyset$ and $I_1 \cup I_2 = 1{:}\ell$. Then $V_t^{[\ell],T}$, given by (3.1), satisfies the Polchinski equation (2.20).*

**Proof.** A straightforward computation; see Appendix B.2. □

Our goal is now to carry out the program of Section 2.7, i.e. to propagate an inductive bound on the kernels through the equation (3.3) by a proper choice of boundary conditions. To do this, we have to define the norm $|||\cdot|||$ in which we estimate.

Given $x_I := \{x_i\}_{i \in I} \subseteq \mathbb{R}^d$ for a finite index set $I \subseteq \mathbb{N}$, denote the Steiner distance by

$$\mathrm{St}(x_I) := \min_{x_I \subseteq x_J} \min_{\tau \in c(x_J)} L(\tau)$$

where $L(\tau)$ is the (Euclidean) length of a tree $\tau$ on $\mathbb{R}^d$, and the second minimum runs over the set of spanning trees $c(x_J)$ of the points $x_J$, i.e. over all trees that connect all points of $x_J$. For a nonnegative function $\Theta$ that might depend on information about the kernel $v^{[\ell]}$, in particular $\ell$, and that we will determine later, we then define the weight

$$\omega_u^{\Theta}(x_I) := \exp(\Theta \sqrt{u}\, \mathrm{St}(x_I)) \tag{3.4}$$

and the associated Banach norm

$$|||v_t^{[\ell]}|||_{t,\Theta,\eta} := \sup_{\xi_1} \int |v_t^{[\ell]}(\xi_{1:\ell})| \omega_t^{\Theta}(\xi_{1:\ell}) \langle \mathrm{St}(x_{1:\ell}) \rangle^{\eta} \mathrm{d}\xi_{2:\ell}. \tag{3.5}$$

for some $\eta \geq 0$ to be specified later.

The Steiner weight $\omega_u^{\Theta}$ forces exponential localization of the kernels as we approach the UV; its role will become apparent later on. For the renormalised potential itself, we will need it for technical reasons when $\beta^2 \geq 4\pi$ in Section 3.2 below; for constructing the sine–Gordon model alone, we could have taken $\Theta := 0$ in the Wick regime $\beta^2 < 4\pi$. However, even in this regime, proving the existence of general observable products require UV localization of their germs (see Lemma 6.5 below), which in turn needed to be inherited from the potential sector.

The scale-independent weight $\langle \mathrm{St}(x_{1:\ell}) \rangle^{\eta}$ is required to control the infinite-volume limit for the observables. We will sometimes suppress the $\eta$ and $\Theta$-dependence in the notation when it does not lead to confusion, and also frequently take $\Theta := \varsigma \in \mathbb{R}_+$ constant.

It is seen inductively from (3.3) that the kernels $v_t^{[\ell]}$ are invariant under permutations under the index set $1{:}\ell$, hence in fact

$$|||v_t^{[\ell]}|||_{t,\Theta,\eta} = \sup_{\xi_j} \int |v_t^{[\ell]}(\xi_{1:\ell})| \omega_t^{\Theta}(\xi_{1:\ell}) \langle \mathrm{St}(x_{(1:\ell)}) \rangle^{\eta} \mathrm{d}\xi_{(1:\ell) \setminus \{j\}}$$

for any $j \in 1{:}\ell$. The following simple lemma shows that these norms are compatible with the structure of the flow equation (3.3).

**Lemma 3.2.** *Let $\varsigma, \eta \in \mathbb{R}_+$ be arbitrary. Then*

$$|||v_t^{[|I_1|]}(\xi_{I_1})\dot{\mathbb{G}}_t(\xi_{I_1},\xi_{I_2})v_t^{[|I_2|]}(\xi_{I_2})|||_{t,\varsigma,\eta} \lesssim_\beta |||v_t^{[|I_1|]}|||_{t,\varsigma,\eta}\, |||v_t^{[|I_2|]}|||_{t,\varsigma,\eta}\, \langle t\rangle^{-2} \tag{3.6}$$

**Proof.** Indeed, given any trees $\tau_i$ connecting $x_{I_i}$ and $\tau$ an edge between a point in $x_{I_1}$ and $x_{I_2}$, we have that $\tau_1 \cup \tau_2 \cup \tau$ connects $x_{I_1 \uplus I_2}$, thus

$$\mathrm{St}(x_{I_1\uplus I_2}) = \min_{x_{I_1\uplus I_2}\subseteq x_J} \min_{\tau\in c(x_J)} L(\tau) \leq \min_{x_{I_1}\subseteq x_J} \min_{\tau\in c(x_J)} L(\tau) + \min_{x_{I_2}\subseteq x_J} \min_{\tau\in c(x_J)} L(\tau) + \tilde{d}(x_{I_1},x_{I_2})$$

where $\tilde{d}(x_{I_1},x_{I_2}) := \inf_{(a,b)\in x_{I_1}\times x_{I_2}} |a-b|_{\ell^2(\mathbb{R}^d)}$. It then follows from the triangle inequality that

$$\begin{array}{rcl} \omega_t^\varsigma(\xi_{I_1\uplus I_2}) & \leq & \omega_t^\varsigma(\xi_{I_1})\omega_t^\varsigma(\xi_{I_2})e^{\varsigma\sqrt{t}\,\tilde{d}(x_{I_1},x_{I_2})} \\ \langle \mathrm{St}(x_{I_1\uplus I_2})\rangle^\eta & \lesssim_\eta & \langle \mathrm{St}(x_{I_1})\rangle^\eta \langle \mathrm{St}(x_{I_2})\rangle^\eta \langle \tilde{d}(x_{I_1},x_{I_2})\rangle^\eta. \end{array} \tag{3.7}$$

Then Young's convolution inequality gives

$$\begin{array}{rl} & |||v_t^{[|I_1|]}(\xi_{I_1})\dot{\mathbb{G}}_t(\xi_{I_1},\xi_{I_2})v_t^{[|I_2|]}(\xi_{I_2})|||_{t,\varsigma,\eta} \\ \leq & \sup\limits_{\xi_1} \displaystyle\int \left| v_t^{[|I_1|]}(\xi_{I_1})\omega_t^\varsigma(\xi_{I_1})\dot{\mathbb{G}}_t(\xi_{I_1},\xi_{I_2})e^{\varsigma\sqrt{t}\,\tilde{d}(x_{I_1},x_{I_2})}v_t^{[|I_2|]}(\xi_{I_2})\omega_t^\varsigma(\xi_{I_2})\langle \mathrm{St}(x_{1:\ell})\rangle^\eta \right| \mathrm{d}\xi_{2:\ell} \\ \lesssim & |||v_t^{[|I_1|]}|||_{t,\varsigma,\eta}|||v_t^{[|I_2|]}|||_{t,\varsigma,\eta} \sup\limits_y \displaystyle\int_{\mathbb{R}^2} |\dot{\mathcal{G}}_t(x-y)e^{\varsigma\sqrt{t}|x-y|}\langle x-y\rangle^\eta|\mathrm{d}x. \end{array}$$

We then conclude from Lemma A.2. □

We impose the boundary conditions

$$\begin{array}{rcl} v_{0,T}^{[1]} & := & \dfrac{\lambda}{2} \\ v_{T,T}^{[\ell]} & := & 0 \quad \text{for } \ell \geq 2. \end{array} \tag{3.8}$$

The motivation for the $\ell = 1$ boundary condition is that at resolution $t = 0$, inserting it into (3.1) we obtain

$$V_0^{[1],T}(\phi) = \sum_{\sigma_1\in\{-1,1\}} \int_{\mathbb{R}^2} \frac{\lambda}{2} e^{i\beta\sigma_1\phi(x_1)}\rho(x_1)\mathrm{d}x_1 = \lambda\int_{\mathbb{R}^2} \cos(\beta\phi(x))\rho(x)\mathrm{d}x$$

the expected form of the sine–Gordon potential with coupling $\lambda$. The boundary conditions at $\ell \geq 2$ follow from the naive power counting argument of Section 2.7, see the proof of the following

**Lemma 3.3.** *Let $\beta^2 < 4\pi \Leftrightarrow \delta > \frac{1}{2}$. Take $\varsigma \in \mathbb{R}_+$ arbitrary. With the boundary conditions (3.8), the system of equations (3.3) admits a solution satisfying*

$$|||v_{t,T}^{[\ell]}|||_{t,\varsigma} \lesssim \lambda_t(\lambda_t\langle t\rangle^{-1})^{\ell-1} \asymp \langle t\rangle^{1-\ell\delta} \tag{3.9}$$

*where $\lambda_t := \frac{|\lambda|}{2} e^{\frac{\beta^2}{2}\mathcal{G}_t(0)}$.*

**Proof.** For $\ell$=1, since $v_{\cdot,T}^{[0]} \equiv 0$, the source term of (3.3) vanishes, and plugging in $t = 0$ we obtain

$$\frac{\lambda}{2} = v_{0,T}^{[1]}(\xi_1) = v_{T,T}^{[1]}(\xi_1)e^{\mathbb{G}_{0,T}(\xi_1)} = v_{T,T}^{[1]}(\xi_1)e^{-\frac{\beta^2}{2}\mathcal{G}_{0,T}(0)}$$

which yields $v_{T,T}^{[1]}(\xi_1) = \frac{\lambda}{2} e^{\frac{\beta^2}{2}\mathcal{G}_{0,T}(0)}$; plugging it back to (3.3) gives

$$v_{t,T}^{[1]}(\xi_1) = v_{T,T}^{[1]}(\xi_1)e^{\mathbb{G}_{t,T}(\xi_1)} = \frac{\lambda}{2} e^{\frac{\beta^2}{2}\mathcal{G}_{0,T}(0)} e^{-\frac{\beta^2}{2}\mathcal{G}_{t,T}(0)} = \frac{\lambda}{2} e^{\frac{\beta^2}{2}\mathcal{G}_{0,t}(0)}.$$

It satisfies, for any $\varsigma \geq 0$

$$|||v_{t,T}^{[1]}(\xi_1)|||_{t,\varsigma} = \lambda_t \sup_{x_1} |e^{\varsigma\sqrt{t}0}| = \lambda_t \asymp e^{\frac{\beta^2}{8\pi}\log(1\vee t)} \asymp \langle t\rangle^{1-\delta},$$

verifying the bound (3.9) at $\ell = 1$. Hence from now on we assume $\ell \geq 2$. With the irrelevant boundary condition $v_{T,T}^{[\ell]} = 0$, we estimate (3.3) via Lemma 3.2 and the inductive assumption by

$$\begin{aligned}
|||v_t^{[\ell]}|||_{t,\varsigma} &\lesssim \int_t^\infty \|e^{\mathbb{G}_{t,s}(\xi_{1:\ell})}\|_{L^\infty} \sum_{I_1 \uplus I_2 = 1:\ell} |||v_s^{[|I_1|]}(\xi_{I_1}) \dot{\mathbb{G}}_s(\xi_{I_1}, \xi_{I_2}) v_s^{[|I_2|]}(\xi_{I_2})|||_{t,\varsigma} \mathrm{d}s \\
&\leq \int_t^\infty \|e^{\mathbb{G}_{t,s}(\xi_{1:\ell})}\|_{L^\infty} \sum_{I_1 \uplus I_2 = 1:\ell} |||v_s^{[|I_1|]}(\xi_{I_1}) \dot{\mathbb{G}}_s(\xi_{I_1}, \xi_{I_2}) v_s^{[|I_2|]}(\xi_{I_2})|||_{s,\varsigma} \mathrm{d}s \\
&\lesssim \int_t^\infty \|e^{\mathbb{G}_{t,s}(\xi_{1:\ell})}\|_{L^\infty} \sum_{I_1 \uplus I_2 = 1:\ell} |||v_s^{[|I_1|]}|||_{s,\varsigma} \langle s\rangle^{-2} |||v_s^{[|I_2|]}|||_{s,\varsigma} \mathrm{d}s \\
&\lesssim \int_t^\infty \|e^{\mathbb{G}_{t,s}(\xi_{1:\ell})}\|_{L^\infty} \lambda_s (\lambda_s \langle s\rangle^{-1})^{|I_1|-1} \langle s\rangle^{-2} \lambda_s (\lambda_s \langle s\rangle^{-1})^{|I_2|-1} \mathrm{d}s \\
&\lesssim \int_t^\infty \|e^{\mathbb{G}_{t,s}(\xi_{1:\ell})}\|_{L^\infty} \lambda_s (\lambda_s \langle s\rangle^{-1})^{\ell-1} \langle s\rangle^{-1} \mathrm{d}s.
\end{aligned}$$

where we used that $t \leq s$ to bound the norms. Now as $x \mapsto \dot{\mathcal{G}}_u(x)$ is a positive definite function, with $t \leq s$ it follows that $\|e^{\mathbb{G}_{t,s}(\xi_{1:\ell})}\|_{L^\infty} \lesssim 1$. Moreover, the assumptions $\delta > 1/2$ and $\ell \geq 2$ combined with

$$\lambda_s^\ell \langle s\rangle^{-\ell} \asymp \langle s\rangle^{-\ell\delta}$$

make the above integral convergent, verifying (3.9). □

### 3.2 Extension to $\beta^2 < 6\pi$

When $\beta^2 \geq 4\pi$, the proof of Lemma 3.3 breaks down. Based on the discussion of Section 2.7, one might naively try to change boundary conditions when the backward integral diverges and flow the corresponding kernels *forward*. Ultimately this will turn out to be the solution, however the implementation requires us to be careful and first make two remarks about the structure of the flow equation. First, we note that some kernels are in fact irrelevant in the entire subcritical regime:

**Example 3.4. (2-point case)** Even if $\beta^2 \geq 4\pi$, the $\ell = 1$ case of Lemma 3.3 goes through without modification, so that in the $\ell = 2$ case (3.3) reads, up to constants and the volume cutoff $\rho$, as

$$v_{t,T}^{[2]}(\xi_{1:2}) \asymp v_{T,T}^{[2]}(\xi_{1:2}) e^{\mathbb{G}_{t,T}(\xi_{1:2})} - \int_t^T e^{\mathbb{G}_{t,s}(\xi_{1:2})} \lambda_s^2 \dot{\mathcal{G}}_s(x_1 - x_2) \mathrm{d}s. \tag{3.10}$$

Note that $\mathbb{G}_{t,s}(\xi_{1:2})$ depends on the total charge $\chi(\xi_{1:2}) := \sigma_1 + \sigma_2 \in \{-2, 0, 2\}$ of the contribution:

$$\mathbb{G}_{t,s}(\xi_{1:2}) := -\beta^2(\mathcal{G}_{t,s}(0) + (\mathbb{1}_{\chi(\xi_{1:2}) \neq 0} - \mathbb{1}_{\chi(\xi_{1:2}) = 0}) \mathcal{G}_{t,s}(x_1 - x_2)).$$

For the charged contributions $|\chi(\xi_{1:2})| = 2$, the particles repel each other and the supremum is attained as $x_1 - x_2 \to \infty$ (recall that $t \leq s$ and $x \mapsto \dot{\mathcal{G}}_u(x)$ is positive):

$$\sup_{x_{1:2}} -\beta^2 \mathcal{G}_{t,s}(0) - \beta^2 \mathcal{G}_{t,s}(x_1 - x_2) \leq -\beta^2 \mathcal{G}_{t,s}(0).$$

On the other hand, in the neutral case $|\chi(\xi_{1:2})| = 0$ the particles attract each other and form a dipole, cancelling the $-\beta^2 \mathcal{G}_{t,s}(0)$ gain:

$$\sup_{x_{1:2}} -\beta^2 \mathcal{G}_{t,s}(0) + \beta^2 \mathcal{G}_{t,s}(x_1 - x_2) \leq 0.$$

Thus

$$\sup_{x_{1:2}} e^{\mathbb{G}_{t,s}(\xi_{1:2})} \lesssim (\lambda_t \lambda_s^{-1})^{|\chi(\xi_{1:2})|}. \tag{3.11}$$

Hence for charged contributions we gain from the covariance matrix $\mathbb{G}_{t,s}(\xi_{1:2})$ compared to the naive bound $\|e^{\mathbb{G}_{t,s}(\xi_{1:2})}\|_{L^\infty} \lesssim 1$ we had in the previous proof and this gain is sufficient to restore integrability and flow backward. In the neutral case, the naive bound was in fact sharp, the integral in (3.10) is divergent and the kernel is relevant.

The following analysis on the screening mechanism of $\mathbb{G}_{t,s}$ generalizes the estimate (3.11) to the $\ell$-point case. Since it serves as a base point to the more general case of the external vertex field insertion of Section 5.2, we decide to include a full self-contained proof.

**Lemma 3.5. (Individual screening)** *Let $t \leq s$. Then there exists a constant $C$ such that for all $x_g, x_b, z \in \mathbb{R}^2$ with $|x_g - x_b| \leq |x_b - z|$ we have*

$$\mathcal{G}_{t,s}(x_b - z) - \mathcal{G}_{t,s}(x_g - z) \leq C.$$

**Proof.** See Appendix A.2. □

We can bootstrap the above estimate by recursively extracting neutral pairs from $\mathbb{G}_{t,s}$.

**Lemma 3.6. (Screening bound on $\mathbb{G}_{t,s}$)** *Let $t \leq s$. Then*

$$\sup_{x_{1:\ell}} e^{\mathbb{G}_{t,s}(\xi_{1:\ell})} \lesssim_\ell \lambda_t^{|\chi(\xi_{1:\ell})|} \lambda_s^{-|\chi(\xi_{1:\ell})|}. \tag{3.12}$$

**Proof.** For simplicity denote $\chi := \chi(\xi_{1:\ell})$. First suppose the ensemble is purely charged, i.e. $|\chi| = \ell$. Then $\sigma_j \sigma_k = 1$ for all $j, k \in 1{:}\,\ell$, so

$$e^{\mathbb{G}_{t,s}(\xi_{1:\ell})} = \exp\left(-\ell \frac{\beta^2}{2} \mathcal{G}_{t,s}(0) - \beta^2 \sum_{j<k} \mathcal{G}_{t,s}(x_j - x_k)\right) \leq \exp\left(-\ell \frac{\beta^2}{2} \mathcal{G}_{t,s}(0)\right) = \lambda_t^\ell \lambda_s^{-\ell}, \tag{3.13}$$

verifying (3.12). Suppose from now on that $|\chi| \leq \ell - 1$; we will show that we can remove neutral clusters from $\mathbb{G}_{t,s}(\xi_{1:\ell})$ until we arrive at a purely charged subset of particles $\tilde{\xi}_{1:|\chi|} \in (\{-1, 1\} \times \mathbb{R}^2)^{|\chi|}$ and so that

$$e^{\mathbb{G}_{t,s}(\xi_{1:\ell})} \lesssim e^{\mathbb{G}_{t,s}(\tilde{\xi}_{1:|\chi|})}; \tag{3.14}$$

the statement then follows from applying (3.13) to the purely charged sector $\tilde{\xi}_{1:|\chi|}$.

Since $|\chi| \leq \ell - 1$, there is at least one neutral pair, i.e. $j_+, j_- \in 1{:}\,\ell$ such that $\sigma_{j_+} \sigma_{j_-} = -1$. Denote by $k_+, k_-$ the index of a closest of such neutral pairs, i.e. a pair $\eta := \{(+1, x_{k_+}), (-1, x_{k_-})\} \in (\{-1, 1\} \times \mathbb{R}^2)^2$ such that

$$|x_{k_+} - x_{k_-}| = \min_{j_+, j_- : \sigma_{j_+} \sigma_{j_-} = -1} |x_{j_+} - x_{j_-}|.$$

By possibly relabeling the points, we distinguish the purely charged core and the ensemble of neutral pairs with $\eta$ being the last, i.e. we assume $\sigma_1, \ldots, \sigma_{|\chi|} = \mathrm{sgn}(\chi)$ and $\sigma_{|\chi|+j} = -\sigma_{|\chi|+j+1}$ for all $j \in 1, \ldots, \ell - |\chi| - 1$ and $k_- \neq k_+ \in \{\ell - 1, \ell\}$ The covariance matrix factorizes accordingly:

$$e^{\mathbb{G}_{t,s}(\xi_{1:\ell})} = e^{\mathbb{G}_{t,s}(\xi_{1:(\ell-2)})} e^{2\mathbb{G}_{t,s}(\xi_{1:(\ell-2)}, \eta)} e^{\mathbb{G}_{t,s}(\eta)}.$$

By positive definiteness of $\dot{\mathcal{G}}_u$ we have $e^{\mathbb{G}_{t,s}(\eta)} \lesssim 1$. By construction, for all $j \in 1{:}\,(\ell - 2)$

$$|x_{k_+} - x_{k_-}| \leq \begin{cases} |x_{k_-} - x_j| & \text{if } \sigma_j = 1 \\ |x_{k_+} - x_j| & \text{if } \sigma_j = -1 \end{cases}.$$

Thus we may bound the cross term by

$$2\mathbb{G}_{t,s}(\xi_{1:(\ell-2)},\eta) = -\beta^2\sum_{j=1}^{\ell-2}\sigma_j(\mathcal{G}_{t,s}(x_{k_+}-x_j)-\mathcal{G}_{t,s}(x_{k_-}-x_j)) \lesssim_\ell C$$

by applying Lemma 3.5 with the choice

$$\begin{array}{lll} x_b := x_{k_-}, & x_g := x_{k_+} & z := x_j \ \text{ if } \sigma_j = 1 \\ x_b := x_{k_+} & x_g := x_{k_-} & z := x_j \ \text{ if } \sigma_j = -1. \end{array}$$

Hence we obtain $e^{\mathbb{G}_{t,s}(\xi_{1:\ell})} \lesssim e^{\mathbb{G}_{t,s}(\xi_{1:(\ell-2)})}$, and repeat the above extraction procedure inductively until only the purely charged core remains, yielding (3.14). □

The screening bound ensures that non-neutral configurations can be flown backward. Moving on to the neutral configurations, they are relevant and hence need to be flown forward, cf. Section 2.7. However, as the following example shows, naively changing the boundary condition does not work.

**Example 3.7. (The need for localisation)** Let $\beta^2 \geq 4\pi$. As we have seen in the previous example, the $\ell = 2$ neutral kernel is relevant. Trying to naively change the boundary condition and flow forward yields the equation

$$v_{t,T}^{[2]}(\xi_{1:2}) \asymp \int_0^t e^{\mathbb{G}_{t,s}(\xi_{1:2})}\lambda_s^2\dot{\mathcal{G}}_s(x_1-x_2)\mathrm{d}s.$$

It is clear that the Steiner weight in the norm $|||\cdot|||_{t,\varsigma}$ for $\varsigma > 0$ is incompatible with this equation. But even when $\varsigma = 0$ we face problems: now, as $s \leq t$, the covariance matrix behaves differently: now $\mathcal{G}_{t,s}$ is negative, so in the neutral case the particles escape to infinity and we have

$$\sup_{x_{1:2}} -\beta^2\mathcal{G}_{t,s}(0) + \beta^2\mathcal{G}_{t,s}(x_1-x_2) \leq -\beta^2\mathcal{G}_{t,s}(0),$$

thus $\|e^{\mathbb{G}_{t,s}(\xi_{1:2})}\mathbb{1}_{|\chi(\xi_{1:2})=0|}\|_{L^\infty} \lesssim \lambda_t^2\lambda_s^{-2}$ and trying to flow forward

$$\begin{aligned} |||v_{t,T}^{[2]}(\xi_{1:2})\mathbb{1}_{\chi(\xi_{1:2})=0}|||_{t,0} &\lesssim \int_0^t |||e^{\mathbb{G}_{t,s}(\xi_{1:2})}\mathbb{1}_{\chi(\xi_{1:2})=0}\lambda_s^2\dot{\mathcal{G}}_s(x_1-x_2)|||_{t,0}\mathrm{d}s \\ &\lesssim \int_0^t \|e^{\mathbb{G}_{t,s}(\xi_{1:2})}\mathbb{1}_{\chi(\xi_{1:2})=0}\|_{L^\infty}\lambda_s^2\|\dot{\mathcal{G}}_s(x)\|_{L^1_\eta}\mathrm{d}s \\ &\lesssim \lambda_t^2\int_0^t \langle s\rangle^{-2}\mathrm{d}s \\ &\lesssim \lambda_t^2 \end{aligned}$$

as the integral converges. The bound $|||v_{t,T}^{[2]}(\xi_{1:2})\mathbb{1}_{\chi(\xi_{1:2})=0}|||_{t,0} \lesssim \lambda_t^2$ is worse than the expected inductive estimate $\lambda_t(\lambda_t\langle t\rangle^{-1})$. That is, naively trying to flow forward in the case when the backward integral diverges, the covariance matrix will cause the forward integral to converge, breaking the induction. The problem, as we noted above, is particles escaping to infinity; evaluating $e^{\mathbb{G}_{t,s}(\xi_{1:2})}$ at $x_1 = x_2 = 0$ gives

$$e^{\mathbb{G}_{t,s}(0)}\mathbb{1}_{\chi(\xi_{1:2})=0} = e^{-\beta^2\mathcal{G}_{t,s}(0)+\beta^2\mathcal{G}_{t,s}(0)} = 1$$

and thus formally replacing the $L^\infty$ bound on $e^{\mathbb{G}_{t,s}(\xi_{1:2})}$ with the above diagonal evaluation we obtain, with $\mathfrak{l}_t := (1+\log(1\vee t))$

$$\begin{aligned} \int_0^t e^{\mathbb{G}_{t,s}(0)}\mathbb{1}_{\chi(\xi_{1:2})=0}|||\lambda_s^2\dot{\mathcal{G}}_s(x_1-x_2)|||_{t,0}\mathrm{d}s &\lesssim \int_0^t \lambda_s^2\|\dot{\mathcal{G}}_s(x)\|_{L^1_\eta}\mathrm{d}s \\ &\lesssim \int_0^t \lambda_s^2\langle s\rangle^{-2}\mathrm{d}s \\ &\lesssim \lambda_t^2\langle t\rangle^{-1}\mathfrak{l}_t^{\mathfrak{c}} \end{aligned}$$

with for some $c\in\mathbb{R}_+$ when $\beta^2\geq 4\pi$, which is the expected inductive bound.

To sum up the above example, flowing forward requires forcing $\mathbb{G}_{t,s}(\xi_{1:\ell})$ to be evaluated on its diagonal. Implementing rigorously the idea of forcing the covariance matrix onto the diagonal in a way that also adheres to the weighted norms when $\varsigma>0$ requires the introduction of a localisation scheme.

The localisation scheme we employ is based on the generalisation of the Mayer expansion ansatz (3.1). Thus from now on, we work in the following generalised ansatz:

$$V_t^{[\ell],T}(\phi):=\sum_{\boldsymbol{a}\in\mathfrak{M}}\int v_{t,T}^{[\ell](\boldsymbol{a})}(\xi_{\boldsymbol{a}})\Psi_{\boldsymbol{a}}(\xi_{\boldsymbol{a}})\rho(x_{\boldsymbol{a}})\mathrm{d}\xi_{\boldsymbol{a}} \tag{3.15}$$

where $\Psi_{\boldsymbol{a}}(\xi_{\boldsymbol{a}})$ are functionals of the field configuration $\phi$ to be determined, called *field content*, and the index set $\mathfrak{M}$ is called the *basis of the representation*, again to be determined. We shall employ the formal notation $\langle v_{t,T}^{[\ell](\boldsymbol{a})},\Psi_{\boldsymbol{a}}\rangle$ for terms in the above sum.

The way we determine $\mathfrak{M}$ and $\Psi_{\boldsymbol{a}}(\xi_{\boldsymbol{a}})$ is inductive: we start from the $\ell=1$ case with assuming the naive Mayer ansatz (3.1) as before. Supposing (3.15) is already determined up to a fixed iteration order $\ell-1$, we compute the corresponding bilinear term on the right-hand side of (2.20) using $V^{[0]}\equiv 0$. The so-defined source term then undergoes a post-processing based on the relevance/irrelevance of the associated kernels, involving localisation when it is needed. The basis of the extended representation at level $\ell$ is then determined by adhering to the field contents of the possibly localised source term.

Let us carry out the above general strategy. Starting at the $\ell=1$ order case, we choose (3.15) to be the naive Mayer ansatz (3.1) as before; with the boundary condition $v_{0,T}^{[1]}:=\frac{\lambda}{2}$, we obtain, as in the first half of the proof of Lemma 3.3,

$$|\!|\!|v_{t,T}^{[1]}(\xi_1)|\!|\!|_{t,\varsigma}=\lambda_t\asymp\langle t\rangle^{1-\delta}$$

for any $\varsigma\geq 0$. Then (3.3) gives the $\ell=2$ equation

$$\partial_t\langle v_t^{[2]},\psi_{1:2}\rangle+\frac{1}{2}\Delta_{\dot{G}_t}\langle v_t^{[2]},\psi_{1:2}\rangle=\langle\mathbf{B}_t^{[2]},\psi_{1:2}\rangle \tag{3.16}$$

where we denoted

$$\begin{aligned}\langle\mathbf{B}_t^{[2]}(\xi_{1:2}),\psi_{1:2}\rangle &:= -\frac{\beta^2}{2}\sigma_1\sigma_2\int v_t^{[1]}(\xi_1)v_t^{[1]}(\xi_2)\rho(x_{1:2})\dot{\mathcal{G}}_t(x_1-x_2)\psi_{1:2}\mathrm{d}\xi_{1:2}\\ &\asymp \int\lambda_t^2\rho(x_{1:2})\dot{\mathcal{G}}_t(x_1-x_2)\psi_{1:2}\mathrm{d}\xi_{1:2}\end{aligned}$$

For charged terms we can send this backward via the screening estimates on the covariance matrix, hence in this case the post-processing is to just leave the source term as it is.

For neutral terms, as we have seen above, this is not the way to go, hence we need to “clear” the source term $\langle\mathbf{B}^{[2]},\psi_{1:2}\rangle$ from the right-hand side of (3.16). We achieve this by decomposing it into a local part that will be sent forward, and a more regular remainder that can be flown backward. Concretely, for $\sigma_1=-\sigma_2$, let

$$\psi_{1:2}=\frac{\psi_{1:2}-1}{|x_1-x_2|^{\varkappa}}|x_1-x_2|^{\varkappa}+1=:\psi_{1:2}^{\varkappa}|x_1-x_2|^{\varkappa}+1$$

where $\varkappa>0$ to be determined. Based on this split, the source term decomposes as

$$\begin{aligned}\langle\mathbf{B}_t^{[2]},\psi_{1:2}\rangle &= \frac{\beta^2}{2}\int\lambda_t^2\rho(x_{1:2})\dot{\mathcal{G}}_t(x_1-x_2)\mathrm{d}x_{1:2}+\frac{\beta^2}{2}\int\lambda_t^2\rho(x_{1:2})\psi_{1:2}^{\varkappa}|x_1-x_2|^{\varkappa}\dot{\mathcal{G}}_t(x_1-x_2)\mathrm{d}x_{1:2}\\ &=: \langle\mathbf{B}_t^{[2](0)},1\rangle+\langle\mathbf{B}_t^{[2](\varkappa)}(x_{1:2}),\psi_{1:2}^{\varkappa}\rangle.\end{aligned}$$

Note that $\mathbf{B}_t^{[2](\varkappa)}$ is a function on $(\mathbb{R}^2)^2$, while $\mathbf{B}_t^{[2](0)}$ is a number.

The above localisation thus left us with nonzero source terms corresponding to the field contents $1, \psi_{1:2}^{\varkappa}$. Thus (3.16) is now replaced with the following hierarchy of equations:

$$
\begin{array}{lclcl}
\partial_t \langle v_t^{[2](0)}, 1 \rangle & + & \frac{1}{2}\Delta_{\dot{G}_t} \langle v_t^{[2](0)}, 1 \rangle & = & \mathbf{B}_t^{[2](0)} \\
\partial_t \langle v_t^{[2]}, \psi_{1:2} \rangle & + & \frac{1}{2}\Delta_{\dot{G}_t} \langle v_t^{[2]}, \psi_{1:2} \rangle & = & 0 \\
\partial_t \langle v_t^{[2](\varkappa)}, \psi_{1:2}^{\varkappa} \rangle & + & \frac{1}{2}\Delta_{\dot{G}_t} \langle v_t^{[2](\varkappa)}, \psi_{1:2}^{\varkappa} \rangle & = & \mathbf{B}_t^{[2](\varkappa)}
\end{array}
$$

The point of the localisation was that we eliminated the source term from the equation of the naive term $\langle v_t^{[2]}, \psi_{1:2} \rangle$, which we saw we cannot estimate in the neutral case. To ensure that it stays this way, i.e. that the action of the functional Laplacian $\Delta_{\dot{G}_t}$ on $\psi_{1:2}^{\varkappa}$ does not generate contributions with field content $\psi_{1:2}$ (which would result in the source term of the form $f(v_t^{[2],(\varkappa)})$ in the equation of $v_t^{[2]}$, reintroducing the problem we encountered in the previous example), we proceed with an abstract integration by parts argument: compute

$$
\begin{aligned}
\frac{1}{2}\Delta_{\dot{G}_t} \langle v^{[2](\varkappa)}, \psi_{1:2}^{\varkappa} \rangle &= \frac{1}{2}\Delta_{\dot{G}_t} \left\langle v^{[2](\varkappa)}, \frac{\psi_{1:2} - 1}{|x_1 - x_2|^{\varkappa}} \right\rangle \\
&= \frac{1}{2}\Delta_{\dot{G}_t} \langle v^{[2](\varkappa)} |x_1 - x_2|^{-\varkappa}, \psi_{1:2} - 1 \rangle \\
&\asymp \langle \dot{\mathbb{G}}_t(\xi_{1:2}) v^{[2](\varkappa)} |x_1 - x_2|^{-\varkappa}, \psi_{1:2} \rangle \\
&= \langle \dot{\mathbb{G}}_t(\xi_{1:2}) v^{[2](\varkappa)} |x_1 - x_2|^{-\varkappa}, \psi_{1:2} - 1 + 1 \rangle \\
&= \langle \dot{\mathbb{G}}_t(\xi_{1:2}) v^{[2](\varkappa)}, \psi_{1:2}^{\varkappa} \rangle + \langle \dot{\mathbb{G}}_t(\xi_{1:2}) v^{[2](\varkappa)} |x_1 - x_2|^{-\varkappa}, 1 \rangle.
\end{aligned}
$$

Hence denoting by $\mathbf{H}_t := \partial_t + \dot{\mathbb{G}}_t(\xi_{1:2})$ and

$$
\mathbf{A}(v_t^{[2](\varkappa)}) := \langle \dot{\mathbb{G}}_t(\xi_{1:2}) v_t^{[2](\varkappa)} |x_1 - x_2|^{-\varkappa}, 1 \rangle := \int \dot{\mathbb{G}}_t(\xi_{1:2}) v_t^{[2](\varkappa)}(\xi_{1:2}) |x_1 - x_2|^{-\varkappa} \rho(x_{1:2}) \mathrm{d}x_{1:2},
$$

the $\ell = 2$ equation (3.16) in the naive representation is replaced with the following hierarchy of equations for the kernels, labelled by the field contents on the left:

$$
\chi(\xi_{1:2}) = 0: \quad
\begin{array}{c|lcll}
1 & \partial_t v_t^{[2](0)} & + & \mathbf{A}(v_t^{[2](\varkappa)}) & = \mathbf{B}_t^{[2](0)} \\
\psi_{1:2} & \mathbf{H}_t v_t^{[2]} & & & = 0 \\
\psi_{1:2}^{\varkappa} & \mathbf{H}_t v_t^{[2](\varkappa)} & & & = \mathbf{B}_t^{[2](\varkappa)}
\end{array}
\tag{3.17}
$$

Crucially, as a consequence of our post-processing, the kernel of the naive field content $\psi_{1:2}$ satisfies a linear ODE with vanishing source term and zero boundary conditions, hence is identically zero.

The system (3.17) is the result of the post-processing in the neutral case. For the sake of direct comparison, we formulate the result of the above-discussed trivial post-processing in the charged case (cf. (3.16)):

$$
\chi(\xi_{1:2}) \neq 0: \quad
\begin{array}{c|lcl}
1 & \partial_t v_t^{[2](0)} & = & 0 \\
\psi_{1:2} & \mathbf{H}_t v_t^{[2]} & = & \mathbf{B}_t^{[2]} \\
\psi_{1:2}^{\varkappa} & \mathbf{H}_t v_t^{[2](\varkappa)} & = & 0
\end{array}
\tag{3.18}
$$

We summarize the obtained representation in the following

**Lemma 3.8. (Change of basis)** *Let $V_t^{[1]}$ be given in terms of (3.1), and let*

$$
V_{t,T}^{[2]}(\phi) = \sum_{\sigma_{1:2} \in \{-1,1\}^2} \langle v_{t,T}^{[2](0)}(\xi_{1:2}), 1 \rangle + \langle v_{t,T}^{[2]}(\xi_{1:2}), \psi_{1:2} \rangle + \langle v_{t,T}^{[2](\varkappa)}(\xi_{1:2}), \psi_{1:2}^{\varkappa} \rangle
$$

*with kernels satisfying* (3.17) *when* $\chi(\xi_{1:2})=0$ *and* (3.18) *when* $\chi(\xi_{1:2})\neq 0$. *Then* $\{V_{t,T}^{[\ell]}\}_{\ell\in\{1,2\}}$ *satisfies the Polchinski equation* (2.20) *up to level* $\ell=2$.

In this new basis, circumventing the issues highlighted in Example 3.7, we can propagate the correct inductive estimates and thus extend the result of Lemma 3.3 to the case $\beta^2\geq 4\pi$:

**Lemma 3.9.** *Let* $\beta^2\in(0,8\pi)$, $\varkappa\in(0,2]$ *such that* $\varkappa>\frac{\beta^2-4\pi}{2\pi}=2(1-2\delta)$, *take* $\varsigma,\eta\in\mathbb{R}_+$ *arbitrary. Then imposing the boundary conditions*

$$v_{T,T}^{[2]}:=0,\qquad v_{T,T}^{[2](\varkappa)}:=0,\qquad \text{and}\qquad \begin{array}{rcl} v_{T,T}^{[2](0)} & := & 0 \text{ if } \beta^2<4\pi\\ v_{0,T}^{[2](0)} & := & 0 \text{ if } \beta^2\geq 4\pi\end{array} \tag{3.19}$$

*on the system of equations* (3.17), (3.18), *we have, with* $\mathfrak{l}_t:=(1+\log(1\vee t))$ *and some* $c\in\mathbb{R}_+$, *the upper bounds*

| | $\vert\vert\vert v_{t,T}^{[2](0)}\vert\vert\vert_{t,\varsigma}$ | $\vert\vert\vert v_{t,T}^{[2]}\vert\vert\vert_{t,\varsigma}$ | $\vert\vert\vert v_{t,T}^{[2](\varkappa)}\vert\vert\vert_{t,\varsigma}$ |
|---|---|---|---|
| $\chi(\sigma_{1:2})\neq 0$ | $0$ | $\lambda_t(\lambda_t\langle t\rangle^{-1})\mathfrak{l}_t^c$ | $0$ |
| $\chi(\sigma_{1:2})=0$ | $\Vert\rho\Vert_{L^1}\lambda_t(\lambda_t\langle t\rangle^{-1})\mathfrak{l}_t^c$ | $0$ | $\lambda_t(\lambda_t\langle t\rangle^{-1})\mathfrak{l}_t^c\langle t\rangle^{-\frac{\varkappa}{2}}.$ |

*Moreover, for any* $\alpha\geq 0$ *we have*

$$\vert\vert\vert v_{t,T}^{[2](\varkappa)}(\xi_{1:2})|x_1-x_2|^{\alpha}\vert\vert\vert_{t,\varsigma}\lesssim\lambda_t(\lambda_t\langle t\rangle^{-1})\mathfrak{l}_t^c\langle t\rangle^{-\frac{\varkappa+\alpha}{2}}. \tag{3.20}$$

**Proof.** Whenever a kernel is the solution of a linear ODE with vanishing source term, it vanishes identically.

In the charged case $\chi(\sigma_{1:2})\neq 0$, we estimate similarly to Lemma 3.3, the only difference being the reliance on the screening bound (3.12): indeed, with $|\chi(\sigma_{1:2})|=2$

$$\vert\vert\vert v_{t,T}^{[2]}\vert\vert\vert_{t,\varsigma}\leq\int_t^{\infty}\Vert e^{\mathbb{G}_{t,s}(\xi_{1:2})}\Vert_{L^\infty}\vert\vert\vert v_{s,T}^{[1]}\vert\vert\vert_{s,\varsigma}^2\langle s\rangle^{-2}\mathrm{d}s\lesssim\lambda_t^2\int_t^{\infty}\langle s\rangle^{-2}\mathrm{d}s\lesssim\lambda_t^2\langle t\rangle^{-1}.$$

In the neutral case $\chi(\sigma_{1:2})=0$, for the nonlocal remainder we use the gain from the smoothing kernel $|\cdot|^{\varkappa}$: by Lemma A.2,

$$\vert\vert\vert\mathbf{B}_t^{[2](\varkappa)}\vert\vert\vert_{t,\varsigma}\lesssim\lambda_t^2\Vert|\cdot|^{\varkappa}\dot{\mathcal{G}}_t(\cdot)e^{\varsigma\sqrt{t}|\cdot|}\Vert_{L_\eta^1}\lesssim\lambda_t^2\langle t\rangle^{-2}\langle t\rangle^{-\frac{\varkappa}{2}} \tag{3.21}$$

which allows to send the kernel backward; with $\Vert e^{\mathbb{G}_{t,s}(\xi_{1:2})}\Vert_{L^\infty}\lesssim 1$

$$\vert\vert\vert v_t^{[2](\varkappa)}\vert\vert\vert_{t,\varsigma}\lesssim\int_t^{\infty}\Vert e^{\mathbb{G}_{t,s}(\xi_{1:2})}\Vert_{L^\infty}\vert\vert\vert\mathbf{B}_s^{[2](\varkappa)}\vert\vert\vert_{s,\varsigma}\mathrm{d}s\lesssim\int_t^{\infty}\lambda_s^2\langle s\rangle^{-2}\langle s\rangle^{-\frac{\varkappa}{2}}\mathrm{d}s\lesssim\lambda_t^2\langle t\rangle^{-1}\langle t\rangle^{-\frac{\varkappa}{2}} \tag{3.22}$$

since the integral converges due to the condition on $\varkappa$: indeed, $2(1-\delta)-1-(1-2\delta)-\varepsilon=-\varepsilon$. The source term of the local part is

$$\begin{aligned}\vert\vert\vert\mathbf{B}_t^{[2](0)}\vert\vert\vert_{t,\varsigma} &= \left|\int_{(\mathbb{R}^2)^2}\lambda_t^2\rho(x_{1:2})\dot{\mathcal{G}}_t(x_1-x_2)\mathrm{d}x_{1:2}\right|\\ &\leq \lambda_t^2 t^{-1}e^{-\frac{m^2}{t}}\Vert\rho\Vert_{L^1}\sup_{x_1}\int_{\mathbb{R}^d}\rho(x_2)e^{-\frac{t|x_1-x_2|^2}{4}}\mathrm{d}x_2\\ &\lesssim \lambda_t^2\langle t\rangle^{-2}\Vert\rho\Vert_{L^1}.\end{aligned}$$

To estimate $\mathbf{A}(v_t^{[2](\varkappa)})$, note that in the neutral case

$$\dot{\mathbb{G}}_t(\xi_{1:2})\asymp_\beta-(\dot{\mathcal{G}}_t(0)-\dot{\mathcal{G}}_t(x_1-x_2))$$

so that with $1 \leq \omega_t^{\varsigma}$ it follows from Lemma A.3 that

$$\begin{aligned}
|||\mathbf{A}(v_t^{[2](\varkappa)})|||_{t,\varsigma} &= \left| \int_{(\mathbb{R}^2)^2} \frac{\dot{\mathbb{G}}_t(\xi_{1:2})}{|x_1 - x_2|^{\varkappa}} v_t^{[2](\varkappa)}(\xi_{1:2}) \rho(x_{1:2}) \mathrm{d}x_{1:2} \right| \\
&\lesssim \left\| \frac{\dot{\mathbb{G}}_t(\xi_{1:2})}{|x_1 - x_2|^{\varkappa}} \right\|_{L^\infty} \int_{(\mathbb{R}^2)^2} \left| v_t^{[2](\varkappa)}(\xi_{1:2}) \right| \rho(x_{1:2}) \omega_t^{\varsigma}(x_{1:2}) \langle \mathrm{St}(x_{1:2}) \rangle^{\eta} \mathrm{d}x_{1:2} \\
&\leq \left\| \frac{\dot{\mathbb{G}}_t(\xi_{1:2})}{|x_1 - x_2|^{\varkappa}} \right\|_{L^\infty} |||v_t^{[2](\varkappa)}|||_{t,\varsigma} \|\rho\|_{L^1} \\
&\lesssim \langle t \rangle^{-1+\frac{\varkappa}{2}} \lambda_t^2 \langle t \rangle^{-1} \langle t \rangle^{-\frac{\varkappa}{2}} \|\rho\|_{L^1}
\end{aligned}$$

Hence, since by locality $|v_t^{[2](0)}| = |||v_t^{[2](0)}|||_{t,\varsigma} = |||v_t^{[2](0)}|||_{s,\varsigma}$ for all $s, t \in \mathbb{R}_+$, when $\beta^2 < 4\pi \Leftrightarrow 1 - \delta < 1/2$ the irrelevant boundary condition yields

$$|||v_t^{[2](0)}|||_{t,\varsigma} \lesssim \int_t^{\infty} |||\mathbf{A}(v_s^{[2](\varkappa)})|||_{s,\varsigma} + |||\mathbf{B}_s^{[2](0)}|||_{s,\varsigma} \mathrm{d}s \lesssim \|\rho\|_{L^1} \int_t^{\infty} \lambda_s^2 \langle s \rangle^{-2} \mathrm{d}s \lesssim \|\rho\|_{L^1} \lambda_t^2 \langle t \rangle^{-1}$$

since the integral converges, while if $\beta^2 \geq 4\pi \Leftrightarrow 1 - \delta \geq 1/2$ the relevant boundary condition yields

$$|||v_t^{[2](0)}|||_{t,\varsigma} \lesssim \int_0^t |||\mathbf{A}(v_s^{[2](\varkappa)})|||_{s,\varsigma} + |||\mathbf{B}_s^{[2](0)}|||_{s,\varsigma} \mathrm{d}s \lesssim \|\rho\|_{L^1} \int_0^t \lambda_s^2 \langle s \rangle^{-2} \mathrm{d}s \lesssim \|\rho\|_{L^1} \lambda_t^2 \langle t \rangle^{-1} \mathfrak{l}_t$$

since the integral diverges. Finally, (3.20) follows from (3.21) and (3.22) by replacing $\varkappa$ with $\varkappa + \alpha$. □

**Remark 3.10.** Unifying the notation, with the above lemma, we have proven in the $\ell = 2$ case the inductive bound

$$|||v_{t,T}^{[2](\varkappa)}|||_{t,\varsigma,\eta} \lesssim_{\rho} \lambda_t (\lambda_t \langle t \rangle^{-1})^{\ell-1} \mathfrak{l}_t^c \langle t \rangle^{-\frac{\varkappa}{2}}, \qquad \varkappa \in (0, 2] \cup \{0\}$$

which holds even for $v^{[2]}$ if one insists on it corresponding to $\varkappa = 0$ by an abuse of notation (as we have seen, even though $v^{[2]}$ and its localisation $v^{[2](0)}$ satisfy the same estimates, it is crucial that we treat them as different kernels in order to prove said estimates).
Also note that the estimate depends on $\rho$ only in the case of the vacuum energy term $v^{[2](0)}$.

**Remark 3.11.** For future reference, we record that as $\dot{\mathbb{G}}_t(\xi_{1:2}) = -\beta^2(\dot{\mathcal{G}}_t(0) + \sigma_1\sigma_2\dot{\mathcal{G}}_t(x_1 - x_2))$,

$$v_{t,T}^{[2](\varkappa)}(\xi_{1:2}) = \frac{\beta^2}{2} \int_t^T e^{\mathbb{G}_{t,s}(\xi_{1:2})} \lambda_s^2 \rho(x_{1:2}) |x_1 - x_2|^{\varkappa} \sigma_1\sigma_2 \dot{\mathcal{G}}_s(x_1 - x_2) \mathrm{d}s$$

is symmetric under charge conjugation i.e. under $\sigma_{1:2} \leftrightarrow \overline{\sigma_{1:2}} := -\sigma_{1:2}$.

The previous two Lemmas tell us that in the the charged case $\chi(\sigma_{1:2}) \neq 0$ the field content is $\psi_{1:2}$, hence it propagates the field content $\psi_{1:\ell}$ further down the induction, i.e. for any $V^{[|I|]} = \langle v^{[|I|]}, \psi_{|I|} \rangle$, the bilinear term

$$\langle \mathrm{D}V_t^{[2]} \mathbb{1}_{\chi(\sigma_{1:2}) \neq 0}, \dot{G}_t \mathrm{D}V_t^{[|I|]} \rangle$$

will have associated field content $\psi_{|I|+2}$. On the other hand, the neutral term has the dipole field content $\{1, \psi_{1:2}^{\varkappa}\}$. Compute

$$\mathrm{D}\big(\langle v_t^{[2](0)}, 1 \rangle + \langle v_t^{[2](\varkappa)}, \psi_{1:2}^{\varkappa} \rangle\big) = v_t^{[2](\varkappa)} |x_1 - x_2|^{-\varkappa} \psi_{1:2} \sum_{j=1}^{2} \sigma_j \delta_{x_j} i\beta \tag{3.23}$$

so that for example

$$\langle \mathrm{D}V_t^{[2]} \mathbb{1}_{\chi(\sigma_{1:2}) = 0}, \dot{G}_t \mathrm{D}V_t^{[1]} \rangle \asymp \langle v_t^{[2](\varkappa)} \dot{\mathbb{G}}_t(\xi_{1:2}, \xi_3) v^{[1]}, |x_1 - x_2|^{-\varkappa} \psi_{1:2} \psi_3 \rangle.$$

Now a priori we have two options:

1. we try to absorb the divergence $|x_1-x_2|^{-\varkappa}$ into $\dot{\mathbb{G}}_t(\xi_{1:2},\xi_3)$, leaving us with the field content $\psi_{1:3}$.
2. create a new dipole, i.e. pick a $\tilde{\varkappa}$ so that

$$|x_1-x_2|^{-\varkappa}\psi_{1:2}\psi_3=|x_1-x_2|^{-\varkappa}(\psi_{1:2}-1+1)\psi_3=|x_1-x_2|^{\tilde{\varkappa}-\varkappa}\psi_{1:2}^{\tilde{\varkappa}}\psi_3+\frac{\psi_3}{|x_1-x_2|^{\varkappa}},$$

which yields the decomposition

$$\begin{aligned}\left\langle \mathrm{D}V_t^{[2]}\mathbb{1}_{|\chi(\sigma_{1:2})|=0},\dot{G}\mathrm{D}V_t^{[1]}\right\rangle &\asymp \left\langle v_t^{[2](\varkappa)}\dot{\mathbb{G}}_t(\xi_{1:2},\xi_3)v^{[1]},|x_1-x_2|^{-\varkappa}\psi_{1:2}\psi_3\right\rangle\\ &= \left\langle\int v_t^{[2](\varkappa)}(x_{1:2})\frac{\dot{\mathbb{G}}_t(\xi_{1:2},\xi_3)}{|x_1-x_2|^{\varkappa}}v^{[1]}(x_3)\mathrm{d}x_{1:2},\psi_3\right\rangle\\ &+ \left\langle v_t^{[2](\varkappa)}\frac{\dot{\mathbb{G}}_t(\xi_{1:2},\xi_3)}{|x_1-x_2|^{\varkappa-\tilde{\varkappa}}}v^{[1]},\psi_{1:2}^{\tilde{\varkappa}}\psi_3\right\rangle\\ &=: \langle\mathbf{B}_{\psi_3}(x_3),\psi_3\rangle+\langle\mathbf{B}_{\psi_{1:2}^{\tilde{\varkappa}}\psi_3}(x_{1:3}),\psi_{1:2}^{\tilde{\varkappa}}\psi_3\rangle.\end{aligned}\tag{3.24}$$

This induces a change of basis in the $\ell=3$, $|\chi(\sigma_{1:3})|=1$ case, which would then create further modifications along the induction in $\ell$, making the representation (3.15) differ significantly from the classical Mayer ansatz (3.1).

It turns out that in the case $\beta^2<6\pi$, since we may choose $\varkappa:=1$ (cf. Lemma 3.9), option 1. is viable: it is possible to negate the divergence $|x_1-x_2|^{-1}$ of the dipole, hence returning to the Mayer ansatz from $\ell\geq 3$ onwards. This is shown in the following

**Lemma 3.12.** *Let $\varkappa\in(0,1]$, $\varsigma,\eta\geq 0$, $c>\varsigma$. Suppose $\xi_{1:2}$ is neutral, i.e. $\sigma_1=-\sigma_2$. Then for any kernel $k_t(\xi_I)$, $I\subseteq 3{:}\ell$*

$$\left\|\left|v_t^{[2](\varkappa)}(\xi_{1:2})\frac{\dot{\mathbb{G}}_t(\xi_{1:2},\xi_I)}{|x_1-x_2|^{\varkappa}}k_t(\xi_I)\right\|\right|_{t,\varsigma}\lesssim |||v_t^{[2](\varkappa)}|||_{t,\varsigma+c}\langle t\rangle^{-2+\frac{\varkappa}{2}}|||k_t|||_{t,\varsigma}.\tag{3.25}$$

*Moreover, when $\xi_{3:4}$ is also neutral, i.e. $\sigma_3=-\sigma_4$,*

$$\left\|\left|v_t^{[2](\varkappa)}(\xi_{1:2})\frac{\dot{\mathbb{G}}_t(\xi_{1:2},\xi_{3:4})}{|x_1-x_2|^{\varkappa}|x_3-x_4|^{\varkappa}}v_t^{[2](\varkappa)}(\xi_{3:4})\right\|\right|_{t,\varsigma}\lesssim |||v_t^{[2](\varkappa)}|||_{t,\varsigma+c}\langle t\rangle^{-2+\varkappa}|||v_t^{[2](\varkappa)}|||_{t,\varsigma+c}\tag{3.26}$$

**Proof.** To prove (3.25), fix $j\in I$ and denote by $y:=x_j$. Then using Lemma A.6 with $|\mu|=0$ we have for any $c\geq 0$

$$\frac{|\dot{\mathcal{G}}_t(x_1-y)-\dot{\mathcal{G}}_t(x_2-y)|}{|x_1-x_2|^{\varkappa}}\lesssim|x_1-x_2|^{1-\varkappa}e^{-c\sqrt{t}|x_1-y|}\omega_t^c(x_{1:2})e^{-\frac{m^2}{t}}t^{-\frac{1}{2}}$$

thus writing $\delta_{ij}:=|x_i-x_j|$ it follows from (3.7) and $\tilde{d}(x_{1:2},x_I)\leq|x_1-y|$ that

$$\begin{aligned}&\left\|\left|v_t^{[2](\varkappa)}(\xi_{1:2})\frac{\dot{\mathbb{G}}_t(\xi_{1:2},\xi_I)}{|x_1-x_2|^{\varkappa}}k_t(\xi_I)\right\|\right|_{t,\varsigma}\\ \lesssim\ &\sup_{x_1}\int\left|v_t^{[2](\varkappa)}(\xi_{1:2})\omega_t^{\varsigma}(\xi_{1:2})\right|\sum_{j\in I}\frac{|\dot{\mathcal{G}}_t(x_1-x_j)-\dot{\mathcal{G}}_t(x_2-x_j)|}{|x_1-x_2|^{\varkappa}}e^{\varsigma\sqrt{t}\tilde{d}(x_{1:2},x_I)}|k_t(\xi_I)\omega_t^{\varsigma}(\xi_I)|\mathrm{d}x_2\mathrm{d}x_I\\ \lesssim\ &\sum_{j\in I}|||k_t|||_{t,\varsigma}\sup_{x_1}\int\left|v_t^{[2](\varkappa)}(\xi_{1:2})\omega_t^{\varsigma+c}(\xi_{1:2})\delta_{12}^{1-\varkappa}\right|e^{(\varsigma-c)\sqrt{t}\delta_{1j}}e^{-\frac{m^2}{t}}t^{-\frac{1}{2}}\mathrm{d}x_2\mathrm{d}x_j\\ \lesssim\ &\sum_{j\in I}|||k_t|||_{t,\varsigma}\,|||v_t^{[2](\varkappa)}|\cdot|^{1-\varkappa}|||_{t,\varsigma+c}\sup_{x_1}\int e^{(\varsigma-c)\sqrt{t}\delta_{1j}}e^{-\frac{m^2}{t}}t^{-\frac{1}{2}}\mathrm{d}x_j\\ \lesssim\ &|||k_t|||_{t,\varsigma}\,|||v_t^{[2](\varkappa)}|||_{t,\varsigma+c}\langle t\rangle^{-\frac{1-\varkappa}{2}}\langle t\rangle^{-1-\frac{1}{2}}\end{aligned}$$

where in the last estimate we used (3.20) and that $c>\varsigma$. Thus we proved (3.25) for $\eta=0$, and the same argument extends for any $\eta\geq 0$ after inserting $\langle \mathrm{St}(\cdot)\rangle^{\eta}$.

For (3.26) we proceed similarly: assuming, without loss of generality that $\sigma_1=\sigma_3=1$ and $\sigma_2=\sigma_4=-1$, we have, writing $\vartheta(x_{34}):=x_3+\vartheta(x_4-x_3)$, we have by Lemma A.6 with $|\mu|=1$ that

$$\begin{aligned}
&\frac{\dot{\mathbb{G}}_t(\xi_{1:2},\xi_{3:4})}{\delta_{12}^{\varkappa}\delta_{34}^{\varkappa}}\\
\asymp_\beta\ & \frac{\dot{\mathcal{G}}_t(x_3-x_1)-\dot{\mathcal{G}}_t(x_4-x_1)}{\delta_{12}^{\varkappa}\delta_{34}^{\varkappa}}-\frac{\dot{\mathcal{G}}_t(x_3-x_2)-\dot{\mathcal{G}}_t(x_4-x_2)}{\delta_{12}^{\varkappa}\delta_{34}^{\varkappa}}\\
=\ & \delta_{12}^{-\varkappa}\delta_{34}^{-\varkappa}(x_3-x_4)\cdot\int_0^1 \nabla\dot{\mathcal{G}}_t(\vartheta(x_{34})-x_1)-\nabla\dot{\mathcal{G}}_t(\vartheta(x_{34})-x_2)\mathrm{d}\vartheta\\
\lesssim\ & e^{-\frac{m^2}{t}}\delta_{12}^{1-\varkappa}\delta_{34}^{1-\varkappa}\int_0^1 e^{-c\sqrt{t}|x_1-\vartheta(x_{34})|}e^{c\sqrt{t}\delta_{12}}\mathrm{d}\vartheta\\
\lesssim\ & e^{-\frac{m^2}{t}}\delta_{12}^{1-\varkappa}\delta_{34}^{1-\varkappa}e^{-c\sqrt{t}\delta_{13}}\omega_t^c(x_{1:2})\omega_t^c(x_{3:4})
\end{aligned}$$

where the last bound follows from

$$|x_1-x_3|-\vartheta|x_4-x_3|\leq|x_1-x_3-\vartheta(x_4-x_3)|=|x_1-\vartheta(x_{34})|.$$

Thus we have, with $K_t(i,j):=\left|v_t^{[2](\varkappa)}(\xi_i,\xi_j)\right|\omega_t^{\varsigma+c}(x_i,x_j)\delta_{ij}^{1-\varkappa}$

$$\begin{aligned}
&\left|\!\left|\!\left|v_t^{[2](\varkappa)}(\xi_{1:2})\frac{\dot{\mathbb{G}}_t(\xi_{1:2},\xi_{3:4})}{|x_1-x_2|^{\varkappa}|x_3-x_4|^{\varkappa}}v_t^{[2](\varkappa)}(\xi_{3:4})\right|\!\right|\!\right|_{t,\varsigma}\\
\lesssim\ & \sup_{x_1}\int K_t(1,2)\,e^{(\varsigma-c)\sqrt{t}\delta_{13}}e^{-\frac{m^2}{t}}K_t(3,4)\mathrm{d}x_{2:4}\\
\lesssim\ & |\!|\!|v_t^{[2](\varkappa)}|\cdot|^{1-\varkappa}|\!|\!|_{t,\varsigma+c}\sup_{x_1}\int K_t(1,2)\,e^{(\varsigma-c)\sqrt{t}\delta_{13}}e^{-\frac{m^2}{t}}\mathrm{d}x_{2:3}\\
\lesssim\ & |\!|\!|v_t^{[2](\varkappa)}|\cdot|^{1-\varkappa}|\!|\!|_{t,\varsigma+c}|\!|\!|v_t^{[2](\varkappa)}|\cdot|^{1-\varkappa}|\!|\!|_{t,\varsigma+c}\langle t\rangle^{-1}\\
\lesssim\ & |\!|\!|v_t^{[2](\varkappa)}|\!|\!|_{t,\varsigma+c}\langle t\rangle^{-\frac{1-\varkappa}{2}}|\!|\!|v_t^{[2](\varkappa)}|\!|\!|_{t,\varsigma+c}\langle t\rangle^{-\frac{1-\varkappa}{2}}\langle t\rangle^{-1}
\end{aligned}$$

by (3.20) and $c>\varsigma$. Again, the argument generalizes to any $\eta\geq 0$ after inserting $\langle \mathrm{St}(\cdot)\rangle^{\eta}$. □

**Remark 3.13. (Above 6$\pi$)** When $\beta^2>6\pi$, i.e. we need to choose $\varkappa\in(1,2]$, Lemma 3.12 cannot be used. We then have to choose option 2 in the discussion preceding the lemma, creating a "tripole" in (3.24). The point is that the remainder bilinear term $\mathbf{B}_{\psi_{1:2}^{\tilde{\varkappa}}\psi_3}(x_{1:3})$ has a lowered divergence $|x_1-x_2|^{-(\varkappa-\tilde{\varkappa})}$; it turns out that we may pick $\tilde{\varkappa}:=1$ so that $\varkappa-\tilde{\varkappa}\in(0,1]$ and Lemma 3.12 can be used to estimate $\mathbf{B}_{\psi_{1:2}^{\tilde{\varkappa}}\psi_3}$; on the other hand, $\mathbf{B}_{\psi_3}(x_3)$ has integration over $\mathrm{d}x_{1:2}$ which provides the additional smoothing effect to negate the divergence $|x_1-x_2|^{-\varkappa}$. As we have mentioned, the price to pay is the proliferating change of bases further along the flow, which results in a representation for $V^{[\ell]}$ significantly different from the naive Mayer ansatz. We refer to [36] for further discussion on this issue.

Note that we introduced in (3.5) the function $\Theta$ in the weight, but so far we only took norms with $\Theta:=\varsigma\in\mathbb{R}_+$. Since the above lemma changes the weight in the norm, this is no longer feasible. Hence from now on, we take (3.5) with

$$\Theta:=\Theta(\ell):=4^{-\ell}.$$

Lemma 3.12 allows us to finally determine the form of the generalised representation (3.15) for all $\ell\geq 3$, and subsequently to propagate the induction through the corresponding kernels, i.e. to extend Lemma 3.3 from $\beta^2<4\pi$ to $\beta^2<6\pi$. Define the basis of representation

$$\mathfrak{M}:=\{(\ell,\sigma_{1:\ell},\varkappa):\ell\in\mathbb{N}_0,\sigma_{1:\ell}\in\{-1,1\}^{\ell},\varkappa\in\{0,1\}\}.$$

Given an index $\boldsymbol{a} := (\ell, \sigma_{1:\ell}, \varkappa) \in \mathfrak{M}$, let us introduce the notation

$$\ell(\boldsymbol{a}) := \ell, \qquad \sigma(\boldsymbol{a}) := \sigma_{1:\ell}, \qquad \varkappa(\boldsymbol{a}) := \varkappa.$$

We decompose the basis as $\mathfrak{M} = \mathfrak{M}_p \uplus (\mathfrak{M} \setminus \mathfrak{M}_p)$, where $\mathfrak{M}_p$ is the subset of indices that needed nontrivial post-processing, namely

$$\mathfrak{M}_p := \{(2, \sigma_{1:2}, \varkappa) \colon \chi(\sigma_{1:2}) = 0\}$$

The associated field content of the representation is then

$$\Psi_{\boldsymbol{a}} := \begin{cases} \psi_{1:\ell} & \text{if } \ \boldsymbol{a} \in \mathfrak{M} \setminus \mathfrak{M}_p, \varkappa(\boldsymbol{a}) = 0 \\ 0 & \text{if } \ \boldsymbol{a} \in \mathfrak{M} \setminus \mathfrak{M}_p, \varkappa(\boldsymbol{a}) = 1 \\ \psi_{1:2}^{\varkappa} & \text{if } \ \boldsymbol{a} \in \mathfrak{M}_p \end{cases} \qquad \text{where} \quad \psi_{1:2}^{\varkappa} := \begin{cases} \dfrac{\psi_{1:2} - 1}{|x_1 - x_2|^{\varkappa}} & \text{if } \ \varkappa > 0 \\ 1 & \text{if } \ \varkappa = 0. \end{cases} \tag{3.27}$$

We note that this is a redefinition of our previous notation; we have already started moving in this direction in Remark 3.10.

The boundary conditions are already determined for almost all kernels of the first 2 iteration orders, cf. (3.19) and $v_{0,T}^{[1](\boldsymbol{a})} := \frac{\lambda}{2}$ when $\varkappa(\boldsymbol{a}) = 0$; what remains is setting

$$v_{0,T}^{[1](\boldsymbol{a})} := 0 \quad \text{if} \quad \varkappa(\boldsymbol{a}) = 1.$$

For all kernels with $\ell \geq 3$ we impose the irrelevant boundary condition, i.e. we set

$$v_{T,T}^{[\ell](\boldsymbol{a})} := 0 \quad \text{if} \quad \ell \geq 3, \quad \boldsymbol{a} \in \mathfrak{M} \setminus \mathfrak{M}_p. \tag{3.28}$$

With this, the boundary conditions have been fixed for all kernels of the index set $\mathfrak{M}$.

Having determined the generalised ansatz (3.15), we can now translate the Polchinski equation (2.20) from $V_t^{[\ell],T}$ to a hierarchy of equations for the kernels $v_{t,T}^{[\ell](\boldsymbol{a})}$ of the representation and prove the inductive bounds on them.

**Proposition 3.14. (estimates on the potential kernels)** *Let $\beta^2 < 6\pi$. Then there exist kernels $\{v_{t,T}^{[\ell](\boldsymbol{a})}\}_{\boldsymbol{a} \in \mathfrak{M}}$ obeying boundary conditions $v_{0,T}^{[1](\boldsymbol{a})} := \frac{\lambda}{2}$, (3.19) and (3.28), such that $\{V_t^{[\ell],T}\}_{\ell}$, given by (3.15), satisfy (2.20). Moreover, for each $\boldsymbol{a} \in \mathfrak{M}$ and any $\eta \in \mathbb{R}_+$, we have the estimates*

$$|||v_{t,T}^{[\ell](\boldsymbol{a})}|||_{t, \Theta(\ell), \eta} \lesssim \lambda_t (\lambda_t \langle t \rangle^{-1})^{\ell - 1} \mathfrak{l}_t^c \langle t \rangle^{-\frac{\varkappa(\boldsymbol{a})}{2}} \tag{3.29}$$

*for some $c = c(\boldsymbol{a}, \beta^2) \in \mathbb{R}_+$. All bounds are independent of $\rho$ except that of $\boldsymbol{a} \in \mathfrak{M}_p$ with $\varkappa(\boldsymbol{a}) = 0$.*

**Proof.** First, we set $v_{\cdot,T}^{[\ell](\boldsymbol{a})} \equiv 0$ whenever $\boldsymbol{a} \in \mathfrak{M} \setminus \mathfrak{M}_p$ and $\varkappa(\boldsymbol{a}) = 1$. These trivially satisfy (3.29) and also consistent with the Polchinski equation because of the 0 field content. Thus from now on we only consider the remaining kernels.

The starting level $\ell = 1$ is done just as in the proof of Lemma 3.3. The case $\ell = 2$ has been treated in Lemma 3.9 (cf. also Remark 3.10 for notation). Thus from now on, suppose $\ell \geq 3$; we will check (3.29) inductively with the understanding that $c \in \mathbb{R}_+$ in $\mathfrak{l}_t^c$ is a constant we do not track and may change from line to line.

Based on Lemma 3.8, the bilinear part of the Polchinski equation (2.20) may be decomposed as

$$\begin{aligned} &\sum_{\ell' + \ell'' = \ell} \langle \mathrm{D} V_s^{[\ell']}, \dot{G}_s \mathrm{D} V_s^{[\ell'']} \rangle \\ &= \sum_{\ell' + \ell'' = \ell} \Big\langle \sum_{\boldsymbol{b} \in \mathfrak{M}} C_{\boldsymbol{b}} \mathrm{D} \langle v_s^{[\ell'](\boldsymbol{b})}, \Psi_{\boldsymbol{b}} \rangle, \dot{G}_s \sum_{\boldsymbol{c} \in \mathfrak{M}} C_{\boldsymbol{c}} \langle \mathrm{D} \langle v_s^{[\ell''](\boldsymbol{c})}, \Psi_{\boldsymbol{c}} \rangle \rangle \Big\rangle \\ &= \Bigg( \sum_{\substack{\ell' + \ell'' = \ell \\ \boldsymbol{b}, \boldsymbol{c} \in \mathfrak{M} \setminus \mathfrak{M}_p}} + \sum_{\substack{\ell' + \ell'' = \ell \\ \boldsymbol{b} \in \mathfrak{M}_p, \boldsymbol{c} \in \mathfrak{M} \setminus \mathfrak{M}_p}} + \sum_{\substack{\ell' + \ell'' = \ell \\ \boldsymbol{b}, \boldsymbol{c} \in \mathfrak{M}_p}} \Bigg) C_{\boldsymbol{b}, \boldsymbol{c}} \langle \mathrm{D} \langle v_s^{[\ell'](\boldsymbol{b})}, \Psi_{\boldsymbol{b}} \rangle, \dot{G}_s \mathrm{D} \langle v_s^{[\ell''](\boldsymbol{c})}, \Psi_{\boldsymbol{c}} \rangle \rangle \end{aligned}$$

with $C_{\boldsymbol{b}}, C_{\boldsymbol{c}}, C_{\boldsymbol{b},\boldsymbol{c}}$ nonnegative combinatorial constants. We treat the three parts in separate cases.

1. Case $\boldsymbol{b},\boldsymbol{c}\in\mathfrak{M}\setminus\mathfrak{M}_p$: by (3.27) we compute, just as in the proof of Lemma 3.1, that

$$\langle \mathrm{D}\langle v_s^{[|I_1|](\boldsymbol{b})},\Psi_{\boldsymbol{b}}\rangle,\dot{G}_s\mathrm{D}\langle v_s^{[|I_2|](\boldsymbol{c})},\Psi_{\boldsymbol{c}}\rangle\rangle = C_{\boldsymbol{b},\boldsymbol{c}}\langle v_s^{[|I_1|](\boldsymbol{b})}(\xi_{I_1})\dot{\mathbb{G}}_s(\xi_{I_1},\xi_{I_2})v_s^{[|I_2|](\boldsymbol{c})}(\xi_{I_2}),\psi_{1:\ell}\rangle$$

   for some constant $C_{\boldsymbol{b},\boldsymbol{c}}$, so that

$$\begin{aligned}&|||v_s^{[|I_1|](\boldsymbol{b})}(\xi_{I_1})\dot{\mathbb{G}}_s(\xi_{I_1},\xi_{I_2})v_s^{[|I_2|](\boldsymbol{c})}(\xi_{I_2})|||_{s,\Theta(\ell)}\\ \lesssim\ & |||v_s^{[|I_1|](\boldsymbol{b})}|||_{s,4^{-\ell}}\langle s\rangle^{-2}|||v_s^{[|I_2|](\boldsymbol{c})}|||_{s,4^{-\ell}}\\ \lesssim\ & |||v_s^{[|I_1|](\boldsymbol{b})}|||_{s,\Theta(|I_1|)}\langle s\rangle^{-2}|||v_s^{[|I_2|](\boldsymbol{c})}|||_{s,\Theta(|I_2|)}\\ \lesssim\ & \lambda_s(\lambda_s\langle s\rangle^{-1})^{\ell-1}\langle s\rangle^{-1}\mathfrak{l}_s^c\end{aligned}$$

2. Case $\boldsymbol{b}\in\mathfrak{M}_p,\boldsymbol{c}\in\mathfrak{M}\setminus\mathfrak{M}_p$: computing as in (3.23),

$$\langle \mathrm{D}\langle v_s^{[|I_1|](\boldsymbol{b})},\Psi_{\boldsymbol{b}}\rangle,\dot{G}_s\mathrm{D}\langle v_s^{[|I_2|](\boldsymbol{c})},\Psi_{\boldsymbol{c}}\rangle\rangle = C_{\boldsymbol{b},\boldsymbol{c}}\left\langle v_s^{[2](\boldsymbol{b})}(\xi_{1:2})\frac{\dot{\mathbb{G}}_s(\xi_{1:2},\xi_{I_2})}{|x_1-x_2|^{\varkappa}}v_s^{[|I_2|](\boldsymbol{c})}(\xi_{I_2}),\psi_{1:\ell}\right\rangle$$

   where $\varkappa:=\varkappa(\boldsymbol{b})=1$. Hence estimating by (3.25) with $c:=2\cdot 4^{-\ell}$ and using the bounds of Lemma 3.9 we get, since $\ell\geq 3$

$$\begin{aligned}&\left|\left|\left|v_s^{[2](\boldsymbol{b})}(\xi_{1:2})\frac{\dot{\mathbb{G}}_s(\xi_{1:2},\xi_{I_2})}{|x_1-x_2|^{\varkappa}}v_s^{[|I_2|](\boldsymbol{c})}(\xi_{I_2})\right|\right|\right|_{s,\Theta(\ell)}\\ \lesssim\ & |||v_s^{[2](\boldsymbol{b})}|||_{s,3\cdot 4^{-\ell}}\langle s\rangle^{-2+\frac{\varkappa}{2}}|||v_s^{[|I_2|](\boldsymbol{c})}|||_{s,4^{-\ell}}\\ \lesssim\ & |||v_s^{[2](\boldsymbol{b})}|||_{s,4^{-(\ell-1)}}\langle s\rangle^{-2+\frac{\varkappa}{2}}|||v_s^{[|I_2|](\boldsymbol{c})}|||_{s,4^{-\ell}}\\ \lesssim\ & |||v_s^{[2](\boldsymbol{b})}|||_{s,\Theta(2)}\langle s\rangle^{-2+\frac{\varkappa}{2}}|||v_s^{[|I_2|](\boldsymbol{c})}|||_{s,\Theta(|I_2|)}\\ \lesssim\ & \lambda_s(\lambda_s\langle s\rangle^{-1})^{\ell-1}\langle s\rangle^{-1}\mathfrak{l}_s^c.\end{aligned}$$

3. Case $\boldsymbol{b},\boldsymbol{c}\in\mathfrak{M}_p$: it follows from the definition of $\mathfrak{M}_p$ that this only happens in when $\ell=4$. Similarly to (3.23) we obtain

$$\begin{aligned}&\langle \mathrm{D}\langle v_s^{[|I_1|](\boldsymbol{b})},\Psi_{\boldsymbol{b}}\rangle,\dot{G}_s\mathrm{D}\langle v_s^{[|I_2|](\boldsymbol{c})},\Psi_{\boldsymbol{c}}\rangle\rangle\\ =\ & C_{\boldsymbol{b},\boldsymbol{c}}\left\langle v_s^{[2](\boldsymbol{b})}(\xi_{1:2})\frac{\dot{\mathbb{G}}_s(\xi_{1:2},\xi_{3:4})}{|x_1-x_2|^{\varkappa}|x_3-x_4|^{\varkappa}}v_s^{[2](\boldsymbol{c})}(\xi_{3:4}),\psi_{1:4}\right\rangle\end{aligned}$$

   where $\varkappa:=\varkappa(\boldsymbol{b})=\varkappa(\boldsymbol{c})=1$. Thus estimating by (3.26) with $c:=2\cdot 4^{-\ell}$ and Lemma 3.9

$$\begin{aligned}&\left|\left|\left|v_s^{[2](\boldsymbol{b})}(\xi_{1:2})\frac{\dot{\mathbb{G}}_s(\xi_{1:2},\xi_{3:4})}{|x_1-x_2|^{\varkappa}|x_3-x_4|^{\varkappa}}v_s^{[2](\boldsymbol{c})}(\xi_{3:4})\right|\right|\right|_{s,\Theta(\ell)}\\ \lesssim\ & |||v_s^{[2](\boldsymbol{b})}|||_{s,4^{-(\ell-1)}}\langle s\rangle^{-2+\varkappa}|||v_s^{[2](\boldsymbol{c})}|||_{s,4^{-(\ell-1)}}\\ \lesssim\ & |||v_s^{[2](\boldsymbol{b})}|||_{s,\Theta(2)}\langle s\rangle^{-2+\varkappa}|||v_s^{[2](\boldsymbol{c})}|||_{s,\Theta(2)}\\ \lesssim\ & \lambda_s(\lambda_s\langle s\rangle^{-1})^{\ell-1}\langle s\rangle^{-1}\mathfrak{l}_s^c.\end{aligned}$$

Hence combining the above cases and using the screening bound (3.12) we get

$$\begin{aligned}|||v_t^{[\ell](\boldsymbol{a})}|||_{t,\Theta(\ell)} &\lesssim \int_t^\infty \|e^{\mathbb{G}_{t,s}(\xi_{1:\ell})}\|_{L^\infty}\lambda_s(\lambda_s\langle s\rangle^{-1})^{\ell-1}\langle s\rangle^{-1}\mathfrak{l}_s^c\mathrm{d}s\\ &\lesssim \int_t^\infty (\lambda_t\lambda_s^{-1})^{|\chi(\xi_{1:\ell})|}\lambda_s(\lambda_s\langle s\rangle^{-1})^{\ell-1}\langle s\rangle^{-1}\mathfrak{l}_s^c\mathrm{d}s\\ &\lesssim \lambda_t(\lambda_t\langle t\rangle^{-1})^{\ell-1}\mathfrak{l}_t^c.\end{aligned}$$

where the integral converges because the assumption $\ell\geq 3$ implies either $\ell=3$, in which case $|\chi(\sigma_{1:\ell})|\geq 1$, or $\ell\geq 4$, in which case $\lambda_s(\lambda_s\langle s\rangle^{-1})^{\ell-1}\asymp\langle s\rangle^{1-4\delta}\lesssim\langle s\rangle^{-\varepsilon}$ for some $\varepsilon>0$ as $\beta^2<6\pi\Leftrightarrow\delta>1/4$. □

A straightforward consequence of the above uniform bounds is the removal of the UV cutoff $T$. To this end, we extend the kernels trivially from $[0,T]$ to $\mathbb{R}_+$; namely, put

$$\tilde{v}_{t,T}^{[\ell](\boldsymbol{a})} := v_{t\wedge T,T}^{[\ell](\boldsymbol{a})}, \tag{3.30}$$

for any $t \in \mathbb{R}_+$; we abuse notation by simply writing $v_{t,T}^{[\ell](\boldsymbol{a})}$ instead of $\tilde{v}_{t,T}^{[\ell](\boldsymbol{a})}$.

**Corollary 3.15.** *For any $\eta \geq 0$, the $|||\cdot|||_{t,\Theta,\eta}$-limits*

$$v_{t,\infty}^{[\ell](\boldsymbol{a})} := \lim_{T\to\infty} v_{t,T}^{[\ell](\boldsymbol{a})}$$

*exist, are independent of $\eta$, and satisfy the bounds (3.29). Moreover, $V^{[\ell],\infty,\rho}(\phi)$, given by*

$$V_t^{[\ell],\infty,\rho} := \sum_{\substack{\boldsymbol{a}\in\mathfrak{M} \\ \ell(\boldsymbol{a})=\ell}} V_t^{[\ell](\boldsymbol{a}),\infty,\rho} \qquad \text{where} \qquad V_t^{[\ell](\boldsymbol{a}),\infty,\rho}(\phi) := \int v_{t,\infty}^{[\ell](\boldsymbol{a})}(\xi_{\boldsymbol{a}})\Psi_{\boldsymbol{a}}(\xi_{\boldsymbol{a}})\rho(x_{\boldsymbol{a}})\mathrm{d}\xi_{\boldsymbol{a}} \tag{3.31}$$

*satisfies (2.20).*

**Proof.** See Appendix B.3 □

To simplify notation, we frequently suppress reference to $\rho$, and also write $v_t^{[\ell](\boldsymbol{a})} := v_{t,\infty}^{[\ell](\boldsymbol{a})}$ and $V_t^{[\ell]} := V_t^{[\ell],\infty}$.

# 4 The effective force observable

With the inductive bounds on the kernels of the representation proven, we may translate them to estimates on the approximate effective force. We tacitly assume that every object with superscript $T$ is extended from $[0,T]$ to $\mathbb{R}_+$ trivially, i.e. just as in (3.30) above.

**Proposition 4.1. (estimates on the approximate force)** *Let $\beta^2 < 6\pi$. Denote by $E := L^\infty$ or $L^2(\langle x\rangle^{-n})$ for $n > 2$. Then with $F_t^{[\ell],T,\rho} := \mathrm{D}V_t^{[\ell],T,\rho}$, for all $X \in \mathcal{S}'$, we have*

$$\sup_T \sup_\rho \sup_{\phi\in\mathcal{H}} \left\| Q_t F_t^{[\ell],T,\rho}(\phi) \right\|_E \lesssim \lambda_t (\lambda_t \langle t\rangle^{-1})^{\ell-1} \mathfrak{l}_t^c \langle t\rangle^{-1}. \tag{4.1}$$

$$\sup_T \sup_\rho \sup_{\phi\in\mathcal{H}} \left\| \langle \mathrm{D}F_t^{[\ell],T,\rho}(\phi), \dot{G}_t X\rangle \right\|_E \lesssim \lambda_t (\lambda_t \langle t\rangle^{-1})^{\ell-1} \mathfrak{l}_t^c \langle t\rangle^{-1} \|Q_t X\|_E. \tag{4.2}$$

*Moreover, denoting $E := L^\infty$ or $L^2(\langle x\rangle^{-n})$ for any $n \in \mathbb{Z}$, we have*

$$\sup_T \sup_\rho \left\| Q_t\big(F_t^{[\ell],T,\rho}(\phi) - F_t^{[\ell],T,\rho}(\phi')\big) \right\|_E \lesssim \lambda_t (\lambda_t \langle t\rangle^{-1})^{\ell-1} \mathfrak{l}_t^c \langle t\rangle^{-1} \|\phi - \phi'\|_E \tag{4.3}$$

$$\sup_T \sup_\rho \left\| \langle \mathrm{D}\big(F_t^{[\ell],T,\rho}(\phi) - F_t^{[\ell],T,\rho}(\phi')\big), \dot{G}_t X\rangle \right\|_E \lesssim \lambda_t (\lambda_t \langle t\rangle^{-1})^{\ell-1} \mathfrak{l}_t^c \langle t\rangle^{-1} \|Q_t X\|_{L^\infty} \|\phi - \phi'\|_E \tag{4.4}$$

*Finally, if $\|Q_t X\|_{L^\infty} < \infty$, then as $T \to \infty$,*

$$\sup_\rho \sup_{\phi\in\mathcal{H}} \left\| Q_t\big(F_t^{[\ell],\infty,\rho} - F_t^{[\ell],T,\rho}\big)(\phi) \right\|_{L^\infty} + \left\| \langle \mathrm{D}\big(F_t^{[\ell],\infty,\rho} - F_t^{[\ell],T,\rho}\big)(\phi), \dot{G}_t X\rangle \right\|_{L^\infty} \to 0.$$

*and for $n > 2$, denoting by $L^{2,-n} := L^2(\langle x\rangle^{-n})$, as $\rho \to \mathbb{1}_{\mathbb{R}^2}$,*

$$\sup_T \sup_{\phi\in\mathcal{H}} \left\| Q_t\big(F_t^{[\ell],T,\mathbb{1}_{\mathbb{R}^2}} - F_t^{[\ell],T,\rho}\big)(\phi) \right\|_{L^{2,-n}} + \left\| \langle \mathrm{D}\big(F_t^{[\ell],T,\mathbb{1}_{\mathbb{R}^2}} - F_t^{[\ell],T,\rho}\big)(\phi), \dot{G}_t X\rangle \right\|_{L^{2,-n}} \to 0.$$

**Proof.** The bounds follow from the estimates of the previous section in a straightforward way.

1. **The force bound (4.1):**
The vacuum energy contribution $v_{t,T}^{[2](0)}$ is annihilated by the derivative, thus denoting by
$$M^{(\boldsymbol{b})}(x_I) := \mathbb{1}_{\boldsymbol{b}\in\mathfrak{M}\setminus\mathfrak{M}_p} + \mathbb{1}_{\boldsymbol{b}\in\mathfrak{M}_p}\mathbb{1}_{\varkappa(\boldsymbol{b})>0}|x_{I(1)} - x_{I(2)}|^{-\varkappa(\boldsymbol{b})}$$
we compute
$$\begin{aligned} Q_t F_t^{[\ell],T}(\phi)(y) &= Q_t \mathrm{D}V_t^{[\ell],T}(\phi)(y) \\ &= i\beta \sum_{\substack{\boldsymbol{a}\in\mathfrak{M}\\ \ell(\boldsymbol{a})=\ell}} \int \rho(x_{1:\ell}) v_{t,T}^{[\ell](\boldsymbol{a})} M^{(\boldsymbol{a})}(\xi_{1:\ell}) \sum_{j=1}^{\ell} \sigma_j k_t^Q(x_j - y)\psi_{1:\ell}\mathrm{d}\xi_{1:\ell}. \end{aligned}$$
Hence for $\boldsymbol{a}\in\mathfrak{M}\setminus\mathfrak{M}_p$, we simply have, uniformly in $\|\rho\|_{L^\infty}\leq 1$, that
$$\begin{aligned} \big\|Q_t F_t^{[\ell](\boldsymbol{a}),T}(\phi)\big\|_{L^\infty} &\lesssim \sup_y \int \big|v_{t,T}^{[\ell](\boldsymbol{a})}(\xi_{1:\ell})\big| \sum_{j=1}^{\ell} k_t^Q(x_j - y)\mathrm{d}\xi_{1:\ell} \\ &\leq |||v_{t,T}^{[\ell](\boldsymbol{a})}|||_{t,\Theta,\eta} \sup_y \int_{\mathbb{R}^d} k_t^Q(x-y)\mathrm{d}x \\ &\lesssim |||v_{t,T}^{[\ell](\boldsymbol{a})}|||_{t,\Theta,\eta}\langle t\rangle^{-1}, \end{aligned}$$
On the other hand, for $\boldsymbol{a}\in\mathfrak{M}_p$ with $\varkappa := \varkappa(\boldsymbol{a}) = 1$, by neutrality $\sigma_1 = -\sigma_2$, using Lemma A.6, and (3.20) yields
$$\begin{aligned} \big|Q_t F_t^{[\ell](\boldsymbol{a}),T}(\phi)(y)\big| &\lesssim \int v_{t,T}^{[\ell](\boldsymbol{a})}(\xi_{1:2}) \frac{\sigma_1 k_t^Q(x_1 - y) + \sigma_2 k_t^Q(x_2 - y)}{|x_1 - x_2|^{\varkappa}} \psi_{1:2}\mathrm{d}\xi_{1:2} \\ &\lesssim \int v_{t,T}^{[\ell](\boldsymbol{a})}(\xi_{1:2}) \omega_t^c(x_{1:2}) e^{-c\sqrt{t}|x_1 - y|} e^{-\frac{m^2}{2t}} t^{\frac{1}{2}} \mathrm{d}x_{1:2} \\ &\lesssim |||v_{t,T}^{[\ell](\boldsymbol{a})}|||_{t,\Theta,\eta}\langle t\rangle^{-\frac{1}{2}} \end{aligned}$$
and the $L^\infty$ statement follows from Proposition 3.14 and Corollary 3.15. The $L^2(\langle x\rangle^{-n})$ statement is done verbatim.

2. **The contraction bound (4.2):**
Since $\dot{G}_t X = \dot{\mathcal{G}}_t * X \in \mathcal{C}^\infty(\mathbb{R}^d)$, point evaluations make sense and we may compute
$$\begin{aligned} &\big\langle \mathrm{D}F_t^{[\ell](\boldsymbol{a}),T}(\phi), \dot{G}_t X\big\rangle(y) \\ =\ & \mathrm{D}F_t^{[\ell](\boldsymbol{a}),T}(\phi)[\dot{G}_t X](y) \\ =\ & -\beta^2 \int \rho(x_{1:\ell}) v_{t,T}^{[\ell](\boldsymbol{a})} M^{(\boldsymbol{a})}(\xi_{1:\ell}) \sum_{j,k=1}^{\ell} \sigma_j\sigma_k (\dot{G}_t X)(x_j)\delta(x_k - y)\psi_{1:\ell}\mathrm{d}\xi_{1:\ell}. \end{aligned}$$
When $\boldsymbol{a}\in\mathfrak{M}\setminus\mathfrak{M}_p$, we simply bound it by
$$\begin{aligned} &\lesssim \|\dot{G}_t X\|_{L^\infty} \int \big|v_{t,T}^{[\ell](\boldsymbol{a})}\big|(\xi_{1:\ell}) \sum_{k=1}^{\ell} \delta(x_k - y)\mathrm{d}\xi_{1:\ell} \\ &\lesssim \langle t\rangle^{-1}\|Q_t X\|_{L^\infty} |||v_{t,T}^{[\ell](\boldsymbol{a})}|||_{t,\Theta,\eta}. \end{aligned}$$
and we conclude from Proposition 3.14 and Corollary 3.15. When $\boldsymbol{a}\in\mathfrak{M}_p$ with $\varkappa := \varkappa(\boldsymbol{a}) = 1$, by neutrality $\sigma_1 = -\sigma_2$ it is estimated via
$$\begin{aligned} &\lesssim \int \big|v_{t,T}^{[\ell](\boldsymbol{a})}\big|(\xi_{1:2}) \frac{|(\dot{G}_t X)(x_1) - (\dot{G}_t X)(x_2)|}{|x_1 - x_2|} \sum_{k=1}^{2} \delta(x_k - y)\mathrm{d}\xi_{1:2} \\ &\lesssim \|\nabla \dot{G}_t X\|_{L^\infty} |||v_{t,T}^{[\ell](\boldsymbol{a})}|||_{t,\Theta,\eta} \\ &\lesssim \langle t\rangle^{-\frac{1}{2}} \|Q_t X\|_{L^\infty} |||v_{t,T}^{[\ell](\boldsymbol{a})}|||_{t,\Theta,\eta} \end{aligned}$$

by (A.2).

3. **The Lipschitz bounds (4.3) and (4.4):**

Denoting by $\psi_j' := e^{i\beta\sigma_j\phi'(x_j)}$, we have

$$|\psi_{1:\ell} - \psi_{1:\ell}'| = \sum_{j\in 1:\ell} |e^{i\beta\sigma_j\phi(x_j)} - e^{i\beta\sigma_j\phi'(x_j)}| \leq \beta \sum_{j\in 1:\ell} |\phi(x_j) - \phi'(x_j)|$$

from which the $L^\infty$ bounds follow immediately by repeating the argument of (4.1) and (4.2). Thus from now on, we focus only on the weighted $L^2$ bounds. First let $\boldsymbol{a} \in \mathfrak{M} \setminus \mathfrak{M}_p$. Fix a charge configuration $\sigma_{1:\ell}$ and pick any $j, k \in 1{:}\,\ell$; we have, with the change of variables $r = x_j - y$, $z_j = 0$, $z_m = x_m - x_j$ $(m \neq j)$

$$\begin{aligned} I_{j,k}(y) &= \int |v_{t,T}^{[\ell](\boldsymbol{a})}(x_{1:\ell})| k_t^Q(x_j - y)|\phi(x_k) - \phi'(x_k)| \mathrm{d}x_{1:\ell} \\ &= \int |v_{t,T}^{[\ell](\boldsymbol{a})}(z_{1:\ell} + y + r)| k_t^Q(r)|\phi(z_k + y + r) - \phi'(z_k + y + r)| \mathrm{d}z_{1:\ell\setminus\{j\}} \mathrm{d}r \\ &= \int |v_{t,T}^{[\ell](\boldsymbol{a})}(z_{1:\ell})| k_t^Q(r)|\phi(z_k + y + r) - \phi'(z_k + y + r)| \mathrm{d}z_{1:\ell\setminus\{j\}} \mathrm{d}r \end{aligned}$$

where the last equality follows from translation invariance of the kernels $v_{t,T}^{[\ell](\boldsymbol{a})}$. Then by Minkowski's integral inequality, the estimate $\|u(\cdot + h)\|_{L^{2,-n}} \lesssim_h \langle h\rangle^{|n|} \|u\|_{L^{2,-n}}$ and

$$\|\langle\cdot\rangle^{-n} I_{j,k}\|_{L^2} \lesssim |||v_{t,T}^{[\ell](\boldsymbol{a})}|||_{t,\Theta,\eta} \|k_t^Q\|_{L^1} \|\phi - \phi'\|_{L^{2,-n}}$$

where $\langle r\rangle^{|n|}$ is eaten by $\|k_t^Q\|_{L^1}$ and $\langle z_k\rangle^{|n|}$ by the parameter $\eta$ in the kernel norm. This shows (4.3) in the case $\boldsymbol{a} \in \mathfrak{M} \setminus \mathfrak{M}_p$. For $\boldsymbol{a} \in \mathfrak{M}_p$, we combine the above argument with how this case was handled in proving (4.1). Finally, (4.4) follows along similar lines from the proof of (4.2).

4. **The $\boldsymbol{T}$ and $\boldsymbol{\rho}$ limits:**

The claim about the $T \to \infty$ limits follows from Corollary 3.15.

For the claim about the $\rho \to \mathbb{1}_{\mathbb{R}^2}$, limit, estimate, when $\boldsymbol{a} \in \mathfrak{M} \setminus \mathfrak{M}_p$,

$$\begin{aligned} &\|Q_t(F_t^{[\ell](\boldsymbol{a}),\infty,\rho} - F_t^{[\ell](\boldsymbol{a}),T,\rho})(\phi)\|_{L^{2,-n}}^2 \\ \asymp &\int \left| \int (1 - \rho(x_{1:\ell})) v_t^{[\ell](\boldsymbol{a})}(\xi_{1:\ell}) \sum_{j\in 1:\ell} \sigma_j k_t^Q(x_j - y) \psi_{1:\ell} \mathrm{d}\xi_{1:\ell} \right|^2 \langle y\rangle^{-n} \mathrm{d}y \\ \lesssim &\sup_{x_{1:\ell}} |1 - \rho(x_{1:\ell})|^2 |||v_t^{[\ell](\boldsymbol{a})}|||_{t,\Theta,\eta}^2 |||k_t^Q|||_{L^1}^2 \int \langle y\rangle^{-n} \mathrm{d}y \end{aligned}$$

and apply dominated convergence. On the other hand, when $\boldsymbol{a} \in \mathfrak{M}_p$ with $\varkappa := \varkappa(\boldsymbol{a}) = 1$, we use neutrality and Lemma A.6 to get

$$\begin{aligned} &\int \left| \int (1 - \rho(x_{1:2})) v_t^{[\ell](\boldsymbol{a})}(\xi_{1:2}) \frac{k_t^Q(x_1 - y) - k_t^Q(x_2 - y)}{|x_1 - x_2|} \mathrm{d}\xi_{1:2} \right|^2 \langle y\rangle^{-n} \mathrm{d}y \\ \lesssim &\int \left| \int (1 - \rho(x_{1:2})) v_t^{[\ell](\boldsymbol{a})}(\xi_{1:2}) e^{-c\sqrt{t}|x_1 - y|} e^{c\sqrt{t}|x_1 - x_2|} e^{-\frac{m^2}{2t}} t^{\frac{1}{2}} \mathrm{d}x_{1:2} \right|^2 \langle y\rangle^{-n} \mathrm{d}y \\ \lesssim &\; t^{\frac{1}{2}} \sup_{x_{1:2}} |1 - \rho(x_{1:2})|^2 |||v_t^{[\ell](\boldsymbol{a})}|||_{t,\Theta,\eta}^2 \left\| e^{-c\sqrt{t}|\cdot|} e^{-\frac{m^2}{2t}} \right\|_{L^1} \int \langle y\rangle^{-n} \mathrm{d}y \end{aligned}$$

and again conclude from dominated convergence.

The case of $\langle \mathrm{D}F_t^{[\ell],T}(\phi), \dot{G}_t X\rangle$ follows similarly, by combining the above argument with the proof of (4.2). □

Combining the bounds (4.1) and (4.2), we obtain that for $\ell', \ell''$ such that $\ell' + \ell'' \geq \ell$, we have

$$\begin{aligned} \| \langle \mathrm{D}F_t^{[\ell'],T}(\phi), \dot{G}_t F_t^{[\ell''],T}(\phi) \rangle \|_{L^\infty} &\lesssim \lambda_t (\lambda_t \langle t \rangle^{-1})^{\ell'-1} \mathfrak{l}_t^c \langle t \rangle^{-1} \lambda_t (\lambda_t \langle t \rangle^{-1})^{\ell''-1} \mathfrak{l}_t^c \langle t \rangle^{-1} \\ &\lesssim (\lambda_t \langle t \rangle^{-1})^{\ell} \mathfrak{l}_t^c \end{aligned} \tag{4.5}$$

again uniformly in $\rho$, $\phi$, and $T \in [0, \infty]$.

We are in now position to prove the existence of the sine–Gordon scale decomposition process $\varphi$. We start with recalling the following

**Theorem 4.2. ([29])** *Let $\beta^2 < 6\pi$ and denote by $\tilde{F}^T := F^{[\leq 3],T}$. Then for any $\rho \in \mathcal{C}_c^\infty$ and $T \in [0, \infty)$ there exists a unique solution $(\varphi^{T,\rho}, \tilde{R}^{T,\rho})$ to the equation*

$$\begin{cases} \varphi_t^{T,\rho} = -\displaystyle\int_0^t \dot{G}_s (\tilde{F}_s^{T,\rho}(\varphi_s^{T,\rho}) + \tilde{R}_s^{T,\rho}) \mathrm{d}s + B_t \\ \tilde{R}_t^{T,\rho} = \mathbb{E}_t^B \left[ \displaystyle\int_t^T \left( \partial_s + \frac{1}{2} \Delta_{\dot{G}_s} \right) \tilde{F}_s^{T,\rho}(\varphi_s^{T,\rho}) - \langle \mathrm{D}\tilde{F}_s^{T,\rho}(\varphi_s^{T,\rho}), \dot{G}_s (\tilde{F}_s^{T,\rho}(\varphi_s^{T,\rho}) + \tilde{R}_s^{T,\rho}) \rangle \mathrm{d}s \right] \end{cases} \tag{4.6}$$

*satisfying*

$$\sup_\rho \sup_T \| \tilde{R}_t^{T,\rho} \|_{L^\infty(\Omega; L^\infty)} \lesssim |\lambda| + \lambda_t^4 \langle t \rangle^{-3} \tag{4.7}$$

**Proof.** Existence and uniqueness is Lemma 2.6 and Theorem 2.2 in [29]. The uniform boundedness estimate (4.7) follows from Lemma 4.1. □

**Corollary 4.3.** *Let $n > 2$. Then*

1. *The laws of $(\varphi_\cdot^{T,\rho}, B_\cdot)$ are tight on $\mathcal{C}_{\mathrm{loc}}(\mathbb{R}_+; L^2(\langle x \rangle^{-n}))^2$*
2. *The laws of $\tilde{R}_\cdot^{T,\rho}$ are tight on $L^2_{\mathrm{loc}}(\mathbb{R}_+; L^2(\langle x \rangle^{-n}))$, equipped with the weak topology of $L^2([0, S]; L^2(\langle x \rangle^{-n}))$ for each $S \in \mathbb{N}$.*

**Proof.** Fix an $S < \infty$. Denote by $Z := \varphi - B$. It follows from (4.1), (4.7) and (A.2) that uniformly in $T$ and $\rho$,

$$\sup_{t \in [0,S]} (\| Z_t^{T,\rho} \|_{L^\infty(\Omega; W^{1,\infty})} + \| \dot{Z}_t^{T,\rho} \|_{L^\infty(\Omega; L^\infty)}) \leq C_{\lambda,S}.$$

It is clear that $W^{1,\infty}(\mathbb{R}^2)$ is compactly embedded in $L^2(\langle x \rangle^{-n})$. Indeed, pick a bounded $(u_k) \subseteq W^{1,\infty}(\mathbb{R}^2)$; on every closed ball of $\mathbb{R}^2$ we may apply Arzelà–Ascoli, and a diagonal argument allows to extract a subsequence $(u_{k_j})$ converging locally uniformly on $\mathbb{R}^2$ to a $u \in \mathcal{C}_b$. By integrability of $\langle x \rangle^{-n}$, dominated convergence yields $u_{k_j} \to u$ in $L^2(\langle x \rangle^{-n})$.

Then by Lipschitzness obtained from the bound above on $\dot{Z}_t^T$, we may use Arzelà–Ascoli to see that $\{Z_\cdot^{T,\rho}\}_{T,\rho}$ lives almost surely in a compact subset of $\mathcal{C}([0, S]; L^2(\langle x \rangle^{-n}))$, hence its laws are tight.

On the other hand, it follows from (4.7) that uniformly in $T$ and $\rho$, we have almost surely

$$\| \tilde{R}^{T,\rho} \|^2_{L^2([0,S]; L^2(\langle x \rangle^{-n}))} = \int_0^S \| \tilde{R}_s^{T,\rho} \|^2_{L^2(\langle x \rangle^{-n})} \mathrm{d}s \lesssim \int_0^S |\lambda| + \lambda_s^4 \langle s \rangle^{-3} \mathrm{d}s,$$

hence by the Banach–Alaoglu theorem and reflexivity, $\{\tilde{R}^{T,\rho}\}_{T,\rho}$ lives almost surely in a weakly compact subset of $L^2([0, S]; L^2(\langle x \rangle^{-n}))$. □

By Prokhorov's theorem, the sequence of the above measures admit a weakly convergent subsequence. The next proposition identifies such limits as the laws of the solution to an FBSDE corresponding to the cutoff-removed effective force coefficients.

**Theorem 4.4.** *For every accumulation point $(\nu_\varphi, \nu_R)$ of* $\mathrm{Law}(\varphi^{T,\rho}, \tilde{R}^{T,\rho})_{T,\rho}$ *in Theorem 4.2, there exists a filtered probability space $(\Omega, \mathscr{A}, \mathscr{F}, \mathbb{P})$ with admissible filtration $\mathscr{F}$ such that $(\varphi, \tilde{R})$ has law $(\nu_\varphi, \nu_R)$ and solves*

$$\begin{cases} \varphi_t &= -\int_0^t \dot{G}_s(\tilde{F}_s^{\infty,\mathbb{1}_{\mathbb{R}^2}}(\varphi_s) + \tilde{R}_s)\mathrm{d}s + B_t \\ \tilde{R}_t &= \mathbb{E}_t^{\mathscr{F}}\Big[\int_t^\infty \Big(\partial_s + \frac{1}{2}\Delta_{\dot{G}_s}\Big)\tilde{F}_s^{\infty,\mathbb{1}_{\mathbb{R}^2}}(\varphi_s) - \langle \mathrm{D}\tilde{F}_s^{\infty,\mathbb{1}_{\mathbb{R}^2}}(\varphi_s), \dot{G}_s(\tilde{F}_s^{\infty,\mathbb{1}_{\mathbb{R}^2}}(\varphi_s) + \tilde{R}_s)\rangle \mathrm{d}s\Big] \end{cases} \tag{4.8}$$

**Proof.** Let $(T, \rho)$ index a convergent subsequence. By the Skorokhod representation theorem, there exists a probability space $(\Omega, \mathscr{A}, \mathbb{P})$ and random variables $(\varphi, \tilde{R})$ such that $\mathbb{P}$-almost surely $\varphi^{T,\rho} \to \varphi$ in $\mathcal{C}([0,S], L^2(\langle x\rangle^{-n}))$ and $\tilde{R}^{T,\rho} \to \tilde{R}$ weakly in $L^2([0,S]; L^2(\langle x\rangle^{-n}))$ for any $S$.

**Convergence of the forward equation**

Write

$$\begin{aligned} &\int_0^t \dot{G}_s(\tilde{F}_s^{\infty,\mathbb{1}_{\mathbb{R}^2}}(\varphi_s) + \tilde{R}_s) - \dot{G}_s(\tilde{F}_s^{T,\rho}(\varphi_s^{T,\rho}) + \tilde{R}_s^{T,\rho})\mathrm{d}s \\ =& \int_0^t \dot{G}_s(\tilde{F}_s^{\infty,\mathbb{1}_{\mathbb{R}^2}}(\varphi_s) - \tilde{F}_s^{\infty,\mathbb{1}_{\mathbb{R}^2}}(\varphi_s^{T,\rho}))\mathrm{d}s + \int_0^t \dot{G}_s(\tilde{F}_s^{\infty,\mathbb{1}_{\mathbb{R}^2}}(\varphi_s^{T,\rho}) - \tilde{F}_s^{T,\rho}(\varphi_s^{T,\rho}))\mathrm{d}s \\ +& \int_0^t \dot{G}_s(\tilde{R}_s - \tilde{R}_s^{T,\rho})\mathrm{d}s \end{aligned}$$

The first term vanishes in $L^2(\langle x\rangle^{-n})$ due to (4.3), the second due to the convergence result of Proposition 4.1, and the third weakly in $L^2(\langle x\rangle^{-n})$ due to $\tilde{R}^{T,\rho} \to \tilde{R}$. Indeed, for the latter, note that for any $h \in L^2(\langle x\rangle^{-n})$ we have

$$\Big\langle \int_0^t \dot{G}_s(\tilde{R}_s - \tilde{R}_s^{T,\rho})\mathrm{d}s, h \Big\rangle_{L^2(\langle x\rangle^{-n})} = \int_0^t \langle (\tilde{R}_s - \tilde{R}_s^{T,\rho}), \dot{G}_s^* h\rangle_{L^2(\langle x\rangle^{-n})}\mathrm{d}s \longrightarrow 0,$$

as $s \mapsto \dot{G}_s^* h$ is in $L^2([0,t]; L^2(\langle x\rangle^{-n}))$ due to Lemma A.9 (recall Notation B.2). Coupled with the convergence $\varphi^{T,\rho} \to \varphi$, we obtain that the forward equation of (4.8) holds in the limit.

**Convergence of the integrand of the backward equation**

We are going to show that the integrands of $\tilde{R}_t^{T,\rho}$ converge to that with the limiting force coefficients. For any $\ell', \ell'' \geq 0$,

$$\begin{aligned} & \langle \mathrm{D}F^{[\ell'],\infty,\mathbb{1}_{\mathbb{R}^2}}(\varphi), \dot{G}F^{[\ell''],\infty,\mathbb{1}_{\mathbb{R}^2}}(\varphi)\rangle - \langle \mathrm{D}F^{[\ell'],T,\rho}(\varphi^{T,\rho}), \dot{G}F^{[\ell''],T,\rho}(\varphi^{T,\rho})\rangle \\ =& \langle \mathrm{D}(F^{[\ell'],\infty,\mathbb{1}_{\mathbb{R}^2}}(\varphi) - F^{[\ell'],T,\rho}(\varphi)), \dot{G}F^{[\ell''],\infty,\mathbb{1}_{\mathbb{R}^2}}(\varphi)\rangle \\ +& \langle \mathrm{D}(F^{[\ell'],T,\rho}(\varphi) - F^{[\ell'],T,\rho}(\varphi^{T,\rho})), \dot{G}F^{[\ell''],\infty,\mathbb{1}_{\mathbb{R}^2}}(\varphi)\rangle \\ +& \langle \mathrm{D}F^{[\ell'],T,\rho}(\varphi^{T,\rho}), \dot{G}(F^{[\ell''],\infty,\mathbb{1}_{\mathbb{R}^2}}(\varphi) - F^{[\ell''],T,\rho}(\varphi))\rangle \\ +& \langle \mathrm{D}F^{[\ell'],T,\rho}(\varphi^{T,\rho}), \dot{G}(F^{[\ell''],T,\rho}(\varphi) - F^{[\ell''],T,\rho}(\varphi^{T,\rho}))\rangle \end{aligned}$$

and the first term vanishes in $\mathcal{C}([0,S], L^2(\langle x\rangle^{-n}))$ due to the convergence result of Proposition 4.1, the second due to (4.4), the third due to (4.2) combined with the convergence result of Proposition 4.1, and the fourth due to (4.2) and (4.3).

On the other hand,

$$\begin{aligned} & \langle \mathrm{D}\tilde{F}^{\infty,\mathbb{1}_{\mathbb{R}^2}}(\varphi), \dot{G}\tilde{R}\rangle - \langle \mathrm{D}\tilde{F}^{T,\rho}(\varphi^{T,\rho}), \dot{G}\tilde{R}^{T,\rho}\rangle \\ =& \langle \mathrm{D}(\tilde{F}^{\infty,\mathbb{1}_{\mathbb{R}^2}}(\varphi) - \tilde{F}^{T,\rho}(\varphi)), \dot{G}\tilde{R}^{T,\rho}\rangle + \langle \mathrm{D}(\tilde{F}^{T,\rho}(\varphi) - \tilde{F}^{T,\rho}(\varphi^{T,\rho})), \dot{G}\tilde{R}^{T,\rho}\rangle \\ +& \langle \mathrm{D}\tilde{F}^{\infty,\mathbb{1}_{\mathbb{R}^2}}(\varphi), \dot{G}(\tilde{R} - \tilde{R}^{T,\rho})\rangle \end{aligned}$$

and the first term vanishes in strongly in $L^2([0,S];L^2(\langle x\rangle^{-n}))$ due to the convergence result of Proposition 4.1 combined with (4.7), the second due to (4.4) and (4.7). Finally, the third vanishes weakly in $L^2([0,S];L^2(\langle x\rangle^{-n}))$. To see this, denote by

$$A_s\colon L^2(\langle x\rangle^{-n})\to L^2(\langle x\rangle^{-n}),\qquad A_s(X):=\left\langle \mathrm{D}\tilde{F}_s^{\infty,\mathbb{1}_{\mathbb{R}^2}}(\varphi_s),\dot{G}_sX\right\rangle.$$

Then it follows from (4.2) that

$$\sup_{s\in[0,S]}\|A_s\|_{\mathcal{L}(L^2(\langle x\rangle^{-n}))}<\infty.$$

Hence for every $h\in L^2([0,S];L^2(\langle x\rangle^{-n}))$, we have that $A^*h\in L^2([0,S];L^2(\langle x\rangle^{-n}))$, so

$$\int_0^S\left\langle\left\langle \mathrm{D}\tilde{F}_s^{\infty,\mathbb{1}_{\mathbb{R}^2}}(\varphi_s),\dot{G}_s(\tilde{R}_s-\tilde{R}_s^{T,\rho})\right\rangle,h_s\right\rangle_{L^2(\langle x\rangle^{-n})}\mathrm{d}s=\int_0^S\langle\tilde{R}_s-\tilde{R}_s^{T,\rho},A_s^*h_s\rangle_{L^2(\langle x\rangle^{-n})}\mathrm{d}s\longrightarrow 0.$$

Combined with the cancellation induced by the structure of (2.20), the statement follows.

**Convergence of the backward equation**

Take $\mathscr{G}^{T,\rho}$ be the filtration generated by $P^{T,\rho}:=(\varphi_\cdot^{T,\rho},B_\cdot,\tilde{R}_\cdot^{T,\rho})$. Denote by $\mathscr{F}^{T,\rho}$ its right-continuous augmentation.

Fix $S<\infty$ and let $a<b\in[0,S]$. Denoting by $I^{T,\rho}$ the integrand of the backward equation, for any $G$ bounded continuous function and $h\in L^2(\langle x\rangle^{-n})$, by Fubini and the tower rule we have

$$\begin{aligned}
&\mathbb{E}\left[\int_a^b G(P^{T,\rho}|_{[0,a]})\left\langle\tilde{R}_t^{T,\rho}-\int_t^T I_s^{T,\rho}\mathrm{d}s,h\right\rangle_{L^2(\langle x\rangle^{-n})}\mathrm{d}t\right]\\
=&\int_a^b\mathbb{E}\left[\mathbb{E}_t^B\left\langle\mathbb{E}_t^B\left[\int_t^T I_s^{T,\rho}\mathrm{d}s\right]-\int_t^T I_s^{T,\rho}\mathrm{d}s,G(P^{T,\rho}|_{[0,a]})h\right\rangle_{L^2(\langle x\rangle^{-n})}\right]\mathrm{d}t\\
=&\ 0.
\end{aligned}$$

By the convergence of $I^{T,\rho}$, $\tilde{R}^{T,\rho}$ and $P^{T,\rho}$, and the uniform tail bounds (4.7) and (4.5), we may pass to the limit in the above equation. Thus approximating in $L^2(\Omega)$ yields

$$\mathbb{E}\left[X\int_a^b\left\langle\tilde{R}_t-\int_t^\infty I_s^{\infty,\mathbb{1}_{\mathbb{R}^2}}\mathrm{d}s,h\right\rangle_{L^2(\langle x\rangle^{-n})}\mathrm{d}t\right]=0$$

for any bounded $\mathscr{G}_a$-measurable $X$. Now note that the same holds with $\tilde{R}_t$ replaced by $\mathbb{E}_t^{\mathscr{F}}\left[\int_t^\infty I_u^{\infty,\mathbb{1}_{\mathbb{R}^2}}\mathrm{d}u\right]$. Thus considering the process

$$M_t:=\int_0^t\tilde{R}_s-\mathbb{E}_s^{\mathscr{F}}\left[\int_s^\infty I_u^{\infty,\mathbb{1}_{\mathbb{R}^2}}\mathrm{d}u\right]\mathrm{d}s,$$

we have, for any bounded $\mathscr{G}_a$-measurable $X$, that

$$\mathbb{E}[X\langle M_b-M_a,h\rangle_{L^2(\langle x\rangle^{-n})}]=0.$$

Moreover, $M$ is continuous and $\mathscr{F}$-adapted, so also $\mathscr{G}$-adapted. Hence $M$ is a martingale, and since it of bounded variation, it is identically zero. Thus $\tilde{R}_s=\mathbb{E}_s^{\mathscr{F}}\left[\int_s^\infty I_u^{\infty,\mathbb{1}_{\mathbb{R}^2}}\mathrm{d}u\right]$ $\mathrm{d}t\otimes\mathrm{d}\mathbb{P}$-almost everywhere.

**Admissibility of the filtration**

Recall that $\mathscr{F}$ was defined to be the usual augmentation of $\mathscr{G}$, the filtration generated by $P:=(\varphi_\cdot,B_\cdot,\tilde{R}_\cdot)$, so we only need to show that $B$ is an $\mathscr{F}$-martingale. Before removing the cutoffs, since $P^{T,\rho}|_{[0,s]}$ is $\mathscr{F}_s^B$-measurable, we have for any $s<t$

$$\mathrm{Law}(P^{T,\rho}|_{[0,s]},B_t-B_s)=\mathrm{Law}(P^{T,\rho}|_{[0,s]})\otimes\mathrm{Law}(B_t-B_s),$$

which passes to the limit. Now let $A \in \mathscr{F}_s$. For any $s < v < t$, we have $A \in \mathscr{G}_v$, hence

$$\mathbb{E}[\mathbb{1}_A(B_t - B_v)] = 0.$$

We then conclude from the fact that $B_v \to B_s$ in $L^2(\Omega; L^2(\langle x \rangle^{-n}))$ as $v \to s$. $\square$

**Corollary 4.5.** *Let $\varphi$ be as in Theorem 4.4 with associated convergent subsequence indexed by $(T, \rho)$. Then for any $n > 2$ and $\varepsilon > 0$, we have*

$$\mathrm{Law}(\varphi_T^{T,\rho}) \to \mathrm{Law}(\varphi_\infty)$$

*weakly in $H^{-\varepsilon}(\langle x \rangle^{-n})$.*

**Proof.** Write $E := H^{-\varepsilon}(\langle x \rangle^{-n})$ and $Z := \varphi - B$, From (4.1), (4.7) and (A.2) and the embedding $L^2(\langle x \rangle^{-n}) \hookrightarrow E$ we have

$$\sup_{T,\rho} \|\dot{Z}_s^{T,\rho}\|_{L^\infty(\Omega; E)} \lesssim \langle s \rangle^{-1-\delta} \mathfrak{l}_s^c$$

so that $Z_\infty := \lim_{t \to \infty} Z_t$ exists almost surely. Fixing $S < \infty$, we have

$$\begin{aligned} \|Z_T^{T,\rho} - Z_\infty\|_E &\leq \|Z_T^{T,\rho} - Z_S^{T,\rho}\|_E + \|Z_S^{T,\rho} - Z_S\|_E + \|Z_S - Z_\infty\|_E \\ &\lesssim \int_S^T \langle s \rangle^{-1-\delta} \mathfrak{l}_s^c \mathrm{d}s + \|Z_S^{T,\rho} - Z_S\|_E + \int_S^\infty \langle s \rangle^{-1-\delta} \mathfrak{l}_s^c \mathrm{d}s. \end{aligned}$$

Since $Z_S^{T,\rho} \to Z_S$ in $E$ almost surely, the second term vanishes. Then letting $S \to \infty$ shows $\|Z_T^{T,\rho} - Z_\infty\|_E \to 0$ almost surely. The statement then follows from convergence in probability of $B_T$. $\square$

A process $\varphi$ in the above theorem is the (UV-complete) scale decomposition process of the massive sine–Gordon model and $\mathrm{Law}(\varphi_\infty)$ is the sine–Gordon measure. Let us remark that when the coupling constant $\lambda$ is small or the volume is finite, the measure is unique.

From now on, we are going to omit the reference to the removed cutoffs in (4.8) and write simply $\tilde{F}$ instead of $\tilde{F}^{\infty, \mathbb{1}_{\mathbb{R}^2}}$.

Let us collect some immediate consequences of this theorem, namely, basic probabilistic estimates on local functionals of $\varphi$ that will often be used in what follows. All of these bounds follow straightforwardly from the decomposition $\varphi_t = B_t + Z_t$ to a Gaussian $B$ and a bounded drift $Z$.

**Lemma 4.6.** *Let $\varphi = (\varphi_t)_{t \geq 0}$ be the scale decomposition process of sine–Gordon, $x \in \mathbb{R}^2$, $p \in [1, \infty)$. Then for all $s < \infty$, $n \in \mathbb{N}$, $\alpha \in \mathbb{R}$ and multiindex $\mu \in \mathbb{N}^2$ we have*

$$\|(\partial^\mu \varphi_s)^n(x)\|_{L^p(\Omega)} \lesssim_{\lambda, p, \mu} \langle s \rangle^{\frac{n|\mu|}{2}} \mathfrak{l}_s^n, \tag{4.9}$$

$$\|e^{\alpha \varphi_s(x)}\|_{L^p(\Omega)} \lesssim_{\lambda, p, \alpha} e^{\frac{\alpha^2 p}{2} \mathcal{G}_s(0)} \tag{4.10}$$

*and, for some $c \in \mathbb{R}_+$,*

$$\|\partial^\mu e^{\alpha \varphi_s(x)}\|_{L^1(\Omega)} \lesssim_{\lambda, p, \alpha, \mu} e^{\frac{\alpha^2}{2} \mathcal{G}_s(0)} \langle s \rangle^{\frac{|\mu|}{2}} \mathfrak{l}_s^c \tag{4.11}$$

**Proof.** Denote $\dot{Z}_s := -\dot{G}_s(\tilde{F}_s(\varphi_s) + \tilde{R}_s)$. By Young's convolution inequality, (4.1), (4.7), (A.2), and subcriticality $\delta > 0$ we have

$$\begin{aligned} \|\partial^\mu Z_s\|_{L^\infty(\Omega; L^\infty)} &\leq \int_0^t \| \| \partial^\mu \dot{Z}_s \|_{L^\infty} \|_{L^\infty(\Omega)} \mathrm{d}s \\ &\leq \int_0^t \|\partial^\mu k_s^Q\|_{L^1} \big( \big\| Q_s F_s^{[\leq 3]}(\varphi_s) \big\|_{L^\infty(\Omega; L^\infty)} + \|Q_s \tilde{R}_s\|_{L^\infty(\Omega; L^\infty)} \big) \mathrm{d}s \\ &\lesssim \int_0^t \langle s \rangle^{\frac{|\mu|}{2} - 1} \langle s \rangle^{-1} (\lambda_s \mathfrak{l}_s^c + \lambda + \lambda_s^4 \langle s \rangle^{-3}) \mathrm{d}s \\ &\lesssim \langle t \rangle^{\frac{|\mu|}{2}}. \end{aligned} \tag{4.12}$$

When $|\mu|=0$, for any $q\in[1,\infty)$, by Gaussianity of $B$, hypercontractivity and (4.12) we have

$$\|\varphi_s(x)\|_{L^q(\Omega)}\lesssim_q \|B_s(x)\|_{L^2(\Omega)}+\|\|Z_s\|_{L^\infty}\|_{L^\infty(\Omega)}\lesssim\sqrt{\mathcal{G}_s(0)}+1\lesssim\lambda(1+\log(1\vee s)),$$

and the claim then follows from Hölder with $q:=np$. When $|\mu|\geq 1$, since for the $2k^{\text{th}}$ Hermite polynomial it holds that $\mathrm{He}_{2k}(0,\tau)\asymp\tau^k$, we have

$$\begin{aligned}
\mathbb{E}[(\partial^\mu B_t)^2(x)]=(-1)^{|\mu|}(\partial^{2\mu}\mathcal{G}_t)(0) &\asymp \partial_\xi^{2\mu}\left(\int_0^t s^{-1}e^{-\frac{m^2}{s}}e^{-\frac{s}{4}|\xi|^2}\mathrm{d}s\right)\Bigg|_{\xi=0}\\
&\lesssim \int_0^t s^{-1}e^{-\frac{m^2}{s}}\prod_{i=1}^d s^{2\mu_i}\mathrm{He}_{2\mu_i}(\xi_i,s^{-1})e^{-\frac{s}{4}|\xi|^2}\Bigg|_{\xi=0}\mathrm{d}s\\
&\lesssim \int_0^t e^{-\frac{m^2}{s}}s^{|\mu|-1}\mathrm{d}s\\
&\lesssim \langle t\rangle^{|\mu|}.
\end{aligned}$$

The estimate (4.10) follows simply from

$$\|e^{\alpha\varphi_s(x)}\|_{L^p(\Omega)}\lesssim\|e^{\alpha B_s(x)}\|_{L^p(\Omega)}\big\|e^{\alpha\|Z_s\|_{L^\infty(\mathbb{R}^2)}}\big\|_{L^\infty(\Omega)}\lesssim e^{\frac{\alpha^2 p}{2}\mathcal{G}_s(0)}e^{\alpha C_\lambda}.$$

For (4.11), recall that given jointly Gaussian vectors $(X,Y)$, we have

$$\mathbb{E}[e^{\langle a,X\rangle}e^{\langle u,Y\rangle}]=e^{\frac{1}{2}\langle a,\Sigma_X a\rangle}\mathbb{E}[e^{\langle u,Y+\Sigma_{XY}a\rangle}]$$

for every $u$, where $\Sigma_X$ and $\Sigma_{XY}$ are the corresponding (cross-)covariance matrices. Thus

$$\mathbb{E}[e^{\langle a,X\rangle}F(Y)]=e^{\frac{1}{2}\langle a,\Sigma_X a\rangle}\mathbb{E}[F(Y+\Sigma_{XY}a)]$$

for any polynomial $F$, so that denoting by $P_\mu$ the polynomials obtained from Faà di Bruno,

$$\begin{aligned}
\|\partial^\mu e^{\alpha B_t(x)}\|_{L^1(\Omega)} &= \mathbb{E}[e^{\alpha B_t(y)}|P_\mu(\{B_t^\nu(x)\}_{\nu\leq\mu})|]\\
&= e^{\frac{\alpha^2}{2}\mathcal{G}_t(0)}\|P_\mu(\{\partial^\nu B_t(x)+a\,\mathrm{Cov}(\partial^\nu B_t(x),B_t(x))\}_{\nu\leq\mu})\|_{L^1(\Omega)},
\end{aligned}$$

and since

$$\mathrm{Cov}(\partial^\nu B_t(x),B_t(x))=(-1)^{|\nu|}\partial^\nu\mathcal{G}_t(0)\lesssim\langle t\rangle^{\frac{|\nu|}{2}}\mathfrak{l}_t^c,$$

we get

$$\|\partial^\mu e^{\alpha B_t(x)}\|_{L^1(\Omega)}\lesssim e^{\frac{\alpha^2}{2}\mathcal{G}_t(0)}\langle t\rangle^{\frac{|\mu|}{2}}\mathfrak{l}_t^c.$$

Then using the above $L^\infty$ bounds on $Z$,

$$\begin{aligned}
\|\partial^\mu e^{\alpha\varphi_t(x)}\|_{L^1(\Omega)} &= \|\partial^\mu(e^{\alpha B_t(x)}e^{\alpha Z_t(x)})\|_{L^1(\Omega)}\\
&\lesssim \sum_{\nu\leq\mu}\binom{\mu}{\nu}\|\partial^\nu e^{\alpha B_t(x)}\|_{L^1(\Omega)}\|\partial^{\mu-\nu}e^{\alpha Z_t(x)}\|_{L^\infty(\Omega)}\\
&\lesssim e^{\alpha C_\lambda}e^{\frac{\alpha^2}{2}\mathcal{G}_t(0)}\langle t\rangle^{\frac{|\mu|}{2}}\mathfrak{l}_t^c
\end{aligned}$$

as requested. □

**Proposition 4.7. (existence of the effective force observable and its remainder)** *Let $\beta^2<6\pi$, $y\in\mathbb{R}^2$. Then*

1. *for any $\tilde{\ell}\in\mathbb{N}$ with $\tilde{\ell}>\delta^{-1}-1$, the germ $F^{[\leq\tilde{\ell}]}(y)$ is coherent, and for any $\ell'\geq\tilde{\ell}$, it satisfies (2.21).*
2. *for any $\ell^*\in\mathbb{N}_0$, $\dot{G}R^{[>\ell^*]}$, given by (2.22), is well-defined and satisfies (2.26).*

**Proof.** Using the symmetry of $v_t^{[\ell](\boldsymbol{a})}$ under permutations of indices (which extends to the UV limit by the proof of Corollary (3.15)),

$$\begin{aligned} F_t^{[\ell]}(\phi) &= i\beta \sum_{\substack{\boldsymbol{a}\in\mathfrak{M}\\ \ell(\boldsymbol{a})=\ell}} \int \rho(x_{1:\ell}) v_t^{[\ell](\boldsymbol{a})} M^{(\boldsymbol{a})}(\xi_{1:\ell}) \sum_{j=1}^{\ell} \sigma_j \delta_{x_j} \psi_{1:\ell} \mathrm{d}\xi_{1:\ell} \\ &= i\beta\ell \sum_{\substack{\boldsymbol{a}\in\mathfrak{M}\\ \ell(\boldsymbol{a})=\ell}} \int \rho(x_{1:\ell}) v_t^{[\ell](\boldsymbol{a})} M^{(\boldsymbol{a})}(\xi_{1:\ell}) \sigma_1 \delta_{x_1} \psi_{1:\ell} \mathrm{d}\xi_{1:\ell}. \\ &=: \sum_{\substack{\boldsymbol{a}\in\mathfrak{M}\\ \ell(\boldsymbol{a})=\ell}} F_t^{[\ell](\boldsymbol{a})}(\phi). \end{aligned}$$

We may represent $F_t^{[\ell](\boldsymbol{a})}$ as a function $y \mapsto K_t^{F,[\ell](\boldsymbol{a})}[\phi](y)$ over $\mathbb{R}^d$ via the duality

$$\langle K_t^{F,[\ell](\boldsymbol{a})}[\phi], h \rangle_{L^2(\mathbb{R}^d)} = \langle F_t^{[\ell](\boldsymbol{a})}(\phi), h \rangle_{\mathcal{S},\mathcal{S}'}$$

for $h \in \mathcal{S}$. That is,

$$\begin{aligned} &\int_{\mathbb{R}^d} K_t^{F,[\ell](\boldsymbol{a})}[\phi](y) h(y) \mathrm{d}y \\ = &\int h(x_1)\rho(x_1) \sum_{\sigma_1\in\{-1,1\}} i\beta\ell\sigma_1\psi_1 \int \rho(x_{2:\ell}) v_t^{[\ell](\boldsymbol{a})} M^{(\boldsymbol{a})}(\xi_1, \xi_{2:\ell}) \psi_{2:\ell} \mathrm{d}\xi_{2:\ell} \mathrm{d}x_1 \end{aligned}$$

which yields

$$K_t^{F,[\ell](\boldsymbol{a})}[\phi](y) = \rho(y) i\beta\ell \sum_{\sigma_1\in\{-1,1\}} \sigma_1 \psi_{(\sigma_1,y)} \int \rho(x_{2:\ell}) v_t^{[\ell](\boldsymbol{a})} M^{(\boldsymbol{a})}(\sigma_1, y, \xi_{2:\ell}) \psi_{2:\ell} \mathrm{d}\xi_{2:\ell}.$$

We abuse notation and write $F_t^{[\ell](\boldsymbol{a})}(\phi)$ to mean the function $K_t^{F,[\ell](\boldsymbol{a})}[\phi]$. For $\boldsymbol{a} \in \mathfrak{M} \setminus \mathfrak{M}_p$, we then simply have, uniformly in $\|\rho\|_{L^\infty} \leq 1$, that

$$\left|F_t^{[\ell](\boldsymbol{a})}(\phi)(y)\right| \lesssim \sup_{\xi_1} \int \left|v_t^{[\ell](\boldsymbol{a})}(\xi_{1:\ell})\right| \mathrm{d}\xi_{2:\ell} \lesssim \|\!|v_{t,T}^{[\ell](\boldsymbol{a})}|\!\|_{t,\Theta,\eta},$$

whereas for $\boldsymbol{a} \in \mathfrak{M}_p$, with $\varkappa := \varkappa(\boldsymbol{a}) = 1$, by neutrality and charge conjugation symmetry of $v_{t,T}^{[2](\varkappa)}$ (cf. Remark 3.11 and proof of Corollary (3.15))

$$\begin{aligned} &F_t^{[\ell](\boldsymbol{a})}(\phi)(y) \\ = \;& \rho(y) i\beta\ell e^{+i\beta\phi(y)} \int \rho(x_2) v_t^{[\ell](\boldsymbol{a})} M^{(\boldsymbol{a})}(+, y, -, x_2) e^{-i\beta\phi(x_2)} \mathrm{d}x_2 \\ - \;& \rho(y) i\beta\ell e^{-i\beta\phi(y)} \int \rho(x_2) v_t^{[\ell](\boldsymbol{a})} M^{(\boldsymbol{a})}(-, y, +, x_2) e^{i\beta\phi(x_2)} \mathrm{d}x_2 \\ = \;& \rho(y) i\beta\ell \left( \int \rho(x_2) v_t^{[\ell](\boldsymbol{a})} M^{(\boldsymbol{a})}(+, y, x_2, -) (e^{i\beta\phi(y)} e^{-i\beta\phi(x_2)} - e^{-i\beta\phi(y)} e^{i\beta\phi(x_2)}) \mathrm{d}x_2 \right) \end{aligned}$$

and by the neutrality condition, zeroth order Taylor expansion gives

$$\begin{aligned} g(y, x_2) &:= (e^{i\beta\phi(y)} e^{-i\beta\phi(x_2)} - e^{-i\beta\phi(y)} e^{i\beta\phi(x_2)}) |y - x_2|^{-\varkappa} \\ &= (1 - 1) - i\beta \frac{(x_2 - y)}{|y - x_2|^{\varkappa}} \cdot \int_0^1 \left( e^{i\beta(\phi(y) - \phi(x_2^\vartheta))} + e^{-i\beta(\phi(y) - \phi(x_2^\vartheta))} \right) \nabla\phi(x_2^\vartheta) \mathrm{d}\vartheta \end{aligned}$$

where we denoted $x_2^\vartheta := y + \vartheta(x_2 - y)$.
Thus Minkowski's integral inequality, Proposition 3.14 and Lemma 4.6 yields that for any $p \in [1, \infty)$

$$\left\|F_t^{[\ell](\boldsymbol{a})}(\varphi_t)(y)\right\|_{L^p(\Omega)} \lesssim \|\!|v_t^{[\ell](\boldsymbol{a})}|\!\|_{t,\Theta(\ell),\eta} \sum_{|\nu|=1} \|\partial^\nu \varphi_t(x)\|_{L^p(\Omega)} \lesssim \lambda_t (\lambda_t \langle t \rangle^{-1})^{\ell-1} \mathfrak{l}_t^c.$$

It follows that for any $\ell' \geq \tilde{\ell}$ we have

$$\begin{aligned} \left\| F_t^{[\leq \tilde{\ell}]}(\varphi_t)(y) - F_t^{[\leq \ell']}(\varphi_t)(y) \right\|_{L^1(\Omega)} &\leq \sum_{\ell=\tilde{\ell}+1}^{\ell'} \| F_t^{[\ell]}(\varphi_t)(y) \|_{L^1(\Omega)} \\ &\lesssim \langle t \rangle^{1-\delta} \langle t \rangle^{-\delta \tilde{\ell}} \mathfrak{c}_t^c \end{aligned} \tag{4.13}$$

which vanishes as $t \to \infty$, since $\tilde{\ell} + 1 > \delta^{-1}$.

Recall that it follows from (2.20) that $F^{[\ell]} = \mathrm{D}V^{[\ell]}$ satisfies the flow equation

$$\partial_t F_t^{[\ell]} + \frac{1}{2} \Delta_{\dot{G}_t} F_t^{[\ell]} = \frac{1}{2} \sum_{\ell' + \ell'' = \ell} \mathrm{D} \langle F_t^{[\ell']}, \dot{G}_t F_t^{[\ell'']} \rangle = \sum_{\ell' + \ell'' = \ell} \langle \mathrm{D} F_t^{[\ell']}, \dot{G}_t F_t^{[\ell'']} \rangle.$$

Then for any $\ell^* \in \mathbb{N}$, we have

$$\begin{aligned} &= \left( \partial_s + \frac{1}{2} \Delta_{\dot{G}_s} \right) F_s^{[\leq \ell^*]} - \langle \mathrm{D} F_s^{[\leq \ell^*]}(\varphi_s), \dot{G}_s (F_s^{[\leq 3]}(\varphi_s) + \tilde{R}_s) \rangle (y) \\ &= \sum_{\ell' + \ell'' \leq \ell^*} \langle \mathrm{D} F_s^{[\ell']}(\varphi_s), \dot{G}_s F_s^{[\ell'']}(\varphi_s) \rangle (y) - \sum_{\ell' \leq \ell^*} \sum_{\ell'' \leq 3} \langle \mathrm{D} F_s^{[\ell']}(\varphi_s), \dot{G}_s F_s^{[\ell'']}(\varphi_s) \rangle (y) \\ &- \langle \mathrm{D} F_s^{[\leq \ell^*]}(\varphi_s), \dot{G}_s \tilde{R}_s \rangle (y) \\ &= \sum_{\substack{\ell'' \geq 4 \\ \ell' + \ell'' \leq \ell^*}} \langle \mathrm{D} F_s^{[\ell']}(\varphi_s), \dot{G}_s F_s^{[\ell'']}(\varphi_s) \rangle (y) - \sum_{\substack{\ell' \leq \ell^*, \ell'' \leq 3 \\ \ell' + \ell'' \geq \ell^* + 1}} \langle \mathrm{D} F_s^{[\ell']}(\varphi_s), \dot{G}_s F_s^{[\ell'']}(\varphi_s) \rangle (y) \\ &- \langle \mathrm{D} F_s^{[\leq \ell^*]}(\varphi_s), \dot{G}_s \tilde{R}_s \rangle (y) \end{aligned}$$

and plugging in (4.2), (4.5) and the bound on $\tilde{R}$,

$$\begin{aligned} \sum_{\substack{\ell'' \geq 4 \\ \ell' + \ell'' \leq \ell^*}} \| \langle \mathrm{D} F_s^{[\ell']}(\varphi_s), \dot{G}_s F_s^{[\ell'']}(\varphi_s) \rangle (y) \|_{L^1(\Omega)} &\lesssim \langle t \rangle^{-5\delta} \mathfrak{c}_t^c \\ \sum_{\substack{\ell' \leq \ell^*, \ell'' \leq 3 \\ \ell' + \ell'' \geq \ell^* + 1}} \| \langle \mathrm{D} F_s^{[\ell']}(\varphi_s), \dot{G}_s F_s^{[\ell'']}(\varphi_s) \rangle (y) \|_{L^1(\Omega)} &\lesssim \langle t \rangle^{-(\ell^* + 1)\delta} \mathfrak{c}_t^c \\ \| \langle \mathrm{D} F_s^{[\leq \ell^*]}(\varphi_s), \dot{G}_s \tilde{R}_s \rangle (y) \|_{L^1(\Omega)} &\lesssim \langle t \rangle^{-1-\delta} \mathfrak{c}_t^c \end{aligned}$$

so that $F_s^{[\leq \ell^*]}(y)$ satisfies (2.10) provided $\ell^* \geq \tilde{\ell} > \delta^{-1} - 1$. Conditions 1. and 2. of Definition 2.5 are checked in a straightforward way, see Lemma 6.1 below. This proves item 1.

To prove item 2, note that it follows from item 1, in particular (4.13), that the reconstructed observable of the germ $F^{[\leq \ell^*]}(y)$ is the same as $\tilde{F}_s(\varphi_s) + \tilde{R}_s$. Therefore, for any $\ell^* \geq \tilde{\ell}$, we may use the germ $F^{[\leq \ell^*]}$ in the generator $\mathscr{L}^V$ and compute

$$\begin{aligned} R_t^{[> \ell^*]} &= \mathbb{E}_t \int_t^\infty \mathscr{L}^V F_s^{[\leq \ell^*]}(\varphi_s) \mathrm{d}s \\ &= \mathbb{E}_t \int_t^\infty \left( \partial_s + \frac{1}{2} \Delta_{\dot{G}_s} \right) F_s^{[\leq \ell^*]} - \langle \mathrm{D} F_s^{[\leq \ell^*]}(\varphi_s), \dot{G}_s (F_s^{[\leq \ell^*]}(\varphi_s) + R_s^{[> \ell^*]}) \rangle \mathrm{d}s \\ &= \mathbb{E}_t \int_t^\infty - \sum_{\substack{\ell', \ell'' \leq \ell^* \\ \ell' + \ell'' > \ell^*}} \langle \mathrm{D} F_s^{[\ell']}(\varphi_s), \dot{G}_s F_s^{[\ell'']}(\varphi_s) \rangle - \sum_{\ell \leq \ell^*} \langle \mathrm{D} F_s^{[\ell]}(\varphi_s), \dot{G}_s R_s^{[> \ell^*]} \rangle \mathrm{d}s. \end{aligned} \tag{4.14}$$

Now (4.2) yields

$$\begin{aligned} \sum_{\ell \leq \ell^*} \| \langle \mathrm{D} F_t^{[\ell]}, \dot{G}_t R_t^{[> \ell^*]} \rangle \|_{L^\infty(\Omega; L^\infty)} &\lesssim \sum_{\ell \leq \ell^*} \lambda_t (\lambda_t \langle t \rangle^{-1})^{\ell - 1} \langle t \rangle^{-1} \mathfrak{c}_t^c \| Q_t R_t^{[> \ell^*]} \|_{L^\infty(\Omega; L^\infty)} \\ &\lesssim \lambda_t \langle t \rangle^{-2} \mathfrak{c}_t^c \| R_t^{[> \ell^*]} \|_{L^\infty(\Omega; L^\infty)} \end{aligned}$$

where we used that the smoothing effect of the heat kernel forces the right-hand side of (4.14), hence $R_t^{[>\ell^*]}$, to be function-valued. Plugging this bound and (4.5) into (4.14),

$$\|R_t^{[>\ell^*]}\|_{L^\infty(\Omega;L^\infty)} \lesssim \int_t^\infty \langle s\rangle^{-\delta(\ell^*+1)}\mathfrak{l}_s^c + \langle t\rangle^{-1-\delta}\mathfrak{l}_s^c\|R_s^{[>\ell^*]}\|_{L^\infty(\Omega;L^\infty)}\mathrm{d}s$$

hence by backward Grönwall

$$\begin{aligned}\|R_t^{[>\ell^*]}\|_{L^\infty(\Omega;L^\infty)} &\lesssim \int_t^\infty \langle s\rangle^{-\delta(\ell^*+1)}\mathfrak{l}_s^c e^{\int_t^s \langle u\rangle^{-1-\delta}\mathfrak{l}_u^c \mathrm{d}u}\mathrm{d}s \\ &\lesssim \langle t\rangle^{1-\delta(\ell^*+1)}\mathfrak{l}_t^c\end{aligned} \tag{4.15}$$

and (2.26) follows from Young's inequality and $\|\dot{\mathcal{G}}_t\|_{L^1} \lesssim \langle t\rangle^{-2}$. When $\ell^* < \tilde{\ell}$, combine the estimate just derived with (4.1). □

Similarly to the Taylor expansion for effective force in the previous proposition, we may obtain $L^p(\Omega)$ bounds on $V_t^{[\ell],T}(\varphi_t)$ that, however, depend on the volume cutoff $\rho$. We will not use these estimates in the following and hence invite the reader to skip the remark below.

**Remark 4.8.** It holds for any $p \in [1,\infty)$ uniformly in $T$ that

$$\|V_t^{[\ell],T}(\varphi_t)\|_{L^p(\Omega)} \lesssim \|\rho\|_{L^1}\lambda_t(\lambda_t\langle t\rangle^{-1})^{\ell-1}\mathfrak{l}_t^c. \tag{4.16}$$

Indeed, for $\boldsymbol{a} \in \mathfrak{M}\setminus\mathfrak{M}_p$, we have via $\sup_{x\in\mathbb{R}^d}|\exp(i\beta\phi(x))| \leq 1$, $\rho(x) \in [0,1]$ and $\omega_t^\varsigma \geq 1$ that

$$|V_t^{[\ell](\boldsymbol{a}),T}(\phi)| \leq \int |v_{t,T}^{[\ell](\boldsymbol{a})}(\xi_{1:\ell})|\rho(x_{1:\ell})\mathrm{d}\xi_{1:\ell} \lesssim \|\rho\|_{L^1}|||v_{t,T}^{[\ell](\boldsymbol{a})}|||_{t,\Theta}$$

and the statement follows from Proposition 3.14.

For $\boldsymbol{a} \in \mathfrak{M}_p$ with $\varkappa(\boldsymbol{a}) = 0$ we obtain trivially $|V_t^{[\ell](\boldsymbol{a}),T}(\phi)| = |v_{t,T}^{[2](0)}| = |||v_{t,T}^{[\ell](\boldsymbol{a})}|||_{t,\Theta}$ and again conclude by Proposition 3.14. Given a neutral pair of vertex fields, Taylor expansion around $x_1$ yields by the neutrality condition $\sigma_1 = -\sigma_2$ that

$$\psi_{1:2} = 1 + i\beta\sigma_2\nabla\phi(x_1)\cdot(x_2-x_1) + \psi_1\sum_{|\nu|=2}\frac{2}{\nu!}(x_2-x_1)^\nu\int_0^1(1-\vartheta)\partial^\nu\psi|_{(\sigma_2,x_1+\vartheta(x_2-x_1))}\mathrm{d}\vartheta \tag{4.17}$$

thus for an $\boldsymbol{a} \in \mathfrak{M}_p$, with $\varkappa := \varkappa(\boldsymbol{a}) > 0$, we obtain

$$\begin{aligned}&|V_t^{[\ell](\boldsymbol{a}),T}(\phi)| \\ =\;& \left|\int v_{t,T}^{[2](\varkappa)}(\xi_{1:2})\rho(x_{1:2})\frac{\psi_{\xi_1}\psi_{\xi_2}-1}{|x_1-x_2|^\varkappa}\mathrm{d}\xi_{1:2}\right| \\ \leq\;& \left|\int v_{t,T}^{[2](\varkappa)}(\xi_{1:2})\rho(x_{1:2})\frac{i\beta\sigma_2\nabla\phi(x_1)\cdot(x_2-x_1)}{|x_1-x_2|^\varkappa}\mathrm{d}\xi_{1:2}\right| \\ +\;& \left|\int v_{t,T}^{[2](\varkappa)}(\xi_{1:2})\rho(x_{1:2})\psi_{\xi_1}\sum_{|\nu|=2}\frac{2}{\nu!}\frac{(x_2-x_1)^\nu}{|x_1-x_2|^\varkappa}\int_0^1(1-\vartheta)\partial^\nu\psi|_{(\sigma_2,x_1+\vartheta(x_2-x_1))}\mathrm{d}\vartheta\mathrm{d}\xi_{1:2}\right|\end{aligned}$$

with $\varkappa \leq 2$. The first term is zero by the charge conjugation symmetry of $v_{t,T}^{[2](\varkappa)}$ (cf. Remark 3.11); hence by Minkowski's integral inequality, (3.20), Lemma 4.6

$$\begin{aligned}\|V_t^{[\ell](\boldsymbol{a}),T}(\varphi_t)\|_{L^p(\Omega)} &\lesssim \|\rho\|_{L^1}|||v_{t,T}^{[2](\varkappa)}|x_1-x_2|^{2-\varkappa}|||_{t,\Theta}\sum_{|\nu|=2}\|\partial^\nu\varphi_t\|_{L^p(\Omega)} \\ &\lesssim \|\rho\|_{L^1}|||v_{t,T}^{[2](\varkappa)}|||_{t,\Theta}\langle t\rangle^{-\frac{2-\varkappa}{2}}\langle t\rangle^1\end{aligned}$$

and conclude from Proposition 3.14.

Note that the above argument with the first-order Taylor expansion in (4.17) allows for curing a divergence of order up to $\varkappa \leq 2$. Since we restrict ourselves to $\beta^2 < 6\pi$, we could have replaced (4.17) with zeroth-order Taylor expansion to offset the divergence of order $\varkappa = 1$; in this case we would not have had to rely on Remark 3.11.

# 5 Observable flow equations

## 5.1 Polynomial observables

In this section we construct all polynomial observables of the sine–Gordon field by carrying out the program of Section 2.6. Since the field $\varphi_t(y)$ has only logarithmic blowup in the UV $t \to \infty$, simple Wick ordering is sufficient to make sense of the observable-inserted correlation functions $\langle [\![\phi^n(y)]\!] \star [\![\phi(x_1)]\!] \ldots [\![\phi(x_n)]\!] \rangle$ for any $n \in \mathbb{N}$; in turn, this is already computed by level $\ell = 1$ approximation in the flow equation. However, as we want to prove the prefactorization algebra property of the observables we construct (Theorem 1.2), and in the subsequent sections we will consider observables $\mathfrak{Q}(z)$ of arbitrarily high superficial degree of divergence, we need to develop the perturbative approximation to arbitrary precision even for polynomials in order to conclude the existence of the products $\langle [\![\phi^n(y)]\!] \star \mathfrak{Q}(z) \rangle$.

From now on, we fix $y \in \mathbb{R}^d$, the point where the observable is localized. Using the notations of the previous section, we make the ansatz

$$q_t^{[\ell],T}(\phi) = \int \Gamma_{t,T}^{[\ell]}(\phi, \xi_{1:(\ell-1)}, y) \psi_{1:(\ell-1)}^{\rho} \mathrm{d}\xi_{1:(\ell-1)} \tag{5.1}$$

where we denoted $\psi_{1:(\ell-1)}^{\rho} := \psi_{1:(\ell-1)} \rho(x_{1:(\ell-1)})$, and let

$$\Gamma_{t,T}^{[\ell]}(\phi, \xi_{1:(\ell-1)}, y) = \sum_{k=0}^{n} \gamma_{t,T}^{[\ell](k)}(\xi_{1:(\ell-1)}, y) \phi^k(y) \tag{5.2}$$

with kernels $\gamma_{t,T}^{[\ell](k)}$ to be determined. That is, the start of our induction is the first-order representation

$$q_t^{[1],T}(\phi) = \sum_{k=0}^{n} \gamma_{t,T}^{[1](k)}(y) \phi^k(y)$$

and the periodic structure, i.e. dependence on the vertex fields $\psi_{1:(\ell-1)}$ in $q_t^{[\ell],T}$ enters only at $\ell \geq 2$ iteration order through contact with the Mayer-expanded potential in the flow equation (2.19). It is straightforward to verify that the above representation is compatible with the structure of the flow equation (2.19). To this end, recall the definition of $\mathbb{G}_s$ from (3.2); we now also introduce

$$\dot{\not{\mathbb{G}}}_s(\xi_I, y) := i \frac{\beta}{2} \sum_{j \in I} \sigma_j \dot{\mathcal{G}}_s(y - x_j)$$

for $I \subseteq \mathbb{N}$. Due to the nontrivial post-processing of the potential's basis, we also define

$$M^{(\boldsymbol{b})}(x_I) := \mathbb{1}_{\boldsymbol{b} \in \mathfrak{M} \setminus \mathfrak{M}_p} + \mathbb{1}_{\boldsymbol{b} \in \mathfrak{M}_p} \mathbb{1}_{\varkappa(\boldsymbol{b}) > 0} |x_{I(1)} - x_{I(2)}|^{-\varkappa(\boldsymbol{b})}. \tag{5.3}$$

Finally, given a finite index set $I \subseteq \mathbb{N}$ we denote by $\xi_{I-1}$ for the tuple obtained by omitting one index from $I$. When such notation appears in the formulae below, each summand is understood after consistently relabeling all of its remaining variables.

**Lemma 5.1.** *Let $V^{[\ell]}$ be given by (3.31). Suppose that the kernels $\gamma_{t,T}^{[\ell](k)}$ satisfy the equations*

$$\gamma_{t,T}^{[\ell](k)}(\xi_{1:(\ell-1)}, y) = e^{\mathbb{G}_{t,T}(\xi_{1:(\ell-1)})} \gamma_{T,T}^{[\ell](k)}(\xi_{1:(\ell-1)}, y) + \int_t^T e^{\mathbb{G}_{t,s}(\xi_{1:(\ell-1)})} \mathscr{S}_{s,T}^{[\ell](k)}(\xi_{1:(\ell-1)}, y) \mathrm{d}s \tag{5.4}$$

*where*

$$
\begin{aligned}
& \mathscr{S}^{[\ell](k)}_{s,T}(\xi_{1:(\ell-1)},y) \\
:= \; & \frac{(k+1)(k+2)}{2}\dot{\mathcal{G}}_s(0)\gamma^{[\ell](k+2)}_{s,T}(\xi_{1:(\ell-1)},y) \\
+ \; & 2(k+1)\dot{\not{\mathcal{G}}}_s(\xi_{1:(\ell-1)},y)\gamma^{[\ell](k+1)}_{s,T}(\xi_{1:(\ell-1)},y) \\
- \; & \sum_{I_1\uplus I_2=1:\ell}(k+1)C_{|I_1|,|I_2|}\sum_{\substack{\boldsymbol{b}\in\mathfrak{M}:\\ \ell(\boldsymbol{b})=|I_2|}}\big(v^{[|I_2|](\boldsymbol{b})}_s M^{(\boldsymbol{b})}\big)(\xi_{I_2})\dot{\not{\mathcal{G}}}_s(\xi_{I_2},y)\gamma^{[|I_1|](k+1)}_{s,T}(\xi_{I_1-1},y) \qquad (5.5)\\
- \; & \sum_{I_1\uplus I_2=1:\ell}C_{|I_1|,|I_2|}\sum_{\substack{\boldsymbol{b}\in\mathfrak{M}:\\ \ell(\boldsymbol{b})=|I_2|}}\big(v^{[|I_2|](\boldsymbol{b})}_s M^{(\boldsymbol{b})}\big)(\xi_{I_2})\dot{\mathbb{G}}_s(\xi_{I_2},\xi_{I_1-1})\gamma^{[|I_1|](k)}_{s,T}(\xi_{I_1-1},y)
\end{aligned}
$$

*for some positive combinatorial coefficients* $C_{|I_1|,|I_2|}$. *Then* $q^{[\ell],T}_t$*, given by* (5.1) *and* (5.2)*, satisfies the observable flow equation* (2.19)*.*

**Proof.** A straightforward computation; see Appendix B.2 for details. □

Imposing the boundary conditions

$$
\begin{aligned}
\gamma^{[1](k)}_{0,T} &:= \delta_{n,k} \\
\gamma^{[\ell](k)}_{T,T} &:= 0 \qquad \text{if } \ell\geq 2,
\end{aligned} \tag{5.6}
$$

the above hierarchy of equations allows us to estimate the germ kernels $\gamma^{[\ell](k)}_{t,T}$ inductively using the already established bounds on the potential kernels. We will estimate in the norm

$$
|||\gamma^{[\ell](k)}_{t,T}|||_{t,\Theta,\eta} := \int \big|\gamma^{[\ell](k)}_{t,T}(\xi_{1:(\ell-1)},y)\big|\omega^{\Theta}_t(x_{1:(\ell-1)},y)\langle \mathrm{St}(x_{1:(\ell-1)},y)\rangle^{\eta}\mathrm{d}\xi_{1:(\ell-1)} \tag{5.7}
$$

where $\omega^{\Theta}_u$ is given by (3.4) with the control chosen to be $\Theta := \Theta(\ell) := 4^{-\ell}$ as before, and $\eta\geq 0$ to be chosen later. The purpose of the weights $\langle \mathrm{St}(x_{1:(\ell-1)},y)\rangle^{\eta}$ is to control the infinite volume limit.

**Proposition 5.2.** *Let* $\beta^2<6\pi$. *Then there exist kernels* $\gamma^{[\ell](k)}_{t,T}$ *satisfying* (5.4) *with boundary conditions* (5.6)*, and it holds for any* $\eta\geq 0$ *that*

$$
|||\gamma^{[\ell](k)}_{t,T}|||_{t,\Theta(\ell),\eta} \lesssim (\lambda_t\langle t\rangle^{-1})^{\ell-1}\mathfrak{l}^c_t \tag{5.8}
$$

*for some* $c=c(n,k,\ell,\beta^2)\in\mathbb{R}_+$. *Moreover, the* $|||\cdot|||_{t,\Theta,\eta}$*-limits*

$$
\gamma^{[\ell](k)}_t := \lim_{T\to\infty}\gamma^{[\ell](k)}_{t,T}
$$

*exist, are independent of* $\eta$ *and satisfy the bounds* (5.8)*.*

**Proof.** We proceed inductively in $\ell,k$; the starting point is $(\ell,k)=(1,n)$; from which we descend $k\to k-1$ until $k=0$ is reached, from which we go to $(\ell,k)=(2,n)$ and so forth. As before, we do not track the log-dependence, hence the constant $c\in\mathbb{R}_+$ in $\mathfrak{l}^c_t$ may change from line to line.

Clearly for $\ell=1$ the flow equation computes the Gaussian renormalisation counterterms, hence by Lemma A.1

$$
|||\gamma^{[1](k)}_{t,T}|||_{t,\Theta(1),\eta} = \big|\gamma^{[1](k)}_{t,T}(y)\big| \lesssim_{n,k} 1+\mathcal{G}_t(0)^n \lesssim \mathfrak{l}^c_t.
$$

For the inductive step, we can bound the source term as follows. Of course

$$
\begin{aligned}
|||\dot{\mathcal{G}}_s(0)\gamma^{[\ell](k+2)}_{s,T}(\xi_{1:(\ell-1)},y)|||_{s,\Theta,\eta} &\lesssim \|\dot{\mathcal{G}}_s\|_{L^\infty}|||\gamma^{[\ell](k+2)}_{s,T}|||_{s,\Theta,\eta} \lesssim \langle s\rangle^{-1}(\lambda_s\langle s\rangle^{-1})^{\ell-1}\mathfrak{l}^c_s \\
|||\dot{\not{\mathcal{G}}}_s(\xi_{1:(\ell-1)},y)\gamma^{[\ell](k+1)}_{s,T}(\xi_{1:(\ell-1)},y)|||_{s,\Theta,\eta} &\lesssim \|\dot{\mathcal{G}}_s\|_{L^\infty}|||\gamma^{[\ell](k+1)}_{s,T}|||_{s,\Theta,\eta} \lesssim \langle s\rangle^{-1}(\lambda_s\langle s\rangle^{-1})^{\ell-1}\mathfrak{l}^c_s
\end{aligned}
$$

For the terms involving the potential kernels, note that for nontrivial $M_{(\boldsymbol{b})}(x_{I_2})$ we can use the technology of Lemma 3.12. Indeed, from (3.25) with $c := 2\cdot 4^{-\ell}$ and (3.29) we see

$$\begin{aligned}
& |||v_s^{[|I_2|](\boldsymbol{b})}(\xi_{I_2})M^{(\boldsymbol{b})}(x_{I_2})\dot{\mathbb{G}}_s(\xi_{I_2},\xi_{I_1-1})\gamma_{s,T}^{[|I_1|](k)}(\xi_{I_1-1},y)|||_{s,4^{-\ell},\eta}\\
\lesssim\ & |||v_s^{[|I_2|](\boldsymbol{b})}|||_{s,3\cdot 4^{-\ell},\eta}\langle s\rangle^{-2+\frac{\varkappa(\boldsymbol{b})}{2}}|||\gamma_{s,T}^{[|I_1|](k)}|||_{s,4^{-\ell},\eta}\\
\lesssim\ & |||v_s^{[|I_2|](\boldsymbol{b})}|||_{s,4^{-(\ell-1)},\eta}\langle s\rangle^{-2+\frac{\varkappa(\boldsymbol{b})}{2}}|||\gamma_{s,T}^{[|I_1|](k)}|||_{s,4^{-\ell},\eta}\\
\lesssim\ & \lambda_s(\lambda_s\langle s\rangle^{-1})^{|I_2|-1}\langle s\rangle^{-2}(\lambda_s\langle s\rangle^{-1})^{|I_1|-1}\mathfrak{l}_s^c
\end{aligned}$$

since $|I_2|\leq \ell-1$ apart from the trivial case $|I_2|=\ell$ when $|I_1|=0$ so that $q^{[0]}\equiv 0$. The remaining source term is bounded similarly via a straightforward modification of (3.25) to $\dot{\not{\mathbb{G}}}_s(\xi_{1:2},y)$:

$$\begin{aligned}
& |||v_s^{[|I_2|](\boldsymbol{b})}(\xi_{I_2})M^{(\boldsymbol{b})}(x_{I_2})\dot{\not{\mathbb{G}}}_s(\xi_{I_2},y)\gamma_{s,T}^{[|I_1|](k+1)}(\xi_{I_1-1},y)|||_{s,4^{-\ell},\eta}\\
\lesssim\ & |||v_s^{[|I_2|](\boldsymbol{b})}|||_{s,4^{-(\ell-1)},\eta}\langle s\rangle^{-2+\frac{\varkappa(\boldsymbol{b})}{2}}|||\gamma_{s,T}^{[|I_1|](k+1)}|||_{s,4^{-\ell},\eta}\\
\lesssim\ & \lambda_s(\lambda_s\langle s\rangle^{-1})^{|I_2|-1}\langle s\rangle^{-2}(\lambda_s\langle s\rangle^{-1})^{|I_1|-1}\mathfrak{l}_s^c.
\end{aligned}$$

Combining these estimates and inserting them to (5.5) and (5.4), we obtain

$$\begin{aligned}
|||\gamma_{t,T}^{[\ell](k)}|||_{t,\Theta,\eta} &\leq \int_t^\infty \|e^{\mathbb{G}_{t,s}(\xi_{1:(\ell-1)})}\|_{L^\infty}|||\mathscr{S}_{s,T}^{[\ell](k)}|||_{t,\Theta,\eta}\mathrm{d}s\\
&\lesssim \int_t^\infty |||\mathscr{S}_{s,T}^{[\ell](k)}|||_{s,\Theta,\eta}\mathrm{d}s\\
&\lesssim \int_t^\infty (\lambda_s\langle s\rangle^{-1})^{\ell-1}\mathfrak{l}_s^c\langle s\rangle^{-1}\mathrm{d}s\\
&\lesssim (\lambda_t\langle t\rangle^{-1})^{\ell-1}\mathfrak{l}_t^c
\end{aligned}$$

by the positive definiteness bound $\|e^{\mathbb{G}_{t,s}(\xi_{1:(\ell-1)})}\|_{L^\infty}\lesssim 1$ that holds when $t\leq s$. This completes the induction.

The claim about the $|||\cdot|||_{t,\Theta}$-limits follows along the lines of Corollary 3.15. Indeed, we claim that the limiting kernels are solutions to the equations

$$\begin{aligned}
\gamma_t^{[\ell](n)}(\xi_{1:(\ell-1)}) &= 1 && \text{if } \ell=1, k=n\\
\gamma_t^{[\ell](k)}(\xi_{1:(\ell-1)}) &= -\int_0^t e^{\mathbb{G}_{t,s}(\xi_{1:(\ell-1)})}\mathscr{S}_{s,\infty}^{[\ell](k)}(\xi_{1:(\ell-1)},y)\mathrm{d}s && \text{if } \ell=1, k<n\\
\gamma_t^{[\ell](k)}(\xi_{1:(\ell-1)}) &= \int_t^\infty e^{\mathbb{G}_{t,s}(\xi_{1:(\ell-1)})}\mathscr{S}_{s,\infty}^{[\ell](k)}(\xi_{1:(\ell-1)},y)\mathrm{d}s && \text{if } \ell\geq 2
\end{aligned}\tag{5.9}$$

where $\mathscr{S}_{s,\infty}^{[\ell](k)}$ are source terms obtained by replacing the $\gamma_{t,T}^{[\ell](k)}$ in (5.5) with the limits $\gamma_t^{[\ell](k)}$; this is well-defined since the system is triangular in $(\ell,k)$. One then may see inductively that there exists a sufficiently small $\varepsilon>0$, depending only on $\beta$, such that

$$|||\gamma_t^{[\ell](k)}-\gamma_{t,T}^{[\ell](k)}|||_{t,\Theta(\ell),\eta}\lesssim(\lambda_t\langle t\rangle^{-1})^{\ell-1}\mathfrak{l}_t^c\langle t\rangle^{\varepsilon}\langle T\rangle^{-\varepsilon}.\qquad\square$$

Finally, we can show that the UV-complete germ is indeed a good approximate solution of the observable flow equation, i.e. we verify rigorously the estimate (2.27). Define

$$\begin{aligned}
q_t^{[\ell]}(\phi) &:= \int\Gamma_t^{[\ell]}(\phi,\xi_{1:(\ell-1)},y)\psi_{1:(\ell-1)}^{\rho}\mathrm{d}\xi_{1:(\ell-1)}\\
\Gamma_t^{[\ell]}(\phi,\xi_{1:(\ell-1)},y) &:= \sum_{k=0}^{n}\gamma_t^{[\ell](k)}(\xi_{1:(\ell-1)},y)\phi^k(y)
\end{aligned}\tag{5.10}$$

where the UV limit kernels $\gamma_t^{[\ell](k)}$ were constructed in Proposition 5.2.

**Proposition 5.3.** *Suppose $\beta^2<6\pi$, and let $\gamma_{t,T}^{[\ell](k)}$ be as in Proposition 5.2. $q_t^{[\leq\ell^*]}:=\sum_{\ell\leq\ell^*}q_t^{[\ell]}$, where $q_t^{[\ell]}$ is given by (5.10), is in the space $\hat{\mathscr{O}}_{\infty-,0,y}^{0,\ell^*\delta}$.*

**Proof.** The boundedness of the complex exponential and $1 \le \omega_t^{\Theta}(x_I)$ yields

$$
\begin{aligned}
q_t^{[\ell]}(\phi) &= \int \sum_{k=0}^{n} \gamma_t^{[\ell](k)}(\xi_{1:(\ell-1)}, y)\phi^k(y)\psi^{\rho}_{1:(\ell-1)} \mathrm{d}\xi_{1:(\ell-1)} \\
&\le \sum_{k=0}^{n} |||\gamma_t^{[\ell](k)}|||_{t,\Theta(\ell),\eta} |\phi^k(y)|,
\end{aligned}
$$

so we deduce from (4.9) and (5.8) that for all $\ell \ge 1$ there exists a $c = c(n, k, \ell, \beta^2) \in \mathbb{R}_+$ such that for any $p \in [1, \infty)$,

$$
\left\| q_t^{[\ell]}(\varphi_t) \right\|_{L^p(\Omega)} \lesssim (\lambda_t \langle t \rangle^{-1})^{\ell-1} \mathfrak{l}_t^c. \tag{5.11}
$$

Since the UV limiting kernels satisfy (5.9), adapting the proof of Lemma 5.1 yields that $q^{[\ell]}$ satisfy the observable flow equation (2.19). Thus Lemma 2.19 applies, and computing

$$
\begin{aligned}
&\langle \mathrm{D}q_t^{[\ell]}(\phi), \dot{G}R_t^{[>(\ell^*-\ell)]} \rangle \\
= &\int \sum_{k=0}^{n} \gamma_t^{[\ell](k)}(\xi_{1:(\ell-1)}, y) k\phi^{k-1}(y)\psi^{\rho}_{1:(\ell-1)} \dot{G}R_t^{[>(\ell^*-\ell)]}(y) \mathrm{d}\xi_{1:(\ell-1)} \\
+ &\int \sum_{k=0}^{n} \gamma_t^{[\ell](k)}(\xi_{1:(\ell-1)}, y)\phi^k(y)\psi^{\rho}_{1:(\ell-1)} i\beta \sum_{j \in 1:(\ell-1)} \sigma_j \dot{G}R_t^{[>(\ell^*-\ell)]}(x_j) \mathrm{d}\xi_{1:(\ell-1)} \\
\lesssim &\sum_{k=0}^{n} |||\gamma_t^{[\ell](k)}|||_{t,\Theta(\ell),\eta} (|\phi^{k-1}(y)| + |\phi^k(y)|) \left\| \dot{G}R_t^{[>(\ell^*-\ell)]} \right\|_{L^\infty(\Omega; L^\infty)},
\end{aligned}
$$

we obtain from the estimates (5.8), (4.9) and (2.26) that

$$
\left\| \mathscr{L}^V q_t^{[\le \ell^*]}(\varphi_t) \right\|_{L^p(\Omega)} \lesssim \langle t \rangle^{-\ell^*\delta - 1} \mathfrak{l}_t^c. \tag{5.12}
$$

Combining (5.11) and (5.12) then yields the statement. □

## 5.2 Vertex fields in the optimal Wick regime

Let $y \in \mathbb{R}^d$ be the point the observable is localized; let $\alpha \in \mathbb{C}$ and denote by

$$
\tilde{\psi}_y := \exp(i\alpha\phi(y)).
$$

As we will see at the end of this section (Corollary 5.8), the real part of $\alpha$ contains all difficulties and extending the results in the imaginary direction is straightforward. Since treating the general complex case from the beginning would only clutter the notation, we temporarily restrict ourselves to $\alpha \in \mathbb{R}$.

Our ansatz for the vertex field germ is

$$
q_t^{[\ell],T}(\phi) := \int \gamma_{t,T}^{[\ell]}(y, \xi_{1:(\ell-1)}) \tilde{\psi}_y \psi^{\rho}_{1:(\ell-1)} \mathrm{d}\xi_{1:(\ell-1)} \tag{5.13}
$$

The key to controlling the flow equation for the vertex field is understanding the action of the functional Laplacian $\Delta_{\dot{G}_t}$ on the field content $\tilde{\psi}_y \psi_{1:(\ell-1)}$. To motivate the discussion below, let us consider the example of the $\ell = 2$ case explicitly. Then

$$
\frac{1}{2}\Delta_{\dot{G}_t}[\tilde{\psi}_y \psi_1] = -\left( \frac{\alpha^2}{2}\dot{\mathcal{G}}_t(0) + \frac{\beta^2}{2}\dot{\mathcal{G}}_t(0) + \alpha\beta\sigma_1 \dot{\mathcal{G}}_t(x_1 - y) \right) \tilde{\psi}_y \psi_1
$$

It is clear that the sign of $\alpha\beta\sigma_1$ influences the behaviour crucially. When $t \le s$, it holds that $-\mathcal{G}_{t,s}(x_1 - y) \le 0$, hence the off-diagonal part has $L^\infty$ bounds

$$
\sup_{x_1} e^{-\alpha\beta\sigma_1 \mathcal{G}_{t,s}(x_1 - y)} \lesssim \begin{cases} 1 & \text{if} \quad \alpha\beta\sigma_1 \ge 0 \\ e^{-\alpha\beta\sigma_1 \mathcal{G}_{t,s}(0)} & \text{if} \quad \alpha\beta\sigma_1 \le 0 \end{cases}.
$$

That is, when $\mathrm{sgn}(\alpha)=\sigma_1$, the external vertex particle $\tilde{\psi}_y$ and the internal $\psi_1$ repel each other, behaving like the charged 2-point case $\psi_{1:2}$ of the potential, yielding much better bounds than the case $\mathrm{sgn}(\alpha)\neq\sigma_1$, in which case $\tilde{\psi}_y$ and $\psi_1$ attract each other and form a dipole. As such, to simplify notation, we assume from now on, without loss of generality, that $\alpha\geq 0$ (so if $\alpha\leq 0$ everything can be done verbatim, with the straightforward modification of some signs).

Evaluating the functional Laplacian on the field content $\tilde{\psi}_y\psi_{1:(\ell-1)}$ of (5.13) results in

$$\dot{\mathbf{G}}_t(y,\xi_{1:(\ell-1)}):=\dot{\mathbb{G}}_t(\xi_{1:(\ell-1)})+\left(-\frac{\alpha^2}{2}\dot{\mathcal{G}}_t(0)\right)+\left(-\alpha\beta\sum_{j=1}^{\ell-1}\sigma_j\dot{\mathcal{G}}_t(x_j-y)\right). \tag{5.14}$$

We may think of $\dot{\mathbf{G}}_t(y,\xi_I)$ as the entrywise sum of a “covariance matrix” $[\dot{\mathbf{G}}_t(y,\xi_I)]$ which factorizes as

$$[\dot{\mathbf{G}}_t(y,\xi_I)]=\begin{bmatrix} -\frac{\alpha^2}{2}\dot{\mathcal{G}}_t(0) & \left(\dot{\mathbb{G}}_t^{(\alpha)}(\xi_I)\right)^\top \\ \dot{\mathbb{G}}_t^{(\alpha)}(\xi_I) & [\dot{\mathbb{G}}_t(\xi_I)] \end{bmatrix}\in\mathbb{C}^{(1+|I|)\times(1+|I|)}$$

where $[\dot{\mathbb{G}}_t(\xi_I)]\in\mathbb{C}^{|I|\times|I|}$ is the usual covariance matrix of sine–Gordon, i.e. $[\dot{\mathbb{G}}_s(\xi_I)]_{j,k}:=-\frac{\beta^2}{2}\sigma_j\sigma_k\dot{\mathcal{G}}_s(x_j-x_k)$, and the off-diagonal vector $\dot{\mathbb{G}}_t^{(\alpha)}(\xi_I)$, given by $\left[\dot{\mathbb{G}}_t^{(\alpha)}(\xi_I)\right]_j:=-\frac{\alpha\beta}{2}\dot{\mathcal{G}}_t\sigma_j(x_j-y)$, encodes the interaction effect of perturbing the charge structure of sine–Gordon with the external vertex field $\exp(i\alpha\phi(y))$. The strength of this perturbation depends on the relation $\alpha\lesseqgtr\beta$ and we expect the two cases to require different ideas to control.

When $\alpha>\beta$ we expect that adding $\exp(i\alpha\phi(y))$ destroys the usual charge structure, so that the interaction is dominated by the off-diagonal perturbation term $\dot{\mathbb{G}}_t^{(\alpha)}(\xi_I)$. On the other hand, when $\alpha<\beta$, we expect the perturbation to be irrelevant and $\dot{\mathbf{G}}_t(y,\xi_I)$ to behave similarly to the sine–Gordon potential, i.e. governed by the bulk covariance $[\dot{\mathbb{G}}_t(\xi_I)]$. This observation is leveraged in the following

**Lemma 5.4. (external vertex screening)** *Let $0\leq\alpha,\beta$, suppose $t\leq s$. Then denoting by $\chi:=\chi(\xi_{1:\ell})$ the global charge, we have*

$$\sup_{x_{1:\ell}} e^{\mathbf{G}_{t,s}(y,\xi_{1:\ell})}\lesssim_\ell e^{-\frac{\alpha^2}{2}\mathcal{G}_{t,s}(0)}(\lambda_t\lambda_s^{-1})^{K_\ell}e^{\left(|\chi|+\frac{K_\ell-|\chi|}{2}\right)\alpha\beta\mathcal{G}_{t,s}(0)} \tag{5.15}$$

*with $K_\ell:=|\chi|$ if $\alpha\leq\beta$ and $K_\ell:=\ell$ if $\beta\leq\alpha$.*

**Proof.** First suppose that the ensemble $\xi_{1:\ell}$ is purely charged, i.e. $|\chi|=\ell$. By positivity of $\dot{\mathcal{G}}_u$

$$\begin{aligned} &\mathbf{G}_{t,s}(y,\xi_{1:\ell}) \\ =\;& -\frac{\alpha^2}{2}\mathcal{G}_{t,s}(0)+\mathbb{G}_{t,s}(\xi_{1:\ell})-\alpha\beta\sum_{j\in 1:\ell}\sigma_j\mathcal{G}_{t,s}(x_j-y) \\ =\;& -\frac{\alpha^2}{2}\mathcal{G}_{t,s}(0)-\frac{\beta^2}{2}\sum_{j\in 1:\ell}\mathcal{G}_{t,s}(0)-\frac{\beta^2}{2}\sum_{j\neq k\in 1:\ell}\mathcal{G}_{t,s}(x_j-x_k)-\alpha\beta\sum_{j\in 1:\ell}\sigma_j\mathcal{G}_{t,s}(x_j-y) \\ \leq\;& -\frac{\alpha^2}{2}\mathcal{G}_{t,s}(0)-\frac{\ell}{2}\beta^2\mathcal{G}_{t,s}(0)+\ell\alpha\beta\mathcal{G}_{t,s}(0). \end{aligned} \tag{5.16}$$

yielding the statement. Hence from now on we assume that $|\chi|\leq\ell-1$ so that there exists at least one neutral pair; denote by $k_+,k_-$ the index of a closest of such pairs, i.e. a pair $\eta:=\{(+1,x_{k_+}),(-1,x_{k_-})\}\in(\{-1,1\}\times\mathbb{R}^2)^2$ such that

$$|x_{k_+}-x_{k_-}|=\min_{j_+,j_-:\sigma_{j_+}\sigma_{j_-}=-1}|x_{j_+}-x_{j_-}|.$$

By possibly relabeling the points, we distinguish the purely charged core and the ensemble of neutral pairs with $\eta$ being the last, i.e. we assume $\sigma_1,\ldots,\sigma_{|\chi|}=\mathrm{sgn}(\chi)$ and $\sigma_{|\chi|+j}=-\sigma_{|\chi|+j+1}$ for all $j\in 1,\ldots,\ell-|\chi|-1$ and $(k_-,k_+)=(\ell-1,\ell)$. Then split

$$\mathbf{G}_{t,s}(y,\xi_{1:\ell})=\mathbf{G}_{t,s}(y,\xi_{1:(\ell-2)})+2\mathbb{G}_{t,s}(\xi_{1:(\ell-2)},\eta)+\mathscr{C}_{t,s}(\eta,y) \tag{5.17}$$

where

$$\mathscr{C}_{t,s}(\eta,y):=\beta^2\mathcal{G}_{t,s}(x_{k_+}-x_{k_-})-\beta^2\mathcal{G}_{t,s}(0)+\alpha\beta(\mathcal{G}_{t,s}(x_{k_-}-y)-\mathcal{G}_{t,s}(x_{k_+}-y))$$

We show that the coloured term, governed by the triangle $x_{k_+}, x_{k_-}, y$, is bounded uniformly by

$$\mathscr{C}_{t,s}(\eta,y)\lesssim\begin{cases}C & \text{if } \alpha\leq\beta\\ \alpha\beta\mathcal{G}_{t,s}(0)-\beta^2\mathcal{G}_{t,s}(0) & \text{if } \beta\leq\alpha.\end{cases}$$

1. If the edge $x_{k_+}-y$ is a shortest, then the red term is nonpositive, and the blue term is nonpositive as usual.
2. If the edge $x_{k_+}-x_{k_-}$ is a shortest, then we can bound the red term with Lemma 3.5, while the blue term is nonpositive as usual.
3. If the edge $x_{k_-}-y$ is a shortest, then
   a. when $\alpha\leq\beta$, since the red term in $\mathscr{C}_{t,s}$ is nonnegative,

   $$\mathscr{C}_{t,s}(\eta,y)\leq\beta^2(\mathcal{G}_{t,s}(x_{k_-}-y)-\mathcal{G}_{t,s}(0))+\beta^2(\mathcal{G}_{t,s}(x_{k_+}-x_{k_-})-\mathcal{G}_{t,s}(x_{k_+}-y))$$

   where now the first difference is nonpositive, and the second difference is estimated via Lemma 3.5 with the prescription $z:=x_{k_+}$, $x_b:=x_{k_-}$, $x_g:=y$.

   b. when $\beta\leq\alpha$, then $-\alpha\beta\leq-\beta^2$, so

   $$\mathscr{C}_{t,s}(\eta,y)\leq\alpha\beta\mathcal{G}_{t,s}(x_{k_-}-y)-\beta^2\mathcal{G}_{t,s}(0)+\beta^2(\mathcal{G}_{t,s}(x_{k_+}-x_{k_-})-\mathcal{G}_{t,s}(x_{k_+}-y))$$

   and the first difference is bounded by

   $$\alpha\beta\mathcal{G}_{t,s}(0)-\beta^2\mathcal{G}_{t,s}(0),$$

   while the second difference is estimated via Lemma 3.5 with the prescription $z:=x_{k_+}$, $x_b:=x_{k_-}$, $x_g:=y$.

The term $\mathbb{G}_{t,s}(\xi_{1:(\ell-2)},\eta)$ is treated as in the proof of Lemma 3.6: since by construction

$$|x_{k_+}-x_{k_-}|\leq|\mathbb{1}_{\sigma_j=1}x_{k_-}+\mathbb{1}_{\sigma_j=-1}x_{k_+}-x_j|$$

for all $j\in 1:(\ell-2)$, we may apply Lemma 3.5 to estimate

$$2\mathbb{G}_{t,s}(\xi_{1:(\ell-2)},\eta)=-\beta^2\sum_{j=1}^{\ell-2}\sigma_j(\mathcal{G}_{t,s}(x_{k_+}-x_j)-\mathcal{G}_{t,s}(x_{k_-}-x_j))\lesssim_\ell C.$$

Thus from (5.17) we obtain

$$\begin{aligned}e^{\mathbf{G}_{t,s}(y,\xi_{1:\ell})}&\lesssim_\ell e^{\mathbf{G}_{t,s}(y,\xi_{1:(\ell-2)})+\mathbb{1}_{\beta\leq\alpha}(\alpha\beta\mathcal{G}_{t,s}(0)-\beta^2\mathcal{G}_{t,s}(0))}\\&= e^{\mathbf{G}_{t,s}(y,\xi_{1:(\ell-2)})}((\lambda_t\lambda_s^{-1})^2e^{\alpha\beta\mathcal{G}_{t,s}(0)})^{\mathbb{1}_{\beta\leq\alpha}}\end{aligned}$$

and repeat the above extraction procedure inductively until only the purely charged core $\tilde{\xi}_{1:|\chi|}$ remains for which we have seen in (5.16) the bound

$$e^{\mathbf{G}_{t,s}(y,\tilde{\xi}_{1:|\chi|})}\lesssim e^{-\frac{\alpha^2}{2}\mathcal{G}_{t,s}(0)}(\lambda_t\lambda_s^{-1})^{|\chi|}e^{|\chi|\alpha\beta\mathcal{G}_{t,s}(0)}.$$ □

Note that the above bounds are not the sharpest in the sense that only the "bad" negative charges $\sigma_j=-1$ contribute to divergence in the $\alpha\beta\sigma_j\mathcal{G}_{t,s}(x_j-y)$ term; it remains bounded for the "good" positive ones $\sigma_j=1$.

In light of the previous sections it is now a standard exercise to derive the observable flow equation of the kernels of the ansatz (5.13).

**Lemma 5.5.** *Let $V^{[\ell]}$ be given by (3.31). Suppose that the kernels $\gamma_{t,T}^{[\ell]}$ satisfy the equations*

$$\gamma_{t,T}^{[\ell]}(\xi_{1:(\ell-1)},y)=e^{\mathbf{G}_{t,T}(y,\xi_{1:(\ell-1)})}\gamma_{T,T}^{[\ell]}(\xi_{1:(\ell-1)},y)+\int_t^T e^{\mathbf{G}_{t,s}(y,\xi_{1:(\ell-1)})}\mathscr{S}_{s,T}^{[\ell]}(\xi_{1:(\ell-1)},y)\mathrm{d}s \qquad (5.18)$$

*where* $\dot{\mathbf{G}}_t$ *is defined in* (5.14) *and*

$$
\begin{aligned}
- & \ \mathscr{S}_{s,T}^{[\ell]}(\xi_{1:(\ell-1)}, y) \\
= & \sum_{I_1 \uplus I_2 = 1:\ell} C_{|I_1|,|I_2|} \sum_{\substack{\boldsymbol{b} \in \mathfrak{M}: \\ \ell(\boldsymbol{b}) = |I_2|}} \big(v_s^{[|I_2|](\boldsymbol{b})} M^{(\boldsymbol{b})}\big)(\xi_{I_2}) i\alpha \not{\dot{\mathbb{G}}}_s(\xi_{I_2}, y) \gamma_{s,T}^{[|I_1|]}(y, \xi_{I_1 - 1}) \\
+ & \sum_{I_1 \uplus I_2 = 1:\ell} C_{|I_1|,|I_2|} \sum_{\substack{\boldsymbol{b} \in \mathfrak{M}: \\ \ell(\boldsymbol{b}) = |I_2|}} \big(v_s^{[|I_2|](\boldsymbol{b})} M^{(\boldsymbol{b})}\big)(\xi_{I_2}) \dot{\mathbb{G}}_s(\xi_{I_2}, \xi_{I_1 - 1}) \gamma_{s,T}^{[|I_1|]}(y, \xi_{I_1 - 1})
\end{aligned}
\tag{5.19}
$$

*for some positive combinatorial coefficients* $C_{|I_1|,|I_2|}$. *Then* $q_t^{[\ell],T}$, *given by* (5.13), *satisfies the observable flow equation* (2.19).

**Proof.** Straightforward; similar to Lemma 5.1. □

Imposing the boundary conditions

$$
\begin{aligned}
\gamma_{0,T}^{[1]} &:= 1 \\
\gamma_{T,T}^{[\ell]} &:= 0 \quad \text{if } \ell \geq 2,
\end{aligned}
\tag{5.20}
$$

allows us to estimate inductively the germ kernels in the norms (5.7):

**Proposition 5.6.** *Let* $\beta^2 < 6\pi$ *and* $0 \leq \alpha < 4\pi / \beta$. *Then there exist kernels* $\gamma_{t,T}^{[\ell]}$ *satisfying* (5.18) *with boundary conditions* (5.20), *and for any* $\eta \geq 0$ *it holds that*

$$
|||\gamma_{t,T}^{[\ell]}|||_{t,\Theta(\ell),\eta} \lesssim \langle t \rangle^{\frac{\alpha^2}{8\pi}} (\lambda_t \langle t \rangle^{-1})^{\ell-1} \mathfrak{l}_t^c \tag{5.21}
$$

*for some* $c = c(\ell, \beta^2) \in \mathbb{R}_+$. *Moreover, the* $|||\cdot|||_{t,\Theta,\eta}$*-limits*

$$
\gamma_t^{[\ell]} := \lim_{T \to \infty} \gamma_{t,T}^{[\ell]}
$$

*exist, are independent of* $\eta$ *and satisfy the bounds* (5.21).

**Proof.** The boundary condition (5.20) at $\ell = 1$ computes Wick ordering, so that by Lemma A.1

$$
|||\gamma_{t,T}^{[1]}(y)|||_{t,\Theta(1),\eta} = e^{\frac{\alpha^2}{2} \mathcal{G}_t(0)} \sup_y |e^{\varsigma \sqrt{t} 0}| \asymp e^{\frac{\alpha^2}{8\pi} \log(1 \vee t)} \asymp \langle t \rangle^{\frac{\alpha^2}{8\pi}},
$$

verifying (5.21). When $\ell \geq 2$, to estimate the source term, we use the inductive hypothesis, the potential bound (3.29), and for nontrivial $M^{(\boldsymbol{b})}(x_{I_2})$ the technology of Lemma 3.12. Indeed, from (3.25) with $c := 2 \cdot 4^{-\ell}$

$$
\begin{aligned}
& |||v_s^{[|I_2|](\boldsymbol{b})}(\xi_{I_2}) M^{(\boldsymbol{b})}(x_{I_2}) \dot{\mathbb{G}}_s(\xi_{I_2}, \xi_{I_1 - 1}) \gamma_{s,T}^{[|I_1|]}(\xi_{I_1 - 1}, y)|||_{s,4^{-\ell},\eta} \\
\lesssim & \ |||v_s^{[|I_2|](\boldsymbol{b})}|||_{s,3 \cdot 4^{-\ell},\eta} \langle s \rangle^{-2 + \frac{\varkappa(\boldsymbol{b})}{2}} |||\gamma_{s,T}^{[|I_1|]}|||_{s,4^{-\ell},\eta} \\
\lesssim & \ |||v_s^{[|I_2|](\boldsymbol{b})}|||_{s,4^{-(\ell-1)},\eta} \langle s \rangle^{-2 + \frac{\varkappa(\boldsymbol{b})}{2}} |||\gamma_{s,T}^{[|I_1|]}|||_{s,4^{-\ell},\eta} \\
\lesssim & \ \lambda_s (\lambda_s \langle s \rangle^{-1})^{|I_2| - 1} \mathfrak{l}_s^c \langle s \rangle^{-2} \langle s \rangle^{\frac{\alpha^2}{8\pi}} (\lambda_s \langle s \rangle^{-1})^{|I_1| - 1}
\end{aligned}
$$

since $|I_2| \leq \ell - 1$ apart from the trivial case $|I_2| = \ell$ when $|I_1| = 0$ so that $q^{[0]} \equiv 0$. Similarly, via a straightforward modification of (3.25) to $\not{\dot{\mathbb{G}}}_s(\xi_{1:2}, y)$

$$
\begin{aligned}
& |||v_s^{[|I_2|](\boldsymbol{b})}(\xi_{I_2}) M^{(\boldsymbol{b})}(x_{I_2}) i\alpha \not{\dot{\mathbb{G}}}_s(\xi_{I_2}, y) \gamma_{s,T}^{[|I_1|]}(\xi_{I_1 - 1}, y)|||_{s,4^{-\ell},\eta} \\
\lesssim & \ |||v_s^{[|I_2|](\boldsymbol{b})}|||_{s,4^{-(\ell-1)},\eta} \langle s \rangle^{-2 + \frac{\varkappa(\boldsymbol{b})}{2}} |||\gamma_{s,T}^{[|I_1|]}|||_{s,4^{-\ell},\eta} \\
\lesssim & \ \lambda_s (\lambda_s \langle s \rangle^{-1})^{|I_2| - 1} \mathfrak{l}_s^c \langle s \rangle^{-2} \langle s \rangle^{\frac{\alpha^2}{8\pi}} (\lambda_s \langle s \rangle^{-1})^{|I_1| - 1}.
\end{aligned}
$$

It follows from (5.19) that the source term satisfies

$$|||\mathscr{S}^{[\ell]}_{s,T}(\xi_{1:(\ell-1)},y)|||_{s,\Theta,\eta} \lesssim \langle s\rangle^{\frac{\alpha^2}{8\pi}}(\lambda_s\langle s\rangle^{-1})^{\ell-1}\langle s\rangle^{-1}\mathfrak{l}^c_s. \tag{5.22}$$

Inserting this and the screening bound (5.15) into (5.18) with the boundary condition (5.20) and denoting by $\chi := \chi(\xi_{1:(\ell-1)})$ the global charge yields

$$\begin{aligned} &|||\gamma^{[\ell]}_{t,T}|||_{t,\Theta,\eta} \\ \lesssim\ & \int_t^\infty \|e^{\mathbf{G}_{t,s}(y,\xi_{1:(\ell-1)})}\|_{L^\infty}|||\mathscr{S}^{[\ell]}_{s,T}(\xi_{1:(\ell-1)},y)|||_{s,\Theta,\eta}\mathrm{d}s \\ \lesssim\ & \langle t\rangle^{\frac{\alpha^2}{8\pi}}\int_t^\infty \langle s\rangle^{-\frac{\alpha^2}{8\pi}}(\lambda_t\lambda_s^{-1})^{K_{\ell-1}}e^{\left(|\chi|+\frac{K_{\ell-1}-|\chi|}{2}\right)\alpha\beta\mathcal{G}_{t,s}(0)}\langle s\rangle^{\frac{\alpha^2}{8\pi}}(\lambda_s\langle s\rangle^{-1})^{\ell-1}\langle s\rangle^{-1}\mathfrak{l}^c_s\mathrm{d}s \\ \lesssim\ & \langle t\rangle^{\frac{\alpha^2}{8\pi}}(\lambda_t\langle t\rangle^{-1})^{\ell-1}\mathfrak{l}^c_t \end{aligned}$$

provided the integral converges. To see that it does, we need to show (cf. Lemma A.1) that

$$\Sigma_\sigma := -(1-\delta)K_{\ell-1} + \left(|\chi| + \frac{K_{\ell-1}-|\chi|}{2}\right)\frac{\alpha\beta}{4\pi} - \delta(\ell-1) < 0.$$

Now by Lemma 5.4, either $K_{\ell-1} = \ell-1$, in which case $|\chi| + \frac{K_{\ell-1}-|\chi|}{2} \leq \ell-1$, so by the Wick regime assumption $\alpha\beta/4\pi < 1$

$$\Sigma_\sigma < -(1-\delta)(\ell-1) + \ell - 1 - \delta(\ell-1) = 0;$$

or $K_{\ell-1} = |\chi|$; if so, then $|\chi| = 0$ is trivial, and otherwise

$$\Sigma_\sigma < -(1-\delta)|\chi| + |\chi| - \delta(\ell-1) = \delta(|\chi| - (\ell-1)) \leq 0.$$

Taking the $T\to\infty$ limit is done similarly to the proof of Proposition 5.2: it is seen inductively that the limits solve

$$\begin{aligned} \gamma^{[1]}_t(y) &= e^{\frac{\alpha^2}{2}\mathcal{G}_t(0)} \\ \gamma^{[\ell]}_t(\xi_{1:(\ell-1)},y) &= \int_t^\infty e^{\mathbf{G}_{t,s}(y,\xi_{1:(\ell-1)})}\mathscr{S}^{[\ell]}_{s,\infty}(\xi_{1:(\ell-1)},y)\mathrm{d}s \quad \text{if } \ell\geq 2, \end{aligned} \tag{5.23}$$

where $\mathscr{S}^{[\ell]}_{s,\infty}$ are obtained by replacing the $\gamma^{[\ell]}_{t,T}$ in (5.19) with the limits $\gamma^{[\ell]}_t$, by propagating inductively the bound

$$|||\gamma^{[\ell]}_t - \gamma^{[\ell]}_{t,T}|||_{t,\Theta(\ell),\eta} \lesssim \langle t\rangle^{\frac{\alpha^2}{8\pi}}(\lambda_t\langle t\rangle^{-1})^{\ell-1}\mathfrak{l}^c_t\langle t\rangle^{\varepsilon}\langle T\rangle^{-\varepsilon}$$

for some $\varepsilon > 0$ small enough. □

As we have mentioned, in the case $\beta \leq \alpha$ the strong influence of the vertex field insertion destroys the usual charge structure of the covariance matrix. This observation allows for an alternative proof of the above proposition when the sine–Gordon potential is in its Wick regime $\beta^2 < 4\pi$, which does not rely on the external screening mechanism of Lemma 5.4 and operates only on a naive estimate on $\mathbf{G}_{t,s}(y,\xi_{1:\ell})$. This alternative proof is elaborated in the following side remark; we do not use the ideas presented there in the subsequent sections and hence invite the reader to skip it.

**Remark 5.7. (Alternative construction when $\beta^2 < 4\pi$)** As promised, we will not use Lemma 5.4 here. Instead, denoting by

$$\sharp_- := \sum_{j\in I}\mathbb{1}_{\sigma_j=-1}, \qquad \sharp_+ := \sum_{j\in I}\mathbb{1}_{\sigma_j=+1},$$

from (5.14) we deduce the following naive estimate: when $t\leq s$, and $\alpha,\beta\geq 0$, we have

$$\begin{aligned} \sup_{x_{1:\ell}} e^{\mathbf{G}_{t,s}(y,\xi_{1:\ell})} &\leq e^{-\frac{\alpha^2}{2}\mathcal{G}_{t,s}(0)}\left(\sup_{x_{1:\ell}} e^{\mathbb{G}_{t,s}(\xi_{1:\ell})}\right)\left(\sup_{x_{1:\ell}} e^{-\alpha\beta\sum_{j=1}^{\ell}\sigma_j\mathcal{G}_{t,s}(x_j-y)}\right) \\ &\lesssim e^{-\frac{\alpha^2}{2}\mathcal{G}_{t,s}(0)}(\lambda_t\lambda_s^{-1})^{|\chi(\xi_{1:\ell})|}e^{\sharp_-\alpha\beta\mathcal{G}_{t,s}(0)}; \end{aligned} \tag{5.24}$$

where the bound on the $\mathbb{G}_{t,s}$ sector has been proven in (3.12), and for the bound on the last term note that as $\dot{\mathcal{G}}_u \geq 0$, we have

$$\begin{aligned}
\sup_x \exp(-\alpha\beta\sigma\mathcal{G}_{t,s}(x-y)) &= \sup_x \exp(-\alpha\beta\mathcal{G}_{t,s}(x-y)) = 1 && \text{if } \sigma = 1\\
\sup_x \exp(-\alpha\beta\sigma\mathcal{G}_{t,s}(x-y)) &= \sup_x \exp(\alpha\beta\mathcal{G}_{t,s}(x-y)) = e^{\alpha\beta\mathcal{G}_{t,s}(0)} && \text{if } \sigma = -1
\end{aligned}$$

so that only the terms with $\sigma_j = -1$ contribute. Inserting (5.24) into (5.18), denoting by $\chi := \chi(\xi_{1:(\ell-1)})$ the global charge and relying on the inductive estimates (5.22) on the source term derived in the proof of the above Proposition, we arrive at

$$\begin{aligned}
|||\gamma_{t,T}^{[\ell]}|||_{t,\Theta,\eta} &\lesssim \int_t^\infty \|e^{\mathbf{G}_{t,s}(y,\xi_{1:(\ell-1)})}\|_{L^\infty} |||\mathscr{S}_{s,T}^{[\ell]}(\xi_{1:(\ell-1)},y)|||_{s,\Theta,\eta}\mathrm{d}s\\
&\lesssim \lambda_t^{|\chi|}\langle t\rangle^{-\sharp_- \frac{\alpha\beta}{4\pi}}\langle t\rangle^{\frac{\alpha^2}{8\pi}}\int_t^\infty \lambda_s^{-|\chi|}\langle s\rangle^{\sharp_-\frac{\alpha\beta}{4\pi}}(\lambda_s\langle s\rangle^{-1})^{\ell-1}\langle s\rangle^{-1}\mathfrak{l}_s^c\mathrm{d}s.
\end{aligned}$$

To finish the induction, it remains to show that this integral converges, for which we need

$$\Sigma_{\ell,\delta} := -(1-\delta)|\chi| + \sharp_- \frac{\alpha\beta}{4\pi} - \delta(\ell-1) < 0$$

for all $\ell \geq 2$. When $\beta^2 < 4\pi$, this in turn follows from simple charge algebra. Indeed, since $\sharp_+ + \sharp_- = (\ell-1)$ and $\sharp_+ - \sharp_- = \chi$, we have

$$\sharp_- = (\ell-1) - \sharp_+ = (\ell-1) - \sharp_- - \chi,$$

so

$$2\sharp_- = (\ell-1) - \chi \leq (\ell-1) + |\chi|.$$

Then denoting by $R := \frac{\alpha\beta}{4\pi} < 1$,

$$\begin{aligned}
\Sigma_{\ell,\delta} &\leq -(1-\delta)|\chi| + R\frac{|\chi| + (\ell-1)}{2} - \delta(\ell-1)\\
&= |\chi|\left(\frac{R}{2} - 1 + \delta\right) + (\ell-1)\left(\frac{R}{2} - \delta\right)\\
&= |\chi|\left(\frac{R}{2} - \frac{1}{2} - \frac{1}{2} + \delta\right) + (\ell-1)\left(\frac{R}{2} - \frac{1}{2} + \frac{1}{2} - \delta\right)\\
&= ((\ell-1) + |\chi|)\left(\frac{R-1}{2}\right) + ((\ell-1) - |\chi|)\left(\frac{1}{2} - \delta\right)
\end{aligned}$$

and the statement $\Sigma_{\ell,\delta} < 0$ follows from $\beta^2 < 4\pi \Leftrightarrow \delta > 1/2$ and the fact that $(\ell-1) - |\chi|$ and $(\ell-1) + |\chi|$ cannot be both zero for $\ell \geq 2$. This finishes the alternative construction.

As we have noted in the beginning of the section, it is straightforward to extend the construction from real $\alpha$ to complex.

**Corollary 5.8.** *Let $\beta^2 < 6\pi$ and $\alpha \in \mathbb{C}$ with real part $|\Re(\alpha)| < 4\pi/\beta$. Then the analogue of Proposition 5.6 holds, namely we have*

$$|||\gamma_t^{[\ell]}|||_{t,\Theta(\ell),\eta} \lesssim e^{\frac{\Re(\alpha)^2}{2}\mathcal{G}_t(0)} e^{-\frac{\Im(\alpha)^2}{2}\mathcal{G}_t(0)} (\lambda_t\langle t\rangle^{-1})^{\ell-1}\mathfrak{l}_t^c \tag{5.25}$$

*for some $c = c(\ell, \beta^2) \in \mathbb{R}_+$.*

**Proof.** For $\alpha = \Re(\alpha) + i\Im(\alpha) \in \mathbb{C}$, we may write

$$\tilde{\psi}_y = e^{i\Re(\alpha)\phi(y)} e^{-\Im(\alpha)\phi(y)}$$

so that a simple computation gives that Wick ordering is now

$$\gamma_{t,T}^{[1]} = e^{\frac{\Re(\alpha)^2}{2}\mathcal{G}_t(0)} e^{-\frac{\Im(\alpha)^2}{2}\mathcal{G}_t(0)} e^{i\Re(\alpha)\Im(\alpha)\mathcal{G}_t(0)}$$

verifying the inductive bound (5.25) for $\ell = 1$. The action of the functional Laplacian is now

$$\begin{aligned} \dot{\mathbf{G}}_t^{\mathbb{C}}(y, \xi_{1:(\ell-1)}) &= \dot{\mathbb{G}}_t(\xi_{1:(\ell-1)}) + \left(-\frac{\Re(\alpha)^2}{2}\dot{\mathcal{G}}_t(0)\right) + \left(-\Re(\alpha)\beta\sum_{j=1}^{\ell-1} \sigma_j \dot{\mathcal{G}}_t(x_j - y)\right) \\ &+ \frac{\Im(\alpha)^2}{2}\dot{\mathcal{G}}_t(0) + \left(-i\Im(\alpha)\beta\sum_{j=1}^{\ell-1} \sigma_j \dot{\mathcal{G}}_t(x_j - y)\right) + (-i\Im(\alpha)\Re(\alpha)\dot{\mathcal{G}}_t(0)) \end{aligned} \tag{5.26}$$

so that bounding trivially

$$\left|\exp\left(-i\Im(\alpha)\left(\beta\sum_j \sigma_j \mathcal{G}_{t,s}(x_j - y) + \Re(\alpha)\dot{\mathcal{G}}_t(0)\right)\right)\right| \leq 1$$

we get

$$\|e^{\mathbf{G}_{t,s}^{\mathbb{C}}}\|_{L^\infty} \leq \|e^{\mathbf{G}_{t,s}}\|_{L^\infty} e^{\frac{\Im(\alpha)^2}{2}\mathcal{G}_{t,s}(0)}$$

and adapting the end of the proof of Proposition 5.6 propagates the inductive hypothesis (5.25) via the screening estimate (5.15) on $e^{\mathbf{G}_{t,s}}$ and source term bound obtained just as (5.22). □

Having constructed the UV-complete kernels $\gamma_t^{[\ell]}$, we may define

$$q_t^{[\ell]}(\phi) := \int \gamma_t^{[\ell]}(y, \xi_{1:(\ell-1)}) \tilde{\psi}_y \psi_{1:(\ell-1)}^{\rho} \mathrm{d}\xi_{1:(\ell-1)} \tag{5.27}$$

**Proposition 5.9.** *Suppose $\beta^2 < 6\pi$ and $\alpha \in \mathbb{C}$ with real part $|\Re(\alpha)| < \frac{4\pi}{\beta}$. Then $q_t^{[\leq \ell^*]} := \sum_{\ell \leq \ell^*} q_t^{[\ell]}$, where $q_t^{[\ell]}$ is given by (5.27), is in the space $\hat{\mathscr{O}}_{\infty-, i\alpha, y}^{\kappa, \sigma}$ with $\kappa := \frac{\Re(\alpha)^2 - \Im(\alpha)^2}{8\pi}$, $\sigma := \ell^*\delta - \frac{\Re(\alpha)^2 - \Im(\alpha)^2}{8\pi}$.*

**Proof.** The boundedness of the complex exponential and $1 \leq \omega_t^{\Theta}(x_I)$ yields

$$\begin{aligned} e^{-\langle i\alpha, \phi(y)\rangle} q_t^{[\ell]}(\phi) &= \int \gamma_t^{[\ell]}(y, \xi_{1:(\ell-1)}) e^{-i\alpha\phi(y)} e^{i\alpha\phi(y)} \psi_{1:(\ell-1)}^{\rho} \mathrm{d}\xi_{1:(\ell-1)} \\ &\leq |||\gamma_t^{[\ell]}|||_{t, \Theta(\ell), \eta} \end{aligned}$$

hence (5.25) yields that for all $\ell \geq 1$ there exists a $c = c(\ell, \beta^2) \in \mathbb{R}_+$ such that for any $p \geq 1$

$$\|e^{-\langle i\alpha, \varphi_t(y)\rangle} q_t^{[\ell]}(\varphi_t)\|_{L^p(\Omega)} \lesssim \langle t\rangle^{\frac{\Re(\alpha)^2 - \Im(\alpha)^2}{8\pi}} (\lambda_t \langle t\rangle^{-1})^{\ell-1} \mathfrak{l}_t^c. \tag{5.28}$$

Since the UV limiting kernels satisfy (5.23), it follows just as Lemma 5.5 that $q^{[\ell]}$ satisfy the observable flow equation (2.19). Thus Lemma 2.19 applies, and computing

$$\begin{aligned} & e^{-\langle i\alpha, \phi(y)\rangle} \langle \mathrm{D}q_t^{[\ell]}(\phi), \dot{G}R_t^{[>(\ell^* - \ell)]}\rangle \\ = \ & i\int \gamma_t^{[\ell]}(y, \xi_{1:(\ell-1)}) \psi_{1:(\ell-1)}^{\rho} \left(\alpha \dot{G}R_t^{[>(\ell^* - \ell)]}(y) + \beta \sum_{j \in 1:(\ell-1)} \sigma_j \dot{G}R_t^{[>(\ell^* - \ell)]}(x_j)\right) \mathrm{d}\xi_{1:(\ell-1)} \\ \lesssim \ & |||\gamma_t^{[\ell]}|||_{t, \Theta(\ell), \eta} \|\dot{G}R_t^{[>(\ell^* - \ell)]}\|_{L^\infty}, \end{aligned}$$

we obtain from the estimates (5.25) and (2.26) that

$$\|e^{-\langle i\alpha, \varphi_t(y)\rangle} \mathscr{L}^V q_t^{[\leq \ell^*]}(\varphi_t)\|_{L^p(\Omega)} \lesssim \langle t\rangle^{\frac{\Re(\alpha)^2 - \Im(\alpha)^2}{8\pi} - \ell^*\delta - 1} \mathfrak{l}_t^c. \tag{5.29}$$

Combining (5.28) and (5.29) yields the statement. □

## 5.3 Derivative fields

The polynomial local observables have very mild superficial degree of divergence, hence we could flow backwards at all levels $\ell \geq 2$ in the representation (5.1). This is not the case for sufficiently high derivatives of the field anymore and new ideas are needed to solve the flow equation. To see why, let us compute the source term with the Wick ordered first-order germ $q_t^{[1]}(\phi) := \partial^\mu \phi(y)$ with $y \in \mathbb{R}^d$ and $\mu \in \mathbb{N}_0^d$ multiindex. It is

$$\langle \mathrm{D}q_t^{[1]}, \dot{G}_t \mathrm{D}V_t^{[1]} \rangle = i\beta \int \sigma_1 \partial_y^\mu \dot{\mathcal{G}}_t(y - x_1) v_t^{[1]} \rho(x_1) \psi_1 \mathrm{d}\xi_1.$$

Then it is clear from the bound $\|\partial^\mu \dot{\mathcal{G}}_t\|_{L_\eta^1} \lesssim \langle s \rangle^{-2+\frac{|\mu|}{2}}$ that repeating the strategy of the previous section does not work here for arbitrary multiindex $\mu$, since we would have

$$|||\gamma_{t,T}^{[2]}|||_{t,\Theta,\eta} \lesssim \int_t^\infty \|e^{\mathbb{G}_{t,s}(\xi_1)}\|_{L^\infty} |||\mathscr{S}_{s,T}^{[2]}|||_{s,\Theta,\eta} \mathrm{d}s \lesssim \lambda_t \int_t^\infty \langle s \rangle^{-1+\frac{|\mu|}{2}} \langle s \rangle^{-1} \mathrm{d}s$$

which does not converge when $|\mu| \geq 2$.

The idea then is to change the representation of the germ: the field content in (2.28) will not be products of vertex fields $\psi_{1:\ell}$ as in (5.1); instead, we make a new ansatz. The goal is to shift the divergence caused by the derivative $\partial^\mu$ from the kernels to the field content, which is not flown.

Given an index set $I \subseteq 1{:}\,\ell$, define the operator $\partial_I$ by

$$\partial_I f := \sum_{j \in I} \partial_{x_j} f(x_1, \ldots, x_\ell)$$

where $x_j \in \mathbb{R}^d$. Then postulate

$$q_t^{[\ell],T}(\phi) = \begin{cases} \gamma_{t,T}^{[1]}(y) \partial_y^\mu \phi(y) & \text{if } \ell = 1 \\ \int \gamma_{t,T}^{[\ell]}(y, \xi_{1:(\ell-1)}) \partial_{1:(\ell-1)}^\mu \psi_{1:(\ell-1)}^\rho \mathrm{d}\xi_{1:(\ell-1)} & \text{if } \ell \geq 2 \end{cases} \tag{5.30}$$

We have to show that this ansatz is compatible with the structure of the equation (2.19). This is not immediately obvious, since the representation we obtained before for $V^{[\ell]}$ involves the vertex fields $\psi_{1:\ell}$ without derivatives on them. However, the following lemma shows that translation invariance and exponential decay of the fluctuation free field covariance $\dot{\mathcal{G}}_t$ ensures that the Picard scheme (2.19) propagates the original field insertion $\partial_y^\mu \phi$ through the induction in such a way that is compatible with (5.30).

**Lemma 5.10.** *Let $q^{[\ell],T}$ be given by (5.30) and $V^{[\ell]}$ by (3.31). Then for any $I_1 \uplus I_2 = 1{:}\,\ell$ we have*

$$\langle \mathrm{D}q_t^{[|I_1|],T}, \dot{G}_t \mathrm{D}V_t^{[|I_2|]} \rangle = \int \mathbf{B}_t^{[|I_1|,|I_2|]}(y, \xi_{1:(\ell-1)}) \partial_{1:(\ell-1)}^\mu \psi_{1:(\ell-1)}^\rho \mathrm{d}\xi_{1:(\ell-1)} \tag{5.31}$$

*for some function $\mathbf{B}_t^{[|I_1|,|I_2|]}$.*

**Proof.** Since $\dot{\mathcal{G}}_t(x_1 - x_2)$ is invariant under simultaneous translations, it is straightforward to see from the induction of Section 3.2 that the same is true for the potential kernels $v_t^{[\ell](\boldsymbol{a})}$ (understood in the degenerate sense for the first-order kernel $v_t^{[1]} := \lambda_t$ and the vacuum energy contribution $v_t^{[2](0)}$), and also for the divergence-corrected kernels $v_t^{[\ell](\boldsymbol{a})} M^{(\boldsymbol{a})}$ where the correction $M^{(\boldsymbol{a})}$ was defined in (5.3). That is, for any $z \in \mathbb{R}^d$ we have

$$\big(v_t^{[\ell](\boldsymbol{a})} M^{(\boldsymbol{a})}\big)(\sigma_{\boldsymbol{a}}, x_1 + z, \ldots, x_\ell + z) = \big(v_t^{[\ell](\boldsymbol{a})} M^{(\boldsymbol{a})}\big)(\sigma_{\boldsymbol{a}}, x_1, \ldots, x_\ell).$$

When $|I_1|=1$, compute

$$\left\langle \mathrm{D}q_t^{[1],T}, \dot{G}_t \mathrm{D}V_t^{[|I_2|]}\right\rangle = C_{|I_2|} i\beta\gamma_{t,T}^{[1]}(y) \sum_{\substack{\boldsymbol{b}\in\mathfrak{M}:\\ \ell(\boldsymbol{b})=|I_2|}} \int \sum_{j\in I_2} \sigma_j \partial_y^\mu \dot{\mathcal{G}}_t(y-x_j)\big(v_t^{[|I_2|](\boldsymbol{b})} M^{(\boldsymbol{b})}\big)(\xi_{I_2})\psi_{I_2}^\rho \mathrm{d}\xi_{I_2}.$$

Fix a charge configuration $\sigma_{I_2}$. For any $a\in\mathbb{R}^d$, denote by $\psi_{I_2}^\rho(x_{I_2}-a)$ the usual vertex field product evaluated at $\phi(x_k-a)$, $k\in I_2$. Consider

$$\begin{aligned} H(a) &:= \int \sum_{j\in I_2} \sigma_j \dot{\mathcal{G}}_t(y-x_j)\big(v_t^{[|I_2|](\boldsymbol{b})} M^{(\boldsymbol{b})}\big)(x_{I_2})\psi_{I_2}^\rho(x_{I_2}-a)\mathrm{d}x_{I_2} \\ &= \int \sum_{j\in I_2} \sigma_j \dot{\mathcal{G}}_t(y-a-z_j)\big(v_t^{[|I_2|](\boldsymbol{b})} M^{(\boldsymbol{b})}\big)(z_{I_2})\psi_{I_2}^\rho(z_{I_2})\mathrm{d}z_{I_2} \end{aligned}$$

where in the second line we used the unit Jacobian change of variables $z_k:=x_k-a$, $k\in I_2$ and translation invariance of $v_t^{[|I_2|](\boldsymbol{b})} M^{(\boldsymbol{b})}$. We compute $\partial_a^\mu$ in both of these representations. On one hand, since $\partial_a^\mu \psi_{I_2}^\rho(x_{I_2}-a)=(-1)^{|\mu|}\partial_{I_2}^\mu \psi_{I_2}^\rho(x_{I_2}-a)$, provided we can differentiate under the integral (which we postpone to the end of the proof),

$$\partial_a^\mu H = (-1)^{|\mu|}\int \sum_{j\in I_2} \sigma_j \dot{\mathcal{G}}_t(y-x_j)\big(v_t^{[|I_2|](\boldsymbol{b})} M^{(\boldsymbol{b})}\big)(x_{I_2})\partial_{I_2}^\mu \psi_{I_2}^\rho(x_{I_2}-a)\mathrm{d}x_{I_2}.$$

while on the other, again provided we can differentiate under the integral,

$$\begin{aligned} \partial_a^\mu H &= \int \sum_{j\in I_2} \sigma_j \partial_a^\mu \dot{\mathcal{G}}_t(y-a-z_j)\big(v_t^{[|I_2|](\boldsymbol{b})} M^{(\boldsymbol{b})}\big)(z_{I_2})\psi_{I_2}^\rho(z_{I_2})\mathrm{d}z_{I_2} \\ &= (-1)^{|\mu|}\int \sum_{j\in I_2} \sigma_j \partial_y^\mu \dot{\mathcal{G}}_t(y-a-z_j)\big(v_t^{[|I_2|](\boldsymbol{b})} M^{(\boldsymbol{b})}\big)(z_{I_2})\psi_{I_2}^\rho(z_{I_2})\mathrm{d}z_{I_2}. \end{aligned}$$

Evaluating at $a=0$, we obtain

$$\int \sum_{j\in I_2} \sigma_j \partial_y^\mu \dot{\mathcal{G}}_t(y-z_j)\big(v_t^{[|I_2|](\boldsymbol{b})} M^{(\boldsymbol{b})}\psi_{I_2}^\rho\big)(z_{I_2})\mathrm{d}z_{I_2} = \int \sum_{j\in I_2} \sigma_j \dot{\mathcal{G}}_t(y-x_j)\big(v_t^{[|I_2|](\boldsymbol{b})} M^{(\boldsymbol{b})}\partial_{I_2}^\mu \psi_{I_2}^\rho\big)(x_{I_2})\mathrm{d}x_{I_2}$$

which means that $\left\langle \mathrm{D}q_t^{[1],T}, \dot{G}_t \mathrm{D}V_t^{[|I_2|]}\right\rangle$ is of the form (5.31), since $|I_1|=1$ implies $|I_2|=\ell-1$.

When $|I_1|\geq 2$, since

$$\mathrm{D}\partial_I^\mu \psi_I = \partial_I^\mu \mathrm{D}\psi_I = \partial_I^\mu\left(i\beta\sum_{j\in I}\sigma_j\delta_{x_j}\psi_I\right) = i\beta\sum_{j\in I}\sigma_j\sum_{\nu\leq\mu}\binom{\mu}{\nu}\partial_I^\nu\delta_{x_j}\partial_I^{\mu-\nu}\psi_I^\rho, \tag{5.32}$$

writing $\mathbb{G}_t(\xi_j,\xi_k):=-\frac{\beta^2}{2}\sigma_j\sigma_k\dot{\mathcal{G}}_t(x_j-x_k)$ and $I:=1{:}(\ell-1)$, we have

$$\begin{aligned} &\left\langle \mathrm{D}q_t^{[|I_1|],T}, \dot{G}_t \mathrm{D}V_t^{[|I_2|]}\right\rangle \\ &= C_{|I_1|,|I_2|}\sum_{\substack{\boldsymbol{b}\in\mathfrak{M}:\\ \ell(\boldsymbol{b})=|I_2|}} \int \gamma_{t,T}^{[|I_1|]}(y,\xi_{I_1-1}) \\ &\times \sum_{\substack{j\in I_1-1\\ k\in I_2}}\sum_{\nu\leq\mu}\binom{\mu}{\nu}\partial_{I_1-1}^\nu \mathbb{G}_t(\xi_j,\xi_k)\big(v_t^{[|I_2|](\boldsymbol{b})} M^{(\boldsymbol{b})}\big)(\xi_{I_2})\psi_{I_2}^\rho \partial_{I_1}^{\mu-\nu}\psi_{I_1}^\rho \mathrm{d}\xi_I \end{aligned}$$

By similar argument as above, for any $j\in I_1-1$, for each $x_j\in\mathbb{R}^d$, we obtain

$$\int \sum_{k\in I_2}\sigma_k \partial_{x_j}^\nu \dot{\mathcal{G}}_t(x_j-x_k) v_t^{[|I_2|](\boldsymbol{b})} M^{(\boldsymbol{b})}\psi_{I_2}^\rho(x_{I_2})\mathrm{d}x_{I_2} = \int \sum_{k\in I_2}\sigma_k \dot{\mathcal{G}}_t(x_j-x_k) v_t^{[|I_2|](\boldsymbol{b})} M^{(\boldsymbol{b})}\partial_{I_2}^\nu \psi_{I_2}^\rho(x_{I_2})\mathrm{d}x_{I_2}$$

which, combined with the formula for $\langle \mathrm{D}q_t^{[|I_1|],T}, \dot{G}_t \mathrm{D}V_t^{[|I_2|]}\rangle$ yields the statement, since

$$\sum_{\nu\leq\mu}\binom{\mu}{\nu}\partial_{I_2}^{\nu}\psi_{I_2}^{\rho}\partial_{I_1}^{\mu-\nu}\psi_{I_1}^{\rho}=\partial_{1:(\ell-1)}^{\mu}\psi_{1:(\ell-1)}^{\rho}$$

by the binomial theorem.

It then remains to show we can apply the Leibniz integral rule. Let $j\in I_2$. For any fixed $y\in\mathbb{R}^d$

$$\int\left|\sum_{j\in I_2}\sigma_j\partial_y^{\mu}\dot{\mathcal{G}}_t(y-z_j)\big(v_t^{[|I_2|](\boldsymbol{b})}M^{(\boldsymbol{b})}\big)(z_{I_2})\psi_{I_2}^{\rho}(z_{I_2})\right|\mathrm{d}z_{I_2}\lesssim|||v_t^{[|I_2|](\boldsymbol{b})}|||_{t,\varsigma,\eta}\langle t\rangle^{\frac{|\mu|+\varkappa(\boldsymbol{b})}{2}-2}$$

for any $\varsigma>0$, where in the case $\varkappa(\boldsymbol{b})>0$ so that $M^{(\boldsymbol{b})}(z_{I_2})=|z_1-z_2|^{-\varkappa(\boldsymbol{b})}$, we proceed as in Lemma 3.12, i.e. leverage neutrality and feed the divergence to $\sum_{j\in I_2}\sigma_j\partial_y^{\mu}\dot{\mathcal{G}}_t$ via Lemma A.6.

For the other integral, the function that is needed to be dominated involves $\partial_{I_2}^{\mu}\psi_{I_2}^{\rho}$, so that the boundedness of the complex exponential is not sufficient. To control the infinite volume limit, we need to rely on the Steiner weights in our kernels. Note that for a fixed $j\in I_2$, we have

$$\prod_{k\in I_2}\langle x_k\rangle\lesssim\langle x_j\rangle\prod_{k\in I_2\backslash\{j\}}\langle x_k-x_j\rangle\langle x_j\rangle\lesssim\langle\mathrm{St}(x_{I_2})\rangle^{|I_2|-1}\langle x_j-y\rangle^{|I_2|}\langle y\rangle^{|I_2|}.$$

Thus, when $\varkappa(\boldsymbol{b})=0$, applying Faà di Bruno to the vertex fields, there exists polynomials $P^{\mu}$ of the field and its derivatives such that

$$\begin{aligned}
&\int\sum_{j\in I_2}\sigma_j\dot{\mathcal{G}}_t(y-x_j)\big(v_t^{[|I_2|](\boldsymbol{b})}M^{(\boldsymbol{b})}\big)(x_{I_2})\partial_{I_2}^{\mu}\psi_{I_2}^{\rho}(x_{I_2}-a)\mathrm{d}x_{I_2}\\
&\lesssim\sum_{j\in I_2}|||v_t^{[|I_2|](\boldsymbol{b})}|||_{t,\varsigma,\eta}\|\dot{\mathcal{G}}_t\|_{L_\eta^1}\sup_{x_{I_2}}\left|\langle x_j-y\rangle^{-\eta}\langle\mathrm{St}(x_{I_2})\rangle^{-\eta}\prod_{j\in I_2}P^{\mu}(\phi(x_j-a))\right|\\
&\lesssim|||v_t^{[|I_2|](\boldsymbol{b})}|||_{t,\varsigma,\eta}\langle t\rangle^{-2}\prod_{j\in I_2}\left\|\langle\cdot\rangle^{-\frac{\eta}{|I_2|}}P^{\mu}(\phi(\cdot-a))\right\|_{L^\infty}\\
&\lesssim|||v_t^{[|I_2|](\boldsymbol{b})}|||_{t,\varsigma,\eta}\langle t\rangle^{-2}\|\phi\|_{\mathcal{H}}^{c}
\end{aligned}$$

for some $c$, where in the last line we used Besov embedding and chose $\eta$ large enough. When $\varkappa>0$, we proceed similarly, the only difference being the use of Lemma A.6:

$$\begin{aligned}
&\int\frac{\dot{\mathcal{G}}_t(y-x_1)-\dot{\mathcal{G}}_t(y-x_2)}{|x_1-x_2|^{\varkappa}}\big(v_t^{[|I_2|](\boldsymbol{b})}\big)(x_{1:2})\partial_{I_2}^{\mu}\psi_{I_2}^{\rho}(x_{1:2}-a)\mathrm{d}x_{1:2}\\
&\lesssim|||v_t^{[|I_2|](\boldsymbol{b})}|x_1-x_2|^{1-\varkappa}|||_{t,\varsigma,\eta}e^{-\frac{m^2}{t}}t^{-\frac{1}{2}}\int e^{-\varsigma\sqrt{t}|x_1-y|}\langle x_1-y\rangle^{\eta}\mathrm{d}x_1\\
&\times\sup_{x_{1:2}}\left|\langle x_1-y\rangle^{-\eta}\langle\mathrm{St}(x_{1:2})\rangle^{-\eta}\prod_{j\in1:2}P^{\mu}(\phi(x_j-a))\right|.
\end{aligned}$$

The justification for the Leibniz integral rule when $|I_1|\geq2$ follows verbatim. □

The point of the above lemma was not just to shift the divergence from the kernels to the field, but also to do it in such a way that the resulting basis of the representation consists solely of eigenfunctions of the functional Laplacian $\Delta_{\dot{G}_t}$. Indeed, translation invariance of $\dot{\mathcal{G}}_t$ ensures that there are no "trickle down" maps when evaluating $\Delta_{\dot{G}_t}$ on the field content $\partial_{1:\ell}^{\mu}\psi_{1:\ell}$ (as opposed to the case of the potential, cf. the map $\mathbf{A}\big(v_t^{[2](\varkappa)}\big)$ around (3.17)). This observation leads to

**Lemma 5.11.** *Fix $\mu\in\mathbb{N}_0^d$ multiindex, and let $V^{[\ell]}$ be given by (3.31). Suppose that the kernels $\gamma_{t,T}^{[\ell]}$ satisfy the equations*

$$\gamma_{t,T}^{[\ell]}(\xi_{1:(\ell-1)},y)=e^{\mathbb{G}_{t,T}(\xi_{1:(\ell-1)})}\gamma_{T,T}^{[\ell]}(\xi_{1:(\ell-1)},y)+\int_t^T e^{\mathbb{G}_{t,s}(\xi_{1:(\ell-1)})}\mathscr{S}_{s,T}^{[\ell]}(\xi_{1:(\ell-1)},y)\mathrm{d}s \tag{5.33}$$

*where*

$$
\begin{aligned}
&-\mathscr{S}^{[\ell]}_{s,T}(\xi_{1:(\ell-1)}, y)\\
= &\sum_{\substack{I_1 \uplus I_2 = 1:\ell\\ |I_1|=1}} C_{|I_1|,|I_2|} \sum_{\substack{\boldsymbol{b}\in\mathfrak{M}:\\ \ell(\boldsymbol{b})=|I_2|}} \gamma^{[1]}_{s,T}(y)\dot{\mathcal{G}}_s(\xi_{I_2}, y)\big(v_s^{[|I_2|](\boldsymbol{b})}M^{(\boldsymbol{b})}\big)(\xi_{I_2})\\
+ &\sum_{\substack{I_1 \uplus I_2 = 1:\ell\\ |I_1|\geq 2}} C_{|I_1|,|I_2|} \sum_{\substack{\boldsymbol{b}\in\mathfrak{M}:\\ \ell(\boldsymbol{b})=|I_2|}} \gamma^{[|I_1|]}_{s,T}(y, \xi_{I_1-1})\dot{\mathbb{G}}_s(\xi_{I_1-1}, \xi_{I_2})\big(v_s^{[|I_2|](\boldsymbol{b})}M^{(\boldsymbol{b})}\big)(\xi_{I_2})
\end{aligned}
\tag{5.34}
$$

*for some positive combinatorial coefficients* $C_{|I_1|,|I_2|}$. *Then* $q_t^{[\ell],T}$, *given by* (5.30), *satisfies the observable flow equation* (2.19).

**Proof.** Let $I\subseteq\mathbb{N}$ be a finite index set. Compute

$$
\begin{aligned}
&\frac{1}{2}\Delta_{\dot{G}_t}[\partial_I^\mu \psi_I^\rho]\\
= &-\frac{\beta^2}{2}\int_{(\mathbb{R}^2)^2}\dot{\mathcal{G}}_t(z_1-z_2)\partial^\mu_{x_I}\left(\sum_{j_1,j_2\in I}\sigma_{j_1}\sigma_{j_2}\delta(x_{j_1}-z_1)\delta(x_{j_2}-z_2)\psi_I^\rho\right)\mathrm{d}z_{1:2}\\
= &-\frac{\beta^2}{2}\int_{(\mathbb{R}^2)^2}\dot{\mathcal{G}}_t(z_1-z_2)\sum_{\nu\leq\mu}\binom{\mu}{\nu}\partial^\nu_{x_I}\left(\sum_{j_1,j_2\in I}\sigma_{j_1}\sigma_{j_2}\delta(x_{j_1}-z_1)\delta(x_{j_2}-z_2)\right)\partial_I^{\mu-\nu}\psi_I^\rho\mathrm{d}z_{1:2}\\
= &\sum_{\nu\leq\mu}\binom{\mu}{\nu}\partial^\nu_{x_I}\left(-\frac{\beta^2}{2}\sum_{j_1,j_2\in I}\sigma_{j_1}\sigma_{j_2}\int\dot{\mathcal{G}}_t(z_1-z_2)\delta(x_{j_1}-z_1)\delta(x_{j_2}-z_2)\mathrm{d}z_{1:2}\right)\partial_I^{\mu-\nu}\psi_I^\rho\\
= &\sum_{\nu\leq\mu}\binom{\mu}{\nu}\partial^\nu_{x_I}\dot{\mathbb{G}}_t(\xi_I)\partial_I^{\mu-\nu}\psi_I^\rho\\
= &\;\dot{\mathbb{G}}_t(\xi_I)\partial_I^\mu\psi_I^\rho
\end{aligned}
$$

since for any nontrivial multiindex, the translation invariant covariance matrix $\dot{\mathbb{G}}_t(\xi_I)$ is annihilated by the diagonal derivative $\partial_I=\sum_{j\in I}\partial_{x_j}$. The statement then follows from Lemma 5.10. □

Imposing the boundary conditions

$$
\begin{aligned}
\gamma^{[1]}_{0,T} &:= 1\\
\gamma^{[\ell]}_{T,T} &:= 0 \quad \text{if } \ell\geq 2,
\end{aligned}
\tag{5.35}
$$

allows us to estimate inductively the germ kernels in the norms (5.7) and transfer these bounds onto the germ and its approximate generator.

**Proposition 5.12.** *Let* $\beta^2<6\pi$ *and* $\mu\in\mathbb{N}_0^d$ *multiindex. Then there exist kernels* $\gamma^{[\ell]}_{t,T}$ *satisfying* (5.33) *with boundary conditions* (5.35), *and for any* $\eta\geq 0$ *it holds that*

$$
|||\gamma^{[\ell]}_{t,T}|||_{t,\Theta(\ell),\eta}\lesssim(\lambda_t\langle t\rangle^{-1})^{\ell-1}\mathfrak{l}_t^c \tag{5.36}
$$

*for some* $c=c(\ell,\beta^2)\in\mathbb{R}_+$. *Moreover, the* $|||\cdot|||_{t,\Theta,\eta}$*-limits* $\gamma_t^{[\ell]}:=\lim_{T\to\infty}\gamma^{[\ell]}_{t,T}$ *exist, are independent of* $\eta$ *and satisfy the bounds* (5.36). *Finally, defining*

$$
q_t^{[\ell]}(\phi):=\begin{cases}\gamma_t^{[1]}(y)\partial^\mu_y\phi(y) & \text{if } \ell=1\\ \displaystyle\int\gamma_t^{[\ell]}(y,\xi_{1:(\ell-1)})\partial^\mu_{1:(\ell-1)}\psi^\rho_{1:(\ell-1)}\mathrm{d}\xi_{1:(\ell-1)} & \text{if } \ell\geq 2\end{cases}, \tag{5.37}
$$

*it holds that* $q_t^{[\leq\ell^*]}:=\sum_{\ell\leq\ell^*}q_t^{[\ell]}$ *is in the space* $\hat{\mathscr{O}}^{\frac{|\mu|}{2},\ell^*\delta-\frac{|\mu|}{2}}_{\infty-,0,y}$.

**Proof.** The bound (5.36) is verified similarly to Proposition 5.2 by repeating standard estimates on the source term (5.34). The $T\to\infty$ limit is also obtained via a straightforward modification: defining $\mathscr{S}^{[\ell]}_{s,\infty}$ by replacing the $\gamma^{[\ell]}_{t,T}$ in (5.34) with the limits $\gamma^{[\ell]}_t$, it is seen inductively that the limits solve

$$\begin{aligned}\gamma^{[1]}_t(y) &= 1\\ \gamma^{[\ell]}_t(\xi_{1:(\ell-1)},y) &= \int_t^\infty e^{\mathbb{G}_{t,s}(\xi_{1:(\ell-1)})}\mathscr{S}^{[\ell]}_{s,\infty}(\xi_{1:(\ell-1)},y)\mathrm{d}s \quad \text{if } \ell\geq 2,\end{aligned} \tag{5.38}$$

by propagating the bound

$$|||\gamma^{[\ell]}_t-\gamma^{[\ell]}_{t,T}|||_{t,\Theta(\ell),\eta}\lesssim(\lambda_t\langle t\rangle^{-1})^{\ell-1}\mathfrak{l}^c_t\langle t\rangle^{\varepsilon}\langle T\rangle^{-\varepsilon}$$

for some $\varepsilon>0$ small enough.

We claim that for any $\ell\geq 1$ there exists $c=c(\ell,\beta^2)\in\mathbb{R}_+$ such that for every $p\geq 1$

$$\left\|q^{[\ell]}_t(\varphi_t)\right\|_{L^p(\Omega)}\lesssim\langle t\rangle^{\frac{|\mu|}{2}}(\lambda_t\langle t\rangle^{-1})^{\ell-1}\mathfrak{l}^c_t. \tag{5.39}$$

Indeed, when $\ell=1$ simply apply (4.9). For $\ell\geq 2$, first note that Faà di Bruno and (4.9) yields for any $b\in\mathbb{R}$, $x\in\mathbb{R}^2$ that

$$\|\partial^\mu_x e^{ib\varphi_t(x)}\|_{L^p(\Omega)}\lesssim\langle t\rangle^{\frac{|\mu|}{2}}\mathfrak{l}^c_t,$$

so that denoting by $\psi_j(t):=e^{i\beta\sigma_j\varphi_t(x_j)}$, using Minkowski's integral inequality, $1\leq\omega^\Theta_t$ and the multinomial theorem we obtain

$$\begin{aligned}&\left\|q^{[\ell]}_t(\varphi_t)\right\|_{L^p(\Omega)}\\ \leq&\int\left|\gamma^{[\ell]}_t(y,\xi_{1:(\ell-1)})\right|\|\partial^\mu_{1:(\ell-1)}\psi^\rho_{1:(\ell-1)}(t)\|_{L^p(\Omega)}\mathrm{d}\xi_{1:(\ell-1)}\\ \leq&\int\left|\gamma^{[\ell]}_t(y,\xi_{1:(\ell-1)})\right|\sum_{\sum_{j=1}^{\ell-1}\mu^{(j)}=\mu}\binom{\mu}{\mu^{(1)},\ldots,\mu^{(\ell-1)}}\prod_{j=1}^{\ell-1}\left\|\partial^{\mu^{(j)}}_{x_j}\psi^\rho_j(t)\right\|_{L^{(\ell-1)p}(\Omega)}\mathrm{d}\xi_{1:(\ell-1)}\\ \lesssim&\ |||\gamma^{[\ell]}_t|||_{t,\Theta(\ell),\eta}\langle t\rangle^{\frac{|\mu|}{2}}\end{aligned}$$

and inserting (5.36) yields (5.39). On the other hand, we claim

$$\left\|\mathscr{L}^V q^{[\leq\ell^*]}_t(\varphi_t)\right\|_{L^p(\Omega)}\lesssim\langle t\rangle^{\frac{|\mu|}{2}}\langle t\rangle^{-\ell^*\delta-1}\mathfrak{l}^c_t. \tag{5.40}$$

Indeed, this follows in the case $\ell=1$ simply from (2.26); for $\ell\geq 2$, we proceed similarly as above: compute with (5.32)

$$\begin{aligned}&\left\|\langle\mathrm{D}q^{[\ell]}_t(\varphi_t),\dot{G}R^{[>(\ell^*-\ell)]}_t\rangle\right\|_{L^p(\Omega)}\\ \lesssim&\int\left|\gamma^{[\ell]}_t(y,\xi_{1:(\ell-1)})\right|\left\|\sum_{j\in 1:(\ell-1)}\sigma_j\sum_{\nu\leq\mu}\binom{\mu}{\nu}\partial^\nu\dot{G}R^{[>(\ell^*-\ell)]}_t(x_j)\partial^{\mu-\nu}_{1:(\ell-1)}\psi_{1:(\ell-1)}\right\|_{L^p(\Omega)}\mathrm{d}\xi_{1:(\ell-1)}\\ \lesssim&\int\left|\gamma^{[\ell]}_t(y,\xi_{1:(\ell-1)})\right|\sum_{\nu\leq\mu}\|\partial^\nu\dot{G}R^{[>(\ell^*-\ell)]}_t\|_{L^\infty(\Omega;L^\infty)}\|\partial^{\mu-\nu}_{1:(\ell-1)}\psi_{1:(\ell-1)}\|_{L^p(\Omega)}\mathrm{d}\xi_{1:(\ell-1)}\\ \lesssim&\ |||\gamma^{[\ell]}_t|||_{t,\Theta(\ell),\eta}\sum_{\nu\leq\mu}\langle t\rangle^{\frac{|\nu|}{2}-1-\delta-(\ell^*-\ell)\delta}\langle t\rangle^{\frac{|\mu|-|\nu|}{2}}\end{aligned}$$

where in the last line we used (2.26) to estimate $\partial^\nu\dot{G}R^{[>(\ell^*-\ell)]}_t$ (more precisely, we estimate with (4.15), Young's inequality, (A.2) and (4.1)). It follows (5.38) that the proof of Lemma 5.11 applies and thus $q^{[\ell]}$ satisfy the observable flow equation (2.19). Thus inserting (5.36) in the above estimate and using Lemma 2.19, we conclude (5.40). Combining (5.39) and (5.40) yields the statement. □

Similar considerations hold for the derivatives of the vertex field $\tilde{\psi}_y := \exp(i\alpha\phi(y))$ with $\alpha$ as in the previous section. Make the ansatz

$$q_t^{[\ell],T}(\phi) = \int \gamma_{t,T}^{[\ell]}(y, \xi_{1:(\ell-1)}) \partial_{y,1:(\ell-1)}^{\mu} \tilde{\psi}_y \psi_{1:(\ell-1)}^{\rho} \mathrm{d}\xi_{1:(\ell-1)} \tag{5.41}$$

where the operator $\partial_{y,I}$ is defined by

$$\partial_{y,I} f(y, x_1, \ldots, x_\ell) := (\partial_y + \partial_I) f(y, x_1, \ldots, x_\ell).$$

Similar ideas to those in Lemma 5.10 lead us to the following

**Lemma 5.13.** *Let $q_t^{[\ell]}$ be given by (5.41) and $V^{[\ell]}$ by (3.31). Then for all $I_1 \uplus I_2 = 1{:}\ell$ we have*

$$\left\langle \mathrm{D}q_t^{[|I_1|]}, \dot{G}_t \mathrm{D}V_t^{[|I_2|]} \right\rangle = \int \mathbf{B}_t^{[|I_1|,|I_2|]}(y, \xi_{1:(\ell-1)}) \partial_{y,1:(\ell-1)}^{\mu} \tilde{\psi}_y \psi_{1:(\ell-1)}^{\rho} \mathrm{d}\xi_{1:(\ell-1)} \tag{5.42}$$

*for some function $\mathbf{B}_t^{[|I_1|,|I_2|]}$.*

**Proof.** Compute

$$\begin{aligned} \mathrm{D}_\phi \partial_{y,I}^{\mu} \tilde{\psi}_y \psi_I^{\rho} &= \partial_{y,I}^{\mu} \mathrm{D}_\phi \tilde{\psi}_y \psi_I^{\rho} \\ &= \partial_{y,I}^{\mu} \left( \left( i\alpha\delta_y + i\beta \sum_{j \in I} \sigma_j \delta_{x_j} \right) \tilde{\psi}_y \psi_I^{\rho} \right) \\ &= i \sum_{\nu \leq \mu} \binom{\mu}{\nu} \left( \alpha \partial_y^{\nu} \delta_y + \beta \sum_{j \in I} \sigma_j \partial_{x_j}^{\nu} \delta_{x_j} \right) \partial_{y,I}^{\mu-\nu} (\tilde{\psi}_y \psi_I^{\rho}) \end{aligned}$$

so that

$$\begin{aligned} &\left\langle \mathrm{D}q_t^{[|I_1|]}, \dot{G}_t \mathrm{D}V_t^{[|I_2|]} \right\rangle \\ = \;& C_{|I_1|,|I_2|} \sum_{\substack{\boldsymbol{b} \in \mathfrak{M}: \\ \ell(\boldsymbol{b}) = |I_2|}} \int \gamma_t^{[|I_1|]}(y, \xi_{I_1 - 1}) \sum_{k \in I_2} \sum_{\nu \leq \mu} \binom{\mu}{\nu} \left( i\alpha \partial_y^{\nu} \dot{\not{\mathbb{G}}}_t(y, \xi_k) + \sum_{j \in I_1 - 1} \partial_{x_j}^{\nu} \dot{\mathbb{G}}_t(\xi_j, \xi_k) \right) \\ \times \;& \left( v_t^{[|I_2|](\boldsymbol{b})} M^{(\boldsymbol{b})} \right)(\xi_{I_2}) \psi_{I_2}^{\rho} \partial_{y,I_1-1}^{\mu-\nu} (\tilde{\psi}_y \psi_{I_1-1}^{\rho}) \mathrm{d}\xi_{1:(\ell-1)} \end{aligned}$$

Similarly to the proof of Lemma 5.10, we obtain

$$\int \sum_{k \in I_2} i\alpha \partial_y^{\nu} \dot{\not{\mathbb{G}}}_t(y, \xi_k) \left( v_t^{[|I_2|](\boldsymbol{b})} M^{(\boldsymbol{b})} \right)(\xi_{I_2}) \psi_{I_2}^{\rho} \mathrm{d}\xi_{I_2} = \int \sum_{k \in I_2} i\alpha \dot{\not{\mathbb{G}}}_t(y, \xi_k) \left( v_t^{[|I_2|](\boldsymbol{b})} M^{(\boldsymbol{b})} \right)(\xi_{I_2}) \partial_{I_2}^{\nu} \psi_{I_2}^{\rho} \mathrm{d}x_{I_2}$$

and for each fixed $x_j \in \mathbb{R}^d$, $j \in I_1 - 1$,

$$\int \sum_{k \in I_2} \partial_{x_j}^{\nu} \dot{\mathbb{G}}_t(\xi_j, \xi_k) \left( v_t^{[|I_2|](\boldsymbol{b})} M^{(\boldsymbol{b})} \right)(\xi_{I_2}) \psi_{I_2}^{\rho} \mathrm{d}\xi_{I_2} = \int \sum_{k \in I_2} \dot{\mathbb{G}}_t(\xi_j, \xi_k) \left( v_t^{[|I_2|](\boldsymbol{b})} M^{(\boldsymbol{b})} \right)(\xi_{I_2}) \partial_{I_2}^{\nu} \psi_{I_2}^{\rho} \mathrm{d}x_{I_2}.$$

The statement then follows from

$$\sum_{\nu \leq \mu} \binom{\mu}{\nu} \partial_{I_2}^{\nu} \psi_{I_2}^{\rho} \partial_{y,I_1-1}^{\mu-\nu} (\tilde{\psi}_y \psi_{I_1-1}^{\rho}) = \partial_{y,1:(\ell-1)}^{\mu} \tilde{\psi}_y \psi_{1:(\ell-1)}^{\rho}. \qquad \square$$

We then have the analogue of Lemma 5.5 and 5.11:

**Lemma 5.14.** *Fix $\mu \in \mathbb{N}_0^d$ multiindex, and let $V_t^{[\ell]}$ be given by (3.31). Suppose that the kernels $\gamma_{t,T}^{[\ell]}$ satisfy the equations*

$$\gamma_{t,T}^{[\ell]}(\xi_{1:(\ell-1)}, y) = e^{\mathbf{G}_{t,T}^{\mathbb{C}}(y, \xi_{1:(\ell-1)})} \gamma_{T,T}^{[\ell]}(\xi_{1:(\ell-1)}, y) + \int_t^T e^{\mathbf{G}_{t,s}^{\mathbb{C}}(y, \xi_{1:(\ell-1)})} \mathscr{S}_{s,T}^{[\ell]}(\xi_{1:(\ell-1)}, y) \mathrm{d}s \tag{5.43}$$

*where* $\dot{\mathbf{G}}_t^{\mathbb{C}}$ *is defined in* (5.26) *and* $\mathscr{S}_{s,T}^{[\ell]}$ *in* (5.19). *Then* $q_t^{[\ell],T}$, *given by* (5.41), *satisfies the observable flow equation* (2.19).

**Proof.** Similarly to the proof of Lemma 5.11, we have for any $I \subseteq \mathbb{N}$ finite index set

$$\frac{1}{2}\Delta_{\dot{G}_t}[\partial_{y,I}^{\mu}\tilde{\psi}_y\psi_I^{\rho}] = \dot{\mathbf{G}}_t^{\mathbb{C}}(y,\xi_I)\partial_{y,I}^{\mu}\tilde{\psi}_y\psi_I^{\rho}.$$

The statement then follows from the previous lemma. □

Imposing the boundary conditions

$$\begin{array}{rcl}\gamma_{0,T}^{[1]} & := & 1\\ \gamma_{T,T}^{[\ell]} & := & 0 \quad \text{if } \ell \geq 2,\end{array} \tag{5.44}$$

finishes the argument:

**Proposition 5.15.** *Let* $\beta^2 < 6\pi$, $\mu \in \mathbb{N}_0^d$ *multiindex and* $\alpha \in \mathbb{C}$ *with real part* $|\Re(\alpha)| < 4\pi/\beta$. *Then there exist kernels* $\gamma_{t,T}^{[\ell]}$ *satisfying* (5.43) *with boundary conditions* (5.44), *and for any* $\eta \geq 0$ *it holds that*

$$|||\gamma_{t,T}^{[\ell]}|||_{t,\Theta(\ell),\eta} \lesssim e^{\frac{\Re(\alpha)^2}{2}\mathcal{G}_t(0)} e^{-\frac{\Im(\alpha)^2}{2}\mathcal{G}_t(0)} (\lambda_t\langle t\rangle^{-1})^{\ell-1}\mathfrak{l}_t^c \tag{5.45}$$

*for some* $c = c(\ell,\beta^2) \in \mathbb{R}_+$. *Moreover, the* $|||\cdot|||_{t,\Theta,\eta}$*-limits* $\gamma_t^{[\ell]} := \lim_{T\to\infty}\gamma_{t,T}^{[\ell]}$ *exist, are independent of* $\eta$ *and satisfy the bounds* (5.45). *Finally, defining*

$$q_t^{[\ell]}(\phi) := \int \gamma_t^{[\ell]}(y,\xi_{1:(\ell-1)})\partial_{y,1:(\ell-1)}^{\mu}\tilde{\psi}_y\psi_{1:(\ell-1)}^{\rho}\mathrm{d}\xi_{1:(\ell-1)}, \tag{5.46}$$

*it holds that* $q_t^{[\leq \ell^*]} := \sum_{\ell \leq \ell^*} q_t^{[\ell]}$ *is in the space* $\hat{\mathscr{O}}_{\infty-,i\alpha,y}^{\kappa,\sigma}$ *with* $\kappa := \frac{|\mu|}{2} + \frac{\Re(\alpha)^2 - \Im(\alpha)^2}{8\pi}$, $\sigma := \ell^*\delta - \frac{\Re(\alpha)^2 - \Im(\alpha)^2}{8\pi} - \frac{|\mu|}{2}$.

**Proof.** The inductive bound (5.45) follows from comparing Lemma 5.14 with Lemma 5.5 and applying Proposition 5.6 and Corollary 5.8. Sending $T \to \infty$ is again routine: defining $\mathscr{S}_{s,\infty}^{[\ell]}$ by replacing the $\gamma_{t,T}^{[\ell]}$ in (5.19) with the limits $\gamma_t^{[\ell]}$, it is seen inductively that the limits solve

$$\begin{array}{rcl}\gamma_t^{[1]}(y) & = & e^{\frac{\Re(\alpha)^2}{2}\mathcal{G}_t(0)} e^{-\frac{\Im(\alpha)^2}{2}\mathcal{G}_t(0)} e^{i\Re(\alpha)\Im(\alpha)\mathcal{G}_t(0)}\\ \gamma_t^{[\ell]}(\xi_{1:(\ell-1)},y) & = & \displaystyle\int_t^{\infty} e^{\mathbf{G}_{t,s}^{\mathbb{C}}(y,\xi_{1:(\ell-1)})}\mathscr{S}_{s,\infty}^{[\ell]}(\xi_{1:(\ell-1)},y)\mathrm{d}s \quad \text{if } \ell \geq 2,\end{array} \tag{5.47}$$

by propagating the bound

$$|||\gamma_t^{[\ell]} - \gamma_{t,T}^{[\ell]}|||_{t,\Theta(\ell),\eta} \lesssim e^{\frac{\Re(\alpha)^2}{2}\mathcal{G}_t(0)} e^{-\frac{\Im(\alpha)^2}{2}\mathcal{G}_t(0)} (\lambda_t\langle t\rangle^{-1})^{\ell-1}\mathfrak{l}_t^c\langle t\rangle^{\varepsilon}\langle T\rangle^{-\varepsilon}$$

for some $\varepsilon > 0$ small enough.

By the multinomial theorem and Faà di Bruno,

$$\begin{aligned} & e^{-\langle i\alpha,\phi(y)\rangle} q_t^{[\ell]} \\ = \ & e^{-i\alpha\phi(y)} \int \gamma_t^{[\ell]}(y,\xi_{1:(\ell-1)})\partial_{y,1:(\ell-1)}^{\mu} e^{i\alpha\phi(y)}\psi_{1:(\ell-1)}^{\rho}\mathrm{d}\xi_{1:(\ell-1)} \\ = \ & e^{-i\alpha\phi(y)} \int \gamma_t^{[\ell]}(y,\xi_{1:(\ell-1)}) \sum_{\sum_{j=0}^{\ell-1}\mu^{(j)}=\mu} \binom{\mu}{\mu^{(0)},\ldots,\mu^{(\ell-1)}} \partial_y^{\mu^{(0)}} e^{i\alpha\phi(y)} \prod_{j=1}^{\ell-1} \partial_{x_j}^{\mu^{(j)}}\psi_j^{\rho}\mathrm{d}\xi_{1:(\ell-1)} \\ = \ & \int \gamma_t^{[\ell]}(y,\xi_{1:(\ell-1)}) \sum_{\sum_{j=0}^{\ell-1}\mu^{(j)}=\mu} \binom{\mu}{\mu^{(0)},\ldots,\mu^{(\ell-1)}} P_{\mu^{(0)}}(\{\partial^{\nu}\phi(y)\}_{\nu\leq\mu^{(0)}}) \prod_{j=1}^{\ell-1} \partial_{x_j}^{\mu^{(j)}}\psi_j^{\rho}\mathrm{d}\xi_{1:(\ell-1)} \end{aligned}$$

for polynomials $P_{\mu^{(0)}}$ with total scaling degree $|\mu^{(0)}|$, and consequently by Minkowski's integral inequality, Hölder, (4.9), $1 \le \omega_t^{\Theta}$, and (5.45) we obtain

$$\left\| e^{-\langle i\alpha, \varphi_t(y)\rangle} q_t^{[\ell]}(\varphi_t) \right\|_{L^p(\Omega)} \lesssim \langle t\rangle^{\frac{|\mu|}{2} + \frac{\Re(\alpha)^2 - \Im(\alpha)^2}{8\pi}} (\lambda_t \langle t\rangle^{-1})^{\ell-1} \mathfrak{c}_t^c. \tag{5.48}$$

On the other hand, computing

$$\begin{aligned} & e^{-\langle i\alpha, \phi(y)\rangle} \left\langle \mathrm{D} q_t^{[\ell]}(\phi), \dot{G}R_t^{[>(\ell^*-\ell)]} \right\rangle \\ = \; & i \int \gamma_t^{[\ell]}(y, \xi_{1:(\ell-1)}) \sum_{\nu \le \mu} \binom{\mu}{\nu} \partial^{\nu}_{y,1:(\ell-1)} \left( \alpha \dot{G}R_t^{[>(\ell^*-\ell)]}(y) + \beta \sum_{j=1}^{\ell-1} \sigma_j \dot{G}R_t^{[>(\ell^*-\ell)]}(x_j) \right) \\ \times \; & e^{-i\alpha\phi(y)} \partial^{\mu-\nu}_{y,1:(\ell-1)} e^{i\alpha\phi(y)} \psi^{\rho}_{1:(\ell-1)} \mathrm{d}\xi_{1:(\ell-1)} \end{aligned}$$

and proceeding similarly as above, by (2.26) (more precisely, by (4.15), Young's inequality, (A.2) and (4.1)), we get

$$\left\| e^{-\langle i\alpha, \varphi_t(y)\rangle} \left\langle \mathrm{D} q_t^{[\ell]}(\varphi_t), \dot{G}R_t^{[>(\ell^*-\ell)]} \right\rangle \right\|_{L^p(\Omega)} \lesssim |||\gamma_t^{[\ell]}|||_{t,\Theta(\ell),\eta} \sum_{\nu \le \mu} \langle t\rangle^{\frac{|\nu|}{2} - 1 - \delta - (\ell^*-\ell)\delta} \langle t\rangle^{\frac{|\mu| - |\nu|}{2}}.$$

It follows from (5.47) that the proof of Lemma 5.1 applies and thus $q^{[\ell]}$ satisfy the observable flow equation (2.19). Thus inserting (5.45) in the above estimate and using Lemma 2.19, we obtain

$$\left\| e^{-\langle i\alpha, \varphi_t(y)\rangle} \mathscr{L}^V q_t^{[\le \ell^*]}(\varphi_t) \right\|_{L^p(\Omega)} \lesssim \langle t\rangle^{\frac{|\mu|}{2} + \frac{\Re(\alpha)^2 - \Im(\alpha)^2}{8\pi}} \langle t\rangle^{-\ell^*\delta - 1} \mathfrak{c}_t^c. \tag{5.49}$$

The statement follows from combining (5.48) and (5.49). □

# 6 Observable products, prefactorization algebra property

In this section we assemble the ideas developed throughout the paper to finish the proof of Theorems 1.1 and 1.2.

We begin with unifying the notation of the previous section: given a $y \in \mathbb{R}^d$ where the germ is localized, $\mu \in \mathbb{N}_0^d$, $\alpha \in \mathbb{C}$, we write

$$q_t^{[\ell]}(\phi) = \int \sum_{k=0}^{n} \gamma_t^{[\ell](\alpha,k)}(\xi_{1:(\ell-1)}, y) \partial^{\mu}_{y,1:(\ell-1)} \left( \phi^k(y) \tilde{\psi}_y^{(\alpha)} \psi^{\rho}_{1:(\ell-1)} \right) \mathrm{d}\xi_{1:(\ell-1)} \tag{6.1}$$

where $\tilde{\psi}_y^{(\alpha)} := e^{i\alpha\phi(y)}$. This ansatz includes everything we have constructed in the previous section as special cases, cf. (5.10), (5.27), (5.37) and (5.46).

Recall that to prove existence of observables, in light of the Ito reconstruction (Lemma 2.6), we need to construct coherent germs. It is intuitively clear that the conditions 1 and 2 in Definition 2.5 are only technical, and the actually constraining part is the coherence estimate (2.10). In the lemma below we verify these conditions to the ansatz (6.1).

**Lemma 6.1.** *Let $q^{\langle 1\rangle}$ and $q^{\langle 2\rangle}$ be a sum of functionals of the form (6.1). Then $\alpha q^{\langle 1\rangle} + q^{\langle 2\rangle}$ and $\alpha q^{\langle 1\rangle} q^{\langle 2\rangle}$, $\alpha \in \mathbb{R}$ satisfy Conditions 1 and 2 of Definition 2.5.*

**Proof.** Let $m > d/2$. Fix a $T \in \mathbb{R}_+$ arbitrary, and consider the restriction of $q^{\langle j\rangle}$ to $[0, T]$. Take $\mathcal{H} := B^s_{2,2}(\langle x\rangle^{-m})$ with $s > |\mu| + \frac{d}{2}$, so that Morrey's inequality yields that $\partial^{\mu}\delta_x$ are continuous functionals on $\mathcal{H}$. In particular, given a sequence $\phi_r \xrightarrow{\|\cdot\|_{\mathcal{H}}} \phi$, we have for every $x \in \mathbb{R}^d$ and $\zeta \in \mathbb{C}$ that

$$(\partial^{\mu}\phi_r(x))^n e^{\zeta\phi_r(x)} \to (\partial^{\mu}\phi(x))^n e^{\zeta\phi(x)}$$

and similarly for all $k\leq n$ and $\nu$ multiindex with $|\nu|\leq|\mu|$.

As $|x_j|\leq|y|+|x_j-y|\leq|y|+\mathrm{St}(x_{1:(\ell-1)},y)$, we obtain

$$\prod_{j=1}^{\ell-1}\langle x_j\rangle\lesssim\langle y\rangle^{(\ell-1)}\langle\mathrm{St}(x_{1:(\ell-1)},y)\rangle^{(\ell-1)}.$$

Applying Faà di Bruno to the vertex fields in (6.1) we get $P_k^\mu$ polynomials $\phi(x_j)$ and their derivatives. Then choosing $\eta$ large enough, by Besov embedding we have,

$$\begin{aligned}\left|q_t^{[\ell]}(\phi)\right| &\lesssim \sum_{k=0}^{n}|||\gamma_t^{[\ell](\alpha,k)}|||_{0,\Theta,\eta}\sup_{x_{1:(\ell-1)}}\left|\langle\mathrm{St}(x_{1:(\ell-1)},y)\rangle^{-\eta}P_k^\mu(\phi(y))\prod_{j=1}^{\ell-1}P_k^\mu(\phi(x_j))e^{-\Im(\alpha)\phi(y)}\right|\\ &\lesssim \sum_{k=0}^{n}|||\gamma_t^{[\ell](\alpha,k)}|||_{0,\Theta,\eta}|\langle y\rangle^{\eta}P_k^\mu(\phi(y))e^{-\Im(\alpha)\phi(y)}|\prod_{j=1}^{\ell-1}\left\|\langle\cdot\rangle^{-\frac{\eta}{(\ell-1)}}P_k^\mu(\phi(\cdot))\right\|_{L^\infty}\\ &\lesssim \sum_{k=0}^{n}|||\gamma_t^{[\ell](\alpha,k)}|||_{0,\Theta,\eta}|\langle y\rangle^{\eta}P_k^\mu(\phi(y))e^{-\Im(\alpha)\phi(y)}|||\phi\|_{\mathcal{H}}^c\end{aligned}$$

for some $c$. Combining this with the pointwise convergence above, by dominated convergence we obtain $q_t^{[\ell]}(\phi_r)\to q_t^{[\ell]}(\phi)$.

On the other hand, we see from the induction on the kernels that $t\mapsto\gamma_t^{[\ell](\alpha,k)}$ and $t\mapsto\partial_t\gamma_t^{[\ell](\alpha,k)}$ are continuous in $|||\cdot|||_{0,\Theta,\eta}$, since $\gamma_t^{[\ell](\alpha,k)}$ satisfies classically a differential equation with continuous coefficients. Indeed, for relevant kernels continuity is trivial, and for irrelevant, note that

$$|||\gamma_t^{[\ell](\alpha,k)}-\gamma_{t_r}^{[\ell](\alpha,k)}|||_{0,\Theta}\lesssim\int_{t\wedge t_r}^{t\vee t_r}\|e^{\mathbb{U}_{t\wedge t_r,s}}\|_{L^\infty}|||\mathscr{S}_s|||_{s,\Theta}\mathrm{d}s+\int_{t\vee t_r}^{\infty}\|e^{\mathbb{U}_{t,s}}-e^{\mathbb{U}_{t_r,s}}\|_{L^\infty}|||\mathscr{S}_s|||_{s,\Theta}\mathrm{d}s$$

for some $\mathbb{U}_{\cdot,s}=\int_\cdot^s\dot{\mathbb{U}}_v\mathrm{d}v$ with $\dot{\mathbb{U}}_v$ a linear combination of $\mathcal{G}_{\cdot,s}$ (see the preceding section), hence continuous in $L^\infty$. The case of time derivatives follow by induction, as

$$\partial_t\gamma_t^{[\ell](\alpha,k)}=\dot{\mathbb{U}}_t\gamma_t^{[\ell](\alpha,k)}+\mathscr{S}_t$$

and we have $\|\dot{\mathbb{U}}_t\|_{L^\infty}\lesssim1$.

Combined with the above argument for the bound on $\left|q_t^{[\ell]}(\phi)\right|$, we obtain joint continuity on $[0,T]\times\mathcal{H}$ of $q^{[\ell]}$ and $\partial_t q^{[\ell]}$. Boundedness and continuity of $\mathrm{D}_\phi q^{[\ell]}$ and $\mathrm{D}_\phi^2 q^{[\ell]}$ on bounded subsets of $[0,T]\times\mathcal{H}$ follow along similar lines by noting that (6.1) ensures that they all depend on $\phi$ through smooth gradients and Hessians applied to the point evaluations $\partial^\mu\delta_x$.

For condition 2, note that similarly to the above argument for the bound on $\left|q_t^{[\ell]}(\phi)\right|$, there exists an $N\in\mathbb{N}$ and polynomial $P$ of $\phi(y)$ and its derivatives such that

$$\|\mathrm{D}q^{[\ell]}|_{\varphi_s}^\sharp\|_{\mathcal{H}}=\sup_{\|v\|_{\mathcal{H}}\leq1}|\mathrm{D}q^{[\ell]}|_{\varphi_s}(v)|\lesssim|||\gamma_t^{[\ell](\alpha,k)}|||_{0,\Theta,\eta}|P(\varphi_s(y))e^{-\Im(\alpha)\varphi_s(y)}|||\varphi_s\|_{\mathcal{H}}^N.$$

Then plugging in the bounds on the kernels $|||\gamma^{[\ell](\alpha,k)}|||$ of the previous section, (4.12) and Lemmas 4.6 and A.9, we obtain that there exists a $\tilde{N}\in\mathbb{Z}$ such that

$$\mathbb{E}\int_0^T\|Q_s^*\mathrm{D}q^{[\ell]}|_{\varphi_s}^\sharp\|_{L^2}^2\mathrm{d}s\leq\int_0^T\|Q_s^*\|_{\mathcal{L}^2(\mathcal{H};L^2)}^2\mathbb{E}\|\mathrm{D}q^{[\ell]}|_{\varphi_s}^\sharp\|_{\mathcal{H}}^2\mathrm{d}s\lesssim\int_0^T\langle s\rangle^{\tilde{N}}\mathrm{d}s<\infty.$$

We thus have proven that $\alpha q^{\langle1\rangle}+q^{\langle2\rangle}$ satisfy Conditions 1 and 2 of Definition 2.5. Extension to $\alpha q^{\langle1\rangle}q^{\langle2\rangle}$ is straightforward. □

**Proof. (of Theorem 1.1)** Consider the germs $q_t^{[\leq \ell^*]}$ constructed in the previous section. As we have seen in Propositions 5.3, 5.9, 5.12 and 5.15, they are in the spaces $\hat{\mathscr{O}}^{\kappa,\sigma}_{\infty-,i\alpha,y}$ for $\kappa := \frac{|\mu|}{2} + \frac{\Re(\alpha)^2 - \Im(\alpha)^2}{8\pi}$ and $\sigma := \ell^*\delta - \frac{\Re(\alpha)^2 - \Im(\alpha)^2}{8\pi} - \frac{|\mu|}{2}$. Pick $p > 1$ arbitrary. Then (4.10) yields

$$\begin{aligned}\left\|\mathscr{L}^V q_t^{[\leq \ell^*]}(\varphi_t)\right\|_{L^1(\Omega)} &\leq \|e^{i\alpha\phi(y)}\|_{L^p(\Omega)} \left\|e^{-i\alpha\varphi_t(y)}\mathscr{L}^V q_t^{[\leq \ell^*]}(\varphi_t)\right\|_{L^{p'}(\Omega)} \\ &\lesssim \langle t\rangle^{\frac{|\mu|}{2} + \frac{\Re(\alpha)^2}{8\pi} + (p-1)\frac{\Im(\alpha)^2}{8\pi} - \ell^*\delta - 1} \mathfrak{l}_t^c \|q^{[\leq \ell^*]}\|_{\hat{\mathscr{O}}^{\kappa,\sigma}_{\infty-,i\alpha,y}}\end{aligned}$$

so as $\delta > 0$, choosing

$$\ell^* := \left\lceil \frac{\frac{|\mu|}{2} + \frac{\Re(\alpha)^2}{8\pi} + (p-1)\frac{\Im(\alpha)^2}{8\pi} + 1}{\delta} \right\rceil < \infty \tag{6.2}$$

ensures that the germs $q_t^{[\leq \ell^*]}$ satisfy (2.10). Coupled with the previous lemma, this means that they are coherent; the statement about existence then follows from Lemma 2.6. The existence of $\varphi$ and the convergence of measures is the content of Theorem 4.4 and Corollary 4.5. □

A couple of remarks are in order.

**Remark 6.2. (Sharp estimates)** Obviously the Hölder inequality we use in the above proof does not give the optimal scaling for exponentials. Instead, recalling the proofs of Propositions 5.9 and 5.15, we see that the field content of $q_t^{[\leq \ell^*]}(\phi)$ and $\mathscr{L}^V q_t^{[\leq \ell^*]}(\phi)$ is of the form $e^{i\alpha\phi(y)}P_\mu(\phi)$ with $P_\mu$ polynomial in $\phi(y), \psi_{1:(\ell-1)}$ and their derivatives. Thus proceeding similarly to the proof of (4.11), the sharp estimates

$$\left\|q_t^{[\leq \ell^*]}(\varphi_t)\right\|_{L^1(\Omega)} \lesssim \langle t\rangle^{\frac{|\mu|}{2} + \frac{\Re(\alpha)^2}{8\pi}} \mathfrak{l}_t^c, \qquad \left\|\mathscr{L}^V q_t^{[\leq \ell^*]}(\varphi_t)\right\|_{L^1(\Omega)} \lesssim \langle t\rangle^{\frac{|\mu|}{2} + \frac{\Re(\alpha)^2}{8\pi} - \ell^*\delta - 1} \mathfrak{l}_t^c$$

follow. We chose to include the suboptimal bounds in the above proof and the subsequent results, since they suffice for our purposes while having simpler presentation.

**Remark 6.3. (Required order of expansion)** Clearly, the $+1$ in (6.2) is not necessary apart from the edge cases where the fraction, without it, would be an integer. Then using the sharp estimates in the above remark, recalling that $\delta := 1 - \frac{\beta^2}{8\pi}$, $\beta^2 \in (0, 8\pi)$ and plugging in the saturated optimal condition $|\Re(\alpha)|\beta = 4\pi$ into (6.2), we obtain that the existence of the observable $[\![\partial^\mu e^{i\alpha\phi(y)}]\!]$ in the optimal regime requires the flow equation to be solved up to iteration

$$\ell^*(\beta^2, \mu) = \left\lceil \frac{4\pi(\beta^2|\mu| + 4\pi)}{\beta^2(8\pi - \beta^2)} \right\rceil. \tag{6.3}$$

This blows up as approaching criticality, $\beta^2 \to 8\pi$, indicating the breakdown of perturbation theory, and also as $\beta^2 \to 0$ due to the allowed superficial degree of divergence blowing up along the saturation of the condition $|\Re(\alpha)| < \frac{4\pi}{\beta}$.

The function (in $\beta^2$) inside the ceiling in (6.3) has global minimum at

$$\beta^2_{\min} = \frac{8\pi}{1 + \sqrt{1 + 2|\mu|}} \qquad \text{with value} \qquad \frac{(1 + \sqrt{1 + 2|\mu|})^2}{4}$$

When $|\mu| = 0$, i.e. in the case of the observable $[\![e^{i\alpha\phi(y)}]\!]$, the function $\beta^2 \mapsto \frac{16\pi^2}{\beta^2(8\pi - \beta^2)}$ is symmetric around its global minimum at the free fermion point, $\beta^2_{\min} = 4\pi$, with value 1. That is, exactly at the free fermion point, even first-order perturbation theory is sufficient to construct the observable in the full optimal regime $|\Re(\alpha)| < \frac{4\pi}{\beta}$, and none of the screening theory developed in Section 5.2 is required.

As we will soon see, the above considerations are true not only for the 1-point function $\langle [\![e^{i\alpha\phi(y)}]\!] \rangle$, but also for general $n$-point functions $\langle [\![e^{i\alpha\phi(y)}]\!] \star [\![\phi^{k_1}(y_1)]\!] \star \ldots [\![\phi^{k_n}(y_n)]\!] \rangle$ with $y \neq y_1 \ldots \neq y_n$, since polynomials of the field have $\kappa = 0$ superficial degree of divergence owing to their log-blowup. For products with other vertex fields or derivatives, higher $\ell^*$ is needed, so that the decay $\langle t\rangle^{-\sigma}$ coming from the accuracy $\sigma$ of the germ can offset the other germs' divergence on top of its own.

**Remark 6.4. (Independence of order of expansion)** The previous remark suggests the need for discussion on whether the observables we construct depend on the order of expansion. Suppose there exists an $\ell^*$ (not necessarily given by (6.2)) for which the construction of Theorem 1.1 goes through, i.e. such that

$$\frac{|\mu|}{2}+\frac{\Re(\alpha)^2}{8\pi}+(p-1)\frac{\Im(\alpha)^2}{8\pi}-\ell^*\delta<0 \tag{6.4}$$

for some $p>1$. Then for any $\tilde{\ell}>\ell^*$, by the proofs of the previous section, in particular (5.11), (5.28), (5.39), (5.48), we have

$$\begin{aligned}\left\|q_t^{[\le\tilde{\ell}]}(\varphi_t)-q_t^{[\le\ell^*]}(\varphi_t)\right\|_{L^1(\Omega)} &\le \sum_{\ell=\ell^*+1}^{\tilde{\ell}}\|e^{\langle i\alpha,\varphi_t(y)\rangle}\|_{L^p(\Omega)}\left\|e^{-\langle i\alpha,\varphi_t(y)\rangle}q_t^{[\ell]}(\varphi_t)\right\|_{L^{p'}(\Omega)}\\ &\lesssim \sum_{\ell=\ell^*+1}^{\tilde{\ell}}\langle t\rangle^{\frac{|\mu|}{2}+\frac{\Re(\alpha)^2}{8\pi}+(p-1)\frac{\Im(\alpha)^2}{8\pi}}(\lambda_t\langle t\rangle^{-1})^{\ell-1}\mathfrak{r}_t^c\\ &\lesssim \langle t\rangle^{\frac{|\mu|}{2}+\frac{\Re(\alpha)^2}{8\pi}+(p-1)\frac{\Im(\alpha)^2}{8\pi}-\ell^*\delta}\mathfrak{r}_t^c\end{aligned}$$

which goes to 0 as $t\to\infty$ by (6.4). Thus $q^{[\le\tilde{\ell}]}\sim_\varphi q^{[\le\ell^*]}$, and hence by (2.5) the germ $q^{[\le\tilde{\ell}]}$ reconstructs the same martingale as $q^{[\le\ell^*]}$. In particular, the correlation functions we constructed are independent of the order of perturbative approximation.

To construct the observable product, we exploit that the germs obtained as approximate solution on the observable flow equation (2.18) localize deep in the UV. This is the content of the following

**Lemma 6.5.** *Let $q_t^{[\ell_j]\langle j\rangle}$, $j\in\{1,2\}$ be given by (6.1), localized at $y_1\neq y_2\in\mathbb{R}^2$, respectively. Then there exists a $K=K(\mu_j,n_j,\alpha_j)$ and $C=C(\ell_1,\ell_2)$ such that for all $p\in[1,\infty)$,*

$$\left\|e^{-i\langle\boldsymbol{\alpha},\varphi_t(\boldsymbol{y})\rangle}\mathrm{Tr}\left[\dot{G}_t\mathrm{D}q_t^{[\ell_1]\langle 1\rangle}(\varphi_t)\mathrm{D}q_t^{[\ell_2]\langle 2\rangle}(\varphi_t)\right]\right\|_{L^p(\Omega)}\lesssim\langle t\rangle^K e^{-\frac{\sqrt{t}}{C}|y_1-y_2|}\prod_{j\in 1:2}\sum_{k_j\le n_j}|||\gamma_{t,T}^{[\ell_j](\alpha_j,k_j)}|||_{t,\Theta(\ell_j)}$$

*with $\boldsymbol{\alpha}:=\alpha_1\veebar\alpha_2:=(\alpha_1,\alpha_2)\in\mathbb{C}^2$ and $\boldsymbol{y}:=y_1\veebar y_2:=(y_1,y_2)\in(\mathbb{R}^d)^2$.*

**Proof.** Evaluating the functional derivative on (6.1), we see that after contraction, $\mathrm{Tr}\left[\mathrm{D}q_t^{[\ell_1]\langle 1\rangle}(\phi)\mathrm{D}q_t^{[\ell_2]\langle 2\rangle}(\phi)\right]$ are linear combination of terms of the form

$$\iint\prod_{j\in 1:2}\gamma_t^{[\ell_j](\alpha_j,k_j)}(\xi_{1:(\ell_j-1)},y_j)\mathfrak{G}(\boldsymbol{\xi},\boldsymbol{y})P^\mu_{k_j,\alpha_j}(\phi)\mathrm{d}\xi_{1:(\ell_1-1)}\mathrm{d}\xi_{1:(\ell_2-1)} \tag{6.5}$$

for some $k_j\le n_j$, where $P^\mu_{k_j,\alpha_j}(\phi)$ are polynomials of $\phi(y_j)$, $\psi_{y_j}^{(\alpha_j)}$, $\psi_{1:(\ell_j-1)}$ and their derivatives and where $\mathfrak{G}(\boldsymbol{\xi},\boldsymbol{y})$ is linear combination of $\dot{\mathcal{G}}$ and its derivatives, contracting between the sets $(x_{1:(\ell_1-1)},y_1)$ and $(x_{1:(\ell_2-1)},y_2)$. Inserting

$$1=\prod_{j\in 1:2}\omega_t^\Theta(y_j,\xi_{1:(\ell_j-1)})\omega_{-t}^\Theta(y_j,\xi_{1:(\ell_j-1)})$$

with weight $\omega_t^\Theta$ defined in (3.4), for any $k_j\in 1:(\ell_j-1)$ and multiindices $\nu_1,\nu_2$, bound

$$\begin{aligned}&e^{-\Theta(\ell_1)\sqrt{t}\,\mathrm{St}(x_{1:(\ell_j-1)},y_j)}\partial_{x_{k_1}}^{\nu_1}\partial_{x_{k_2}}^{\nu_2}\dot{\mathcal{G}}_t(x_{k_1}-x_{k_2})e^{-\Theta(\ell_2)\sqrt{t}\,\mathrm{St}(x_{1:(\ell_j-1)},y_j)}\\ \lesssim\;&\langle t\rangle^{\frac{|\nu_1|+|\nu_2|}{2}-1}e^{-\frac{\sqrt{t}}{4^{\ell_1}}|x_{k_1}-y_1|}e^{-\frac{t}{8}|x_{k_1}-x_{k_2}|^2}e^{-\frac{\sqrt{t}}{4^{\ell_2}}|x_{k_2}-y_2|}\\ \le\;&\langle t\rangle^{\frac{|\nu_1|+|\nu_2|}{2}-1}\exp\left(-\frac{\sqrt{t}}{4^{\ell_1\vee\ell_2}}|y_1-y_2|\right)e^{\frac{2}{16^{\ell_1\vee\ell_2}}}\end{aligned}$$

where in the second line we absorbed into the Gaussian the polynomial obtained from differentiation, and the last inequality follows, with $C := 4^{\ell_1 \vee \ell_2}$, from

$$\begin{aligned} \frac{\sqrt{t}}{C}|y_1 - y_2| - \frac{2}{C^2} &\leq \frac{1}{8}\left(\sqrt{t}|x_{k_1} - x_{k_2}| - \frac{4}{C}\right)^2 + \frac{\sqrt{t}}{C}|y_1 - y_2| - \frac{2}{C^2} \\ &= \frac{t}{8}|x_{k_1} - x_{k_2}|^2 + \frac{\sqrt{t}}{C}(|y_1 - y_2| - |x_{k_1} - x_{k_2}|). \end{aligned}$$

Similar estimates hold when the contraction is $\dot{\mathcal{G}}_t(x_{k_1} - y_2)$ and $\dot{\mathcal{G}}_t(y_1 - y_2)$, with the corresponding weights $\omega_{-t}^{\Theta}$ estimated trivially by 1. The statement then follows from Minkowski's integral inequality, the multinomial theorem, Faà di Bruno, Hölder and Lemma 4.9. □

**Proof. (of Theorem 1.2)** Given $\mathfrak{Q}^{(j)} \in \mathcal{O}_W^{\text{loc}}$, by results of the above section we may pick germs $q^{[\leq \ell_j^*]\langle j \rangle} \in \hat{\mathscr{O}}^{\kappa_j, \sigma_j}_{\infty-, i\alpha_j, y_j}$ with $\kappa_j := \frac{|\mu_j|}{2} + \frac{\Re(\alpha_j)^2 - \Im(\alpha_j)^2}{8\pi}$, $\sigma_j := \ell_j^* \delta - \kappa_j$. Choose

$$\ell_j^* := \left\lceil \frac{1 + \sum_{m=1}^n \kappa_m + \frac{n}{8\pi}\sum_{m=1}^n \Im(\alpha_m)^2}{\delta} \right\rceil \tag{6.6}$$

for any $j \in 1\colon n$, so that

$$\sigma_{1:n} := \bigwedge_{j \in 1:n} \left(\sigma_j - \sum_{m \in (1:n) \setminus \{j\}} \kappa_m\right) = \bigwedge_{j \in 1:n} \left(\ell_j^* \delta - \sum_{m=1}^n \kappa_m\right) \geq \frac{n}{8\pi} \sum_{m=1}^n \Im(\alpha_m)^2 + 1. \tag{6.7}$$

Combining Lemma 6.5 with estimates on the kernels $|||\gamma_{t,T}^{[\ell_j](\alpha_j, k_j)}|||_{t, \Theta(\ell_j)}$, of the previous section and the fact that $q^{[\leq \ell_j^*]\langle j \rangle} \in \hat{\mathscr{O}}^{\kappa_j, \sigma_j}_{\infty-, i\alpha_j, y_j}$, for any $p \in [1, \infty)$ we obtain that for any $j_1 \neq j_2 \in 1\colon n$,

$$\left\| e^{-i\langle \boldsymbol{\alpha}_{1:n}, \varphi_t(\boldsymbol{y}_{1:n})\rangle} \mathrm{Tr}\Big[\dot{G}_t \mathrm{D} q_t^{[\leq \ell_{j_1}^*]\langle j_1 \rangle}(\varphi_t) \mathrm{D} q_t^{[\leq \ell_{j_2}^*]\langle j_2 \rangle}(\varphi_t)\Big] \prod_{j \in 1:n \setminus \{j_1, j_2\}} q_t^{[\leq \ell_j^*]\langle j \rangle}(\varphi_t) \right\|_{L^p(\Omega)} \lesssim \langle t \rangle^{\tilde{K}} e^{-C\sqrt{t}}$$

with $\boldsymbol{\alpha}_{1:n} := \alpha_1 \veebar \ldots \veebar \alpha_n$, $\boldsymbol{y}_{1:n} := y_1 \veebar \ldots \veebar y_n$ for some $\tilde{K}$ and $C = C(\ell_{j_1}^*, \ell_{j_2}^*, y_{j_1}, y_{j_2})$. Hence the conditions of Lemma 2.17 are satisfied and thus we have that the product germ is

$$q^{[\leq \boldsymbol{\ell}_{1:n}^*]\langle 1:n \rangle} := \prod_{j \in 1:n} q^{[\leq \ell_j^*]\langle j \rangle} \in \hat{\mathscr{O}}^{\kappa_{1:n}, \sigma_{1:n}}_{\infty-, i\boldsymbol{\alpha}_{1:n}, \boldsymbol{y}_{1:n}}$$

for $\kappa_{1:n} := \sum_{j=1}^n \kappa_j$ and $\sigma_{1:n}$ defined above in (6.7). It is a consequence of (4.10) and $|\mathcal{G}_\infty(x)| \lesssim 1$ for all $x \neq 0$ that for any $p \in [1, \infty)$, $\boldsymbol{\alpha} \in \mathbb{C}^n$, and $\boldsymbol{y} \in (\mathbb{R}^2)^n$ pairwise distinct, we have

$$\|e^{i\langle \boldsymbol{\alpha}, \varphi_t(\boldsymbol{y})\rangle}\|_{L^p(\Omega)} \lesssim e^{\frac{p}{2}\sum_{j=1}^n \Im(\alpha_j)^2 \mathcal{G}_t(0)}.$$

Indeed, $\langle \boldsymbol{\alpha}, \varphi_t(\boldsymbol{y})\rangle = \langle \boldsymbol{\alpha}, B_t(\boldsymbol{y})\rangle + \langle \boldsymbol{\alpha}, Z_t(\boldsymbol{y})\rangle$, $Z_t$ is bounded in $L^\infty$, and $\langle \boldsymbol{\alpha}, B_t(\boldsymbol{y})\rangle$ is a centered Gaussian process with variance $\sum_{j,k=1}^m \alpha_j \alpha_k \mathcal{G}_t(y_j - y_k)$.

Thus choosing $p > 1$ small enough, by Hölder the germ $q^{[\leq \boldsymbol{\ell}_{1:n}^*]\langle 1:n \rangle}$ satisfies (2.10). Indeed, by (6.7)

$$\begin{aligned} \left\|\mathscr{L}^V q_t^{[\leq \boldsymbol{\ell}_{1:n}^*]\langle 1:n \rangle}(\varphi_t)\right\|_{L^1(\Omega)} &\leq \|e^{i\langle \boldsymbol{\alpha}_{1:n}, \varphi_t(\boldsymbol{y}_{1:n})\rangle}\|_{L^p(\Omega)} \left\|e^{-i\langle \boldsymbol{\alpha}_{1:n}, \varphi_t(\boldsymbol{y}_{1:n})\rangle} \mathscr{L}^V q_t^{[\leq \boldsymbol{\ell}_{1:n}^*]\langle 1:n \rangle}(\varphi_t)\right\|_{L^{p'}(\Omega)} \\ &\lesssim \langle t \rangle^{\frac{p}{8\pi}\sum_{j=1}^n \Im(\alpha_j)^2} \langle t \rangle^{-\sigma_{1:n} - 1} \mathbb{I}_t^c \|q^{[\leq \boldsymbol{\ell}_{1:n}^*]\langle 1:n \rangle}\|_{\hat{\mathscr{O}}^{\kappa_{1:n}, \sigma_{1:n}}_{\infty-, i\boldsymbol{\alpha}_{1:n}, \boldsymbol{y}_{1:n}}} \\ &\lesssim \langle t \rangle^{\frac{p-n}{8\pi}\sum_{j=1}^n \Im(\alpha_j)^2 - 2} \mathbb{I}_t^c \end{aligned}$$

and the statement follows. Coupled with Lemma 6.1, this shows that $q^{[\leq \boldsymbol{\ell}_{1:n}^*]\langle 1:n \rangle}$ is coherent.

For the existence of the observable product, it then remains to show that (2.13) holds. Denoting by $r^{\langle j\rangle}$ the canonical Ito remainders (cf. (2.11)) corresponding to the germ $q^{[\leq \ell_j^*]\langle j\rangle}$, and by $r^{\langle 1:n\rangle}$ the remainder corresponding to that of the product, we have

$$\begin{aligned}
&\left\| \mathfrak{R}\left( \prod_{j\in 1:n} q^{[\leq \ell_j^*]\langle j\rangle} \right)_t - \prod_{j\in 1:n} \mathfrak{Q}_t^{(j)} \right\|_{L^1(\Omega)} \\
&= \left\| \prod_{j\in 1:n} q_t^{[\leq \ell_j^*]\langle j\rangle}(\varphi_t) + r_{t,\infty}^{\langle 1:n\rangle} - \prod_{j\in 1:n} \left( q_t^{[\leq \ell_j^*]\langle j\rangle}(\varphi_t) + r_{t,\infty}^{\langle j\rangle} \right) \right\|_{L^1(\Omega)} \\
&\leq \| r_{t,\infty}^{\langle 1:n\rangle} \|_{L^1(\Omega)} + \sum_{\substack{I_1 \uplus I_2 = 1:n \\ I_2 \neq \emptyset}} \left\| \prod_{j_1\in I_1} q_t^{[\leq \ell_{j_1}^*]\langle j_1\rangle}(\varphi_t) \prod_{j_2\in I_2} r_{t,\infty}^{\langle j_2\rangle} \right\|_{L^1(\Omega)}
\end{aligned}$$

where each term in the sum has at least one remainder $r_{t,\infty}^{\langle j_2\rangle}$ due to $\prod_{j\in 1:n} q^{[\leq \ell_j^*]\langle j\rangle}$ being cancelled. The fact that $\| r_{t,\infty}^{\langle 1:n\rangle} \|_{L^1(\Omega)}$ vanishes as $t\to\infty$ is clear from (6.7):

$$\begin{aligned}
&\| r_{t,\infty}^{\langle 1:n\rangle} \|_{L^1(\Omega)} \\
&\leq \int_t^\infty \langle s\rangle^{-\sigma_{1:n}-1} \| e^{i\langle \boldsymbol{\alpha}_{1:n}, \varphi_s(\boldsymbol{y}_{1:n})\rangle} \|_{L^p(\Omega)} \langle s\rangle^{\sigma_{1:n}+1} \left\| e^{-i\langle \boldsymbol{\alpha}_{1:n}, \varphi_s(\boldsymbol{y}_{1:n})\rangle} \mathscr{L}^V q_s^{[\leq \boldsymbol{\ell}_{1:n}^*]\langle 1:n\rangle}(\varphi_s) \right\|_{L^{p'}(\Omega)} \mathrm{d}s \\
&\lesssim \int_t^\infty \langle s\rangle^{\frac{p-n}{8\pi}\sum_{j=1}^n \Im(\alpha_j)^2 - 2} \mathfrak{l}_s^c \mathrm{d}s \, \| q^{[\leq \boldsymbol{\ell}_{1:n}^*]\langle 1:n\rangle} \|_{\hat{\mathscr{O}}_{\infty-, i\boldsymbol{\alpha}_{1:n}, \boldsymbol{y}_{1:n}}^{\kappa_{1:n}, \sigma_{1:n}}}
\end{aligned}$$

and one must choose $p>1$ small enough. Now let $\varepsilon>0$ and denote by $q\geq 0$ for which $(n+\varepsilon)^{-1} + q^{-1} = n^{-1}$. For any $I_1 \uplus I_2 = 1:n$ with $I_2 \neq \emptyset$, by (6.6)

$$\begin{aligned}
&\left\| \prod_{j_1\in I_1} q_t^{[\leq \ell_{j_1}^*]\langle j_1\rangle}(\varphi_t) \prod_{j_2\in I_2} r_{t,\infty}^{\langle j_2\rangle} \right\|_{L^1(\Omega)} \\
&\leq \prod_{j_1\in I_1} \left\| q_t^{[\leq \ell_{j_1}^*]\langle j_1\rangle}(\varphi_t) \right\|_{L^n(\Omega)} \prod_{j_2\in I_2} \| r_{t,\infty}^{\langle j_2\rangle} \|_{L^n(\Omega)} \\
&\leq \langle t\rangle^{\sum_{j_1\in I_1}\kappa_{j_1}} \prod_{j_1\in I_1} \| e^{i\langle \alpha_{j_1}, \varphi_t(y_{j_1})\rangle} \|_{L^{n+\varepsilon}(\Omega)} \langle t\rangle^{-\kappa_{j_1}} \left\| e^{-i\langle \alpha_{j_1}, \varphi_t(y_{j_1})\rangle} q_t^{[\leq \ell_j^*]\langle j\rangle}(\varphi_t) \right\|_{L^q(\Omega)} \\
&\times \prod_{j_2\in I_2} \int_t^\infty \langle s\rangle^{-\sigma_{j_2}-1} \| e^{-\langle \alpha_{j_2}, \varphi_s(y_{j_2})\rangle} \|_{L^{n+\varepsilon}(\Omega)} \langle s\rangle^{\sigma_{j_2}+1} \left\| e^{-\langle \alpha_{j_2}, \varphi_s(y_{j_2})\rangle} \mathscr{L}^V q_s^{[\leq \ell_{j_2}^*]\langle j_2\rangle}(\varphi_s) \right\|_{L^q(\Omega)} \mathrm{d}s \\
&\lesssim \prod_{j\in 1:n} \| q^{[\leq \ell_j^*](j)} \|_{\hat{\mathscr{O}}_{\infty-, \alpha_j, y_j}^{\kappa_j, \sigma_j}} \langle t\rangle^{\sum_{j_1\in I_1}\kappa_{j_1} + \frac{n+\varepsilon}{8\pi}\Im(\alpha_{j_1})^2} \prod_{j_2\in I_2} \int_t^\infty \langle s\rangle^{\frac{n+\varepsilon}{8\pi}\Im(\alpha_{j_2})^2 - \sigma_{j_2} - 1} \mathrm{d}s \\
&\lesssim \langle t\rangle^{\sum_{j\in 1:n}\frac{n+\varepsilon}{8\pi}\Im(\alpha_j)^2 + \sum_{j_1\in I_1}\kappa_j - \sum_{j_2\in I_2}\sigma_{j_2}} \\
&= \langle t\rangle^{\sum_{j\in 1:n}\kappa_j + \frac{n+\varepsilon}{8\pi}\Im(\alpha_j)^2 - \sum_{j_2\in I_2}\ell_{j_2}^*\delta} \\
&\leq \langle t\rangle^{\frac{\varepsilon}{8\pi}\sum_{j\in 1:n}\Im(\alpha_j)^2 - 1}.
\end{aligned}$$

and choosing $\varepsilon>0$ small enough we obtain that the germs $\{q^{[\leq \ell_j^*]\langle j\rangle}\}_{j=1}^n$ are compatible.

It remains to prove associativity, which, in light of Proposition 2.14, entails proving the existence of product with germs grouped together. Let $P$ be an arbitrary partition of $1{:}n$. For any $J\in P$, repeating the argument in the beginning of the proof, we obtain

$$q^{[\leq \boldsymbol{\ell}_J^*]\langle J\rangle} := \prod_{j\in J} q^{[\leq \ell_j^*]\langle j\rangle} \in \hat{\mathscr{O}}_{\infty-, i\boldsymbol{\alpha}_J, \boldsymbol{y}_J}^{\kappa_J, \sigma_J}$$

for $\kappa_J := \sum_{j\in J}\kappa_j$, $\boldsymbol{\alpha}_J := \veebar_{j\in J}\alpha_j$, $\boldsymbol{y}_{1:n} := \veebar_{j\in J} y_j$, $\sigma_j := \ell_j^*\delta - \kappa_j$

$$\sigma_J := \bigwedge_{j\in J}\left( \sigma_j - \sum_{m\in J\setminus\{j\}} \kappa_m \right) = \bigwedge_{j\in J}\left( \ell_j^*\delta - \sum_{m\in J}\kappa_m \right) \geq \frac{n}{8\pi}\sum_{m=1}^n \Im(\alpha_m)^2 + \sum_{m\in(1:n)\setminus J}\kappa_m + 1, \tag{6.8}$$

and that the product $\mathfrak{Q}^{(J)} := \bigstar_{j\in J}\mathfrak{Q}^{(j)} := \mathfrak{R}(q^{[\leq \boldsymbol{\ell}_J^*]\langle J\rangle})$ exist. Denoting by $r^{\langle J\rangle}$ the corresponding Ito remainders,

$$\begin{aligned}
&\left\| \mathfrak{R}\left( \prod_{J\in P} q^{[\leq \boldsymbol{\ell}_J^*]\langle J\rangle} \right)_t - \prod_{J\in P} \mathfrak{R}(q^{[\leq \boldsymbol{\ell}_J^*]\langle J\rangle})_t \right\|_{L^1(\Omega)} \\
= &\left\| \prod_{j\in 1:n} q_t^{[\leq \ell_j^*]\langle j\rangle}(\varphi_t) + r_{t,\infty}^{\langle 1:n\rangle} - \prod_{J\in P} \left( q^{[\leq \boldsymbol{\ell}_J^*]\langle J\rangle}(\varphi_t) + r_{t,\infty}^{\langle J\rangle} \right) \right\|_{L^1(\Omega)} \\
\leq &\ \| r_{t,\infty}^{\langle 1:n\rangle} \|_{L^1(\Omega)} + \sum_{\substack{\mathscr{I}_1 \uplus \mathscr{I}_2 = P \\ \mathscr{I}_2 \neq \emptyset}} \left\| \prod_{J_1\in\mathscr{I}_1} q^{[\leq \boldsymbol{\ell}_{J_1}^*]\langle J_1\rangle}(\varphi_t) \prod_{J_2\in\mathscr{I}_2} r_{t,\infty}^{\langle J_2\rangle} \right\|_{L^1(\Omega)}
\end{aligned}$$

again $\mathscr{I}_2$ nonempty because of the cancellation of $\prod_{J\in P} q^{[\leq \boldsymbol{\ell}_J^*]\langle J\rangle}$. Denoting by $\Im(\boldsymbol{\alpha}_J)^2 := \sum_{j\in J} \Im(\alpha_j)^2$ and adapting the argument above, for all $\mathscr{I}_1 \uplus \mathscr{I}_2 = P$ with $\mathscr{I}_2 \neq \emptyset$ we obtain

$$\begin{aligned}
&\left\| \prod_{J_1\in\mathscr{I}_1} q^{[\leq \boldsymbol{\ell}_{J_1}^*]\langle J_1\rangle}(\varphi_t) \prod_{J_2\in\mathscr{I}_2} r_{t,\infty}^{\langle J_2\rangle} \right\|_{L^1(\Omega)} \\
\leq &\prod_{J_1\in\mathscr{I}_1} \| q^{[\leq \boldsymbol{\ell}_{J_1}^*]\langle J_1\rangle}(\varphi_t) \|_{L^{|P|}(\Omega)} \prod_{J_2\in\mathscr{I}_2} \| r_{t,\infty}^{\langle J_2\rangle} \|_{L^{|P|}(\Omega)} \\
\lesssim &\prod_{J\in P} \| q^{[\leq \boldsymbol{\ell}_J^*]\langle J\rangle} \|_{\hat{\mathscr{O}}_{\infty-, i\boldsymbol{\alpha}_J, \boldsymbol{y}_J}^{\kappa_J, \sigma_J}} \langle t\rangle^{\sum_{J_1\in\mathscr{I}_1} \kappa_{J_1} + \frac{|P|+\varepsilon}{8\pi}\Im(\boldsymbol{\alpha}_{J_1})^2} \prod_{J_2\in\mathscr{I}_2} \int_t^\infty \langle s\rangle^{\frac{|P|+\varepsilon}{8\pi}\Im(\boldsymbol{\alpha}_{J_2})^2 - \sigma_{J_2} - 1} \mathrm{d}s \\
\lesssim &\ \langle t\rangle^{\frac{|P|+\varepsilon}{8\pi}\sum_{J\in P}\Im(\boldsymbol{\alpha}_J)^2 + \sum_{J_1\in\mathscr{I}_1}\kappa_{J_1} - \sum_{J_2\in\mathscr{I}_2}\sigma_{J_2}} \\
\lesssim &\ \langle t\rangle^{\frac{|P|-n+\varepsilon}{8\pi}\Im(\boldsymbol{\alpha}_{1:n})^2 - 1}
\end{aligned}$$

where the last line follows from (6.8) and $\mathscr{I}_1 \cap \mathscr{I}_2 = \emptyset$. Since $|P| \leq n$, choosing $\varepsilon > 0$ small enough we obtain that the product $\bigstar_{J\in P}\mathfrak{Q}^{(J)}$ exist; it then follows from Proposition 2.14 that

$$\bigstar_{j\in 1:n}\mathfrak{Q}^{(j)} = \bigstar_{J\in P}(\bigstar_{j\in J}\mathfrak{Q}^{(j)})$$

for any partition $P$ of $1{:}\,n$, as requested.

The claimed algebraic properties of the product then follow from Proposition 2.14. □

**Remark 6.6.** Combining the compatibility of the germs with Remark 6.4, it follows from Proposition 2.13 that the observable products do not depend on the order of perturbative approximation.

**Corollary 6.7.** *The space of multilocal observables*

$$\mathscr{O}_{\mathrm{mloc}}^W := \mathrm{span}\{\bigstar_{j\in 1:n}\mathfrak{Q}^{(j)}(y_j) \colon n\in\mathbb{N}, \mathfrak{Q}^{(j)} \in \mathscr{O}_{\mathrm{loc}}^W, y_i \neq y_j \text{ if } i\neq j\}$$

*is well-defined.*

# 7 Smoothness of correlation functions

This section is concerned with the proof of Theorem 1.3.

In the above section we saw that there exist $q^{\langle j\rangle} := q^{[\leq \ell_j^*]\langle j\rangle}$ coherent germs for $\mathfrak{Q}^{(j)}$, given by the ansatz (6.1), such that $(q^{\langle 1\rangle}, \ldots, q^{\langle n\rangle})$ are compatible (cf. Theorem 1.2). By Ito reconstruction, we have

$$\langle \mathfrak{Q}^{(1)}(y_1) \star \ldots \star \mathfrak{Q}^{(n)}(y_n) \rangle = q_0^{\langle 1:n\rangle}(\varphi_0) + \mathbb{E}\int_0^\infty \mathscr{L}^V q_s^{\langle 1:n\rangle}(\varphi_s)\mathrm{d}s \tag{7.1}$$

with $q^{\langle 1:n\rangle} := \prod_{j\in 1:n} q^{\langle j\rangle}$ the product germ.

Both of the above terms contain reference to the germ kernels $\gamma_t^{[\ell](\alpha,k)}$. The difficulty of the proof is that such kernels are *not* smooth in the localization point $y$. Indeed, differentiating the kernels in $y$ may destroy the inductions of Section 5: looking e.g. at (5.18), one obtains a term with $\partial_y^\nu \mathbf{G}_{t,s}(y, \xi_{1:(\ell-1)})$. The following simple example shows that there indeed exist non-smooth germ kernels.

**Example 7.1.** Let $\alpha > 0$ real and consider the kernel $\gamma_t := \gamma_t^{[2]}(\xi_1, y)$ for the vertex field with $\sigma_1 := -1$. With $r := x_1 - y$, the covariance is

$$\dot{\mathbf{G}}_t(y, \xi_1) = -\left(\frac{\alpha^2}{2}\dot{\mathcal{G}}_t(0) + \frac{\beta^2}{2}\dot{\mathcal{G}}_t(0) - \alpha\beta\dot{\mathcal{G}}_t(r)\right).$$

Referring back to Lemma 5.5 and (5.20), we have, for a nonzero constant $C(t)$ depending on $t$,

$$\gamma_t(r) = C(t)\alpha\beta\int_t^\infty e^{\alpha\beta\mathcal{G}_{t,s}(r)}\dot{\mathcal{G}}_s(r)\mathrm{d}s = C(t)(e^{\alpha\beta\mathcal{G}_{t,\infty}(r)} - 1).$$

As $|r| \to 0$, we have

$$\mathcal{G}_{t,\infty}(r) = \frac{1}{4\pi}\int_t^\infty s^{-1}e^{-\frac{m^2}{s}}e^{-\frac{s}{4}|r|^2}\mathrm{d}s = \frac{1}{4\pi}\log(|r|^{-2}) + g_t(r)$$

for $g_t(\cdot)$ uniformly bounded. Hence $\gamma_t(r) \asymp_t |r|^{-\frac{\alpha\beta}{2\pi}}$.

To circumvent this problem, we have to shift all derivatives away from the germ kernels. This is possible by the following simple but crucial

**Lemma 7.2.** *Let $\gamma_t^{[\ell](\alpha,k)}$ be a kernel obtained in Section 5. Then for any $h \in \mathbb{R}^2$,*

$$\gamma_t^{[\ell](\alpha,k)}(\sigma_{1:(\ell-1)}, x_1, \ldots, x_{\ell-1}, y) = \gamma_t^{[\ell](\alpha,k)}(\sigma_{1:(\ell-1)}, x_1 + h, \ldots, x_{\ell-1} + h, y + h).$$

**Proof.** It is straightforward to see inductively from the flow equations for the germ kernels (Lemmas 5.1 5.5 5.11, 5.14) that they inherit translation invariance from the fluctuation covariance $\dot{\mathcal{G}}_t$ and the potential kernels $v_t^{[\ell]}$. □

**Proof. (of Theorem 1.3)**
To simplify notation, write $I := 1\colon(\ell - 1)$.

**Smoothness of $q_0^{\langle 1:n\rangle}(\varphi_0)$**
Fix an arbitrary charge configuration $\sigma_I$. Denoting by

$$\Phi(y, x_I) := \phi^k(y)e^{i\alpha\phi(y)}\prod_{j\in I}\rho(x_j)e^{i\beta\sigma_j\phi(x_j)},$$

the field content, by translation invariance

$$\int\gamma_t^{[\ell](\alpha,k)}(x_I, y)\partial_{y,I}^\mu\Phi(y, x_I)\mathrm{d}x_I = \int\gamma_t^{[\ell](\alpha,k)}(z_I, 0)[\partial_{y,I}^\mu\Phi(y, x_I)]_{x_I = y + z_I}\mathrm{d}z_I.$$

Since $\varphi_0 = 0$, we have that $q_0(\varphi_0)$ is a sum of terms of the form

$$\int\gamma_0^{[\ell](\alpha,0)}(z_I, 0)[\partial_I^\mu\rho(x_I)]_{x_I = y + z_I}\mathrm{d}z_I.$$

Since $\rho \in \mathcal{C}_b^\infty$ and the kernel estimates imply $\gamma_0^{[\ell](\alpha,k)}(\cdot, 0) \in L^1$, by the Leibniz integral rule we obtain that $y \mapsto q_0^{\langle j\rangle}(\varphi_0)$ is smooth.

**Smoothness of $\mathbb{E}\int_0^\infty\mathscr{L}^V q_s^{\langle 1:n\rangle}(\varphi_s)\mathrm{d}s$**
The generator of the product is

$$\mathscr{L}^V q^{\langle 1:n\rangle} = \sum_j \mathscr{L}^V q^{\langle j\rangle}\prod_{k\neq j} q^{\langle k\rangle} + \sum_{j<k}\mathrm{Tr}[\dot{G}\mathrm{D}q^{\langle j\rangle}\mathrm{D}q^{\langle k\rangle}]\prod_{r\neq j,k} q^{\langle r\rangle}.$$

Recall from the proofs of Section 5 the formulae for $\langle \mathrm{D}q_t^{[\ell]}(\phi), \dot{G}_t R_t^{[>(\ell^*-\ell)]}\rangle$: we have

$$\langle \mathrm{D}q_t^{[\ell]}(\phi), \dot{G}_t R_t^{[>(\ell^*-\ell)]}\rangle = \sum \int \gamma_t^{[\ell](\alpha,k)}(\xi_I, y)\Upsilon_t(\phi, \xi_I, y)\mathrm{d}\xi_I$$

with $\Upsilon_t$ functions depending on point evaluations of $\dot{G}_t R_t^{[>(\ell^*-\ell)]}$. By translation invariance,

$$\int \gamma_t^{[\ell](\alpha,k)}(x_I, y)\Upsilon_t(\phi, x_I, y)\mathrm{d}x_I \;=\; \int \gamma_t^{[\ell](\alpha,k)}(z_I, 0)\Upsilon_t(\phi, z_I + y, y)\mathrm{d}z_I.$$

Hence we have, for any $\nu$ multiindex, that

$$\partial_y^\nu \langle \mathrm{D}q_t^{[\ell]}(\phi), \dot{G}_t R_t^{[>(\ell^*-\ell)]}\rangle = \sum \int \gamma_t^{[\ell](\alpha,k)}(\xi_I, 0)\partial_y^\nu \Upsilon_t(\phi, \xi_I + y, y)\mathrm{d}\xi_I$$

where the validity of the Leibniz integral rule is checked similarly to the end of the proof of Lemma 5.10. Then

$$\begin{aligned} \sup_y \left\| e^{-i\alpha\varphi_t(y)} \partial_y^\nu \mathscr{L}^V q_t^{[\leq \ell^*]}(\varphi_t)\right\|_{L^p(\Omega)} &\leq \sup_y \sum_{\ell=1}^{\ell^*} \left\| e^{-i\alpha\varphi_t(y)} \partial_y^\nu \langle \mathrm{D}q_t^{[\ell]}, \dot{G}_t R_t^{[>(\ell^*-\ell)]}\rangle \right\|_{L^p(\Omega)} \\ &\lesssim \sup_y \sum_{\ell=1}^{\ell^*} \| |\gamma_t^{[\ell]}| \|_{t,\Theta(\ell)} \langle t\rangle^{\frac{|\mu|+|\nu|}{2} - 1 - \delta - (\ell^*-\ell)\delta} \mathfrak{l}_t^c \\ &\lesssim \langle t\rangle^{\frac{|\nu|}{2} + \kappa - \ell^*\delta - 1} \mathfrak{l}_t^c \end{aligned} \tag{7.2}$$

for any $p \geq 1$ with the standard $\kappa$ of the previous section (i.e. $\kappa := \frac{|\mu|}{2} + \frac{\Re(\alpha)^2 - \Im(\alpha)^2}{8\pi}$) with constants independent of $y$.

Similarly, we may shift the derivatives away from the kernels in (6.1) to the field content to obtain

$$\sup_y \left\| e^{-i\alpha\varphi_t(y)} \partial_y^\nu q_t^{[\leq \ell^*]}(\varphi_t)\right\|_{L^p(\Omega)} \lesssim \langle t\rangle^{\frac{|\nu|}{2} + \kappa} \mathfrak{l}_t^c. \tag{7.3}$$

Fix a compact $K \subseteq (\mathbb{R}^2)^n \setminus \mathrm{Diag}_n$ and let $j_1 \neq j_2 \in 1{:}\,n$. For any multiindices $\nu_{j_1}, \nu_{j_2}$, we may shift the derivatives away from the kernels in $\partial_{y_{j_1}}^{\nu_{j_1}} \partial_{y_{j_2}}^{\nu_{j_2}} \mathrm{Tr}\Big[\dot{G}_t \mathrm{D}q_t^{[\leq \ell_{j_1}^*]\langle j_1\rangle}(\varphi_t) \mathrm{D}q_t^{[\leq \ell_{j_2}^*]\langle j_2\rangle}(\varphi_t)\Big]$ to the field content and the heat kernel contraction $\dot{\mathcal{G}}_t$ (see e.g. (6.5)). Thus denoting by $\boldsymbol{\alpha} := \alpha_{j_1} \veebar \alpha_{j_2}$, $\boldsymbol{y} := y_{j_1} \veebar y_{j_2}$, the proof of Lemma 6.5 goes through and we obtain

$$\begin{aligned} &\sup_{y_{1:n} \in K} \left\| e^{-i\langle \boldsymbol{\alpha}, \varphi_t(\boldsymbol{y})\rangle} \partial_{y_{j_1}}^{\nu_{j_1}} \partial_{y_{j_2}}^{\nu_{j_2}} \mathrm{Tr}\Big[\dot{G}_t \mathrm{D}q_t^{[\leq \ell_{j_1}^*]\langle j_1\rangle}(\varphi_t) \mathrm{D}q_t^{[\leq \ell_{j_2}^*]\langle j_2\rangle}(\varphi_t)\Big] \right\|_{L^p(\Omega)} \\ \lesssim\; &\sup_{y_{1:n} \in K} \langle t\rangle^M e^{-C\sqrt{t}|y_{j_1} - y_{j_2}|} \\ \lesssim\; &\langle t\rangle^M e^{-\tilde{C}\sqrt{t} d_K} \end{aligned} \tag{7.4}$$

for some $M$ depending on $\nu_{j_1}, \nu_{j_2}$, where

$$d_K := \inf_{y_{1:n} \in K} \min_{j<k} |y_j - y_k| > 0$$

by the assumption $K \subseteq (\mathbb{R}^2)^n \setminus \mathrm{Diag}_n$.

Fix $m \in \mathbb{N}$. By choosing $\ell_j^* = \ell_j^*(m)$ large enough, (7.2), (7.3), (7.4) and Lemma 4.6 yields that for any multiindindex tuple $\boldsymbol{\nu} = \nu_1, \ldots, \nu_n$ with $|\boldsymbol{\nu}| \leq m$, we have

$$\int_0^\infty \sup_{y_{1:n} \in K} \left\| \partial_{y_1}^{\nu_1} \ldots \partial_{y_n}^{\nu_n} \mathscr{L}^V q_s^{\langle 1:n\rangle}(\varphi_s)\right\|_{L^1(\Omega)} \mathrm{d}s < \infty$$

(cf. the proof of Theorem 1.2; see also Remark 6.4). For each fixed $t>0$, the bounds above also imply continuity of $(y_1,\dots,y_n)\mapsto\partial_{y_1}^{\nu_1}\dots\partial_{y_n}^{\nu_n}\mathscr{L}^V q_s^{\langle 1:n\rangle}(\varphi_s)$ in $L^1(\Omega)$. Thus we may differentiate under the integrals in (7.1), and the statement follows from dominated convergence. □

# A Heat kernel interpolation

## A.1 Heat kernel estimates

Here we collect some standard heat kernel estimates. Recall that $\dot{\mathcal{G}}_s$ and $k_s^Q$ denote the $L^2$ integral kernels of $\dot{G}_s$ and $Q_s$, respectively, cf. (2.7). Since the bounds are needed mostly for the part of the paper concerned with the special case of the sine–Gordon model, we take $d=2$ unless explicitly stated otherwise.

**Lemma A.1.** *There is a uniformly bounded function $h$ such that*

$$\mathcal{G}_t(0)=\frac{1}{4\pi}\log(1\vee t)+h(t)$$

**Proof.** Simply split

$$\begin{aligned}\mathcal{G}_t(0) &= \frac{1}{4\pi}\int_0^{1\wedge t}s^{-1}e^{-\frac{m^2}{s}}\mathrm{d}s+\frac{\mathbb{1}_{t>1}}{4\pi}\int_1^t s^{-1}e^{-\frac{m^2}{s}}\mathrm{d}s\\ &= \frac{1}{4\pi}\int_0^{1\wedge t}s^{-1}e^{-\frac{m^2}{s}}\mathrm{d}s+\frac{\mathbb{1}_{t>1}}{4\pi}\int_1^t s^{-1}\Big(e^{-\frac{m^2}{s}}-1\Big)\mathrm{d}s+\frac{\mathbb{1}_{t>1}}{4\pi}\log(t)\\ &=: h(t)+\frac{1}{4\pi}\log(1\vee t)\end{aligned}$$

□

**Lemma A.2.** *For any $\varsigma,\eta,\alpha\geq 0$, we have*

$$\int_{\mathbb{R}^2}\dot{\mathcal{G}}_t(x)|\cdot|^{\alpha}e^{\varsigma\sqrt{t}|x|}\langle x\rangle^{\eta}\mathrm{d}x \;\lesssim\; \langle t\rangle^{-2-\frac{\alpha}{2}}$$

**Proof.** Indeed, since $\varsigma\geq 0$, going to spherical coordinates yields

$$\begin{aligned}\int_{\mathbb{R}^2}\dot{\mathcal{G}}_t(x)e^{\varsigma\sqrt{t}|x|}|\cdot|^{\alpha}\langle x\rangle^{\eta}\mathrm{d}x &\lesssim t^{-1}e^{-\frac{m^2}{t}}\int_{\mathbb{R}^2}e^{-\frac{t}{4}|x|^2+\varsigma\sqrt{t}|x|}|\cdot|^{\alpha}\langle x\rangle^{\eta}\mathrm{d}x\\ &\lesssim \langle t\rangle^{-2-\frac{\alpha}{2}}.\end{aligned}$$

□

**Lemma A.3.** *For any $\varkappa\in[0,2]$ we have*

$$\sup_{x\in\mathbb{R}^2}\left|\frac{\dot{\mathcal{G}}_t(0)-\dot{\mathcal{G}}_t(x)}{|x|^{\varkappa}}\right|\lesssim\langle t\rangle^{\frac{\varkappa}{2}-1}$$

**Proof.** Using $1-e^{-u}\leq 1\wedge u$ for $u\geq 0$, we have

$$\frac{\dot{\mathcal{G}}_t(0)-\dot{\mathcal{G}}_t(x)}{|x|^{\varkappa}}\asymp t^{-1}e^{-\frac{m^2}{t}}\left(\frac{1-e^{-\frac{t}{4}|x|^2}}{|x|^{\varkappa}}\right)\leq t^{-1}e^{-\frac{m^2}{t}}\Big(|x|^{-\varkappa}\wedge\frac{t}{4}|x|^{2-\varkappa}\Big)$$

and the statement follows from splitting the cases $|x|\leq 2/\sqrt{t}$ and $|x|\geq 2/\sqrt{t}$. □

**Lemma A.4.** *Fix an arbitrary $\alpha \in (0,1]$. For any $x \in \mathbb{R}^2$ and any $t \in \mathbb{R}_+$ we have*

$$|\dot{\mathcal{G}}_t(x) - \dot{\mathcal{G}}_t(0)| \lesssim |x|^\alpha \langle t \rangle^{\frac{\alpha}{2}-1}. \tag{A.1}$$

**Proof.** Let $f_t^\alpha: \mathbb{R}_+ \to \mathbb{R}$ be given by $f_t^\alpha(y) := t^{-1}\exp\Big(-\frac{m^2}{t} - \frac{t}{4}y^{2/\alpha}\Big)$. Then

$$|\nabla_y f_t^\alpha| \asymp y^{\frac{2}{\alpha}-1} e^{-\frac{m^2}{t}} e^{-\frac{t}{4}y^{2/\alpha}}$$

attains its maximum at $y^* = Ct^{-\alpha/2}$, so for all $y \in \mathbb{R}_+$

$$|f_t^\alpha(y) - f_t^\alpha(0)| \lesssim y \sup_{y \in \mathbb{R}_+} |\nabla_y f_t^\alpha| \lesssim y e^{-\frac{m^2}{t}} t^{\frac{\alpha}{2}-1}.$$

We then conclude from $\dot{\mathcal{G}}_t$ being radial and noting that with $r = |x|$ we have $\dot{\mathcal{G}}_t(x) = f_t^\alpha(r^\alpha)$. □

**Lemma A.5.** *For any multiindex $\mu \in \mathbb{N}_0^2$ we have*

$$\|\partial^\mu k_s^Q\|_{L^1} \lesssim \langle s \rangle^{\frac{|\mu|}{2}-1} \tag{A.2}$$

**Proof.** Denote $g_s := e^{-\frac{s}{2}|\cdot|^2}$. Note that $g_s = g_1 \circ f$ where $f(x) := \sqrt{s}x$, thus by chain rule

$$\partial^\mu g_s(x) = s^{\frac{|\mu|}{2}} (\partial^\mu g_1)(\sqrt{s}x). \tag{A.3}$$

Hence

$$\|\partial^\mu g_s\|_{L^1} = \int_{\mathbb{R}^d} |\partial^\mu g_s(x)| \mathrm{d}x = s^{\frac{|\mu|}{2}} \int_{\mathbb{R}^d} |(\partial^\mu g_1)(\sqrt{s}x)| \mathrm{d}x = s^{\frac{|\mu|-d}{2}} \int_{\mathbb{R}^d} |(\partial^\mu g_1)(\tilde{x})| \mathrm{d}\tilde{x} \lesssim s^{\frac{|\mu|-d}{2}}.$$ □

**Lemma A.6.** *Let $\mu \in \mathbb{N}_0^2$ multiindex. We have for any $s \in \mathbb{R}_+$, $c \geq 0$ that*

$$\begin{aligned} |\partial^\mu \dot{\mathcal{G}}_s(x_1 - y) - \partial^\mu \dot{\mathcal{G}}_s(x_2 - y)| &\lesssim |x_1 - x_2| e^{-c\sqrt{s}|x_1 - y|} e^{c\sqrt{s}|x_1 - x_2|} e^{-\frac{m^2}{s}} s^{\frac{|\mu|-1}{2}} \\ |\partial^\mu k_s^Q(x_1 - y) - \partial^\mu k_s^Q(x_2 - y)| &\lesssim |x_1 - x_2| e^{-c\sqrt{s}|x_1 - y|} e^{c\sqrt{s}|x_1 - x_2|} e^{-\frac{m^2}{2s}} s^{\frac{|\mu|+1}{2}} \end{aligned}$$

**Proof.** Denote by $z_i := x_i - y$. Taylor expansion yields with $z(\vartheta) := z_1 - \vartheta(z_1 - z_2)$, $\vartheta \in [0,1]$

$$|\partial^\mu \dot{\mathcal{G}}_s(z_1) - \partial^\mu \dot{\mathcal{G}}_s(z_2)| \leq |z_1 - z_2| \int_0^1 |\nabla \partial^\mu \dot{\mathcal{G}}_s(z(\vartheta))| \mathrm{d}\vartheta.$$

We have $|z_1| = |z_1 - \vartheta(z_1 - z_2) + \vartheta(z_1 - z_2)| \leq |z_1 - \vartheta(z_1 - z_2)| + |z_1 - z_2|$, so as $\vartheta \in [0,1]$

$$-|z_1 - \vartheta(z_1 - z_2)| \leq |z_1 - z_2| - |z_1|. \tag{A.4}$$

On the other hand, since for all $b > 0$, $r \in \mathbb{R}$ it holds that $\exp(-br^2) \leq \exp(c^2/4b)\exp(-cr)$, choosing $r := \sqrt{s}|\xi|$, $b := 1/8$ we have for all $\xi$ that

$$|\nabla \partial^\mu \dot{\mathcal{G}}_s(\xi)| \lesssim_\mu s^{\frac{|\mu|+1}{2}} e^{-\frac{s}{8}|\xi|^2} e^{-\frac{m^2}{s}} s^{-1} \lesssim_c s^{\frac{|\mu|-1}{2}} e^{-c\sqrt{s}|\xi|} e^{-\frac{m^2}{s}}$$

Combining this with (A.4), taking $\xi := z_1 - \vartheta(z_1 - z_2)$ we obtain

$$|\nabla \partial^\mu \dot{\mathcal{G}}_s(z(\vartheta))| \lesssim s^{\frac{|\mu|-1}{2}} e^{c\sqrt{s}|z_1 - z_2|} e^{-c\sqrt{s}|z_1|} e^{-\frac{m^2}{s}}$$

proving the estimate for $\dot{\mathcal{G}}_s$. For $k_s^Q$, we proceed the same way. □

## A.2 Proof of Lemma 3.5

**Proof. (of Lemma 3.5)** Recall that we want to prove the following: given $t \le s$, there exists a constant $C$ such that for all $x_1, x_2, z \in \mathbb{R}^2$ with $|x_g - x_b| \le |x_b - z|$ we have

$$\mathcal{G}_{t,s}(x_b - z) - \mathcal{G}_{t,s}(x_g - z) \le C.$$

If $|x_g - z| \le |x_b - z|$, then

$$\mathcal{G}_{t,s}(x_b - z) - \mathcal{G}_{t,s}(x_g - z) \le 0.$$

Thus from now on assume $|x_b - z| \le |x_g - z|$. Denoting by $r := |x_g - x_b|$, if $r = 0$, then

$$\mathcal{G}_{t,s}(x_b - z) - \mathcal{G}_{t,s}(x_g - z) = 0$$

so that we may assume $r > 0$. Then we split the scales $[t, s]$ according to $r$:

$$\begin{aligned}
& \mathcal{G}_{t,s}(x_b - z) - \mathcal{G}_{t,s}(x_g - z) \\
= & \int_t^s \mathbb{1}_{\{u \le r^{-2}\}} (\dot{\mathcal{G}}_u(x_b - z) - \dot{\mathcal{G}}_u(x_g - z)) \mathrm{d}u + \int_t^s \mathbb{1}_{\{u \ge r^{-2}\}} (\dot{\mathcal{G}}_u(x_b - z) - \dot{\mathcal{G}}_u(x_g - z)) \mathrm{d}u \\
= & \mathcal{I}_{\text{low}} + \mathcal{I}_{\text{high}}.
\end{aligned}$$

In the low region, the Gaussian is flat and the mean value theorem is a good approximation: using (A.1)

$$\begin{aligned}
|\mathcal{I}_{\text{low}}| & \lesssim r \int_t^s \mathbb{1}_{\{u \le r^{-2}\}} \langle u \rangle^{-\frac{1}{2}} \mathrm{d}u \\
& \lesssim r (s \wedge r^{-2})^{\frac{1}{2}} \\
& \lesssim 1.
\end{aligned}$$

In the high region we exploit the exponential decay of the Gaussian: since $r \le |x_g - z|, |x_b - z|$

$$\begin{aligned}
|\mathcal{I}_{\text{high}}| & \le 2 \int_t^s \mathbb{1}_{\{u \ge r^{-2}\}} |\dot{\mathcal{G}}_u(r)| \mathrm{d}u \\
& \lesssim \int_{r^{-2}}^{\infty} u^{-1} e^{-\frac{m^2}{u}} e^{-u\frac{r^2}{4}} \mathrm{d}u \\
& \lesssim 1 + \int_{r^{-2} \vee 1}^{\infty} u^{-1} e^{-u\frac{r^2}{4}} \mathrm{d}u \\
& \lesssim 1 + \int_{1/4}^{\infty} w^{-1} e^{-w} \mathrm{d}w \\
& = 1 + E_1(1/4)
\end{aligned}$$

yielding the statement. □

## A.3 Scale decomposition of the free field

In this section we prove some standard technical results about the heat kernel decomposition $B_t := \int_0^t Q_s \mathrm{d}W_s$ of the Gaussian free field.

For the Wiener integral $\int_0^t Q_s \mathrm{d}W_s$ to be well-defined, by standard infinite-dimensional stochastic integration theory [20] we have to show the following three properties:

1. $Q_t \in \mathcal{L}_2(L^2; \mathcal{H})$ for every $t$, i.e. $Q_t$ is Hilbert–Schmidt when viewed as an operator from $L^2$ to $\mathcal{H}$.

2. The process $(Q_t)_{t \in \mathbb{R}_+}$ is predictable with respect to the canonical filtration of $W$.

3. For every $T\in\mathbb{R}_+$, it has to hold that

$$\mathbb{P}\left[\int_0^T \|Q_t\|^2_{\mathcal{L}_2(L^2;\mathcal{H})}\mathrm{d}t<\infty\right]=1.$$

We take $\mathcal{H}$ to be the weighted Sobolev space $\mathcal{H}:=B^k_{2,2}(\langle x\rangle^{-\ell})$ with $k,\ell>0$.

**Lemma A.7.** *Suppose $\ell>d/2$. Then $Q_t\in\mathcal{L}_2(L^2;\mathcal{H})$.*

**Proof.** By the Fourier multiplier property of the weighted Sobolev spaces (cf. subsequent remark), for an ONB $(e_n)_{n\in\mathbb{N}}$ of $L^2$,

$$\|Q_t\|^2_{\mathcal{L}_2(L^2;\mathcal{H})}=\sum_{n\in\mathbb{N}}\|Q_t e_n\|^2_{\mathcal{H}}\asymp\sum_{n\in\mathbb{N}}\|\langle\cdot\rangle^{-\ell}(-\Delta+1)^{k/2}Q_t e_n\|^2_{L^2}=\|\langle\cdot\rangle^{-\ell}(-\Delta+1)^{k/2}Q_t\|^2_{\mathcal{L}_2(L^2;L^2)}$$

That is, $Q_t$ is Hilbert–Schmidt as an operator viewed as $L^2\to\mathcal{H}$ if and only if $\langle\cdot\rangle^{-\ell}(-\Delta+1)^{k/2}Q_t$, viewed as $L^2\to L^2$, is Hilbert–Schmidt. Now consider the kernel

$$K_t(x,y):=\frac{\langle x\rangle^{-\ell}}{(2\pi)^d}\int_{\mathbb{R}^d}\underbrace{(|\xi|^2+1)^{k/2}\left(\frac{1}{t^2}e^{\frac{-|\xi|^2-m^2}{t}}\right)^{1/2}}_{:=H_t(\xi)}e^{i\langle y-x,\xi\rangle}\mathrm{d}\xi=\langle x\rangle^{-\ell}\mathcal{F}^{-1}(H_t)(x-y)$$

Since $H_t\in\mathcal{S}(\mathbb{R}^d)$, we have that $\mathcal{F}^{-1}(H_t)\in\mathcal{S}(\mathbb{R}^d)$, and a standard computation involving polar coordinates yields $\|\langle x\rangle^{-\ell}\|^2_{L^2(\mathbb{R}^d)}<\infty$ when $\ell>d/2$, so

$$\|K_t\|^2_{L^2(\mathbb{R}^d)}=\|\langle\cdot\rangle^{-\ell}\|^2_{L^2(\mathbb{R}^d)}\|\mathcal{F}^{-1}(H_t)\|^2_{L^2(\mathbb{R}^d)}<\infty.$$

Hence by Theorem 2.11 of [45], it defines a Hilbert–Schmidt integral operator $A_t$: $L^2\to L^2$ with

$$\|A_t\|_{\mathcal{L}_2(L^2;L^2)}=\|K_t\|_{L^2(\mathbb{R}^d\times\mathbb{R}^d)}$$

On the other hand, for any $f\in\mathcal{S}$,

$$\begin{aligned}\langle x\rangle^{-\ell}(Q_t(-\Delta+1)^{k/2}f)(x) &= \langle x\rangle^{-\ell}\mathcal{F}^{-1}\left((|\xi|^2+1)^{\frac{k}{2}}\left(\frac{1}{t^2}e^{\frac{-|\xi|^2-m^2}{t}}\right)^{\frac{1}{2}}\mathcal{F}(f)\right)(x)\\ &= \frac{\langle x\rangle^{-\ell}}{(2\pi)^d}\int_{\mathbb{R}^d}(|\xi|^2+1)^{\frac{k}{2}}\left(\frac{1}{t^2}e^{\frac{-|\xi|^2-m^2}{t}}\right)^{\frac{1}{2}}\left(\int_{\mathbb{R}^d}e^{-i\langle y,\xi\rangle}f(y)\mathrm{d}y\right)e^{i\langle x,\xi\rangle}\mathrm{d}\xi\\ &= \frac{\langle x\rangle^{-\ell}}{(2\pi)^d}\int_{\mathbb{R}^d}\int_{\mathbb{R}^d}(|\xi|^2+1)^{\frac{k}{2}}\left(\frac{1}{t^2}e^{\frac{-|\xi|^2-m^2}{t}}\right)^{\frac{1}{2}}e^{i\langle x-y,\xi\rangle}f(y)\mathrm{d}y\mathrm{d}\xi\\ &= \int_{\mathbb{R}^d}\left(\frac{\langle x\rangle^{-\ell}}{(2\pi)^d}\int_{\mathbb{R}^d}(|\xi|^2+1)^{\frac{k}{2}}\left(\frac{1}{t^2}e^{\frac{-|\xi|^2-m^2}{t}}\right)^{\frac{1}{2}}e^{i\langle x-y,\xi\rangle}\mathrm{d}\xi\right)f(y)\mathrm{d}y\\ &= \int_{\mathbb{R}^d}K_t(x,y)f(y)\mathrm{d}y\end{aligned}$$

so $A\equiv\langle x\rangle^{-\ell}Q_t(-\Delta+1)^{k/2}$. □

**Remark A.8.** Weighted Besov spaces $B^{s,w}_{p,q}$ are defined by multiplying each Littlewood–Paley block with the weight function, "inside" the $L^p$ and $\ell^q$ norms. However, it can be shown (section 4.2.2 of [25]) that the norm $\|wf\|_{B^s_{p,q}}$ is an equivalent norm for $B^{s,w}_{p,q}$. Moreover, it is well-known that unweighted Besov spaces for $p=q=2$ coincide with Sobolev spaces. Hence, in particular, we have

$$\|f\|_{B^{s,w}_{2,2}}\asymp\|wf\|_{B^s_{2,2}}\asymp\|(-\Delta+1)^{s/2}(wf)\|_{L^2}.$$

On the other hand, by Proposition on p.158 of [25], we have, just as in the unweighted case, the Fourier multiplier property

$$\|f\|_{B_{p,q}^{s,w}} \asymp \|(-\Delta+1)^{k/2}(f)\|_{B_{p,q}^{s-k,w}}.$$

Hence

$$\begin{aligned}\|(-\Delta+1)^{s/2}(wf)\|_{L^2} &\asymp \|f\|_{B_{2,2}^{s,w}}\\ &\asymp \|(-\Delta+1)^{s/2}(f)\|_{B_{2,2}^{0,w}}\\ &\asymp \|w(-\Delta+1)^{s/2}(f)\|_{B_{2,2}^{0}}\\ &= \|w(-\Delta+1)^{s/2}(f)\|_{L^2}\end{aligned}$$

**Lemma A.9.** *Let $\ell > d/2$ and $\mathcal{H} := B_{2,2}^{k,\ell} = B_{2,2}^{k}(\langle x\rangle^{-\ell})$. For all $T \in \mathbb{R}_+$, if $k + d/2 > 1$, then*

$$\mathbb{E}[\|B_T\|_{\mathcal{H}}^2] = \int_0^T \|Q_t\|_{\mathcal{L}_2(L^2;\mathcal{H})}^2 \mathrm{d}t \lesssim_{\ell,k,d} 1 + T^{k+d/2-1}; \tag{A.5}$$

*if $d=2$, $k=0$, then*

$$\mathbb{E}[\|B_T\|_{\mathcal{H}}^2] = \int_0^T \|Q_t\|_{\mathcal{L}_2(L^2;\mathcal{H})}^2 \mathrm{d}t \lesssim_{\ell,k,d} 1 + \log(1 \vee T). \tag{A.6}$$

**Proof.** The fact that

$$\mathbb{E}\left[\left\|\int_0^T Q_s \mathrm{d}W_s\right\|_{\mathcal{H}}^2\right] = \int_0^T \|Q_t\|_{\mathcal{L}_2(L^2;\mathcal{H})}^2 \mathrm{d}t$$

is the infinite dimensional version of Ito isometry (see e.g. (4.30) in [20]). To prove the estimates, we use the previous proof, a change of variables $(x,y) \mapsto (x', x'-y')$ and Plancherel to obtain

$$\begin{aligned}\int_0^T \|Q_t\|_{\mathcal{L}_2(L^2;\mathcal{H})}^2 \mathrm{d}t &= \int_0^T \int_{\mathbb{R}^d}\int_{\mathbb{R}^d} (\langle x\rangle^{-\ell}\mathcal{F}^{-1}(H_t)(x-y))^2 \mathrm{d}x\mathrm{d}y\mathrm{d}t\\ &= \|\langle x\rangle^{-\ell}\|_{L^2(\mathbb{R}^d)}^2 \int_0^T \|\mathcal{F}^{-1}(H_t)\|_{L^2(\mathbb{R}^d)}^2 \mathrm{d}t\\ &\lesssim_\ell \int_0^T \frac{1}{t^2} e^{\frac{-m^2}{t}} \Big(\int_{\mathbb{R}^d} (|\xi|^2+1)^k e^{\frac{-|\xi|^2}{t}} \mathrm{d}\xi\Big)\mathrm{d}t.\end{aligned}$$

Now we focus on the $\mathrm{d}\xi$ part of the above integral. Considering the change of variables $\xi \mapsto \sqrt{t}\xi$,

$$\int_{\mathbb{R}^d} (|\xi|^2+1)^k e^{\frac{-|\xi|^2}{t}} \mathrm{d}\xi = t^{d/2} \underbrace{\int_{\mathbb{R}^d} (t|\xi|^2+1)^k e^{-|\xi|^2} \mathrm{d}\xi}_{:=\mathcal{I}}$$

Split

$$\mathcal{I} = \int_{t|\xi|^2 \le 1} (t|\xi|^2+1)^k e^{-|\xi|^2} \mathrm{d}\xi + \int_{t|\xi|^2 \ge 1} (t|\xi|^2+1)^k e^{-|\xi|^2} \mathrm{d}\xi =: \mathbb{I} + \mathbb{II}.$$

Now since $e^{-|\xi|^2} \le 1$, we have

$$\mathbb{I} \le \int_{|\xi|^2 \le \frac{1}{t}} 2^k e^{-|\xi|^2} \mathrm{d}\xi \le 2^k \mathrm{Vol}_d\Big(t^{-\frac{1}{2}}\Big) \lesssim_{k,d} t^{-\frac{d}{2}}$$

thus

$$\int_0^T \frac{1}{t^2} e^{\frac{-m^2}{t}} t^{d/2} \mathbb{I} \mathrm{d}t \lesssim \int_0^T t^{-2} e^{\frac{-m^2}{t}} \mathrm{d}t \lesssim 1.$$

On the other hand,

$$\mathbb{II} \le \int_{t|\xi|^2 \ge 1} 2^k t^k (|\xi|^2)^k e^{-|\xi|^2} \mathrm{d}\xi \le 2^k t^k \int_{\mathbb{R}^d} |\xi|^{2k} e^{-|\xi|^2} \mathrm{d}\xi \lesssim_k t^k$$

since Gaussian random variables are in every $L^p$. Thus when $k+d/2-2>-1$, we have

$$\int_0^T \frac{1}{t^2} e^{\frac{-m^2}{t}} t^{d/2} \mathbb{III} \mathrm{d}t \lesssim \int_0^T t^{k+d/2-2} e^{\frac{-m^2}{t}} \mathrm{d}t \lesssim_{k,d} T^{k+d/2-1}. \tag{A.7}$$

Putting everything together,

$$\int_0^T \|Q_t\|^2_{\mathcal{L}_2(L^2;\mathcal{H})} \mathrm{d}t \lesssim_\ell \int_0^T \|H_t\|^2_{L^2(\mathbb{R}^d)} \mathrm{d}t = \int_0^T \frac{1}{t^2} e^{\frac{-m^2}{t}} t^{d/2} (\mathbb{I} + \mathbb{III}) \mathrm{d}t \lesssim 1 + T^{k+d/2-1}$$

proving (A.5). The bound (A.6) follows similarly with the term $t^{k+d/2-2}$ in (A.7) being replaced by $t^{-1}$. □

Thus, as $(Q_t)_{t\in\mathbb{R}_+}$ is deterministic hence predictable, we obtained that the stochastic integral is well-defined, and also get by standard results the following

**Corollary A.10.** *For any $t<\infty$, $B_t$ is a Gaussian random variable with covariance operator $\int_0^t Q_s^2 \mathrm{d}s$, and $B_t \in \mathcal{H}$ almost surely.*

# B Calculus on the space of germs

## B.1 Computing the trace

Consider $F, G$ Frechet differentiable functionals on $\mathcal{H}$. Recall that by Riesz–Frechet representation, for any $F\in\mathcal{C}^2(\mathcal{H};\mathbb{R})$ we have the identification of $\mathrm{D}^2F|_\phi \in L^{(2)}(\mathcal{H}\times\mathcal{H};\mathbb{R})$ with $\mathrm{D}^2F|_\phi^\sharp \in L(\mathcal{H})$. Similarly, given $F, G\in\mathcal{C}^1(\mathcal{H};\mathbb{R})$, the continuous bilinear form

$$(v,w)\mapsto \mathrm{D}F|_\phi(v)\mathrm{D}G|_\phi(w)$$

is identified with a bounded linear operator $\mathrm{D}F\mathrm{D}G|_\phi^\sharp$ on $\mathcal{H}$:

$$\mathrm{D}F|_\phi(v)\mathrm{D}G|_\phi(w) = \big\langle \big(\mathrm{D}F\mathrm{D}G|_\phi^\sharp\big)(v), w\big\rangle_{\mathcal{H}}.$$

**Lemma B.1.** *For any $F\in\mathcal{C}^2(\mathcal{H};\mathbb{R})$ we have*

$$\mathrm{Tr}_{\mathcal{H}}\big[\mathrm{D}^2F|_\phi^\sharp Q_t Q_t^*\big] = \int_{\mathbb{R}^d} \mathrm{D}^2F|_\phi(k_t^Q(\cdot-z), k_t^Q(\cdot-z))\mathrm{d}z$$

*and for $F, G\in\mathcal{C}^1(\mathcal{H};\mathbb{R})$,*

$$\mathrm{Tr}_{\mathcal{H}}\big[\mathrm{D}F\mathrm{D}G|_\phi^\sharp Q_t Q_t^*\big] = \int_{\mathbb{R}^d} \mathrm{D}F|_\phi(k_t^Q(\cdot-z))\mathrm{D}G|_\phi(k_t^Q(\cdot-z))\mathrm{d}z.$$

**Proof.** The kernel representation of $Q_t$ ensures that for every $f\in L^2$, $h\in\mathcal{H}$, by continuity of the Bochner integral and $\|k_t^Q(\cdot-z)\|_{\mathcal{H}}\in L^2(\mathrm{d}z)$,

$$\langle Q_t^*h, f\rangle_{L^2} = \langle h, Q_t f\rangle_{\mathcal{H}} = \left\langle h, \int f(z) k_t^Q(\cdot-z)\mathrm{d}z\right\rangle_{\mathcal{H}} = \int f(z)\langle h, k_t^Q(\cdot-z)\rangle_{\mathcal{H}}\mathrm{d}z.$$

Thus $Q_t^*h(z) = \langle h, k_t^Q(\cdot-z)\rangle_{\mathcal{H}}$ for almost every $z$, thus applying the kernel representation again,

$$Q_t Q_t^* h = \int k_t^Q(\cdot-z) Q_t^* h(z)\mathrm{d}z = \int k_t^Q(\cdot-z)\langle h, k_t^Q(\cdot-z)\rangle_{\mathcal{H}}\mathrm{d}z.$$

Consider a bounded bilinear form $b\in L^{(2)}(\mathcal{H}\times\mathcal{H};\mathbb{R})$ with Riesz–Frechet representation $A\in L(\mathcal{H})$, i.e. $b(v,w)=\langle A(v),w\rangle_{\mathcal{H}}$. Then for $(e_n)_{n\in\mathbb{N}}$ ONB of $\mathcal{H}$, $u=\sum_n \langle u,e_n\rangle e_n$

$$\begin{aligned}
\mathrm{Tr}_{\mathcal{H}}[AQ_tQ_t^*] &= \sum_{n\in\mathbb{N}} \langle AQ_tQ_t^* e_n, e_n\rangle_{\mathcal{H}} \\
&= \sum_{n\in\mathbb{N}} \left\langle A\int k_t^Q(\cdot-z)\langle e_n, k_t^Q(\cdot-z)\rangle_{\mathcal{H}}\mathrm{d}z, e_n\right\rangle_{\mathcal{H}} \\
&= \int \sum_{n\in\mathbb{N}} \langle Ak_t^Q(\cdot-z), e_n\rangle_{\mathcal{H}}\langle e_n, k_t^Q(\cdot-z)\rangle_{\mathcal{H}}\mathrm{d}z \\
&= \int \langle Ak_t^Q(\cdot-z), k_t^Q(\cdot-z)\rangle_{\mathcal{H}}\mathrm{d}z \\
&= \int b(k_t^Q(\cdot-z), k_t^Q(\cdot-z))\mathrm{d}z.
\end{aligned}$$

Here the interchange follows via Fubini from Cauchy–Schwarz, Parseval, and the aforementioned fact that $\|k_t^Q(\cdot-z)\|_{\mathcal{H}}\in L^2(\mathrm{d}z)$:

$$\begin{aligned}
\sum_{n\in\mathbb{N}} |\langle Ak_t^Q(\cdot-z), e_n\rangle_{\mathcal{H}}\langle e_n, k_t^Q(\cdot-z)\rangle_{\mathcal{H}}| &\leq \sqrt{\sum_{n\in\mathbb{N}} |\langle Ak_t^Q(\cdot-z), e_n\rangle_{\mathcal{H}}|^2}\sqrt{\sum_{n\in\mathbb{N}} |\langle e_n, k_t^Q(\cdot-z)\rangle_{\mathcal{H}}|^2} \\
&\leq \|A\|_{L(\mathcal{H})}\|k_t^Q(\cdot-z)\|_{\mathcal{H}}^2.
\end{aligned}$$

Applying the above computation with $A:=\mathrm{D}^2F\big|_{\phi}^{\sharp}$ and $A:=\mathrm{D}F\mathrm{D}G\big|_{\phi}^{\sharp}$ yields the results. □

**Notation B.2.** *It follows from* (2.7) *that*

$$\int_{\mathbb{R}^d} k_t^Q(x-\xi)k_t^Q(y-\xi)\mathrm{d}\xi = k_t^{Q^2}(x-y) = \dot{\mathcal{G}}_t(x-y),$$

*thus we write $\dot{G}_t$ instead of $Q_tQ_t^*$ to mean the contractions of Lemma B.1.*

Indeed, to highlight the effectiveness of the above abuse of notation, consider $F$ and $G$ to be point evaluations at $x_1,x_2\in\mathbb{R}^d$. Then $\mathrm{D}F\mathrm{D}G|_{\phi}=\delta_{x_1}\otimes\delta_{x_2}$ and thus

$$\mathrm{Tr}_{\mathcal{H}}\big[\mathrm{D}F\mathrm{D}G\big|_{\phi}^{\sharp}Q_tQ_t^*\big] = \int_{\mathbb{R}^d} \delta_{x_1}(k_t^Q(\cdot-z))\delta_{x_2}(k_t^Q(\cdot-z))\mathrm{d}z = \dot{\mathcal{G}}_t(x_1-x_2).$$

More generally, for functionals $F:\mathcal{H}\to\mathbb{R}$ that depend on a finite jet of $\phi$ at a point, we have the following computation.

**Example B.3.** Fix $m\in\mathbb{N}$ and consider $F:\mathcal{H}\to\mathbb{R}$ given by

$$F(\phi) := f\left(\bigotimes_{|\alpha|\leq m} \partial^{\alpha}\phi(x)\right) = f\circ\left(\bigotimes_{|\alpha|\leq m} (\delta_x\circ\partial^{\alpha})\phi\right)$$

for an $f\in\mathcal{C}^2(\mathbb{R}^N;\mathbb{R})$ where $N=N(d,m)\in\mathbb{N}$ is the number of $d$-dimensional multiindices of order at most $m$. Note that this includes, as a particular case, every constant coefficient differential operator $P=\sum_{|\alpha|\leq m}\gamma_{\alpha}\partial^{\alpha}$ with the induced function $f_P(a):=\sum_{i=1}^N \gamma_i a_i$ for all $a\in\mathbb{R}^N$. Here $\mathcal{H}:=H_{\ell}^{k+m}$ where $k=k(d)$ is the smallest integer such that $H_{\ell}^k\hookrightarrow\mathcal{C}$ continuously (Morrey's inequality). Then chain rule and linearity yields that $F\in\mathcal{C}^1$, and $\mathrm{D}F|_{\phi}\in L(\mathcal{H};\mathbb{R})$ is given by

$$\begin{aligned}
\mathrm{D}F|_{\phi} &= \mathrm{D}f|_{(\bigotimes_{|\alpha|\leq m}(\delta_x\circ\partial^{\alpha})\phi)}\circ\left(\bigotimes_{|\alpha|\leq m} \mathrm{D}\delta_x|_{\partial^{\alpha}\phi}\circ\mathrm{D}\partial^{\alpha}|_{\phi}\right) \\
&= \left\langle \nabla f\left(\bigotimes_{|\alpha|\leq m} (\partial^{\alpha}\phi(x))\right), \cdot\right\rangle_{\mathbb{R}^N}\circ\left(\bigotimes_{|\alpha|\leq m} \delta_x\circ\partial^{\alpha}\right),
\end{aligned}$$

that is, denoting by $\Phi_F := \bigotimes_{|\alpha|\leq m}(\partial^\alpha\phi(x)) \in \mathbb{R}^N$, we have for every $v \in \mathcal{H}$ that

$$\mathrm{D}F|_\phi v = \left\langle \nabla f(\Phi_F), \bigotimes_{|\alpha|\leq m} \partial^\alpha v(x) \right\rangle_{\mathbb{R}^N}.$$

Thus

$$\begin{aligned}\mathrm{Tr}_{\mathcal{H}}\big[\mathrm{D}F\mathrm{D}G|_\phi^\sharp Q_t Q_t^*\big] &= \langle \mathrm{D}G|_\phi \otimes \mathrm{D}F|_\phi, \dot{G}_t\rangle \\ &= \left\langle \sum_{|\alpha_f|\leq m_f, |\alpha_g|\leq m_g} c_{\alpha_f,\alpha_g} \partial^{\alpha_f}\delta_x \otimes \partial^{\alpha_g}\delta_y, \dot{G}_t \right\rangle\end{aligned}$$

with

$$c_{\alpha_f,\alpha_g} = (-1)^{|\alpha_f|+|\alpha_g|}[\nabla f(\Phi_F)]_{\alpha_f} \cdot [\nabla g(\Phi_G)]_{\alpha_g}$$

where $[\nabla f(\Phi_F)]_{\alpha_f}$ denotes the partial derivative of $f$ corresponding to the multiindex $\alpha_f$ evaluated at $\Phi_F$.

## B.2 Flow equation computations

**Proof. (of Lemma 3.1)** For simplicity we will sometimes suppress the UV cutoff $T$ in the notation. Of course

$$\partial_s V_s^{[\ell]} = \partial_s \int v_s^{[\ell]}(\xi_{1:\ell})\psi_{1:\ell}\rho(x_{1:\ell})\mathrm{d}\xi_{1:\ell} = \int \partial_s v_s^{[\ell]}(\xi_{1:\ell})\psi_{1:\ell}\rho(x_{1:\ell})\mathrm{d}\xi_{1:\ell}$$

and

$$\begin{aligned}\frac{1}{2}\Delta_{\dot{G}_s} V_s^{[\ell]} &= \frac{1}{2}\mathrm{Tr}\Big[\dot{G}_s \mathrm{D}^2\Big(\int v_s^{[\ell]}(\xi_{1:\ell})\psi_{\xi_{1:\ell}}\rho(x_{1:\ell})\mathrm{d}\xi_{1:\ell}\Big)\Big] \\ &= -\frac{\beta^2}{2}\int v_s^{[\ell]}(\xi_{1:\ell}) \sum_{j,k=1}^{\ell} \sigma_j\sigma_k \dot{\mathcal{G}}_s(x_j - x_k)\psi_{1:\ell}\rho(x_{1:\ell})\mathrm{d}\xi_{1:\ell} \\ &= \int v_s^{[\ell]}(\xi_{1:\ell})\dot{\mathbb{G}}_s(\xi_{1:\ell})\psi_{1:\ell}\rho(x_{1:\ell})\mathrm{d}\xi_{1:\ell}.\end{aligned}$$

Noting that

$$\begin{aligned}\langle \mathrm{D}\psi_{1:\ell_1}, \dot{G}_s \mathrm{D}\psi_{1:\ell_2}\rangle &= \left\langle \sum_{j=1}^{\ell_1} i\beta\sigma_j\delta_{x_j}\psi_{1:\ell_1}, \dot{G}_s \sum_{k=1}^{\ell_2} i\beta\sigma_k\delta_{x_k}\psi_{1:\ell_2} \right\rangle \\ &= -\beta^2 \sum_{j=1}^{\ell_1}\sum_{k=1}^{\ell_2} \sigma_j\sigma_k\dot{\mathcal{G}}_s(x_j - x_k)\psi_{1:\ell} \\ &= 2\dot{\mathbb{G}}_s(\xi_{1:\ell_1}, \xi_{1:\ell_2})\psi_{1:\ell}\end{aligned}$$

we have, via permutation invariance of the kernels, i.e. e.g. $v_s^{[\ell_1]}(\xi_{1:\ell_1}) = v_s^{[\ell_1]}(\xi_{\ell_2:\ell})$ whenever $\ell_1 + \ell_2 = \ell$, that

$$\begin{aligned}&\frac{1}{2}\langle \mathrm{D}V_s^{[\ell_1]}, \dot{G}_s \mathrm{D}V_s^{[\ell_2]}\rangle \\ &= \frac{1}{2}\Big\langle \mathrm{D}\int v_s^{[\ell_1]}(\xi_{1:\ell_1})\psi_{1:\ell_1}\rho(x_{1:\ell_1})\mathrm{d}\xi_{1:\ell_1}, \dot{G}_s \mathrm{D}\int v_s^{[\ell_2]}(\xi_{1:\ell_2})\psi_{1:\ell_2}\rho(x_{1:\ell_2})\mathrm{d}\xi_{1:\ell_2}\Big\rangle \\ &= \iint v_s^{[\ell_1]}(\xi_{1:\ell_1})\psi_{1:\ell_1}\rho(x_{1:\ell_1})\dot{\mathbb{G}}_s(\xi_{1:\ell_1}, \xi_{1:\ell_2}) v_s^{[\ell_2]}(\xi_{1:\ell_2})\psi_{1:\ell_2}\rho(x_{1:\ell_2})\mathrm{d}\xi_{1:\ell'}\mathrm{d}\xi_{1:\ell''} \\ &= \int v_s^{[\ell_1]}(\xi_{1:\ell_1})\dot{\mathbb{G}}_s(\xi_{1:\ell_1}, \xi_{1:\ell_2}) v_s^{[\ell_2]}(\xi_{1:\ell_2})\psi_{1:\ell}\rho(x_{1:\ell})\mathrm{d}\xi_{1:\ell}.\end{aligned}$$

Thus, the equation (2.20) is satisfied if $v_{t,T}^{[\ell]}$ is such that

$$\partial_t v_{t,T}^{[\ell]}(\xi_{1:\ell}) + v_{t,T}^{[\ell]}(\xi_{1:\ell})\dot{\mathbb{G}}_t(\xi_{1:\ell}) = \sum_{I_1 \uplus I_2 = 1:\ell} C_{|I_1|,|I_2|} v_{t,T}^{[|I_1|]}(\xi_{I_1})\dot{\mathbb{G}}_t(\xi_{I_1}, \xi_{I_2}) v_{t,T}^{[|I_2|]}(\xi_{I_2}) \tag{B.1}$$

with $C_{|I_1|,|I_2|}$ combinatorial constants.

It remains to show that if $v^{[\ell]}_{t,T}$ solves (3.3), then it satisfies (B.1). For brevity, we denote the right-hand side of (B.1) by $f_t$, so that (3.3) reads as

$$v^{[\ell]}_{t,T}(\xi_{1:\ell}) = v^{[\ell]}_{T,T}(\xi_{1:\ell}) e^{\mathbb{G}_{t,T}(\xi_{1:\ell})} - \int_t^T f_s e^{\mathbb{G}_{t,s}(\xi_{1:\ell})} \mathrm{d}s.$$

The Leibniz integral rule says

$$\partial_t \int_{a(t)}^{b(t)} h(t,s)\mathrm{d}s = h(t,b(t))\partial_t b(t) - h(t,a(t))\partial_t a(t) + \int_{a(t)}^{b(t)} \partial_t h(t,s)\mathrm{d}s$$

so in particular with $a(t)=t$, $b(t)=T$ we have

$$\begin{aligned}
\partial_t \int_t^T f_s e^{\mathbb{G}_{t,s}(\xi_{1:\ell})}\mathrm{d}s &= -f_t e^{\mathbb{G}_{t,t}(\xi_{1:\ell})} + \int_t^T f_s \partial_t e^{\mathbb{G}_{t,s}(\xi_{1:\ell})}\mathrm{d}s \\
&= -f_t + \int_t^T f_s e^{\mathbb{G}_{t,s}(\xi_{1:\ell})}(-\dot{\mathbb{G}}_t(\xi_{1:\ell}))\mathrm{d}s \\
&= -f_t - \dot{\mathbb{G}}_t(\xi_{1:\ell}) \int_t^T f_s e^{\mathbb{G}_{t,s}(\xi_{1:\ell})}\mathrm{d}s.
\end{aligned}$$

Hence

$$\begin{aligned}
\partial_t v^{[\ell]}_{t,T}(\xi_{1:\ell}) &= \partial_t\big(v^{[\ell]}_{T,T}(\xi_{1:\ell}) e^{\mathbb{G}_{t,T}(\xi_{1:\ell})}\big) - \partial_t \int_t^T f_s e^{\mathbb{G}_{t,s}(\xi_{1:\ell})}\mathrm{d}s \\
&= v^{[\ell]}_{T,T}(\xi_{1:\ell}) e^{\mathbb{G}_{t,T}(\xi_{1:\ell})}(-\dot{\mathbb{G}}_t(\xi_{1:\ell})) + f_t + \dot{\mathbb{G}}_t(\xi_{1:\ell}) \int_t^T f_s e^{\mathbb{G}_{t,s}(\xi_{1:\ell})}\mathrm{d}s \\
&= v^{[\ell]}_{T,T}(\xi_{1:\ell}) e^{\mathbb{G}_{t,T}(\xi_{1:\ell})}(-\dot{\mathbb{G}}_t(\xi_{1:\ell})) + f_t - (-\dot{\mathbb{G}}_t(\xi_{1:\ell}))\left(\int_t^T f_s e^{\mathbb{G}_{t,s}(\xi_{1:\ell})}\mathrm{d}s\right) \\
&= (-\dot{\mathbb{G}}_t(\xi_{1:\ell}))\left(v^{[\ell]}_{T,T}(\xi_{1:\ell}) e^{\mathbb{G}_{t,T}(\xi_{1:\ell})} - \int_t^T f_s e^{\mathbb{G}_{t,s}(\xi_{1:\ell})}\mathrm{d}s\right) + f_t \\
&= (-\dot{\mathbb{G}}_t(\xi_{1:\ell})) v^{[\ell]}_{t,T}(\xi_{1:\ell}) + f_t.
\end{aligned}$$

as requested. □

**Proof. (of Lemma 5.1)** From (5.2) compute

$$\begin{aligned}
\mathrm{D}\Gamma^{[\ell]}_{s,T} &= \sum_{k=0}^n k \gamma^{[\ell](k)}_{s,T}(\xi_{1:(\ell-1)}, y)\phi(y)^{k-1}\delta_y \\
\mathrm{D}^2\Gamma^{[\ell]}_{s,T} &= \sum_{k=0}^n k(k-1)\gamma^{[\ell](k)}_{s,T}(\xi_{1:(\ell-1)}, y)\phi(y)^{k-2}(\delta_y \otimes \delta_y) \\
\mathrm{D}\psi_{1:(\ell-1)} &= i\beta \sum_{j=1}^{\ell-1} \sigma_j \delta_{x_j} \psi_{1:(\ell-1)} \\
\mathrm{D}^2\psi_{1:(\ell-1)} &= -\beta^2 \sum_{j=1}^{\ell-1}\sum_{k=1}^{\ell-1} \sigma_j\sigma_k(\delta_{x_j}\otimes\delta_{x_k})\psi_{1:(\ell-1)}
\end{aligned}$$

Then from (5.1), as

$$\mathrm{D}^2 q^{[\ell],T}_s(\phi) = \int \Gamma^{[\ell]}_{s,T}\mathrm{D}^2\psi^\rho_{1:(\ell-1)} + \psi^\rho_{1:(\ell-1)}\mathrm{D}^2\Gamma^{[\ell]}_{s,T} + 2\mathrm{D}\Gamma^{[\ell]}_{s,T}\mathrm{D}\psi^\rho_{1:(\ell-1)}\mathrm{d}\xi_{1:(\ell-1)}$$

we obtain

$$\begin{aligned}
\frac{1}{2}\Delta_{\dot{G}_s} q_s^{[\ell],T}(\phi) &= \int \psi_{1:(\ell-1)} \sum_{k=0}^{n} \gamma_{s,T}^{[\ell](k)}(\xi_{1:(\ell-1)}, y)\phi(y)^k \dot{\mathbb{G}}_s(\xi_{1:(\ell-1)}) \mathrm{d}\xi_{1:(\ell-1)} \\
&+ \int \psi_{1:(\ell-1)} \sum_{k=0}^{n} \frac{k(k-1)}{2} \gamma_{s,T}^{[\ell](k)}(\xi_{1:(\ell-1)}, y)\phi(y)^{k-2} \dot{\mathcal{G}}_s(0) \mathrm{d}\xi_{1:(\ell-1)} \\
&+ \int \psi_{1:(\ell-1)} \sum_{k=0}^{n} k\gamma_{s,T}^{[\ell](k)}(\xi_{1:(\ell-1)}, y)\phi(y)^{k-1} 2\dot{\not{\mathbb{G}}}_s(\xi_{1:(\ell-1)}, y) \mathrm{d}\xi_{1:(\ell-1)}
\end{aligned}$$

Moreover, as

$$\mathrm{D}q_s^{[\ell']}(\phi) = \int \psi^{\rho}_{1:(\ell'-1)} \mathrm{D}\Gamma_{s,T}^{[\ell']} + \Gamma_{s,T}^{[\ell']} \mathrm{D}\psi^{\rho}_{1:(\ell'-1)} \mathrm{d}\xi_{1:(\ell'-1)}$$

and

$$\mathrm{D}V_s^{[\ell'']} = i\beta \sum_{\substack{\boldsymbol{b}\in\mathfrak{M}:\\ \ell(\boldsymbol{b})=\ell''}} \int v_s^{[\ell''](\boldsymbol{b})}(\xi_{1:\ell''}) \sum_{j=1}^{\ell''} \sigma_j \delta_{x_j} \psi^{\rho}_{1:\ell''} M_{(\boldsymbol{b})}(x_{1:\ell''}) \mathrm{d}\xi_{1:\ell''}$$

we obtain for $\ell', \ell''$ such that $\ell' + \ell'' = \ell$ by the symmetry of kernels that

$$\begin{aligned}
&\big\langle \mathrm{D}q_s^{[\ell']}, \dot{G}_s \mathrm{D}V_s^{[\ell'']} \big\rangle \\
&= \sum_{\substack{\boldsymbol{b}\in\mathfrak{M}:\\ \ell(\boldsymbol{b})=\ell''}} \left\langle v_s^{[\ell''](\boldsymbol{b})}(\xi_{1:\ell''}) M_{(\boldsymbol{b})}(x_{1:\ell''}) 2\dot{\not{\mathbb{G}}}_s(\xi_{1:\ell''}, y) \sum_{k=0}^{n} k\gamma_{s,T}^{[\ell'](k)}(\xi_{1:(\ell'-1)}, y), \phi(y)^{k-1} \psi^{\rho}_{1:(\ell-1)} \right\rangle \\
&+ \sum_{\substack{\boldsymbol{b}\in\mathfrak{M}:\\ \ell(\boldsymbol{b})=\ell''}} \left\langle v_s^{[\ell''](\boldsymbol{b})}(\xi_{1:\ell''}) M_{(\boldsymbol{b})}(x_{1:\ell''}) 2\dot{\mathbb{G}}_s(\xi_{1:(\ell'-1)}, \xi_{1:\ell''}) \sum_{k=0}^{n} \gamma_{s,T}^{[\ell'](k)}(\xi_{1:(\ell'-1)}, y), \phi(y)^{k} \psi^{\rho}_{1:(\ell-1)} \right\rangle
\end{aligned}$$

hence the differential flow equation (2.19) for the coefficient of field content $\psi_{1:(\ell-1)}$ can be written as

$$\begin{aligned}
&\sum_{k=0}^{n} \partial_s \gamma_{s,T}^{[\ell](k)} \phi(y)^k + \gamma_{s,T}^{[\ell](k)} \phi(y)^k \dot{\mathbb{G}}_s(\xi_{1:(\ell-1)}) \\
&+ \sum_{k=0}^{n} \frac{k(k-1)}{2} \gamma_{s,T}^{[\ell](k)} \phi(y)^{k-2} \dot{\mathcal{G}}_s(0) + 2k\gamma_{s,T}^{[\ell](k)} \phi(y)^{k-1} \dot{\not{\mathbb{G}}}_s(\xi_{1:(\ell-1)}, y) \\
&= \sum_{I_1 \uplus I_2 = 1:\ell} \sum_{\substack{\boldsymbol{b}\in\mathfrak{M}:\\ \ell(\boldsymbol{b})=|I_2|}} C_{|I_1|,|I_2|} v_s^{[|I_2|](\boldsymbol{b})}(\xi_{I_2}) M_{(\boldsymbol{b})}(x_{I_2}) \dot{\not{\mathbb{G}}}_s(\xi_{I_2}, y) \sum_{k=0}^{n} k\gamma_{s,T}^{[|I_1|](k)}(\xi_{I_1-1}, y)\phi(y)^{k-1} \\
&+ \sum_{I_1 \uplus I_2 = 1:\ell} \sum_{\substack{\boldsymbol{b}\in\mathfrak{M}:\\ \ell(\boldsymbol{b})=|I_2|}} C_{|I_1|,|I_2|} v_s^{[|I_2|](\boldsymbol{b})}(\xi_{I_2}) M_{(\boldsymbol{b})}(x_{I_2}) \dot{\mathbb{G}}_s(\xi_{I_1-1}, \xi_{I_2}) \sum_{k=0}^{n} \gamma_{s,T}^{[|I_1|](k)}(\xi_{I_1-1}, y)\phi(y)^{k}.
\end{aligned}$$

Matching powers of $\phi(y)$ we obtain a system of equations for the kernels where corresponding to $\phi(y)^k$ we have

$$\partial_s \gamma_{s,T}^{[\ell](k)} + \dot{\mathbb{G}}_s(\xi_{1:(\ell-1)}) \gamma_{s,T}^{[\ell](k)} = -\mathscr{S}_{s,T}^{[\ell](k)}$$

and the statement follows from the Duhamel principle similarly to the proof of Lemma 3.1. □

## B.3 UV stability estimates

**Proof. (of Corollary 3.15)** We are going show that for $\varepsilon > 0$ sufficiently small depending only on $\beta$, it holds for all $t \in \mathbb{R}_+$ that

$$|||v_{t,\infty}^{[\ell](\boldsymbol{a})} - v_{t,T}^{[\ell](\boldsymbol{a})}|||_{t,\Theta(\ell),\eta} \lesssim \lambda_t (\lambda_t \langle t\rangle^{-1})^{\ell-1} \mathfrak{l}_t^c \langle t\rangle^{-\frac{\varkappa(\boldsymbol{a})}{2}} \langle t\rangle^{\varepsilon} \langle T\rangle^{-\varepsilon}. \tag{B.2}$$

The identically zero kernels, i.e. those with $\boldsymbol{a} \in \mathfrak{M} \setminus \mathfrak{M}_p$ and $\varkappa(\boldsymbol{a}) = 1$, are trivial. For the nonzero kernel at $\ell = 1$, the limit is $v_{t,\infty}^{[1](\boldsymbol{a})} := \lambda e^{\frac{\beta^2}{8\pi}\mathcal{G}_t(0)}$: since $t \mapsto \lambda_t$ is increasing,

$$|||v_{t,\infty}^{[1](\boldsymbol{a})} - v_{t,T}^{[1](\boldsymbol{a})}|||_t = |\lambda_t - \lambda_{t\wedge T}| = \mathbb{1}_{\{t\geq T\}}|\lambda_t - \lambda_T| = \lambda_t \mathbb{1}_{\{t\geq T\}}\left(1 - \frac{\lambda_T}{\lambda_t}\right) = \lambda_t \mathbb{1}_{\{t\geq T\}}\left(1 - e^{\frac{\beta^2}{2}\mathcal{G}_{t,T}(0)}\right)$$

and using that $1 - e^{-x} \leq x$ for $x \geq 0$ and Lemma A.1

$$|||v_{t,\infty}^{[1](\boldsymbol{a})} - v_{t,T}^{[1](\boldsymbol{a})}|||_t \leq \lambda_t \mathbb{1}_{\{t\geq T\}} \frac{\beta^2}{2}\mathcal{G}_{T,t}(0) \leq \lambda_t \mathbb{1}_{\{t\geq T\}} \frac{\beta^2}{8\pi}(\log(t \vee 1) + C) \lesssim \lambda_t \mathfrak{l}_t^c \langle t\rangle^{\varepsilon}\langle T\rangle^{-\varepsilon}.$$

Moving on to $\ell \geq 2$, referring back to Section 2.7, we are going to show that the limits $v_t^{[\ell](\boldsymbol{a})}$ satisfy

$$\begin{aligned} v_{t,\infty}^{[\ell](\boldsymbol{a})} &= \int_t^\infty e^{\mathbb{G}_{t,s}}\mathscr{S}_{s,\infty}^{(\boldsymbol{a})}\,\mathrm{d}s \quad \text{if } \boldsymbol{a} \text{ is such that } v_{T,T}^{[\ell](\boldsymbol{a})} = 0 \\ v_{t,\infty}^{[\ell](\boldsymbol{a})} &= -\int_0^t \mathscr{S}_{s,\infty}^{(\boldsymbol{a})}\,\mathrm{d}s \qquad \text{otherwise,} \end{aligned} \tag{B.3}$$

where $\mathscr{S}_{s,\infty}^{(\boldsymbol{a})}$ are source terms obtained from the limiting kernels; since the system is triangular, we may proceed inductively. Since bounds of (3.29) are uniform in $T$, they pass to the limits.
For kernels with irrelevant boundary conditions we have

$$\begin{aligned} v_{t,\infty}^{[\ell](\boldsymbol{a})} - v_{t,T}^{[\ell](\boldsymbol{a})} &= \int_t^\infty e^{\mathbb{G}_{t,s}}\mathscr{S}_{s,\infty}^{(\boldsymbol{a})}\,\mathrm{d}s - \int_{t\wedge T}^T e^{\mathbb{G}_{t,s}}\mathscr{S}_{s,T}^{(\boldsymbol{a})}\,\mathrm{d}s \\ &= \int_{t\vee T}^\infty e^{\mathbb{G}_{t,s}}\mathscr{S}_{s,\infty}^{[\ell](\boldsymbol{a})}\mathrm{d}s + \int_{t\wedge T}^T e^{\mathbb{G}_{t,s}}(\mathscr{S}_{s,\infty}^{[\ell](\boldsymbol{a})} - \mathscr{S}_{s,T}^{[\ell](\boldsymbol{a})})\mathrm{d}s \\ &=: \mathbb{I} + \mathbb{II}. \end{aligned}$$

The first term is bounded via Proposition 3.14:

$$\begin{aligned} |||\mathbb{I}|||_{t,\Theta(\ell)} &\lesssim \int_{t\vee T}^\infty \|e^{\mathbb{G}_{t,s}}\|_{L^\infty}|||\mathscr{S}_{s,\infty}^{[\ell](\boldsymbol{a})}|||_{s,\Theta(\ell)}\mathrm{d}s \\ &\lesssim \lambda_t^{|\chi|}\int_{t\vee T}^\infty \lambda_s^{-|\chi|}\lambda_s(\lambda_s\langle s\rangle^{-1})^{\ell-1}\mathfrak{l}_s^c\langle s\rangle^{-\frac{\varkappa(\boldsymbol{a})}{2}}\langle s\rangle^{-1}\mathrm{d}s \\ &\lesssim \lambda_t^{|\chi|}\lambda_{t\vee T}^{-|\chi|}\lambda_{t\vee T}(\lambda_{t\vee T}\langle t\vee T\rangle^{-1})^{\ell-1}\mathfrak{l}_{t\vee T}^c\langle t\vee T\rangle^{-\frac{\varkappa(\boldsymbol{a})}{2}} \end{aligned}$$

Now choosing $\varepsilon > 0$ to be small enough, depending only on $\delta$, $s \mapsto \lambda_s^{-|\chi|}\lambda_s(\lambda_s\langle s\rangle^{-1})^{\ell-1}\mathfrak{l}_s^c\langle s\rangle^{-\frac{\varkappa(\boldsymbol{a})}{2}}\langle s\rangle^{\varepsilon}$ is decreasing when i

$$1 - \ell\delta - |\chi|(1-\delta) - \frac{\varkappa(\boldsymbol{a})}{2} + \varepsilon < 0 \quad \Longleftarrow \quad (\ell \geq 4) \vee (|\chi| \geq 1) \vee (\varkappa(\boldsymbol{a}) = 1) \vee (\delta > 1/2)$$

since $\delta > 1/4$ and $\ell \geq 2$; all irrelevant kernels satisfy this condition. Thus

$$|||\mathbb{I}|||_{t,\Theta(\ell)} \lesssim \lambda_t(\lambda_t\langle t\rangle^{-1})^{\ell-1}\mathfrak{l}_t^c\langle t\rangle^{-\frac{\varkappa(\boldsymbol{a})}{2}}\langle t\rangle^{\varepsilon}\langle T\rangle^{-\varepsilon}.$$

To estimate $\mathbb{II}$, we streamline notation by defining $\langle\langle v_{s,T}^{[|I_1|](\boldsymbol{b})}, v_{s,T}^{[|I_1|](\boldsymbol{c})}\rangle_{\dot{\mathbb{G}}}, \Psi_{\boldsymbol{a}}\rangle$ via

$$\langle \mathrm{D}\langle v_{s,T}^{[|I_1|](\boldsymbol{b})}, \Psi_{\boldsymbol{b}}\rangle, \dot{G}_s\mathrm{D}\langle v_{s,T}^{[|I_2|](\boldsymbol{c})}, \Psi_{\boldsymbol{c}}\rangle\rangle = \langle\langle v_{s,T}^{[|I_1|](\boldsymbol{b})}, v_{s,T}^{[|I_1|](\boldsymbol{c})}\rangle_{\dot{\mathbb{G}}}, \Psi_{\boldsymbol{a}}\rangle$$

where $\Psi_{\boldsymbol{a}}$ denotes the field content arising from contracting $\Psi_{\boldsymbol{b}}$ and $\Psi_{\boldsymbol{c}}$ in the final representation (that is, possibly after a nontrivial post-processing as in Lemma 3.8 or the lowering of the degree of divergence as in cases 2. and 3. of Proposition 3.14); That is,

$$\langle\mathscr{S}_{s,T}^{[\ell](\boldsymbol{a})}, \Psi_{\boldsymbol{a}}\rangle = \langle \mathbf{A}_{\boldsymbol{d}}^{\boldsymbol{a}}(v_{s,T}^{[\ell](\boldsymbol{d})}), \Psi_{\boldsymbol{a}}\rangle + \sum_{\substack{I_1 \uplus I_2 = 1:\ell \\ \boldsymbol{b},\boldsymbol{c}\in\mathfrak{M}}} \langle\langle v_{s,T}^{[|I_1|](\boldsymbol{b})}, v_{s,T}^{[|I_1|](\boldsymbol{c})}\rangle_{\dot{\mathbb{G}}}, \Psi_{\boldsymbol{a}}\rangle$$

where $\mathbf{A}^{\boldsymbol{a}}_{\boldsymbol{d}}\big(v^{[\ell](\boldsymbol{d})}_{s,T}\big) := \mathbf{A}\big(v^{[2](\varkappa)}_{s,T}\big)$ introduced in (3.17).

Then using Proposition 3.14 and the induction on (B.2), by the linearity of $\mathbf{A}^{\boldsymbol{a}}_{\boldsymbol{d}}$ we obtain

$$|||\mathbf{A}^{\boldsymbol{a}}_{\boldsymbol{d}}\big(v^{[\ell](\boldsymbol{d})}_{s,\infty}\big) - \mathbf{A}^{\boldsymbol{a}}_{\boldsymbol{d}}\big(v^{[\ell](\boldsymbol{d})}_{s,T}\big)|||_s \lesssim \lambda_s(\lambda_s\langle s\rangle^{-1})^{\ell-1}\mathfrak{l}^c_s\langle s\rangle^{-\frac{\varkappa(\boldsymbol{a})}{2}}\langle s\rangle^{\varepsilon}\langle T\rangle^{-\varepsilon}\langle s\rangle^{-1} \tag{B.4}$$

and by the bilinearity of $\big\langle v^{[|I_1|](\boldsymbol{b})}_{s,T}, v^{[|I_2|](\boldsymbol{c})}_{s,T}\big\rangle_{\dot{\mathbb{G}}}$

$$\begin{aligned} & |||\big\langle v^{[|I_1|](\boldsymbol{b})}_{s,T}, v^{[|I_2|](\boldsymbol{c})}_{s,T}\big\rangle_{\dot{\mathbb{G}}} - \big\langle v^{[|I_1|](\boldsymbol{b})}_{s,\infty}, v^{[|I_2|](\boldsymbol{c})}_{s,\infty}\big\rangle_{\dot{\mathbb{G}}}|||_s \\ \lesssim\ & |||\big\langle v^{[|I_1|](\boldsymbol{b})}_{s,T} - v^{[|I_1|](\boldsymbol{b})}_{s,\infty}, v^{[|I_2|](\boldsymbol{c})}_{s,T}\big\rangle_{\dot{\mathbb{G}}} - \big\langle v^{[|I_1|](\boldsymbol{b})}_{s,\infty}, v^{[|I_2|](\boldsymbol{c})}_{s,\infty} - v^{[|I_2|](\boldsymbol{c})}_{s,T}\big\rangle_{\dot{\mathbb{G}}}|||_s \\ \lesssim\ & \lambda_s(\lambda_s\langle s\rangle^{-1})^{\ell-1}\mathfrak{l}^c_s\langle s\rangle^{-\frac{\varkappa(\boldsymbol{a})}{2}}\langle s\rangle^{\varepsilon}\langle T\rangle^{-\varepsilon}\langle s\rangle^{-1} \end{aligned} \tag{B.5}$$

so that

$$\begin{aligned} |||\mathbb{III}|||_{t,\Theta(\ell)} &\leq \int_t^T \|e^{\mathbb{G}_{t,s}}\|_{L^\infty}|||\mathscr{S}^{[\ell](\boldsymbol{a})}_{s,\infty} - \mathscr{S}^{[\ell](\boldsymbol{a})}_{s,T}|||_{s,\Theta(\ell)}\mathrm{d}s \\ &\lesssim \lambda_t^{|\chi|}\int_t^\infty \lambda_s^{-|\chi|}\lambda_s(\lambda_s\langle s\rangle^{-1})^{\ell-1}\mathfrak{l}^c_s\langle s\rangle^{-\frac{\varkappa(\boldsymbol{a})}{2}}\langle s\rangle^{\varepsilon}\langle T\rangle^{-\varepsilon}\langle s\rangle^{-1}\mathrm{d}s \\ &\lesssim \lambda_t(\lambda_t\langle t\rangle^{-1})^{\ell-1}\mathfrak{l}^c_t\langle t\rangle^{-\frac{\varkappa(\boldsymbol{a})}{2}}\langle t\rangle^{\varepsilon}\langle T\rangle^{-\varepsilon} \end{aligned}$$

as requested.

It remains to treat kernels with relevant boundary conditions. In the regime we are working in, the only such kernel is the one arising from the $\ell = 2$, $\chi(\xi_{1:2}) = 0$ dipole when $\delta \in \big(\frac{1}{4}, \frac{1}{2}\big]$, cf. the discussion preceding Lemma 3.8. Then

$$\begin{aligned} v^{[\ell](\boldsymbol{a})}_{t,\infty} - v^{[\ell](\boldsymbol{a})}_{t,T} &= \int_0^{t\wedge T}\mathscr{S}^{[\ell](\boldsymbol{a})}_{s,T}\mathrm{d}s - \int_0^t \mathscr{S}^{[\ell](\boldsymbol{a})}_{s,\infty}\mathrm{d}s \\ &= \int_0^{t\wedge T}(\mathscr{S}^{[\ell](\boldsymbol{a})}_{s,T} - \mathscr{S}^{[\ell](\boldsymbol{a})}_{s,\infty})\mathrm{d}s - \int_{t\wedge T}^t \mathscr{S}^{[\ell](\boldsymbol{a})}_{s,\infty}\mathrm{d}s \end{aligned}$$

The first term is bounded via the induction on (B.2) and Proposition 3.14 similarly to (B.4) and (B.5) and using that the source term is local:

$$\begin{aligned} \int_0^{t\wedge T}|||\mathscr{S}^{[\ell](\boldsymbol{a})}_{s,\infty} - \mathscr{S}^{[\ell](\boldsymbol{a})}_{s,T}|||_{t,\Theta(\ell)}\mathrm{d}s &\lesssim \int_0^t \lambda_s(\lambda_s\langle s\rangle^{-1})^{\ell-1}\mathfrak{l}^c_s\langle s\rangle^{-\frac{\varkappa(\boldsymbol{a})}{2}}\langle s\rangle^{\varepsilon}\langle T\rangle^{-\varepsilon}\langle s\rangle^{-1}\mathrm{d}s \\ &\lesssim \lambda_t(\lambda_t\langle t\rangle^{-1})^{\ell-1}\mathfrak{l}^c_t\langle t\rangle^{-\frac{\varkappa(\boldsymbol{a})}{2}}\langle t\rangle^{\varepsilon}\langle T\rangle^{-\varepsilon} \end{aligned} \tag{B.6}$$

since the integral diverges when the kernel is relevant (recall the proof of Lemma 3.9). The second term is trivial: again by locality and Proposition 3.14,

$$\int_{t\wedge T}^t |||\mathscr{S}^{[\ell](\boldsymbol{a})}_{s,\infty}|||_{t,\Theta(\ell)}\mathrm{d}s \lesssim \mathbb{1}_{\{T<t\}}\lambda_t(\lambda_t\langle t\rangle^{-1})^{\ell-1}\mathfrak{l}^c_t\langle t\rangle^{-\frac{\varkappa(\boldsymbol{a})}{2}}$$

and we conclude from $\mathbb{1}_{\{T<t\}} \leq \langle t\rangle^{\varepsilon}\langle T\rangle^{-\varepsilon}$. This completes the proof on the claim about the kernel limits.

To prove the statement about $V^{[\ell],\infty}(\phi)$, note that (B.3) implies that the limiting kernels satisfy the same differential flow equation as the UV-cutoff kernels, and that the symmetries, in particular permutation invariance also extend to the limit. Thus we may conclude from the first part of the proof of Lemma 3.1 and Lemmas 3.8 and 3.14. □